**Michiel Hazewinkel** 1 CWI
Direct line: +31-20-5924204 POBox 94079
Secretary: +31-20-5924233 1090GB Amsterdam
Fax: +31-20-5924166
E-mail: mich@cwi.nl

# Witt vectors. Part 1[1]

by

*Michiel Hazewinkel*
*CWI*
*POBox 94079*
*1090GB Amsterdam*
*The Netherlands*

**Table of contents**

[1] Part 2 (and last) is planned for volume 7 of the Handbook of Algebra.

Part 2 of this survey chapter will include discussions of the following topics: **Symm** vs **BU**; **Symm** vs Exp; **Symm** vs $R_{\mathrm{rat}}(\mathbf{GL}_\infty)$; Formal groups; Classification of affine commutative algebraic groups. Cartier duality; Covectors and bivectors; The $K$-theory of endomorphisms; de Rham-Witt complex. Witt vector cohomology; The lifting game; Witt vectors and deformations; Ramified Witt vectors; Noncommutative Witt vectors; The elementary particle view of **Symm**; Brief comments on more results on and applications of Witt vectors.

## 1. Introduction and delimitation

Let **CRing** be the category of commutative rings with unit element. Most of this chapter is about the big (or large or generalized) Witt vectors; that is about a certain functor

$$W\colon \mathbf{CRing} \longrightarrow \mathbf{CRing} \tag{1.1}$$

that (therefore) assigns to each unital commutative ring $A$ a new unital commutative ring $W(A)$ (and to a unit element preserving ring morphism $A \longrightarrow B$ a unital ring morphism $W(A) \longrightarrow W(B)$). It is also about the many quotient functors of $W$ of which the most important are the $p$-adic Witt vectors $W_{p^\infty}$ , the truncated big Witt vectors $W_n$, and the truncated $p$-adic Witt vectors $W_{p^n}$. Instead of 'truncated' one also finds 'of finite length' in the literature.

The $p$-adic Witt vectors functor gets defined via the socalled $p$-adic Witt polynomials

$$Y_0,\quad Y_0^p + pY_1,\quad Y_0^{p^2} + pY_1^p + p^2Y_2,\quad Y_0^{p^3} + pY_1^{p^2} + p^2Y_2^p + p^3Y_3,\quad \cdots \tag{1.2}$$

and their generalizations, the big Witt polynomials

$$
\begin{aligned}
w_1(X) &= X_1 \\
w_2(X) &= X_1^2 + 2X_2 \\
\cdots \quad & \quad \cdots \\
w_n(X) &= \sum_{d|n} dX_d^{n/d} \\
\cdots \quad & \quad \cdots
\end{aligned}
\tag{1.3}
$$

can [2] be used to define the big Witt vectors. To see that (1.3) specializes to (1.2) relabel $X_{p^i}$ in $w_{p^n}(X)$ as $Y_i$.

Here and there in the published literature one finds the Witt polynomials referred to as "mysterious polynomials that come out of nowhere". It is one of my aims in the present screed to try to argue that that is simply not true. There is something inevitable about the Witt polynomials and they turn up naturally in various contexts and with some frequency.

1.4. Historically, this is also how things started. The setting is that of the investigations by Helmut Hasse and his students and collaborators into the structure of complete discrete valuation rings $A$ with residue field $k$. Oswald Teichmüller discovered that in such a situation there is a multiplicative system of representatives, i.e. a multiplicative section of the natural projection $A \longrightarrow k$. Such a system is unique if $k$ is perfect [3]. These are now called Teichmüller representatives. Assuming that $A$ is unramified, every element of $A$ can be written as a power series in $T$, a generator of the maximal ideal of $A$, with coefficients from any chosen system of representatives. Given the nice properties of the Teichmüller system it was and is natural to use these. As they are already multiplicative the first problem was of course to figure out how Teichmüller representatives should be added in the arithmetic of the ring $A$. In the equal characteristic case, $\text{char}(A) = \text{char}(k) = p > 0$ that is no problem at all. The Teichmüller system is then additive as well. But in the unequal characteristic case, $\text{char}(A) = 0$, $\text{char}(k) = p > 0$ things are very different. In this case one can and does choose $T = p$. In 1936 Teichmüller, [382], found a formula for doing precisely this, i.e. adding Teichmüller representatives. The formula is, see loc. cit. p. 156

$$a + b = \sum_{n=0}^{\infty} c_n p^n \tag{1.5}$$

Here $a$ and $b$ are Teichmüller representatives of, say, $\bar{a}$, $\bar{b}$, and the $c_n$ are the Teichmüller representatives of elements $\bar{c}_n$ in $k$ that satisfy

$$\bar{c}_n^{\,p^n} = r_n(\bar{a}, \bar{b}) \tag{1.6}$$

where the $r_n(X,Y)$ are the integer valued polynomials (recursively) determined by

$$r_0(X,Y)^{p^n} + pr_1(X,Y)^{p^{n-1}} + \cdots + p^n r_n(X,Y) = X^{p^n} + Y^{p^n} \tag{1.7}$$

And here they surface, on the left hand side of (1.7), the $p$-adic Witt polynomials.

1.8. Some 30+ years later, in another context, a very similar picture emerged. As a functor the big Witt vectors are isomorphic to the functor $\Lambda$ that assigns to a commutative unital ring $A$ the Abelian group $\Lambda(A)$ of power series with constant term 1 and coeffients in $A$ (under multiplication of power series). Much of Witt vector theory can be developed from this point of view without ever mentioning the generalized Witt polynomials (1.3). Except that it is far from obvious how to get the $p$-adic Witt vectors in this picture. In the main, this power series treatment is what is given below. But the Witt vectors themselves and the Witt polynomials will not be denied, as will be seen now.

[2] But there are other ways. One of these will be used below (most of the time).

[3] If k is not perfect a multiplicative system of representatives still exists but is not unique.

For any functor one is (or at least should be) interested in the operations on it, i.e. its functorial endomorphisms. In the present case a number of obvious and easy operations are the homothety operators $\langle u\rangle f(t) = f(ut)$ and the Verschiebung operators $\mathbf{V}_n f(t) = f(t^n)$. Rather less obvious are the Frobenius operators $\mathbf{f}_n$ which are as much like raising to the $n$-th power as a morphism of Abelian groups can be. Now what is the sum (pointwise sum of operations) of two homothety operators? The answer is

$$\langle u\rangle + \langle v\rangle = \sum_{n=1}^{\infty} \mathbf{V}_n \langle r_n(u,v)\rangle \mathbf{f}_n \tag{1.9}$$

where this time the polynomials $r_n(X,Y)$ are the integer valued polynomials determined recursively by

$$\sum_{d|n} d r_d(X,Y)^{n/d} = X^n + Y^n \tag{1.10}$$

and there they are, the generalised Witt polynomials (1.3) [4].

To appreciate the near perfect fit of the picture of formulas (1.5)-(1.7) with that of (1.9)-(1.10), reflect that in a field of characteristic $p$ raising to the $p$-th power is the Frobenius morphism and that multiplication by $p$ (in a charateristic zero ring, say, the $p$-adic integers) has a shift built into it. Further, the $\langle u\rangle$, which are multiplicative, can be seen, in a very real sense, as Teichmüller representatives.

1.11. The underlying set of $W(A)$ or, better (for the moment), $\Lambda(A)$, is the set of power series over $A$ with constant term $1$ (the unit element of $A$):

$$\Lambda(A) = \{1 + a_1 t + a_2 t^2 + a_3 t^3 + \cdots : a_i \quad A \text{ for all } i \quad \mathbf{N}\} \tag{1.12}$$

(where, as usual, $\mathbf{N}$ denotes the natural numbers $\mathbf{N} = \{1,2,3,\cdots\}$). This functor (to **Set**, the category of sets), is obviously representable by the ring

$$\mathbf{Symm} = \mathbf{Z}[h_1, h_2, h_3, \cdots] \tag{1.13}$$

of polynomials in an infinity of indeterminates over the integers $\mathbf{Z}$.

As will rapidly become clear one really should see the $h_i$ as the complete symmetric functions (polynomials) [5] in an infinity of indeterminates [6] $\xi_1, \xi_2, \xi_3, \cdots$ . Whence the notation

---

[4] But in this case, I believe, the generalized Witt polynomials were known before formula (1.8) was written down.

[5] Of course **Symm** is also equal to the polynomials in the elementary symmetric functions in the same infinity of indeterminates and one could also work with those (and that is often done). But using the complete

used in (1.13). Also it turns out that many constructions, results, ... for the Witt vectors have their natural counterparts in the theory of symmetric functions (as should be). Thus a chapter on the Witt vectors could very properly include most of symmetric function theory, that very large subject that could fill several volumes. This will, of course, not be done. Also, there will be a separate chapter on the symmetric functions in this Handbook of Algebra.

1.14. **Symm**, the Hopf algebra of the symmetric functions is a truly amazing and rich object. It turns up everywhere and carries more extra structure than one would believe possible. For instance it turns up as the homology of the classifying space **BU** and also as the the cohomology of that space, illustrating its selfduality. It turns up as the direct sdum of the representation spaces of the symmetric group and as the ring of rational representations of the infinite general linear group. This time it is Schur duality that is involved. It is the free $\lambda$-ring on one generator. It has a nondegenerate inner product which makes it selfdual and the associated orthonormal basis of the Schur symmetric functions is such that coproduct and product are positive with respect to these basis functions. In this setting, positive, self dual, Hopf algebras with distinguished basis, there is a uniqueness theorem due to Andrei Zelevinsky, [425]. But this is not yet entirely satisfactory unless it can be shown that the Schur symmetric functions are in some algebraic way canonical (which seems very likely). **Symm** is also the representing ring of the functor of the big Witt vectors and the covariant bialgebra of the formal group of the big Witt vectors (another manisfestation of its autoduality). Most of these things will be at least touched on below.

As the free $\lambda$-ring on one generator it of course carries a $\lambda$-ring structure. In addition it carries ring endomorphisms which define a functorial $\lambda$-ring structure on the rings $W(A) = \mathbf{CRing}(\mathbf{Symm}, A)$ for all unital commutative rings $A$. A sort of higher $\lambda$-ring structure. Being selfdual there are also co-$\lambda$-ring structures and higher co-$\lambda$-ring structures (whatever those may be).

Of course, **Symm** caries still more structure: it has a second multiplication and a second comultiplication (dual to each other) that make it a coring object in the category of algebras and, dually, (almost) a ring object in the category of coalgebras.

The functor represented by **Symm**, i.e. the big Witt vector functor, has a comonad structure and the associated coalgebras are precisely the $\lambda$-rings.

All this by no means exhausts the manifestations of and structures carried by **Symm**. It seems unlikely that there is any object in mathematics richer and/or more beautiful than this one, and many more uniqueness theorems are needed.

In this chapter I will only touch upon the aspects of **Symm** which relate directly to the Witt vector constructions and their properties.

To conclude this introduction let me remark that the Witt vector construction is a very beautiful one. But is takes one out of the more traditional realms of algebra very quickly. For instance the ring of $p$-adic Witt vectors of a field is Noetherian if and only if that field is perfect, [63], chapître IX, page 43, exercise 9. Also, much related, the $p$-adic Witt vector ring $W_{p^\infty}(k[T])$ of one dimensional affine space over a characteristic $p > 0$ field is not Noetherian, which would appear to rule out any kind of systematic use of the Witt vectors in algebraic geometry as currently practiced [7]. It is perhaps because of this that the Witt vector functors and the rings and symmetric functions works out just a bit more elegantly, especially in connection with the autodulaity of the Hopf algebra **Symm**. See below in subsection 11.37 and section 12.

[6] Some of the readers hoped for may not be familiar with working with symmetric polynomials in an infinity of determinates. There is really nothing to it. But for those that do not feel confortable about it (and for the pernickety) there is a short appendix on the matter.

[7] This is not really true; witness cristalline cohomologyy and the de Rham-Witt complex (of which cristalline cohomology is the hyper homology) which is made up of W($k$) modules.

algebras thus arising have not really been much studied[8], except in so far as needed for their applications (which are many and varied).

1.15. *Initial sources of information*. There are a number of (electronic) encyclopedia articles on the Witt vectors, mostly the $p$-adic Witt vectors, see [4, 5, 6, 7]. These can serve to obtain a first sketchy impression. There are also a number of (introductory) lecture notes and chapters in books on the subject, see e.g. [57, 74, 86, 188], [192 Chapter 3], [212, 263], [271 §26], [326 Lectures 6-12], [336]. However, the full story of the Witt vectors is a long and varied one. Here, in the present chapter, I try to present a first outline.

## 2. Terminology.

The big and $p$-adic and truncated Witt vectors carry ring and algebra structures, and hence, naturally, are sometimes referred to as a ring or algebra of Witt vectors and then (by erosion) as Witt ring and Witt algebra. This is a bit unfortunate and potentially confusing because these phrases mostly carry other totallly unrelated meanings.

Mostly, 'Witt ring' refers to a ring of equivalence classes of quadratic forms (with addition and multiplication induced by direct product and tensor product respectively), see e.g. [96, 296].

Further 'Witt algebra' mostly refers to something like the Lie algebra of differential operators spanned by the $d_i = x^{i+1}\frac{d}{dx}$, $i \in \mathbf{Z}$ under the commutator product $[d_i, d_j] = (j-i)d_{i+j}$, i.e. the centerless Virasoro algebra, much studied in physics; see e.g. [182, 410, 426]. The fact that the sub Lie algebra spanned by the $d_i$ for $i \geq -1$ is usually denoted $W_1$ adds a bit more potential confusion, and this becomes worse when one encounters $W_n$ for the more variable version of this algebra of differential operators.

On the other hand the term 'Witt group' mostly indeed refers to a group of Witt vectors[9].

## 3. The *p*-adic Witt vectors. More historical motivation.

It has become customary to 'motivate' the introduction of the Witt vectors by looking at unramified complete discrete valuation fields, or, more precisely their rings of integers. To start with, consider the ring of $p$-adic integers $\mathbf{Z}_p$, i.e. the completion of the integers with respect to the norm $\|n\| = p^{-v_p(n)}$ where the valution $v_p(n)$ of an integer $n$ is the largest power of $p$ that divides $n$. Here $p$ is a prime number.

Every $p$-adic integer $\alpha$, i.e. element of $\mathbf{Z}_p$, can be uniquely written as a convergent sum

$$\alpha = a_0 + a_1 p + a_2 p^2 + a_3 p^3 + \cdots, \; a_i \in \{0,1,\cdots,p-1\} \tag{3.1}$$

More generally, instead of the set $\{0,1,\cdots,p-1\}$ one can choose any set of representatives of the residue field $\mathbf{Z}/(p) = \mathbf{F}_p = \mathbf{GF}(p)$[10]. Now what happens if two such expressions are added or multiplied; i.e what are the coefficients of a sum or product of $p$-adic integers? As in the case of the familiar decimal notation for arithmetic with integers this involves carry-overs. That is perfectly acceptable for a calculating machine but irritating to algebraists who would like universal formulas which always work. Things get much worse if the residue field is not a prime field where it is not even clear (a priori) what set of representatives to choose. Of course at this stage it is far from obvious whether a universal formula can exist; that even seems unlikely at first sight.

[8] Should a reader inadvertently get really interested in the Witt vector ring functors he/she is recommended to work through the 57 exercises on the subject in [Bourbaki, 1983 #180], 26 pp worth for the statements only, mostly contributed, I have been told, by Pierre Cartier.

[9] But not always, cf e.g. [82], which is again about quadratic form related matters.

[10] To use three of the standard notations for the prime field of $p$ elements.

As it turns out there are universal formulas: the addition and multiplication polynomials of the $p$-adic Witt vectors.

However, judging from the introduction in Witt's seminal paper [420] the arithmetic of $p$-adic fields was not Witt's primary motivation. This seems also indicated by the fact that this arithmetic is not mentioned in the title of the paper, but only in the subtitle. Instead Witt seems to have been mainly motivated by the desire to obtain a theory of cyclic Galois extensions for general charateristic $p$ fields similar to the theory he had himself obtained for power series fields in [415, 416] and, especially, to understand some mysterious formulas of Hermann Ludwig Schmid, [355] concerning central cyclic algebras in characteristic $p$. All this more or less in the framework of a class field theory for function fields that was being vigorously developed at the time by H Hasse, F K Schmidt, O Teichmüller, H L Schmid, C Chevalley, E Witt and others. See [338] for a thorough and very readable account of this part of the history of class field theory. Some of the important original papers are [16, 23, 87, 191, 190, 355, 356, 357, 382, 415, 416, 420].

The situation as regards these cyclic central algebras was as follows. For $\alpha,\beta \in k$ (of characteristic $p>0$) let $(\alpha,\beta]$ denote the simple central cyclic [11] algebra of degree $p$ with two generators $u,\theta$ subject to the relations

$$u^p=\alpha,\ \ \theta^p-\theta=\beta,\ \ u\theta u^{-1}=\theta+1$$

In the Brauer group there are the relations

$$(\alpha,\beta]\otimes(\alpha',\beta]=(\alpha\alpha',\beta]\ \text{ and }\ (\alpha,\beta]\otimes(\alpha,\beta']=(\alpha,\beta+\beta']$$

H L Schmid in [355] succeeded in defining simple central cyclic algebras of degree $p^n$ denoted $(\alpha|\beta]=(\alpha|\beta_0,\beta_1,\cdots,\beta_{n-1}]$ and by dint of some heroic formula manipulation found rules of the form

$$(\alpha|\beta]\otimes(\alpha'|\beta]=(\alpha\alpha'|\beta]\ \text{ and }\ (\alpha|\beta]\otimes(\alpha|\beta']=(\alpha|s_0(\beta,\beta'),\cdots,s_{n-1}(\beta,\beta')]$$

with, as Witt writes, “certain polynomials $s_i(x,y)$ that are initially defined in characteristic 0 and only after they have been proved to be integral are taken modulo $p$”. In his fundamental paper [420] Witt proves that these formulas are in fact the addition formulas of the $p$-adic Witt vectors; see also section 5 below. I will also briefly return to these algebra in section 8 below.

So, here is already a third way in which the Witt vector polynomials turn up naturally and unavoidably.

## 4. Teichmüller representatives

Let $A$ be a ring with an ideal $\mathfrak{m}$ such that the quotient ring $k=A/\mathfrak{m}$ of characteristic $p>0$ is perfect, which means that the ring morphism $\mathbf{f}_p\colon k\longrightarrow k,\ x\mapsto x^p$, the Frobenius morphism, is supposed to be bijective [12]. Suppose,, moreover, that $A$ is complete in the $\mathfrak{m}$-adic topology. For instance $A$ can be the ring of $p$-adic numbers $\mathbf{Z}_p$,, the quotient ring $\mathbf{Z}/(p^m)$, or the ring of power series $k((T))$ over a perfect field $k$ of characteristic $p$. Here the corresponding ideals $\mathfrak{m}$ are repectively the principal ideals $p\mathbf{Z}$, $p\mathbf{Z}/(p^m)$, and $(T)$.

[11] The word ‘cyclic’ here refers to the fact that the algebras considered are cross products involving a cyclic Galois group; see [96], chapter 7, especially §7.2 and §7.5 for more detail.. See also section 8 below.

There is now the following simple observation:

$$\text{For all } a,b \quad A, \text{ if } a \equiv b \mod \mathfrak{m}^r,\ r \geq 0, \text{ then } a^p \equiv b^p \mod \mathfrak{m}^{r+1} \tag{4.1}$$

Now for any $x \quad k$ take lifts $y_r \quad A$ of $\sigma^{-1}(x)$, $r = 0,1,2,\cdots$; i.e. $q(y_r) = \sigma^{-r}(x)$ where $q: A \to k$ is the canonical projection, so that all the $y_r^{p^r}$ are lifts of $x$. Now consider the sequence $y_0, y_1^p, y_2^{p^2}, \cdots, y_r^{p^r}, \cdots$. It follows from (4.1) and the completeness hypothesis on $A$ that this sequence converges to a limit that is a lift of $x$. Moreover , again by (4.1), this limit does not depend on the choice of the $y_r$. This limit is denoted $t(x)$ or $t_A(x)$ and called the Teichmüller representative of $x$. This Teichmüller system of representatives has the following properties

$$t(0) = 0,\ t(1) = 1,\ t(xx') = t(x)t(x') \tag{4.2}$$

i.e. it is multiplicative. (And if $A$ is also of characteristic $p$ it is also additive: $t_A(x + x') = t_A(x) + t_A(x')$. But in general this is not the case.)

The Teichmüller system is also the unique multiplicative one and the unique one which commutes with $p$-th powers.

## 5. Construction of the functor of the $p$-adic Witt vectors

Let $p$ be a prime number and consider the following polynomials in a countable infinity of commuting indeterminates $X_0, X_1, X_2, \cdots$

$$\begin{aligned} w_0(X) &= X_0 \\ w_1(X) &= X_0^p + pX_1 \\ \cdots & \quad \cdots \\ w_n(X) &= X_0^{p^n} + pX_1^{p^{n-1}} + \cdots + p^{n-1}X_{n-1}^p + p^n X_n \\ \cdots & \quad \cdots \end{aligned} \tag{5.1}$$

These are called the ($p$-adic) Witt polynomials. Now there occurs what has been called 'the miracle of the Witt polynomials'. This is the following integrality statement.

5.2. *Theorem.* Let $\varphi(X,Y,Z)$ be a polynomial over the integers in three (or less, or more) commuting indeterminates. then there are unique polynomials $\varphi_n(X_0, X_1, \cdots, X_n; Y_0, Y_1, \cdots, Y_n; Z_0, Z_1, \cdots, Z_n) = \varphi_n(X;Y;Z)$, $n = 0,1,2,\cdots$ such that for all $n$

$$w_n(\varphi_0(X;Y;Z), \cdots, \varphi_n(X;Y;Z)) = \varphi(w_n(X), w_n(Y), w_n(Z)) \tag{5.3}$$

---

[12] For instance the finite fields $\mathbf{GF}(p^f)$ are perfect. There are also perfect rings that are not fields, for instance, the ring $k[T^{p^{-i}}, i = 0,1,2,\cdots]$ for a perfect field of constants $k$.

The proof is really quite simple. For instance by induction using the following simple observation.

5.4. *Lemma.* Let $\psi(X) = \psi(X_i,\ i \quad I)$ be a polynomial (or power series for that matter) in any set of commuting indeterminates with integer coefficients. Write $\psi(X^p)$ for the polynomial obtained from $\psi(X)$ by replacing each indeterminate by its $p$-th power. Then

$$\psi(X^p) \equiv \psi(X)^p \bmod p \quad \text{and} \quad (\psi(X^p)^{p^j} \equiv \psi(X)^{p^{j+1}} \bmod p^{j+1} \tag{5.5}$$

Note the similarity with (4.1).

5.6. *Proof of theorem* 5.2. Obviously the polynomials $\varphi_n$ are unique and can be recursively calculated from (5.3) over the rationals, starting with $\varphi_0(X_0;Y_0;Z_0) = \varphi(X_0;Y_0;Z_0)$. This also provides the start of the induction. So suppose with induction that the $\varphi_i$ have been proved integral for $i = 0,1,\cdots,n-1$. Now observe that

$$w_n(X) \equiv w_{n-1}(X^p) \bmod p$$

and hence

$$\begin{aligned} \varphi(w_n(X),w_n(Y),w_n(Z)) &\equiv \varphi(w_{n-1}(X^p),w_{n-1}(Y^p),w_{n-1}(Z^p)) \\ &= w_{n-1}(\varphi_0(X^p;Y^p;Z^p),\cdots,\varphi_{n-1}(X^p;Y^p;Z^p)) \end{aligned} \tag{5.7}$$

Using (5.5) one has

$$p^{n-i}\varphi_{n-i}^{p^i}(X;Y;Z) \equiv p^{n-i}\varphi_{n-i}^{p^{i-1}}(X^p;Y^p;Z^p) \bmod p^n$$

and so the $n$ terms of the last expression in (5.7) are term for term equal $\bmod p^n$ to the first $n$ terms of

$$w_n(\varphi_0(X;Y;Z),\cdots,\varphi_n(X;Y;Z)) = \varphi_0^{p^n} + p\varphi_1^{p^{n-1}} + \cdots + p^{n-1}\varphi_{n-1}^p + p^n\varphi_n.$$

Hence $p^n\varphi_n(X;Y;Z) \equiv 0 \bmod p^n$, i.e. $\varphi_n(X;Y;Z)$ is integral.

There are also other proofs; for instance a very elegant one due to Lazard, [258], and reproduced in [367].

5.8. It is obvious from the proof that the theorem holds for polynomials $\varphi$ in any number of variables or power series in any number of variables, but things tend to get a bit messy notationally. It is also obvious that there are all kinds of versions for other rings of coefficients; see also subsection 9.77 - 9.98 below for some more remarks on this theme.

Actually, as will be seen later, see section 9 on the big Witt vectors, theorem 5.2 is not at all necessary for the construction of the various functors of Witt vectors and to work with them; nor, for that matter the Witt polynomials. There are several other ways of doing things.

5.9. *The p-adic Witt addition and multiplication polynomials.* The addition polynomials $s_n(X_0,\cdots,X_n;Y_0,\cdots,Y_n)$ and multiplication polynomials $m_n(X_0,\cdots,X_n;Y_0,\cdots,Y_n)$ of the $p$-adic Witt vectors are now defined by

$$w_n(s_0,s_1,\cdots,s_n)=w_n(X)+w_n(Y) \text{ and } w_n(m_0,m_1,\cdots,m_n)=w_n(X)w_n(Y) \qquad (5.10)$$

5.11. And the functor of the $p$-adic Witt vectors itself is defined as

$$W_{p^\infty}(A)=\{(a_0,a_1,\cdots,a_n,\cdots):a_i\in A\}=A^{\mathbf{N}\setminus\{0\}} \qquad (5.12)$$

as a set and with multiplication and addition defined by

$$\begin{aligned} a+_W b &= (a_0,a_1,\cdots,a_n,\cdots)+_W(b_0,b_1,\cdots,b_n,\cdots)=(s_0(a;b),s_1(a;b),\cdots,s_n(a;b_0),\cdots) \\ a\cdot_W b &= (a_0,a_1,\cdots,a_n,\cdots)\cdot_W(b_0,b_1,\cdots,b_n,\cdots)=(m_0(a;b),m_1(a;b),\cdots,m_n(a;b_0),\cdots) \end{aligned} \qquad (5.13)$$

There is a zero element, viz $(0,0,0,\cdots)$ and a unit element, viz $(1,0,0,\cdots)$ and the claim is:

5.14. *Theorem.* The sets $W_{p^\infty}(A)$ together with the addition and multiplication defined by (5.9) and the unit and zero element as specified define a commutative unital ring valued functor, where the ring morphism corresponding to a ring morphism $\alpha: A\longrightarrow B$ is component wise, i.e. $W_{p^\infty}(\alpha)(a_0,a_1,\cdots)=(\alpha(a_0),\alpha(a_1),\cdots)$. Moreover, the Witt polynomials $w_n$ define functor morphisms $w_n:W_{p^\infty}(A)\longrightarrow A$, $a=(a_0,a_1,a_2,\cdots)\mapsto w_n(a)$. Finally for **Q**-algebras $w$: $W_{p^\infty}(A)\longrightarrow A^{\mathbf{N}\setminus\{0\}}$, $a\mapsto w(a)=(w_0(a),w_1(a),\cdots)$ is an isomorphisms onto the ring $A^{\mathbf{N}\setminus\{0\}}$ with componentwise addition and multiplication.

$W_{p^\infty}$ is also (obviously) the unique functor which setwise looks like $A\mapsto A^{\mathbf{N}\setminus\{0\}}$ and for which the $w_n$ are functorial ring morphisms.

The elements $w_n(a)$ for a Witt vector a are often called the ghost components of that Witt vector (originally: 'Nebenkomponente').

5.15. *Proofs of theorem* 5.14. To prove commutativity, associativity, distributivity and that the zero and unit element have the required properties there are several methods. One is to use defining polynomials as in (5.3). For instance the sequences of polynomials $s(m(X;Y);m(X;Z))$ and $m(X;s(Y;Z))$ both satisfy (5.3) for the polynomial $\varphi(X,Y,Z)=X(Y+Z)=XY+XZ$ and so they are equal proving distributivity on the left.

Another way is to rely on functoriality. First, by the last line in the theorem (which is obvious) distributivity etc. hold for **Q**-algebras. Then, because a commutative integral domain embeds injectively into its ring of quotients, the required properties hold for integral domains. Finally for every unital commutative ring there is an integral domain that surjects onto it and so it follows that the required properties hold for every unital commutative ring.

5.16. *Ghost component equations. Universal calculations* (= '*Calculating with universal polynomials*'). The polynomials $s_n(X;Y)$, $m_n(X;Y)$ are solutions of what I call 'ghost component equations'. For instance the ghost component equations for addition can be written

$$w_n(\text{addition}) = w_n(X) + w_n(Y),\ n = 0,1,2,\cdots \tag{5.17}$$

and call for a sequence of polynomials $s(X;Y) = (s_0(X;Y), s_1(X;Y),\cdots)$ such that $w_n(s) = w_n(X) + w_n(Y)$.

For the unit Wiit vector one has the ghost compnent equations

$$w_n(\text{unitvector}) = 1,\ n = 0,1,2,\cdots \tag{5.18}$$

which call for a series of polynomials $u$ (which turn out to be constants) such that $w_n(u) = 1$ for all $n$. (So that $u_0 = 1$, $u_n = 0$ for $n \geq 1$.)

There will be many ‘ghost component equations’ below. They constitute a most useful and elegant tool, though, as has already been remarked on, one can perfectly well do without them. Not all of these ghost component equations fall within the scope of theorem 5.2, see e.g. subsections 5.25 and 5.27 below.

Here are some more general properties of the $p$-adic Witt vector functor.

5.19. *Ideals and topology.* For each $n = 1,2,3,\cdots$

$$\mathfrak{m}_n = \{(0,\cdots,0,a_n,a_{n+1},a_{n+2},\cdots):\ a_j \quad A\} \tag{5.20}$$

is a (functorial) ideal in $W_{p^\infty}(A)$; further $\mathfrak{m}_i\mathfrak{m}_j \quad \mathfrak{m}_{\max(i,j)}$ (obviously, and no more can be expected in general [13]; for instance the multiplication polynomial $m_1$ has the term $pX_1Y_1$ and so in general $\mathfrak{m}_1\mathfrak{m}_1 \quad \mathfrak{m}_2$). The ring $W_{p^\infty}(A)$ is complete and Hausdorff in the topology defined by these ideals. The quotients $W_{p^n}(A) = W_{p^\infty}(A)/m_{n+1}$ are the rings of $p$-adic Witt vectors of finite length $n+1$.

5.21. *Teichmuller representatives.* For each $x \quad A$ let $t(x) = t_{W_{p^\infty}(A)}(x) = (x,0,0,\cdots)$ This is the Teichmüller representative of $x$ (for the natural projection $w_0\colon\ W_{p^\infty}(A) \longrightarrow A = W_{p^\infty}(A)/\mathfrak{m}_1$). This system of representatives is indeed multiplicative because, as is easily checked, $m_0(X_0;Y_0) = X_0Y_0$ and $m_n(X_0,0,\cdots,0;Y_0,0,\cdots,0) = 0$ for $n \geq 1$.

5.22. *Multiplication with a Teichmüller representative.* It is an easy exercise to check that for the multiplication polynomials

$$m_n(X_0,0,\cdots,0;Y_0,Y_1,\cdots,Y_n) = X_0^{p^n}Y_n \tag{5.23}$$

and so in every $W_{p^\infty}(A)$

[13] But things change very much for the case of $p$-adicWitt vectors over a ring of characteristic $p$; see section 6.

$$(a_0,0,\cdots)\ _W\ (b_0,b_1,\cdots,b_n,\cdots) = (a_0b_0, a_0^p b_1,\cdots,a_0^{p^n} b_n,\cdots) \tag{5.24}$$

5.25. *Verschiebung.* Consider the ghost component equations

$$w_0(\mathbf{V}_p) = 0,\ \ w_n(\mathbf{V}_p) = pw_{n-1} \ \text{ for } \ n \geq 1 \tag{5.26}$$

These call for a series of polynomials $v = (v_0, v_1, v_2, \cdots)$ such that $w_0(v) = 0$, $w_n(v) = pw_{n-1}(v)$ for $n \geq 1$ and do not fall within the scope of theorem 5.2. The immediate and obvious solution is

$$v_0(X) = 0,\ \ v_n(X) = X_{n-1} \ \text{ for } \ n \geq 1$$

and so for each $a = (a_0, a_1, \cdots)$ $W_{p^\infty}(A)$ the operation $\mathbf{V}_p$ acts like $\mathbf{V}_p a = (0, a_0, a_1, a_2, \cdots)$. Then because $w_n(\mathbf{V}_p a) = pw_{n-1}(a)$, $w_n(\mathbf{V}_p(a +_W b)) = w_n(\mathbf{V}_p a) + w_n(\mathbf{V}_p b)$, and it follows that $\mathbf{V}_p$ defines a functorial group endomorphism of the Witt vectors [14]. It does not respect the multiplication. It is called Verschiebung. Note that, see (5.19)

$$\mathfrak{m}_i(A) = \mathbf{V}_p^i(W_{p^\infty}(A))$$

5.27. *Frobenius.* The ghost component equations for the the Frobenius operation are

$$w_n(\mathbf{f}_p) = w_{n+1},\quad n = 0,1,2,\cdots \tag{5.28}$$

This calls for a sequence of polynomials $f = (f_0, f_1, f_2, \cdots)$ such that $w_n(f) = w_{n+1}$, a set of equations that also falls outside the scope of theorem 5.2. For one thing the polynomial $f_n$ involves the indeterminate $X_{n+1}$; for instance

$$f_0 = X_0^p + pX_1,\quad f_1 = X_1^p + pX_2 - \sum_{i=0}^{p-1} p^{p-i-1}\binom{p}{i} X_0^{ip} X_1^{p-i} \equiv X_1^p \bmod p \tag{5.29}$$

It is easy to show that the $f_n$ are integral and that, moreover,

$$f_n(X) \equiv X_n^p \bmod p,\quad n = 0,1,2,3,\cdots \tag{5.30}$$

so that for a ring $k$ of characteristic $p$ the Frobenius operations on the $p$-adic Witt vectors over $k$ are given by

$$\mathbf{f}_p(a_0, a_1, a_2, \cdots) = (a_0^p, a_1^p, a_2^p, \cdots) \tag{5.31}$$

It follows from the ghost component equations (5.28) that the Frobenius operation

---

[14] There are pitfalls in calculating with ghost componets as is done here. Such a calculation gives a valid proof of an identity or something else only if it is a universal calculation; that is, makes no use of any properties beyond those that follow from the axioms for a unital commutative ring only. That is the case here. See also the second proof of theorem 5.14 in 5.15.

$$\mathbf{f}_p(a_0, a_1, a_2, \cdots) = (f_0(a), f_1(a), f_2(a), \cdots) \tag{5.32}$$

is a functorial endomorphism of unital rings of the functor of the $p$-adic Witt vectors.

5.33. *Multiplication by $p$ for the $p$-adic Witt vectors.* The operation of taking $p$-fold sums on the $p$-adic Witt vectors, i.e. $a \mapsto p\ _W a = a +_W a +_W \cdots +_W a$ is (of course) given by the polynomials $P_n = s(s(\cdots(s(X;X);X);\cdots;X)$ which satisfy the ghost component equations

$$w_n(P_0, P_1, P_2, \cdots) = pw_n(X) \tag{5.34}$$

The first two polynomials are $P_0 = pX_0$, $P_1 = X_0^p + pX_1 - p^{p-1}X_0^{p^2}$. With induction one sees

$$P_n(X) \equiv X_{n-1}^p \bmod p \quad \text{for} \quad n \geq 1 \quad (\text{and} \quad P_0(X) \equiv 0 \bmod p) \tag{5.35}$$

For a change, here are the details. According to (5.34) the polynomials $P_n$ are recursively given by the formulas

$$P_0^{p^n} + pP_1^{p^{n-1}} + \cdots + p^{n-1}P_{n-1}^p + p^n P_n = pX_0^{p^n} + p^2 X_1^{p^{n-1}} + \cdots + p^n X_{n-1}^p + p^{n+1} X_n \tag{5.36}$$

With induction, assume $P_i \equiv X_{i-1}^p \bmod p$ for $i \geq 1$. This gives $P_i^{p^{n-i}} \equiv X_{i-1}^{p^{n-i+1}} \bmod p^{n-i+1}$ and $p^i P_i^{p^{n-i}} \equiv p^i X_{i-1}^{p^{n-i+1}} \bmod p^{n+1}$. So the middle $n-1$ terms of the left hand side of (5.36) are term for term congruent to the first $n-1$ terms of the right hand side leaving

$$P_0^{p^n} + p^n P_n \equiv p^n X_{n-1}^p + p^{n+1} X_n \bmod p^{n+1}$$

so that indeed $P_n \equiv X_{n-1}^p \bmod p$ because $P_0 = pX_0$ and $p^n > n$.

5.37. *Taking $p$-th powers in the Wiit vectors.* This operation is governed by the ghost component equations

$$w_n(M_0, M_1, M_2, \cdots) = w_n(X)^p \tag{5.38}$$

These polynomials $M_i$ are of course integral as they can be obtained by repeated substitution of the $p$-adic Witt vector multiplication polynomials into themselves. For the first few polynomials $M_i$ one finds (by direct calculation)

$$M_0(X) = X_0^p,\ M_1(X) \equiv pX_0^{p-1}X_1 \bmod p^2,\ M_2(X) \equiv X_0^{p^2-1}X_1^p \bmod p \tag{5.39}$$

5.40. *Adding ‘disjoint’ Witt vectors.* Finally, suppose for two $p$-adic Witt vectors $a = (a_0, a_1, a_2, \cdots)$ and $b = (b_0, b_1, b_2, \cdots)$ it is the case that for every $n$ at least one of $a_n$ or $b_n$ is zero. Then $(a_i + b_i)^{p^j} = a_i^{p^j} + b_i^{p^j}$ for each $i, j$ and it follows that $w_n(a_0 + b_0, a_1 + b_1, a_2 + b_2, \cdots) = w_n(a) + w_n(b)$ so that in such a case $a +_W b = (a_0 + b_0, a_1 + b_1, a_2 + b_2, \cdots)$. More generally it now follows that each Witt vector $a = (a_0, a_1, a_2, \cdots)$ is equal to the unique convergent sum

$$a = \sum_{i=0}^{\infty} \mathbf{V}_p^i(t(a_i)) \tag{5.41}$$

5.42. *Product formula.* The Frobenius operator, Verschiebung operator, and multiplication for Witt vectors are related in various ways. One of these is the product formula, [192], page 126, formula (17.3.17):

$$\mathbf{V}_p(a\ \ (\mathbf{f}_p b)) = \mathbf{V}_p(a)\ \ b \tag{5.43}$$

This is a kind of formula that shows up in various parts of mathematics, such as when dealing with direct and inverse images of sheaves in algebraic geometry, in (algebraic) $K$-theory and in (abstract representation theory when dealing with restriction and induction of representations [15].

For instance in representation theory such a formula holds with ‘induction’ in place of $\mathbf{V}_p$ and ‘restriction’ in place of $\mathbf{f}_p$. It is important enough in representation theory that it has become an axiom in the part of abstract representation theory known as the theory of Green functors. There it takes the form, [409], page 809:

$$I_K^H(a.R_K^H(b)) = I_K^H(a).b$$

Sometimes this axiom is called the Frobenius axiom. Further in algebraic geometry one has a natural isomorphism $f_*(\mathcal{F}\ \ _{\mathcal{O}_X}\ f^*\mathcal{E})\ \ \ f_*\mathcal{F}\ \ _{O_X}\ \mathcal{E}$ for a morphism of ringed spaces $f\colon (X,\mathcal{O}_X) \longrightarrow (Y,\mathcal{O}_Y)$ and suitable sheaves $\mathcal{E}$ and $\mathcal{F}$ where $f_*$ and $f^*$ stand for taking direct and inverse images under $f$, see [189], exercise 5.1(d) on page 124, while in (étale) cohomology one finds in [295], Ch. VI, §6, proposition 6.5, page 250, a cupproduct formula $i_*(i^*(x)\ \ \ y) = x\ \ \ i_*(y)$. These types of formulae are variously called both product formulae and projection formulae.

It is absolutely not an accident that the same kind of formula turns up in Witt vector theory.

To prove the product formula (5.43) one does a universal calculation. Apply the ghost component morphism $w_n$ to both sides of (5.43) to find respectively (for $n \geq 1$):

$$w_n(\mathbf{V}_p(a\ \ \mathbf{f}_p b) = pw_{n-1}(a\ \ \mathbf{f}_p b) = pw_{n-1}(a)\, w_{n-1}(\mathbf{f}_p b) = pw_{n-1}(a) w_n(b)$$

$$w_n(\mathbf{V}_p a\ \ b) = w_n(\mathbf{V}_p a) w_n(b) = pw_{n-1}(a) w_n(b)$$

and these are equal. As also applying $w_0$ to the two sides gives the same result, the product formula is proved.

5.44. Here is another interrelation between the Frobenius and Verschiebung operators and multiplication:

$$\mathbf{f}_p \mathbf{V}_p = [p] \tag{5.45}$$

---

[15] The terminology ‘product formula’ is not particularly fortunate. There are many things called ‘product formula’ in many part of mathematics. Mostly they refer to formulae that assert that an objct associated to a product is the product of the objects associated to the constituents or to such formulae where a global object for e.g. a global field is the product of corresponding local objects. This happens e.g. for norm residue symbols in class field theory. The terminology ‘projection formula’ (which is also used for formulas like (5.43) is also not particularly fortunate and suffers from similar defects.

where $[p]$ stands for the operation that takes a $p$-adic Witt vector into the $p$-fold sum of itself, $a \mapsto a +_W a +_W \cdots +_W a$.

But in general $\mathbf{V}_p\mathbf{f}_p$ is not equal to $[p]$.

To prove (5.45) one again does a universal calculation: $w_n(\mathbf{f}_p\mathbf{V}_p a) = w_{n+1}(\mathbf{V}_p a) = pw_n(a)$ and also $w_n([p]a) = pw_n(a)$.

## 6. The ring of *p*-adic Witt vectors over a perfect ring of characteristic *p*

In this section $k$ will finally be a perfect ring of characteristic $p > 0$. In this case the ring of $p$-adic Witt vectors over $k$ has a number of additional nice properties including a first nice universality property. Recall that '$k$ is perfect' means that the Frobenius ring morphism $\sigma\colon k \longrightarrow k,\ x \mapsto x^p$ is bijective. One of the aims of this section to show that if $k$ is a perfect field of characteristic $p$ then the $p$-adic Witt vectors form a characteristic zero complete discrete valuation ring with residue field $k$.

6.1. *p-adic Witt vectors over a ring of characteristic p.* First, just assume only that $k$ is of characteristic $p$. Then calculations with $p$-adic Witt vectors simplify quite a bit. For instance:

$$
\begin{aligned}
&\mathbf{f}_p(a_0, a_1.a_2,\cdots) = (a_0^p, a_1^p, a_2^p,\cdots)\\
&[p](a_0,a_1,a_2,\cdots) = p\ (a_0,a_1,a_2,\cdots) = (0,a_0^p,a_1^p,a_2^p,\cdots)\\
&p = (0,1,0,0,0,\cdots)\\
&\mathbf{f}_p\mathbf{V}_p = [p] = \mathbf{V}_p\mathbf{f}_p\\
&\mathbf{V}_p^i a\ \mathbf{V}_p^j b = \mathbf{V}_p^{i+j}(a_0^{p^j}b_0^{p^i},??,??,\cdots)
\end{aligned}
\tag{6.2}
$$

where the ??'s in the last line stand for some not specified polynomial expressions in the the coordinates (components) of $a$ and $b$. The first four of these formulas follow directly from (5.25), (5.30) and (5.35). The last one follows from the fourth one, repeated application of the product formula, and the first formula of (6.2) (and commutativity of the multiplication). Note that the last formula implies that

$$\mathfrak{m}_i\mathfrak{m}_j \subset \mathfrak{m}_{i+j} \tag{6.3}$$

(but equality need not hold). Also there is the corrollary

6.4. *Corollary.* If $k$ is an integral domain of characteristic $p$ then the ring of $p$-adic Witt vectors over $k$ is an integral domain of characteristic zero.

6.5. *Valuation rings.* A (normalized) discrete valuation[16] on a (commutative) field $K$ is a surjective function $v\colon K \longrightarrow \mathbf{Z} \cup \{\infty\}$ such that

$$\begin{aligned} &v(x) = \infty \text{ if and only if } x = 0, \\ &v(xy) = v(x) + v(y), \\ &v(x+y) \geq \min\{v(x), v(y)\} \end{aligned} \tag{6.6}$$

for all $x, y \in K$. Such a valuation can be used to define a norm, and hence a topology, on $K$ by setting $\| x \| = r^{-v(x)}$ where $r$ is e.g. an integer $> 1$. This defines an ultrametric or non-Archimedean metric, $d(x,y) = \| x - y \|$, on $K$ where 'ultra' means that instead of the familiar triangle inequality one has the stronger statement $d(x,y) \leq \max\{d(x,z), d(y,z)\}$.

---

[16] Such valuations as here are also sometimes called 'exponential valuations'. More generally one considers such functions with values in the real numbers (or more general ordered groups) with infinity adjoined, cf e.g. [136], page 20. Then the value group $v(K\setminus\{0\})$ is a subgroup of the additive group of real numbers. Such a group is either discrete (when there is a smallest positive real number in it) or dense. In the discrete case one may as well assume that the value subgroup is the group of integers. Whence the terminology 'normalized discrete'.

The valuation ring of the valued field $(K,v)$ is the ring

$$A = \{x \in K\colon\ v(x) \geq 0\} \tag{6.7}$$

Note that such a ring automatically has two properties

$A$ is local with maximal ideal $\mathfrak{m} = \{x \in K\colon\ v(x) \geq 1\}$ (6.8)

The maximal ideal $\mathfrak{m}$ is principal (6.9)

One can also start with a valuation on an integral domain $A$, which now is a surjective function $v\colon A \longrightarrow \mathbf{N} \cup \{0,\infty\}$ such that the properties (6.6) hold. A valuation on the field of fractions $Q(A)$ is then defined by $v(xy^{-1}) = v(x) - v(y)$ (which is independant of how an element of the field of fractions is written as a fraction). Such an integral domain with a valuation is not necessarily a valuation ring. In particular, one or both of the properties (6.8), (6.9) may fail.

6.10. *Example. The p-adic valuation on the integers*. Take a prime number $p$. Define on the integers the function $v_p(n) = r$ if and only if $p^r$ is the largest power of the prime number $p$ that divides $n$. This is obviously a valuation but $\mathbf{Z}$ is not its corresponding valuation ring. That valuation ring is the localization $\mathbf{Z}_{(p)}$ consisting of all rational numbers that can be written as fractions with a denominator that is prime to $p$.

6.11. However, if $A$ is an integral domain with a valuation and the ring $A$ satisfies properties (6.8), (6.9), then $A$ is a discrete valuation ring. See e.g. [136], page 50, theorem 7.7; [62], Ch. 5, §3, no 6, proposition 9, page 100.

After this intermezzo on discrete valuation ringss let's return to rings of $p$-adic Witt vectors.

6.12. *Units in the ring of p-adic Witt vectors.* Let $k$ be any ring of characteristic $p$, then a $p$-adic Witt vector $a = (a_0, a_1, a_2, \cdots)$ is invertible in $W_{p^\infty}(A)$ if and only if $a_0$ is invertible in $k$.

Indeed, multiplying with $(a_0^{-1}, 0, 0, 0, \cdots)$ it can be assumed by (5.24) that $a_0 = 1$ (multiplication with a Teichmüller representative). Then $a = 1 + \mathbf{V}_p b$ with $b = (a_1, a_2, \cdots)$. The series $1 - \mathbf{V}_p b + (\mathbf{V}_p b)^2 - (\mathbf{V}_p b)^3 + \cdots$ converges in $W_{p^\infty}(k)$ because of (6.2) or (6.3), say, to an element c. Then $a \cdot c = 1$.

It follows in particular that if $k$ is a field of characteristic $p$ then $W_{p^\infty}(k)$ is a local ring with maximal ideal $\mathfrak{m} = \mathfrak{m}_1 = \mathbf{V}_p W_{p^\infty}(k)$.

6.13. *Rings of p-adic Witt vectors over a perfect ring of characteristic p.* Now let $k$ be a ring of characteristic $p$ that is perfect. There are of course such rings that are not fields, for instance a ring $k'[T_1,T_2,\cdots]$ with the relations $T_i^p = T_{i-1}$ for $i \ge 2$ and $k'$ a perfect field. This field can be suggestively written $k'[T,T^{p^{-1}},T^{p^{-2}},\cdots]$.

Now first note that if $k$ is perfect of characteristic $p$ then the ideal $\mathfrak{m} = \mathfrak{m}_1 = \mathbf{V}_p W_{p^\infty}(k)$ is principal and generated by the element $p \in W_{p^\infty}(k)$. This is immediate from the second formula in (6.2) above. Second, observe that in this case $\mathfrak{m}_r = p^r W_{p^\infty}(k) = (\mathfrak{m}_1)^r = \mathbf{V}_p^r W_{p^\infty}(k)$.

In this setting there is the first nice universality property of the *p*-adic Witt vectors.

6.14 *Theorem.* Let $k$ be a ring of characteristic $p$ that is perfect. Let $A$ be a ring with an ideal $\mathfrak{a}$ such that $A/\mathfrak{a} = k$ and such that $A$ is complete and separated in the topology defined by $\mathfrak{a}$. Then there is a unique morphism of rings $\varphi$ such that the following diagram commutes

$$\begin{array}{ccc} W_{p^\infty}(k) & \xrightarrow{\varphi} & A \\ \downarrow w_0 & & \downarrow q \\ k & = & k \end{array}$$

Here the topology defined by $\mathfrak{a}$ is the topology defined by the sequence of ideals $\mathfrak{a}^n$, $n = 1,2,\cdots$ and 'separated' means that $\bigcap_n \mathfrak{a}^n = \{0\}$. Further $q$ is the canonical quotient mapping $A \longrightarrow A/\mathfrak{a} = k$. It is perfectly OK for $A$ to be a ring of characteristic $p$ or for it to be such that $\mathfrak{a}^r = 0$ for some $r$.

6.15. *Proof of theorem* 6.14. First, let's prove uniqueness of $\varphi$ (if it exists at all). To this end observe that by the remarks in (6.13) and (5.41) every element in $W_{p^\infty}(k)$ can be uniqely written in the form of a convergent sum

$$\sum_{n=0}^{\infty} t(a_n^{p^{-n}})p^n, \quad a_n \in k$$

and thus (by the completeness and separatedness of $A$) a ring morphism $\varphi\colon W_{p^\infty}(k) \longrightarrow A$ is uniquely determined by what it does to the Teichmüller representatives. But , as $\varphi$ is a ring morphism and hence multiplicative, the $\varphi(t(a))$, $a \in k$ form a system of Teichmüller representatives for $k$ in $A$. As such systems are unique it follows that there can be at most one $\varphi$ that does the job.

For existence consider the ghost component ring morphism

$$w_n\colon\ W_{p^n}(A) = W_{p^\infty}(A)/\mathfrak{m}_{n+1}(A) \longrightarrow A$$

If $r_0,r_1,\cdots,r_n \in \mathfrak{a}$ then $w_n(r_1,r_2,\cdots,r_n) \in \mathfrak{a}^{n+1}$ and so there is an induced ring morphism $\psi_n$ that makes the following diagram commutative

$$\begin{array}{ccc} W_{p^n}(A) & \xrightarrow{w_n} & A \\ \downarrow W_{p^n}(q) & & \downarrow \\ W_{p^n}(k) & \xrightarrow{\psi_n} & A/\mathfrak{a}^{n+1} \end{array}$$

Now define

$$\varphi_n\colon\ W_{p^n}(k) \longrightarrow A/\mathfrak{a}^{n+1}$$

as the composite

$$W_{p^n}(k) \xrightarrow{W_{p^n}(\sigma^{-n})} W_{p^n}(k) \xrightarrow{\psi_n} A/\mathfrak{a}^{n+1}$$

Then for any $a_0, a_1, \cdots, a_n \in A$ it follows that

$$\varphi_n(q(a_0^{p^n}, q(a_1^{p^n}), \cdots, q(a_n^{p^n})) = a_0^{p^n} + pa_1^{p^{n-1}} + \cdots + p^n a_n \tag{6.16}$$

Now let $x_0, x_1, \cdots, x_n \in k$ and substitute in this formula (6.16) the Teichmüller representatives (in $A$) $t_A(x_i^{p^{-n}})$ to find that

$$\varphi_n(x_0, x_1, \cdots, x_n) = t_A(x_0) + pt_A(x_1^{p^{-1}}) + p^2 t_A(x_2^{p^{-2}}) + \cdots + p^n t_A(x_n^{p^{-n}}) \bmod \mathfrak{a}^{n+1}$$

and so the projective limit $\varphi$ of the $\varphi_n$ exists, is a ring morphism, and is given by the formula

$$\varphi(x_0, x_1, x_2, \cdots) = t_A(x_0) + pt_A(x_1^{p^{-1}}) + p^2 t_A(x_2^{p^{-2}}) + \cdots \tag{6.17}$$

Of course one can write down this formula directly but then one is faced with proving that it is a ring morphism which comes down to much the same calculations as were carried out just now.

Now suppose in addition that $\mathfrak{a} = (p)$ and that $A$ is of characteristic zero. Then the ring morphism

$$W_{p^\infty}(k) \xrightarrow{\varphi} A,\ (x_0, x_1, x_2, \cdots) \mapsto x_0 + pt_A(x_1)^p + p^2 t_A(x_2)^{p^2} + \cdots \tag{6.18}$$

is clearly bijective and hence an isomorphism. So there is a (strong) uniqueness result as follows

6.19. *Theorem.* Let $k$ be a perfect ring of characteristic $p > 0$. Then, up to isomorphism, there is precisely one ring $A$ of characteristic zero with residue ring $A/(p) = k$ and complete in the $p$-adic topology. Moreover $A$ is rigid in the sense that the only ring automorphism of $A$ that induces the identity on $k$ is the identity on $A$.

6.20. *The p-adic Witt vectors over a perfect field of characteristic* $p > 0$. Finallly let $k$ be a perfect field of characteristic $p$. Then by 6.13, $\mathfrak{m} = \mathfrak{m}_1 = \mathbf{V}_p W_{p^\infty}(k)$ is generated by $p$. Define $v$ as follows

$$\begin{aligned} &v\colon\ W_{p^\infty}(k) \longrightarrow \mathbf{N} \cup \{0,\infty\}, \\ &v(0) = \infty,\ v(x_0, x_1, x_2, \cdots) = n \text{ iff } x_n \text{ is the first coordinate unequal to zero} \end{aligned} \tag{6.21}$$

Then $v$ is a valuation on the ring of $p$-adic Witt vectors and makes this ring into a complete discrete valuation ring by (6.11) and (6.12). And this ring is unramified as the maximal ideal is generated by $p$.

Theorem 6.19 now translates into an existence and uniqueness (and rigidity) theorem for complete discrete unramified valuation rings of characteristic zero and residue characteristic $p > 0$.

Thus, for a perfect residue field there are indeed universal formulae that govern the addition and multiplication of $p$-adic expansions as in section 3 above. The treatment here has a bit of a 'deus ex machina' flavour in that the $p$-adic Witt vectors are constructed first and subsequently proved to do the job. Motivationally speaking one can do better. Once it is accepted that it is a good idea to write $p$-adic expansions using Teichmüller representatives one can calculate and come up with addition and multiplication formulae. This is nicely done in [188] . It was also already known in 1936, [382], that there should be integral coefficient formulae that could do the job [17]. Witt's major contribution was finding a way to describe them nicely and recursively.

6.22. *Complete discrete valuation rings with non perfect residue field.* Immediately after Witt's paper, in the same [18] issue of 'Crelle', there is a paper by Teichmüller, [380], in which he proves existence and uniqueness of unramified complete discrete valuation rings with given residue field, thus completing his own arguments from [382]. For the unequal characteristic case he uses the $p$-adic Witt vectors.

Existence and uniqueness of unramifled complete discrete valuation rings had been treated before by Friedrich Carl Schmidt and Helmut Hasse in [191], but that paper had an error, which was later, in 1940, pointed out by Saunders MacLane, [277], and still later corrected by him, [278]. For more see [339].

6.23. *Cohen rings*. The unique characteristic zero complete unramified discrete valuation ring with a given residue field $k$ of characteristic $p > 0$ is nowadays often called the Cohen ring [19] of $k$, see [63], Ch. IX, §2 for a treatment of Cohen rings [20]. If $k$ is perfect the Cohen ring of $k$ is the ring of $p$-adic Witt vectors. If $k$ is not perfect it is a 'much smaller' ring.

The technical definition of Cohen ring in loc.cit. is as follows. A $p$-ring is (by definition) a unital commutative ring such that the ideal $pC$ is maximal and the ring is complete and separate in the $p$-adic topology defined by the powers of this ideal. Given a local separated complete ring $A$ a Cohen subring of it is a $p$-ring that satisfies $A = \mathfrak{m}_A + C$ (where $\mathfrak{m}_A$ is the maximal ideal of $A$.

[17] Indeed in [382], the theorem at the bottom of page 58 says the following. Define polynomials $h_0, h_1, h_2, \cdots$ in two variables by the formulae $x^{p^n} + y^{p^n} = h_0^{p^n} + ph_1^{p^{n-1}} + \cdots + p^{n-1}h_{n-1}^p + p^n h_n$. Then these have integer coefficients. Given $a, b$ in the perfect residue field $k$, let $c_n$, $n = 0,1,2,\cdots$ be given by $c_n^{p^n} = h_n(a,b)$. Then the $p$-adic expansion of the sum to the Teichmüller representatives of $a$ and $b$ is $t(a) + t(b) = t(c_0) + pt(c_1) + p^2 t(c_2) + \cdots$ which of course nicely agrees with the addition of $p$-adic Witt vectors via the isomorphism (6.17). Thus the addition formula for Teichmüller representatives was known, and given by universal formulae, and involved what were to be called Witt polynomials. The multiplication of Teichmüller representaives was of course also known (as these are multiplicative). Distributivity then determines things recursively. There was still much to be done and no predicting that things would come out as nicely as they did, but the seeds were there.

[18] The two papers carry the same date of receipt.

[19] There are other rings that are called 'Cohen rings' in the published literature. First there are certain twisted power series rings which are a kind of noncommutative complete discrete valuation rings, see e.g. [131, 342, 406]. Further in [137], section 18.34B, p. 48, a Cohen ring is defined as a ring $R$ such that $R/P$ is right Artinian for each nonzero prime ideal $P$; see also [138, 413]. This last bit of terminology has to do with the paper [93].

[20] This terminology would appear to come from the so-called Cohen structure theorems for complete local rings, [94].

The existence and uniqueness of such a subring of the $p$-adic Witt vectors is implicit in the work of Teichmüller, [380], and Nagata, [304], Chapter V, §31, p 106ff, but a nice functorial description had to wait till the work of Colette Schoeller, [358] [21].

6.24. *p-basis.* Let $k$ be a field of characteristic $p > 0$. The field is perfect if and only if $k = k^p = \{x^p : x \in k\}$. If it is not perfect there exists a set of elements $\mathcal{B} = \{b_i : i \in I\}$ of $k$ such that the monomials

$$\prod_{i \in I} b_i^{j_i}, \quad j_i \in \{0,1,\cdots,p-1\}, \ j_i \text{ equal to zero for all but finitely many } i$$

form a basis (as a vector space) for $k$ over the subfield $k^p$. The cardinality of $I$, which can be anything, is an invariant of $k$ and sometimes called the imperfection degree of $k$. Such a set is called a $p$-basis. The notion is again due to Teichmüller and first appeared in [381] (for other purposes than valuation theory). For some theory of $p$-bases see [64], Chap. 5, §8, exercise 1, and [304], p. 107ff. If $B$ is a $p$-basis for $k$ the monomials

$$\prod_{i \in I} b_i^{j_i}, \quad j_i \in \{0,1,\cdots,p^n-1\}, \ j_i \text{ equal to zero for all but finitely many } i \tag{6.25}$$

form a vector space basis for $k$ over the subfield $k^{p^n}$.

6.25. *Cohen functor.* Given a characteristic $p > 0$ field $k$ and a chosen $p$-basis let $\mathbf{Z}[\mathcal{B}]$ be the ring of polynomials over the integers in the symbols from $\mathcal{B}$. Then Schoeller defines functors $C_n$ on the category of commutative unital $\mathbf{Z}[\mathcal{B}]$-algebras to the category of unital commutative rings. These functors depend on the choice of the $p$-basis. They define affine group schemes.

Now if the characteristic of the base field $k$ is perfect the finite length $p$-adic Witt vector group schemes serve to classify unipotent affine algebraic groups over $k$, a topic which will get some attention in a later section of this chapter. For a non-perfect base field the Cohen functors defined in [358], see also [248] are a well-working substitute [22].

6.26. *Cohen ring of k.* It follows from (6.25) that $\mathbf{F}_p[\mathcal{B}]$, where $\mathbf{F}_p$ is the field of $p$ elements, naurally embeds in $k$; combined with the canonical projection $\mathbf{Z}[\mathcal{B}] \longrightarrow \mathbf{F}_p[\mathcal{B}]$ this defines a $\mathbf{Z}[\mathcal{B}]$-algebra structure on $k$. The value of the Cohen functors on $k$ can now be described as follows, see [358], section 3.2, p. 260ff; see also [63], Chap. IX, §2, exercise 10, p. 72.

$$\begin{aligned}&C_n(k) \text{ is the subring of } W_{n+1}(k) \text{ generated by } W_{n+1}(k^{p^n}) \subset W_{n+1}(k)\\ &\text{and the Teichmüller representatives } (b,0,0,\cdots,0) \in W_{n+1}(k),\ b \in \mathcal{B}\end{aligned}$$

---

[21] Existence of a multiplicative system of representatives is again crucial. In the general case, i.e. when the residue field $k$ is not necessarily perfect there is still existence (as noticed by Teichmüller, [382]), but no uniqueness. A very nice way of seeing existence is due to Kaplansky , [235], section 26, p. 84ff, who notes that the group of 1-units (Einseinheiten) of $A$, i.e. the invertible elements that are mapped to 1 under the residue mapping, form a direct summand of the group of all invertible elements of $A$. This remark seems to have been of some importance for Schoeller, and reference 4 in her paper, a paper that never appeared, should probably be taken as referring to this.

[22] There may be some difficulties with the constructions in [358] in case the cardinality of the $p$-basis is not finite as is stated in [378], where also an alternative is proposed.

These form a projective system (with surjective morphisms of rings $\mathcal{C}_{n+1}(k) \longrightarrow \mathcal{C}_n(k)$ induced by the projections $W_{n+1}(k) \longrightarrow W_n(k)$) and taking the projective limit one finds a subring $\mathcal{C}(k)$ of the ring of $p$-adic Witt vectors which is a discrete complete unramified valuation ring. It is unique but one loses the rigidity property that holds for the Witt vectors, see theorem 6.19 and 6.20.

## 7. Cyclic Galois extension of degree $p^n$ over a field of characteristic $p$.

As already indicated in the introduction (section 1) for Witt himself one of the most important aspects of the Witt vectors was that they could be used to extend and complete his own results from [416, 415] to obtain a Kummer theory (class field theory) for Abelian extensions of a field of characteristic $p$. This section is a brief outline of this theory. In this whole section $p$ is a fixed prime number.

7.1. *Construction of some Abelian extensions of a field of characteristic* $p > 0$. Consider the functor of $p$-adic Witt vectors $W_{p^\infty}$. It is an easy observation that for each $n \geq 1$, the multiplication and addition polynomials $s_i(X;Y),\ m_i(X;Y),\ i = 0,\cdots,n-1$ only depend on $X_0,\cdots,X_{n-1};Y_0,\cdots,Y_{n-1}$. It follows immediately that the sets

$$\mathbf{V}_p^n W_{p^\infty}(A) = \{(\underbrace{0,\cdots,0}_{n}, a_0, a_1, \cdots, a_m, \cdots)\} \tag{7.2}$$

are (functorial) ideals in the rings $W_{p^\infty}(A)$. And, hence, for each $n$, there is a quotient functor

$$W_{p^{n-1}}(A) = W_{p^\infty}(A) / \mathbf{V}_p^n W_{p^\infty}(A) = \{(a_0, a_1, \cdots, a_{n-1}):\ a_i \quad A\} \tag{7.3}$$

These are called the $p$-adic Witt vectors of length n. The multiplication and addition of these vectors of length $n$ are given by the multipllcation and addition polynomials $s_i(X;Y),\ m_i(X;Y),\ i = 0.\cdots,n-1$.

In other words $W_{p^{n-1}}$ is the functor on **CRing** to **CRing** represented by $\mathbf{Z}[X_0, X_1, \cdots, X_{n-1}]$ provided with the coring object structure given by the polynomials $s_i(X;Y),\ m_i(X;Y),\ i = 0.\cdots,n-1$.

If the ring $A$ is of characteristic $p$ the Frobenius endomorphism on $W_{p^\infty}(A)$ is given by (see (6.2) above)

$$\mathbf{f}_p(a_0, a_1, a_2, \cdots) = (a_0^p, a_1^p, a_2^p, \cdots) \tag{7.4}$$

and thus manifestly takes the ideal $\mathbf{V}_p^n W_{p^\infty}(A)$ into itself. In general this is not the ccase. Thus for rings of characteristic $p$ there is an induced (functorial) Frobenius endomorphism of rings

$$\mathbf{f}_p:\ W_{p^{n-1}}(A) \longrightarrow W_{p^{n-1}}(A) \tag{7.5}$$

From now on in this section $k$ is a field of characteristic $p > 0$. Consider the additive operator

$$\wp:\ W_{p^{n-1}}(k) \longrightarrow W_{p^{n-1}}(k),\quad \wp(\alpha) = \mathbf{f}_p(\alpha) - \alpha \tag{7.6}$$

(Witt vector subtraction of course.) For $n=1$, so that $W_{p^{n-1}}(k)=k$ this is the well-known Artin - Schreier operator $x^p - x$ which governs the theory of cyclic extensions of degree $p$ of the field $k$.

For $n>1$ the operator (7.6) is very similar, particularly if one reflects that the Frobenius operator is 'as much like raising to the power $p$ as an additive operator on the Witt vectors can be' (as has been remarked before).

Moreover, consider the kernel of the operator $\wp$. For a given $n$ his kernel consists of all $p$-adic Witt vectors $\alpha=(a_0,a_1,\cdots,a_{n-1})$ of length $n$ such that $a_i^p=a_i$ for all $i=0,1,\cdots,n-1$. Thus for all these $i$ $a_i$ $\mathbf{F}_p$ $k$ the prime subfield of $p$ elements of $k$. Thus

$$\text{Ker}(\wp)=W_{p^{n-1}}(\mathbf{F}_p)=\mathbf{Z}/p^n\mathbf{Z} \tag{7.7}$$

the cyclic group of order $p^n$. (The last equality in (7.7) is immediate from the considerations of the previous chapter.) This is most encouraging in that it suggests that the study of equations

$$\mathbf{f}_p(\alpha)-\alpha=\wp(\alpha)=\beta \tag{7.8}$$

in $W_{p^{n-1}}(\bar{k})$ with $\beta$ $W_{p^{n-1}}(k)$ and $\bar{k}$ an algebraic closure of $k$ could lead to cyclic extensions of degree $p^n$ just like in standard Artin-Schreier theory as described e.g. in [96], pp 205-206. This is indeed substantially the case and things turn out even better.

For any $\beta$ $W_{p^{n-1}}(k)$ consider a solution $\alpha=(a_0,a_1,\cdots a_{n-1})$ of (7.8). Then $k(\wp^{-1}\beta)$ denotes the extension field of $k$ generated by the components $a_0,a_1,\cdots a_{n-1}$. Because of (7.7) if $\alpha'=(a_0',a_1',\cdots a_{n-1}')$ is another solution $k(a_0,a_1,\cdots,a_{n-1})=k(a_0',a_1',\cdots,a_{n-1}')$ so that the extension $k(\wp^{-1}\beta)$ of $k$ does not depend on the choice of a solution. More generally if $\Omega$ is any subset of $W_{p^{n-1}}(k)$ then $k(\wp^{-1}\Omega)$ denotes the union of all the $k(\wp^{-1}\beta)$ for $\beta$ $\Omega$.

7.9. *The main theorems for Abelian extension of exponent $p^n$ of fields of characteristic $p$.* There is now sufficient notation to formulate the main theorems of Kummer theory for Abelian extensions of exponent $p^n$ of a field of characteristic $p$.

7.10. *Theorem* (*Kummer theory*). Let $\Omega$ $W_{p^{n-1}}(k)$ be a subgroup which contains $\wp W_{p^{n-1}}(k)$ and such that $\Omega/\wp W_{p^{n-1}}(k)$ is finite. Then the Galois group of the extension $k(\wp^{-1}\Omega)$ is isomorphic to $\Omega/\wp W_{p^{n-1}}(k)$

For each Abelian field extension $K/k$ of exponent a divisor of $p^n$ there is precisely one group $\Omega/\wp W_{p^{n-1}}(k)$ such that $K=k(\wp^{-1}\Omega)$.

7.11. *Theorem* (*on cyclic extensions of degree $p^n$*). Let $\beta=(b_0,b_1,\cdots,b_{n-1})$ be a $p$-adic Witt vector of length $n$ over $k$ such that $b_0$ $\wp W_{p^{n-1}}(k)$. Then $k(\wp^{-1}\beta)$ is a cyclic extension of degree $p^n$ of $k$.

A generator of the Galois group is given by $\sigma(\alpha)=\alpha+\mathbf{1}$ (Witt vector addition, and $\mathbf{1}$ the unit of the ring $W_{p^{n-1}}(k)$).

All cyclic extensions of degree $p^n$ can be obtained in this way.

7.12. *On the proofs*. The theorems above already occur in [420]. But the proofs there are very terse. They consists of brief instructions to the reader to first prove a kind of (additve) "Hilbert 90 theorem" for Witt vectors by redoing the proof of theorem I.2 in [415]. This says that a first Galois cohomology group with coefficients in Witt vectors is zero and is 'Satz 11' in [420]. The further instructions are to redo the arguments of [415] using vectors instead of numbers and using 'Satz 11' instead of the usual 'Hilbert 90' as the occasions demand.

For a good complete treatment of this Kummer theory for Abelian extensions of fields of characteristic $p > 0$ see [271], pp 146-150. The statements there are slightly more general and a bit more elegant than in [420] in that the group $\Omega / pW_{p^{n-1}}(k)$ is not required to be finite. The isomorphism statement of theorem 7.10 now becomes a duality statement to the effect that the group $\Omega / pW_{p^{n-1}}(k)$ and the Galois group of the extension are dual to each other under a natural nondegenerate pairing.

## 8. Cyclic central simple algebras of degree $p^n$ over a field of characteristic $p$.

The second main application of the $p$-adic Witt vectors in [420] is to cyclic central simple algebras of prime power degree $p^n$ over a field of charateristic $p > 0$. This is a topic in the theory of simple central algebras over and the Brauer group of a field, [1,2], [96] Chapter 7, [102] Chapter III.

8.1. *Central simple algebras*. A central simple algebra over a field $k$ is what the name indicates: it is a finite dimensional agebra over $k$, it is simple, i.e. no nontrivial ideals, and it is central, i.e. its centre, $\{a \in A: ar = ra \text{ for all } r \in A\}$ coincides with $k$. A central divison algebra over $k$ is a central simple algebra in which every nonzero element is invertible. The classic example is the algebra of quaternions over the reals.

For every central simple algebra $A$ there is a (unique up to isomorphism) central division algebra $D$ such that $A$ is isomorphic to a full matrix ring over $D$.

For every central simple algebra over $k$ there is a field extension[23] $K / k$ such the tensor product $A_K = A \otimes_k K$ is a full matrix ring over $K$. Such a field is called a splitting field.

The tensor product (over $k$) of two central simple algebras over $k$ is again a central simple algebra over $k$.

8.2. *Brauer group.* Two central simple algebras $A$, $B$, over $k$ are called (Brauer) equivalent if for suitable natural numbers $m$ and $n$

$$A \otimes_k M_m(k) \cong B \otimes_k M_n(k) \tag{8.3}$$

Here $M_m(k)$ is the full matrix algebra over $k$ of all $m \times m$ matrices with entries in $k$.

Equivalently $A$ and $B$ are equivalent if they have isomorphic associated division algebras.

The tensor product is compatible with this equivalence notion and defines (hence) a group structure on the equivalence classes. This group is called the Brauer group, $\mathrm{Br}(k)$, of $k$. It can be interpreted as a second Galois cohomology group.

[23] As always the term 'field' implies commutativity.

8.4. *Crossed product central simple algebras*, [102] § III-2, [96] § 7.5.

A central simple algebra over $k$ is called a cross product if it contains a maximal subfield $K$ such that $K/k$ is a Galois extension.

That $K$ is then a splitting field.

This is the 'abstract' definition of a crossed product. There is also an explicit description/construction which is important and explains the terminology. This goes as follows. Let $U$ be the group of invertible elements of a crossed product central simple algebra $A$. Let $N$ be the centralizer of $K^\times = K\setminus\{0\}$ in $U$:

$$N = \{u \in U:\ u^{-1}Ku \subset K\} \tag{8.5}$$

Then $K^\times$ is a normal subgroup of $N$ and hence gives rise to a short exact sequence of groups

$$1 \longrightarrow K^\times \longrightarrow N \longrightarrow \Gamma \longrightarrow 1 \tag{8.6}$$

a group extension of $\Gamma$ by $K^\times$. One easily shows that $\Gamma \cong \mathrm{Gal}(K/k)$.

Being a group externsion, it is determined by what is called a factor system and this same factor system can be used to recover (up to isomorphism) the crossed product central simple algebra $A$.

Not every central simple algebra is a crossed product, but every Brauer equivalence class contains one.

When the Galois extension involved is cyclic, i.e. the group $\Gamma = \mathrm{Gal}(K/k)$ is cyclic, one speaks of a cyclic central simple algebra.

8.7. *Cyclic central algebra associated to a finite length p-adic Witt vector*. From now on in this section $k$ is a fixed field of characteristic $p>0$. Take an element $\alpha \in k$ and a finite length Witt vector $\beta = (b_0, b_1, \cdots, b_{n-1}) \in W_{p^{n-1}}(k)$. To these data associate the algebra generated over $k$ by an indeterminates $u$ and commuting indeterminates $\theta_0, \theta_1, \cdots, \theta_{n-1}$ subject to the relations

$$u^{p^n} = \alpha,\quad \wp\theta = \beta,\quad u\theta u^{-1} = \theta + \mathbf{1} \tag{8.8}$$

Here, $\theta = (\theta_0, \theta_1, \cdots, \theta_{n-1})$ (as a Witt vector) and $u\theta u^{-1} = (u\theta_0 u^{-1}, u\theta_1 u^{-1}, \cdots, u\theta_{n-1}u^{-1})$. This algebra will be denoted $(\alpha|\beta]$. These algebras are central simple algebras of 'degree' $p^n$ (meaning that their dimension over $k$ is $p^{2n}$).

Note that conjugation by $u$ is an operation of order $p^n$. Also, by the results of section 7 above, if $b_0 \notin \wp W_{p^{n-1}}(k)$ the subalgebra $k(\theta_0, \theta_1, \cdots, \theta_{n-1})$ is a cyclic field extension of degree $p^n$ and conjugation by $u$ is the action of a generator of the Galois group. This is just about sufficient to prove that $(\alpha|\beta]$ is a cyclic central simple (crossed product type) algebra. In case $b_0 \in \wp W_{p^{n-1}}(k)$ there is a 'reduction theorem', Satz 15 from [420], which says the following.

8.9. *Theorem*. The algebra $(\alpha|0, b_1, b_2, \cdots, b_{n-1}]$ is Brauer equivalent to $(\alpha|b_1, b_2, \cdots, b_{n-1}]$

The next theorem is rather remarkable and shows once more that there is no escaping the Witt vectors.

8.10. *Theorem* (*Brauer group and Witt vectors*). In the Brauer group of the field $k$ there are the (calculating) rules

$$\begin{aligned} (\alpha|\beta]\,(\alpha'|\beta] &= (\alpha\alpha'|\beta] \\ (\alpha|\beta]\,(\alpha|\beta'] &= (\alpha|\beta+\beta'] \end{aligned} \tag{8.11}$$

(where in the second line of (8.11) the '+' sign means Witt vector addition.

For the time being the sections above conclude the discussions on $p$-adic Witt vectors. There will be more about the $p$-adic Witt vectors in various sections below. But first it is time to say something about that truly universal object the functor of the big Witt vectors; and that will be the subject in the next few sections.

## 9. The functor of the big Witt vectors

In the early 1960's, probably independent of each other, several people, notably Ernst Witt himself and Pierre Cartier, noticed that the $p$-adic Witt polynomials (5.1) are but part of a more general family and that these polynomials can be used in a similar manner to define a ring valued functor $W$: $\mathbf{CRing} \longrightarrow \mathbf{CRing}$ of which the $p$-adic Witt vectors are a quotient, see [105, 149, 192, 255, 419]. This functor $W$ is called the functor of the big Witt vectors [24]. (And also, over $\mathbf{Z}_{(p)}$-algebras the canonical projection $W(A) \longrightarrow W_{p^\infty}(A)$ has an Abelian group section, making in this case the $p$-adic Witt vectors also a sub functor of the big Witt vectors.)

Here, I will first construct the big Witt vectors in another way, before describing those 'generalized' Witt polynomials alluded to, and how the $p$-adic Witt vectors fit with the big ones.

9.1. *The functor of one power series.* For each ring (unital and commutative) let

$$\Lambda(A) = 1 + tA[[t]] = \{f = f(t) = 1 + a_1 t + a_2 t^2 + \cdots :\ a_i \in A\} \tag{9.2}$$

be the set of power series over $A$ with constant term 1. Under the usual multiplication of power series this is an Abelian group (for which the power series 1 acts as the zero element). Thus (9.2) defines an Abelian group valued functor $\Lambda$: $\mathbf{CRing} \longrightarrow \mathbf{AbGroup}$. The morphism of Abelian groups associated to a morphism $\alpha$: $A \longrightarrow B$ is coefficient-wise, i.e.

$$\Lambda(\alpha)(1 + a_1 t + a_2 t^2 + a_3 t^3 + \cdots) = 1 + \alpha(a_1)t + \alpha(a_2)t^2 + \alpha(a_3)t^3 + \cdots \tag{9.3}$$

Giving a power series like (9.2) is of course the same thing as giving an infinite length vector $(a_1, a_2, a_3, \cdots)$ and in turn such a vector is the same as a morphism of commutative rings

$$\mathbf{Symm} = \mathbf{Z}[h_1, h_2, h_3, \cdots] = \mathbf{Z}[h] \longrightarrow A,\quad \varphi_f\colon\ h_i \mapsto a_i \tag{9.4}$$

where, as indicated in the notation, **Symm** is the ring of polynomials in an infinity of commuting indeterminates $h_1, h_2, h_3, \cdots$. Thus the functor $\Lambda$ is representable by **Symm**. The functorial addition on the Abelian group $\Lambda(A)$ then defines a comultiplication on **Symm**

---

[24] In earlier writings also 'generalized Witt vectors'.

$$\mu_S = \mu_{Sum}: \textbf{Symm} \longrightarrow \textbf{Symm} \otimes \textbf{Symm},\ h_n \mapsto \sum_{i+j=n} h_i \otimes h_j \tag{9.5}$$

where by definition and for ease of notation $h_0 = 1$. This makes **Symm** a Hopf algebra. The comultiplication formula of course encodes the ‘universal’ formula (that is ‘recipe’) for multiplying power series which is

$$(1+a_1+a_2+a_3+\cdots)(1+b_1+b_2+b_3+\cdots) = 1+c_1+c_2+c_3+\cdots \quad \Leftrightarrow \quad c_n = \sum_{i+j=n} a_i b_j$$

with, again, $a_0 = b_0 = c_0 = 1$. Which, in turn, is the same as saying that the addition on $\mathbf{CRing}(\mathbf{Symm}, A)$ is given by the convolution product

$$\varphi_{fg} = \textbf{Symm} \xrightarrow{\mu_S} \textbf{Symm} \otimes \textbf{Symm} \xrightarrow{\varphi_f \otimes \varphi_g} A \otimes A \xrightarrow{m_A} A \tag{9.6}$$

where $m_A$ is the multiplication on the ring $A$.

As in the case of the $p$-adic Witt vectors of the previous 7 sections it is often convenient (but never, strictly speaking, necessary) to work with ghost components. These are defined by the formula

$$s(f) = s(a_1, a_2, a_3, \cdots) = s_1 t + s_2 t^2 + s_3 t^3 + \cdots = t\frac{d}{dt}\log(f(t)) = \frac{tf'(t}{f(t)} \tag{9.7}$$

so that for example the first three ghost components are given by the universal formulas

$$s_1(a) = -a_1,\ s_2(a) = a_1^2 - a_2,\ s_3(a) = -a_1^3 + 2a_1a_2 - a_3 \tag{9.8}$$

Because of the properties of ‘log’ as defined formally by

$$\log(1+x) = x - \frac{x^2}{2} + \frac{x^3}{3} - \frac{x^4}{4} + \cdots + (-1)^{n+1}\frac{x^n}{n} + \cdots \tag{9.9}$$

so that $\log(fg) = \log(f) + \log(g)$, this implies that addition in $\Lambda(A)$ (which is multiplication of power series) corresponds to coordinate-wise addition of ghost components (just as in the case of the $p$-adic Witt vectors).

9.10. *The symmetric function point of view.* Now imagine that the $a_i$ are really the complete symmetric functions in a set of elements $\xi_1, \xi_2, \xi_3, \cdots$ (living as it were in some larger ring containing $A$); i.e.

$$a_n = h_n(\xi_1, \cdots, \xi_n, \cdots) \tag{9.11}$$

where the $h_n(X)$ are the familiar complete symmetric functions in the commuting indeterminates $X_1, X_2, X_3, \cdots$ , viz

$$h_n(X) = \sum_{j_1 \le j_2 \le \cdots \le j_n} X_{j_1} X_{j_2} \cdots X_{j_n} \tag{9.12}$$

The relation (9.11) can be conveniently written down in power series terms as

$$\prod_{i=1}^{\infty} \frac{1}{(1-\xi_i)} = 1 + a_1 t + a_2 t^2 + a_3 t^3 + \cdots = f(t) = f \tag{9.13}$$

The ghost components of $f$ are given by (9.7). In terms of the $\xi_i$ that means that

$$\begin{aligned} s_1 t + s_2 t^2 + s_3 t^3 + \cdots &= t\frac{d}{dt} \log f(t) = t\frac{d}{dt}\sum_i \log((1-\xi_i t)^{-1}) \\ &= -t\frac{d}{dt}\sum_i (-\xi_i t - \frac{\xi_i^2 t^2}{2} - \frac{\xi_i^3 t^3}{3} - \cdots) \\ &= \sum_{j=1}^{\infty} (\xi_1^j + \xi_2^j + \xi_3^j + \cdots) t^j \end{aligned} \tag{9.14}$$

so that in terms of the $\xi$ 's the ghost components are given by the power sums.

9.15. *Multiplication on the Abelian groups* $\Lambda(A)$. If there is to be a ring multiplication on $\Lambda(A)$ then in particular there should be a multiplication of the very special power series $(1-xt)^{-1}$ and $(1-yt)^{-1}$. Just about the simplest thing one can imagine is that the product of these very special power series is

$$(1-xt)^{-1} \quad (1-yt)^{-1} = (1-xyt)^{-1} \tag{9.15}$$

Something which fits with what has been seen in the theory of the $p$-adic Witt vectors: the $(1-xt)^{-1}$ are the Teichmüller representatives of the $x \quad A$ which should be multiplicative. Distributivity and functoriality together now force that the multiplication of power series in $\Lambda(A)$ should be given by

$$f(t) = \prod_i (1-\xi_i t)^{-1}, \; g(t) = \prod_i (1-\eta_i t)^{-1} \qquad (f \quad g)(t) = f(t) = \prod_{i,j} (1-\xi_i \eta_j t)^{-1} \tag{9.16}$$

Note that this formula makes perfectly good sense. The right most expression is symmetric in the $\xi$ 's and also symmetric in the $\eta$ 's and so, when written out, gives a power series with coefficients that are complete symmetric functions in the $\xi$ 's and in the $\eta$ 's. And that means that they are universal polynomials in the coefficients of $f$ and $g$.

Here is how this multiplication works out on the ghost components. As in (9.14)

$$\begin{aligned} t\frac{d}{dt}\log(\prod_{i,j} (1-\xi_i \eta_j t)^{-1}) &= \sum_{i,j} ((\xi_i \eta_j) t + (\xi_i \eta_j)^2 t^2 + (\xi_i \eta_j)^3 t^3 + \cdots) \\ &= \sum_{n=1}^{\infty} p_n(\xi) p_n(\eta) t^n \end{aligned} \tag{9.17}$$

Thus, multiplication according to (9.16) translates into component-wise multiplication for the ghost components. Actually because of distributivity it would have been sufficient to do this calculation for (9.15).

Associativity of the multiplication and distributivity follow by functoriality (or directly for that matter).

All in all, what has been defined is a unital-commutative-ring-valued functor

$$\Lambda\colon \mathbf{CRing} \longrightarrow \mathbf{CRing} \tag{9.18}$$

together with a set of 'ghost component' functorial ring morphisms

$$s_n\colon \Lambda(A) \longrightarrow A \tag{9.19}$$

which determine the ring structure on $\Lambda(A)$ by the very requirement that they be functorial ring morphisms.

The unit element of the ring $\Lambda(A)$ is the power series $1+t+t^2+t^3+\cdots = (1-t)^{-1}$.

9.20. *The functor of the big Witt vectors.* Now what has all this to do with Witt-vector-like constructions? As a set define

$$W(A) = A^{\mathbf{N}} = \{(x_1, x_2, x_3, \cdots)\colon\ x_i \in A\} \tag{9.21}$$

and set up a functorial bijection (of sets) between $W(A)$ and $\Lambda(A)$ by

$$e_A\colon (x_1, x_2, x_3, \cdots) \mapsto \prod_{n=1}^{\infty} (1 - x_n t^n)^{-1} \tag{9.22}$$

and transfer the functorial ring structure on the $\Lambda(A)$ to the $W(A)$ by this bijection. This then defines a unital-commutative-ring-valued functor

$$W\colon \mathbf{CRing} \longrightarrow \mathbf{CRing}$$

(that is isomorphic to $\Lambda$, so one might justifiably wonder why one takes the trouble). Let's calculate the ghost components of a Witt vector $(x_1, x_2, x_3, \cdots) \in W(A)$, which by definition means calculating

$$w_n(x) = w_n(x_1, x_2, x_3, \cdots) = s_n(e_A(x)) \tag{9.23}$$

This is easy

$$\begin{aligned} t\frac{d}{dt}\log\Big(\prod_d (1 - x_d t^d)^{-1}\Big) &= t\frac{d}{dt}\sum_d -\Big(\sum_m -m^{-1}(x_d t^d)^m\Big) \\ &= \sum_{d,m} d x_d^m t^{dm} = \sum_{n=1}^{\infty} \sum_{d\mid n} d x_d^{n/d} \end{aligned} \tag{9.24}$$

and thus the ghost components of the ring valued functor $W$ are given by the (generalized) Witt

polynomials

$$w_n(X_1, X_2, \cdots, X_n) = \sum_{d|n} dX_d^{n/d} \tag{9.25}$$

Note that these do indeed generalize the Witt polynomials (5.1) of the previous 8 sections in that, apart from a relabeling the polynomials $w_{p^n}$ of (9.25) are indeed the polynomials $w_n$ of (5.1).

**From now on in this whole chapter, that is in all sections that follow**, $w_n$ **and** $w_n(X)$ **and** $\cdots$ **refer to the polynomials as defined by (9.25)**.

Incidentally, the formula for (second) multiplication of power series, i.e. multiplication in the ring $\Lambda(A)$ of power series written in '$x$-coordinates' as in the right hand side of (9.22) can be written down directly as follows

$$(\prod_m (1-x_d t^d)^{-1}) \quad (\prod_n (1-y_n t^n)^{-1}) = \prod_{r,s\geq 1} (1-x_r^{m/r} y_s^{m/s} t^m)^{-rs/m} \tag{9.27}$$

where on the right hand side $m$ is the least common multiple of $r$ and $s$. This is a formula that was known to Witt, [419]. Thus there is no absolute need to introduce symmetric functions into the game; it is just very convenient.

As far as I know the first substantial treatment of the big Witt vectors from the symmetric functions point of view is due to Pierre Cartier, [84, 86].

9.28. The $a_i$ of (9.2) and (9.13), and the $x_i$ of (9.22) can be best seen as simply being different ways to coordinatize $\Lambda(A)$. As is usual with two different coordinate systems each has its advantages and diadvantages. One good property of the '$x$-coordinates' has already been pointed out. They show that the ring $\Lambda(A)$ is a (far reaching) generalization of the $p$-adic Witt vectors. Another is that the ghost component formulae, i.e the Witt polynomials $w_n(X)$ are very simple. For one thing they do not involve mixed terms in the $X_i$; for another $w_n(X)$ depends only on those $X_d$ for which $d$ is a divisor of $n$. This latter fact makes it clear that there are many interesting quotient functors of $W$ (which is difficult to see in the $\Lambda$ formulation). These quotient functors will be described in the next section, 10.

On the other hand there are good well known formuals both ways (from symmetric function theory) that relate the '$a$-coordinates' (i.e. the power series coordinates) with their ghost components while inverting the Witt polynomials seems to be a bit of a mess.

Quite a number of formulas from symmetric function theory are important (or at least very useful) for Witt vector theory. In the next subsections there are a few.

9.29. *Very partial symmetric function formularium* (1) [25]. Everything takes place in the rings

$$\mathbf{Z}[\xi_i :\ i \quad I] \quad \text{respectively} \quad \mathbf{Z}[[\xi_i :\ i \quad I]]$$

of polynomials, respectively power series, in a countable infinite collection of commuting indeterminates over the integers. (And also in their analogues $\mathbf{Q}[\xi_i :\ i \quad I]$ and $\mathbf{Q}[[\xi_i :\ i \quad I]]$ over the rationals.) For some definitional and foundational remarks on what it means to have a polynomial or power series in an infinite number of variables see the appendix to this chapter.

[25] Not all of the formulas and functions that follow belong to the standarad formularium of the symmetric functions.

The rings just above are graded by giving all the $\xi_i$ degree 1.

The reader who does not like polynomials and power series in an infinite number of variables can just imagine that things are in terms of a finite number of them, large enough for the business at hand.

9.30. *Partitions and monomial symmetric functions.* A partition $\lambda = (\lambda_1, \lambda_2, \cdots, \lambda_n)$ is a finite sequence of elements from $\mathbf{N}\ \ \{0\}$ in non-increasing order, $\lambda_1 \geq \lambda_2 \geq \cdots \geq \lambda_i \geq \cdots$ . Partitions with different numbers of trailing zeros are identified. The length of a partition, $\lg(\lambda)$, is the number of its non zero entries; its weight is $\mathrm{wt}(\lambda) = \lambda_1 + \lambda_2 + \cdots + \lambda_n$. A different notation for a partition is $(1^{m_1(\lambda)} 2^{m_2(\lambda)} \cdots n^{m_n(\lambda)})$ where $m_i(\lambda)$ is the number of entries of it that are equal to $i$. A partition of weight $n$ is said to be a partition of $n$ and its entries are (also) called the parts of that partition.

For any sequence $\alpha = (a_1, a_2, \cdots, a_m)$ of elements from $\mathbf{N}\ \ \{0\}$ let

$$\xi^\alpha = \xi_1^{a_1} \xi_2^{a_2} \cdots \xi_m^{a_m}, \quad \text{e.g. } \xi^{(0,0,0,2,0,3,0,0)} = \xi_4^2 \xi_6^3 \tag{9.31}$$

Given a partition $\lambda$ (strictly speaking with an infinite string of trailing zeros added), the associated monomial symmetric function is the sum

$$m_\lambda = \sum_\alpha \xi^\alpha \tag{9.32}$$

where the sum is over all *distinct* sequences $\alpha$ (infinite with only a finite number of nonzero entries) that are permutations of $\lambda$ (with trailing zeros added). Thus, for example,

$$m_{(2,1)} = \sum_{i \neq j} \xi_i^2 \xi_j = \xi_1^2 \xi_2 + \xi_2^2 \xi_1 + \xi_1^2 \xi_3 + \xi_3^2 \xi_1 + \xi_2^2 \xi_3 + \xi_3^2 \xi_2 + \cdots \tag{9.33}$$

$$m_{(1,1)} = \sum_{i < j} \xi_i \xi_j = \xi_1 \xi_2 + \xi_1 \xi_3 + \xi_2 \xi_3 + \cdots \tag{9.34}$$ [26]

These are the simplest symmetric functions in that if a symmetric function contains one of the monomials from an $m_\lambda$ then it also contains all others (with the same coefficient). The monomial symmetric function $m_\lambda$ is homogeneous of degree $\mathrm{wt}(\lambda) = n$ and they form a basis for the free Abelian group of symmetric functions of weight $n$.

By definition the complete symmetric functions and the elementary symmetric functions are

$$h_n = \sum_{\mathrm{wt}(\lambda) = n} m_\lambda \quad \text{and} \quad e_n = m_{(1,1,\cdots,1)} \quad (n \text{ 1's}) \tag{9.35}$$

and for a partition $\lambda = (\lambda_1, \lambda_2, \cdots, \lambda_m)$ one writes

---

[26] Note the difference between the two formulae which rests on the word ‘distinct’ in the sentence just below (9.32).

$$h_\lambda = h_{\lambda_1} h_{\lambda_2} \cdots h_{\lambda_m} \quad \text{and} \quad e_\lambda = e_{\lambda_1} e_{\lambda_2} \cdots e_{\lambda_m} \tag{9.36}$$ [27]

The complete symmetric functions and the elementary ones are related by the Wronski relations

$$\sum_{i+j=n} (-1)^i h_i e_j \;=\; 0 \quad \text{for } n \geq 1, \text{ where } h_0 = e_0 = 1 \tag{9.37}$$

The $h_\lambda$ (resp. $e_\lambda$) for $\lambda$ of weight $n$ also form a basis for the symmetric functions of degree $n$. By definition

$$\mathbf{Symm} = \quad {}_n\mathbf{Symm}_n \qquad \mathbf{Z}[\xi] \tag{9.38}$$

is the graded ring of symmetric functions in the commuting indeterminates $\xi$. In this whole chapter it will (usually) be seen as the graded ring in the complete symmetric functions

$$\mathbf{Symm} = \mathbf{Z}[h_1, h_2, \cdots, h_m, \cdots] = \mathbf{Z}[h] \tag{9.39}$$

For many purposes one could equally well work with the elementary symmetric functions, but the complete ones just work out better. Actually there are precisely four 'canonical' choices, see subsection 10.18 below on the Liulevicius theorem.

9.40. *Inner product.* One defines a (symmetric positive definite) inner product on **Symm** by declaring the monomial and the complete symmetric functions to form a bi-orthogonal system:

$$< h_\lambda, m_\kappa > \;=\; \delta_{\lambda,\kappa} \quad \text{(Kronecker delta)} \tag{9.41}$$ [28]

There is now a remarkable theorem, see [281], (4.6), p. 63:

9.42. *Theorem.* Let $\{u_\lambda\}$ and $\{v_\lambda\}$ be two sets of symmetric functions indexed by all partitions (including the zero partition) that both form a basis for $\mathbf{Symm}_\mathbf{Q} = \mathbf{Symm} \;\; _\mathbf{Z} \; \mathbf{Q}$ (or **Symm** itself). Then they form a bi-orthonormal system (i.e $< u_\lambda, v_\kappa > \;=\; \delta_{\lambda,\kappa}$) if and only if

$$\sum_\lambda u_\lambda(\xi) v_\lambda(\eta) = \prod_{i,j} (1 - \xi_i \eta_j)^{-1} \tag{9.43}$$

Here $\eta = \{\eta_1, \eta_2, \cdots\}$ is a second set of commuting indeterminates that also commute with the $\xi$ 's, and $v_\lambda(\eta)$ is the same as $v_\lambda(\xi) = v_\lambda$ bur now written as an expression in the $\eta$ 's.

The theorem is especially remarkable to a mathematician with (products of) Witt vectors on his mind [29] and I am far from sure that its consequences in that context have been fully explored and understood.

---

[27] Note that this notation has a rather different meaning than that for $m_\lambda$. Is is definitely not true that $m_\lambda$ is equal to something like $m_{(\lambda_1)} m_{(\lambda_2)} \cdots m_{(\lambda_m)}$.

[28] This inner aproduct is often called 'Hall inner product'.

[29] See (9.16).

Theorem 9.42 fits with formula (9.41) for it is indeed the case that ([281], p. 62)

$$\sum_{\lambda} h_{\lambda}(\xi)m_{\lambda}(\eta) = \prod_{i,\,j}(1-\xi_i\eta_j)^{-1} \tag{9.44}$$

9.45. *Schur functions*. A symmetric function, i.e. an element of **Symm**, is called positive if when expressed as a sum in the monomial basis $\{m_{\lambda}\}$ all its coefficient are non negative. All the specific symmetric functions introduced so far, the $m_{\lambda}$, $h_{\lambda}$, $e_{\lambda}$, are positive.

It now turns out that there is a unique positive orthonormal basis of **Symm**. They are called the Schur functions, and are determined by

$$\sum_{\lambda} s_{\lambda}(\xi)s_{\lambda}(\eta) = \prod_{i,j}(1-\xi_i\eta_j)^{-1} \tag{9.46}$$

and the positivity requirement.

The first few Schur functions are

$$s_{(1)} = m_{(1)} = h_1,\quad s_{(2)} = m_{(2)} + m_{(1,1)} = h_2,\quad s_{(1,1)} = m_{(1,1)} = -h_2 + h_{(1,1)} \tag{9.47}$$

and the ones in weights 3 and 4 are given by the matrices

$$\begin{pmatrix} 1 & 0 & 0 \\ 1 & 1 & 0 \\ 1 & 2 & 1 \end{pmatrix},\ \text{s-m matrix;}\quad \begin{pmatrix} 1 & -1 & 1 \\ 0 & 1 & -2 \\ 0 & 0 & 1 \end{pmatrix},\ \text{s-h matrix} \tag{9.48}$$

$$\begin{pmatrix} 1 & 0 & 0 & 0 & 0 \\ 1 & 1 & 0 & 0 & 0 \\ 1 & 1 & 1 & 0 & 0 \\ 1 & 2 & 1 & 1 & 0 \\ 1 & 3 & 2 & 3 & 1 \end{pmatrix},\ \text{s-m matrix;}\quad \begin{pmatrix} 1 & -1 & 0 & 1 & -1 \\ 0 & 1 & -1 & -1 & 2 \\ 0 & 0 & 1 & -1 & 1 \\ 0 & 0 & 0 & 1 & -3 \\ 0 & 0 & 0 & 0 & 1 \end{pmatrix},\ \text{s-h matrix} \tag{9.49}$$

where the partitions have been ordered lexicographically, (4), (3,1), (2,2), (2,1,1), (1,1,1,1), and the columns give the coefficients of the Schur functions in terms of the monomial symmetric functions (resp. the complete symmetric functions). Thus e.g.

$$\begin{aligned} s_{(3,1)} &= m_{(3,1)} + m_{(2,2)} + 2m_{(2,1,1)} + 3m_{(1,1,1,1)} \\ s_{(2,1,1)} &= h_{(4)} - h_{(3,1)} - h_{(2,2)} + h_{2,1,1)} \end{aligned} \tag{9.50}$$ [30]

[30] Calculating with the monomial symmetric functions directly is a bit messy (and so is calculating with their explicit expressions). For instance $m_{(3,1)}m_{(1)} = m_{(4,1)} + m_{(3,2)} + 2m_{(3,1,1)}$, $m_{(4,1)}m_{(2)} = m_{(6,1)} + m_{(4,3)} + m_{(4,2,1)}$, $m_{(2,1)}m_{(1)} = 2m_{(2,2)} + m_{(3,1)} + 2m_{(2,1,1)}$. It is fairly easy to see what monomial functions should occur in such a product; things are less clear as regards the coefficients with which they occur. For instance, why is there a 2 in front of $m_{(2,2)}$ in the third formula but not one in front of $m_{(3,2)}$ in the first. Things become a good deal better (in my opinion) when one works in the larger ring **QSymm** **Symm** of quasi-symmetric functions where there is a clear easy to use formula for multiplying quasi-symmetric monomial functions. This is the overlapping shuffle

There are many ways to define the Schur functions and to write down formulas for them. The standard (and oldest) definition is the one by Jacobi, see [281], p. 40ff. An explicit expression is given by the Jacobi-Trudy formulas

$$s_\lambda = \det(h_{i} - i + j)_{1\le i,j\le n} \tag{9.51}$$

Other expressions come from the Schur identity

$$\sum_\lambda s_\lambda = \prod_i (1-\xi_i)^{-1} \prod_{i<j} (1-\xi_i\xi_j)^{-1} \tag{9.52}$$

(see [273], 5.4, p. 176). This one has the advantage of showing immediately that Schur symmetric functions are positive (which is far from obvious from (9.51).

9.53. *Forgotten symmetric functions.* The complete symmetric functions $\{h_\lambda\}$ form a basis dual to that of the monomial symmetric functions. There also is of course a dual basis to the basis of elementary symmetric functions $\{e_\lambda\}$. These are obtained by replacing the $h_i$ by the $e_i$ in the formulas for the $m_\lambda$ in the $h_i$. There seem to be no 'nice' formulas for them, [281], p. 22.

9.54. *Power sum syymetric functions.* The power sum symmetric functions are defined by

$$p_n = m_{(n)} = \sum_i \xi_i^n \tag{9.55}$$

Recall that they are the ghost components of the element

$$H(t) = 1 + h_1 t + h_2 t^2 + h_3 t^3 + \cdots \qquad \Lambda(\mathbf{Symm}) = 1 + t\mathbf{Symm}[[t]] \tag{9.56}$$

The power sum symmetric functions are related to the complete symmetric functions by the Newton relations

$$nh_n = p_n + h_1 p_{n-1} + h_2 p_{n-2} + \cdots + h_{n-1} p_1 \tag{9.57}$$

which in terms of the generating functions $H(t)$, see (9.56), and

$$P(t) = \sum_{i\ge1} p_i t^i \quad \text{(so that } p_0 \text{ is set to zero)} \tag{9.58}$$

can be encoded by the 'differental' equation (defining ghost components)

$$P(t) = t\frac{d}{dt}\log H(t) = \frac{tH'(t)}{H(t)} \tag{9.59}$$

From (9.59) one readily obtains a formula for the power sum symmetric functions in terms of the complete symmetric functions, i.e. a universal formula for the ghost components of an element $1 + a_1 t + a_2 t^2 + \cdots$ in terms of its '$a$-coordinates'.

---

product, see 11.26. See also the appendix to this chapter for some more details.

$$p_n = \sum_{r_1+r_2+\cdots+r_k=n} (-1)^{k+1} r_1 h_{r_1} h_{r_2} \cdots h_{r_k}; \quad r_i \quad \mathbf{N} = \{1,2,\cdots\} \tag{9.60}$$

This formula can readily be inversted by ‘solving’ (9.59) using formal exponentials.

For each partition $\lambda$ (with no zero parts) define also

$$p_\lambda = p_{\lambda_1} p_{\lambda_2} \cdots p_{\lambda_n} \tag{9.61}$$

Then, obviously, given the Newton relations, the $p_\lambda$ form a homogeneous basis over the rationals for the vector space $\mathbf{Symm_Q}$. But they are not a basis for **Symm** itself. For a partition $\lambda = (\lambda_1, \lambda_2, \cdots, \lambda_n)$ define the integer $z_\lambda$ by

$$z_\lambda = 1^{m_1(\lambda)} m_1(\lambda)! 2^{m_2(\lambda)} m_2(\lambda)! \cdots n^{m_n(\lambda)} m_n(\lambda)!$$

where, as before, $m_i(\lambda)$ is the number of parts $i$ in $\lambda$. Then ([281], p. 62)

$$\sum_\lambda z_\lambda^{-1} p_\lambda(\xi) p_\lambda(\eta) = \prod_{i,j} (1-\xi_i \eta_j)^{-1} \tag{9.62}$$

showing that suitably normalized $p_\lambda$ form an orthonormal basis for $\mathbf{Symm_Q}$ (but not of course for **Symm** as these normalized $p_\lambda$ are not even elements of **Symm**).

9.63. *On the Witt coordinates*. Now consider the transformation from ‘$a$-coordinates’ of a 1-power-series, i.e. an element of $\Lambda(A)$, to Witt vector coordinates. This is given by the universal formula, see (9.22)

$$\prod_d (1 - x_d t^d)^{-1} = 1 + h_1 t + h_2 t^2 + \cdots \tag{9.64}$$

Assuredly the $x_d$ as defined by (9.64) are symmetric functions and also, obviously, from (9.64) one has

$$x_d = h_d + (\text{homogenous polynomial of weight } d \text{ in the } h_1, \cdots, h_{d-1}) \tag{9.65}$$

and so the

$$x_\lambda = x_{\lambda_1} x_{\lambda_2} x_{\lambda_3} \cdots \tag{9.66}$$

form (yet another) homogeneous basis for **Symm**. It is quite easy to express the $h_n$ in terms of the $x_d$. Indeed rewriting the left hand side of (9.61) as

$$\prod_i (1 + x_i t^i + x_i^2 t^{2i} + \cdots) \tag{9.67}$$

it is immediate that

$$h_n = \sum_{\mathrm{wt}(\lambda)=n} x_\lambda = \sum_{\mathrm{wt}(\lambda)=n} x_1^{m_1(\lambda)} x_2^{m_2(\lambda)} x_3^{m_3(\lambda)} \cdots \tag{9.68}$$

However, finding a formula the other way seems to be rather messy [31]. On the other hand, there is a duality formula. For each symmetric function $f(\xi)$ introduce the shorthand notation

$$f(\xi^i) = f(\xi_1^i, \xi_2^i, \xi_3^i, \cdots) \tag{9.69}$$ [32]

Now for any partition $\lambda = (1^{k_1} 2^{k_2} \cdots n^{k_n})$ introduce

$$r_\lambda = h_{k_1}(\xi^1) h_{k_2}(\xi^2) \cdots h_{k_n}(\xi^n)$$

It is not difficult to check that these also form a basis of **Symm** [33]. There is now the following duality formula

$$< x_\lambda, r_\kappa > \; = \; \delta_{\lambda,\kappa} \tag{9.70}$$

I.e. the $x_\lambda$ form a dual or adjoint basis to the one of the $r_\lambda$. This is due to [350], where there is a quite elegant proof using theorem 9.42.

The symmetric functions $x_n$ also have an interesting, even remarkable, positivity property:

9.71. *Theorem.* The symmetric functions $-x_1, x_2, \cdots, x_n, \cdots$ are Schur positive [34].

Here ‘Schur positive’ means that when expressed as a **Z**-linear combination of the Schur symmetric functions all the coefficients are nonnegative. This is stronger than being positive of course. Theorem 9.71 was conjectured by Christoph Reutenauer, [333, 332], and proved by Thomas Scharf and Jean-Yves Thibon, [350], and also by W Doran. Some more investigations relating to the $x_n$ are reported in loc. cit. Much remains to be done in my opinion, particularly regarding exploiting the positivity result of theorem 9.71 [35].

9.72. *Integrality lemmata. The miracle of the Witt polynomials* (2). For the generalized Witt polynomials $w_n(X)$ there is the same integrality theorem as for the $p$-adic ones.

---

[31] This corresponds to the fact that inverting the Witt vector polynomials appears to be a messy business.

[32] This is the (outer) plethysm with pespect to the power sum symmetric function $p_i$; see subsections 16.76, 18.35 and 18.37 below.

[33] This one has a certain amount of special interest in that for any symmetric function $f$ the integrality conditions $< f, r_\lambda >$ **Z** are equivalent to certain congruences established by Andreas Dress, [118], for testing whether a central function of a symmetric group $S_n$ is a virtual character; see Thomas Scharf and Jean-Yves Thibon, [350].

[34] There is no printing error here; there really is just one minus sign.

[35] Quite generally the systematic use of (partial) orderings and positivy properties in Hopf algebra theory (and hence Witt vector theory) has started only fairly recently.

9.73. *Theorem.* Let $\varphi(X,Y,Z)$ be a polynomial over the integers in three (or less, or more) commuting indeterminates. then there are unique polynomials over the integers $\varphi_n(X_1,X_2,\cdots,X_n;Y_1,Y_2,\cdots,Y_n;Z_1,Z_2,\cdots,Z_n)=\varphi_n(X;Y;Z)$, $n=1,2,3,\cdots$ such that for all $n\geq 1$

$$w_n(\varphi_1(X;Y;Z),\cdots,\varphi_n(X;Y;Z))=\varphi(w_n(X),w_n(Y),w_n(Z)) \qquad (9.74)$$

The proof is (basically) the same as in the $p$-adic case.

The theory of the functor of the big Witt vectors can now be developed starting from this integrality theorem in just the same way as was done for the $p$-adic Witt vectors in section 5 above. That is define addition polynomials $\mu_{S,i}(X_1,\cdots,X_i;Y_1,\cdots,Y_i)$ and multiplication polynomials $\mu_{P,i}(X_1,\cdots,X_i;Y_1,\cdots,Y_i)$ by taking the polynomial $\varphi(X,Y)$ in theorem 9.73 respectively equal to $X+Y$ and $XY$, so that

$$\begin{aligned} w_n(\mu^W_{S,1}(X;Y),\mu^W_{S,n}(X;Y).\cdots,\mu^W_{S,n}(X;Y)) &= w_n(X)+w_n(Y) \\ w_n(\mu^W_{P,1}(X;Y),\mu^W_{P,2}(X;Y).\cdots,\mu^W_{S,n}(X;Y)) &= w_n(X)w_n(Y) \end{aligned} \qquad (9.75)$$

and define an addition and multiplication on Witt vectors $a=(a_1,a_2,a_3,\cdots)$ by

$$\begin{aligned} (a_1,a_2,a_3,\cdots)+_W(b_1,b_2,b_3,\cdots) &= (\mu^W_{S,1}(a;b),\mu^W_{S,2}(a;b),\mu^W_{S,3}(a;b),\cdots) \\ (a_1,a_2,a_3,\cdots)\ {}_W\ (b_1,b_2,b_3,\cdots) &= (\mu^W_{P,1}(a;b),\mu^W_{P,2}(a;b),\mu^W_{P,3}(a;b),\cdots) \end{aligned} \qquad (9.76)$$

This then defines the functorial ring $W(A)$ with multiplication and addition given by (9.76) with zero element $(0,0,0,\cdots)$ and unit element $(1,0,0,0,\cdots)$.

9.77. Theorem 9.73 is only one of a slew of integrality theorems. One of my favourites is the 'functional equation lemma' which is of great use in formal group theory; particularly in the construction of universal formal groups, see [192], Chapter 1, §2.3.

9.78. *Ingredients for the functional equation integrality lemma.* The ingredients for the functional equation lemma are the following

$$A \quad K,\ \sigma: K \to K,\ \mathfrak{p} \quad A,\ p,\ q,\ s_1,s_2,s_3,\cdots \qquad (9.79)$$

Here $A$ is a subring of a ring $K$, $\sigma$ is a ring endomorphism of $K$, $\mathfrak{a}$ is an ideal in $A$, $p$ is a prime number and $q$ is a power of $p$, and the $s_i$ are elements of $K$. These ingredients are supposed to satisfy the following conditions

$$\sigma(A) \quad A,\ \sigma(a)\equiv a^q \bmod \mathfrak{p} \text{ for all } a \quad A,\ p \quad \mathfrak{p},\ s_ib \quad A \text{ for all } b \quad \mathfrak{p} \qquad (9.80)$$

and also

$$\mathfrak{p}^r b \quad \mathfrak{p} \qquad \mathfrak{p}^r\sigma(b) \quad \mathfrak{p} \text{ for all } r \quad \mathbf{N},\ b \quad K \qquad (9.81)$$

a property that is automatically satisfied if the ideal $\mathfrak{p}=(c)$ is principal and $\sigma(c)=uc$ for some unit of $A$.

Here are three examples of such situations. There are many more.

$$A = \mathbf{Z},\ K = \mathbf{Q},\ \sigma = \mathrm{id},\ q = p,\ \mathfrak{p} = p\mathbf{Z},\ s_i \quad p^{-1}\mathbf{Z} \tag{9.82}$$

$$A = \mathbf{Z}_{(p)},\ K = \mathbf{Q},\ \sigma = \mathrm{id},\ q = p,\ \mathfrak{p} = p\mathbf{Z}_{(p)},\ s_i \quad p^{-1}\mathbf{Z}_{(p)} \tag{9.83}$$

$$\begin{gathered} A = \mathbf{Z}[V_1, V_2, \cdots],\ K = \mathbf{Q}[V_1, V_2, \cdots],\ \sigma f(V_1, V_2, \cdots) = f(V_1^p, V_2^p, \cdots), \\ q = p,\ \mathfrak{p} = pA,\ s_i \quad p^{-1}A \end{gathered} \tag{9.84}$$

9.82. *Constructions for the functional equation lemma.* Now, let $g(X) = b_1X + b_2X^2 + \cdots$ be a power series with coefficients in $A$. Using the ingredients (9.78) construct from it a new power series by the recursion formula (or functional equation)

$$f_g(X) = g(X) + \sum_{i=1}^{\infty} s_i \sigma^i f_g(X^{q^i}) \tag{9.85}$$

where $\sigma^i f_g(X)$ is the power series obtained from $f_g(X)$ by applying $\sigma^i$ to the coefficients of $f_g(X)$. Three examples of power series obtained in this way are:

$$X + p^{-1}X^p + p^{-2}X^{p^2} + p^{-3}X^{p^3} + \cdots \tag{9.86}$$

$$\log(1+X) = \sum_{n=1}^{\infty} n^{-1}(-1)^{n+1}X^n \tag{9.87}$$

$$\begin{gathered} X + p^{-1}V_1X^p + (p^{-2}V_1V_1^p + p^{-1}V_2)X^{p^2} \\ +(p^{-3}V_1V_1^pV_1^{p^2} + p^{-2}V_1V_2^p + p^{-2}V_2V_1^{p^2} + p^{-1}V_3)X^{p^3} + \cdots \end{gathered} \tag{9.88}$$

9.89. *Functional equation integrality.*

*Lemma.* Let $A \quad K,\ \sigma,\ \mathfrak{p} \quad A,\ p,\ q,\ s_1, s_2, s_3, \cdots$ be as in (9.77) and be such that (9.80) and (9.81) hold. Let $g(X) = b_1X + b_2X^2 + \cdots$ and $\bar{g}(X) = \bar{b}_1X + \bar{b}_2X^2 + \cdots$ be two power series with coefficients in $A$ and let $b_1$ be invertible. Let this time $f^{-1}(X)$ denote the functional inverse of $f(X)$ (for a power series with zero constant term and nonzero coefficient in degree 1, defined by $f^{-1}(f(X)) = X$. Then

The power series $F_g(X,Y) = f_g^{-1}(f_g(X) + f_g(Y))$ has its coefficients in $A$ (9.90)

The power series $f_g^{-1}(f_{\bar{g}}(X))$ has its coefficients in $A$ (9.91)

If $h(X)$ is a power series with zero constant term over $K$, then

$$h(X) \equiv \bar{g}(X) \bmod \mathfrak{p}^r A[[X]] \quad \Leftrightarrow \quad f_g(h(X)) \equiv f_g(\bar{g}(X)) \bmod \mathfrak{p}^r A[[X]] \tag{9.92}$$

Statement (9.90) says that $F_g(X,Y)$ is a (one dimensional) formal group over $A$ with logarithm $f_g(X)$. Example (9.88) is the logarithm of the universal one dimensional $p$-typical formal group.

There are also more dimensional and infinite dimensional versions and these can be used to construct universal formal groups and also the Witt vectors. For all these statements and the appropriate definitions see [192].

Another application of the functional equation lemma is the following statement, especially useful in the Witt vector context

9.93. *Lemma.* [192], lemma 17.6.1, page 137. Let $A$ be a characteristic zero ring with endomorphisms $\varphi_p$ for all prime numbers $p$ such that $\varphi_p(a) \equiv a^p$ for all $a \in A$. Then for a given sequence $b_1, b_2, \cdots$ in $A$ there exists another sequence $x_1, x_2, \cdots$ in $A$ such that $w_n(x) = b_n$ for all $n$ (i.e. the sequence $b_1, b_2, \cdots$ lies in the image of the ghost mapping) if and only if

$$\varphi_p(b_n) \equiv b_{np} \bmod (p^{v_p(n)+1}) \tag{9.94}$$

Here, as before, see example 6.10, $v_p$ denotes the $p$-adic valuation on $\mathbf{Z}$. A vector with components in $A$ which is in the image of $w\colon W(A) \longrightarrow A^{\mathbf{N}}$, $(a_1, a_2, \cdots) \mapsto (w_1(a), w_2(a), \cdots)$ will be called a ghost-Witt vector over $A$.

As far as I know, this lemma, in this form, first appeared in [86]. Pierre Cartier attributes it to Bernard Dwork and Jean Dieudonné, [107], proposition 1; [132], lemma 1. The formulations there are a bit different.

In e.g. the case $A = \mathbf{Z}$ where one takes all the endomorphisms $\varphi_p$ to be the identity, the integrality condidion (9.94) can be given the following quite elegant formulation. Write down the Mobius inversion transform of the sequence $b_1, b_2, \cdots$, i.e. the sequence

$$c_n = \sum_{d|n} b_{n/d}\mu(d) \tag{9.95}$$

Then the sequence $b_1, b_2, \cdots$ satisfies (9.94) if and only if $c_n \equiv 0 \bmod n$ for all $n$. This remark seems to be due to Albrecht Dold, [113]. Here the Möbius function $\mu\colon \mathbf{N} \longrightarrow \mathbf{Z}$ is defined by $\mu(1) = 1$, $\mu(n) = (-1)^r$ if $n$ is the product of precisely $r$ different primes and $\mu(n) = 0$ otherwise. Still another formulation of the integrality condition is

$$\sum_{i=1}^{n} b_{(i,n)} \equiv 0 \bmod n \tag{9.96}$$

where $(i,n)$ denotes the greatest common divisor of $i$ and $n$.

The (integrality aspects of the) theory of Witt vectors can be developed solely on the basis of this lemma 9.93. This is how things are done in [212].

9.97. *Witt vectors over the integers and fixed points and fixed point indices of iterated mappings.* In [113], Dold proves that a sequence of integers $s = (s_1, s_2, s_3, \cdots)$ is the sequence of fixed point indices of the iterates of a continuous mapping if and only if the formal power series

$$\exp\left(-\sum_{i=1}^{\infty} \frac{s_i}{i} t^i\right)$$

has integral coefficients; that is iff $s$ is the ghost vector of a Witt vector over the integers. This is done using the integrality criterium (9.96).

More concretely, call a sequence of nonnegative integers exactly realizable if there is a set $X$ and a map $f: X \longrightarrow X$ such that the number of fixed points of the $n$-th iterate $f^n$ of $f$ is exactly $a_n$. Then a sequence is exactly realizable if and only if the numbers

$$\sum_{d|n} \mu(d) a_{n/d}$$

are nonnegative and divisible by $n$, [330]. Similar things happen when studying iterates of mappings of the unit interval into itself, i.e. in the part of dynamical system and chaos theory that deals with such mappings, see [70, 292].

The first extensive treatment that I know of of the Witt vector constructions based basically totally on the symmetric function point of view is by Pierre Cartier, [86].

## 10. The Hopf algebra Symm as the representing algebra for the big Witt vectors.

Recall that $\mathbf{Symm} = \mathbf{Z}[h_1, h_2, h_3, \cdots] = \mathbf{Z}[h]$ represents the functor $\Lambda$ and hence the functor $W$ of the big Witt vectors

$$\Lambda(A) = \{a(t) = 1 + a_1 t + a_2 t^2 + a_3 t^3 + \cdots : a_i \in A\} = \mathbf{CRing}(\mathbf{Z}[h], A) \tag{10.1}$$

The fact that this is an Abelian group valued functor means that there is a comultiplication making **Symm** a coalgebra object in the category **CRing** of commutative unital rings, which in turn means that there is a comultiplication, and a co-unit

$$\mu_S: \mathbf{Symm} \longrightarrow \mathbf{Symm} \otimes \mathbf{Symm}, \quad \varepsilon_S: \mathbf{Symm} \longrightarrow \mathbf{Z}, \tag{10.2}$$

that are ring morphisms and that make $(\mathbf{Symm}, \mu_S, \varepsilon_S)$ [36] a co-associative, cocommutative co-unital coalgebra. See [293] for a lot of material on coalgebras.

Because the addition on the set $\Lambda(A)$ is defined as multiplication of power series and the unit element is the power series 1, the morphisms $\mu_S, \varepsilon_S$ are given by

$$\mu_S(h_n) = 1 \otimes h_n + \sum_{i=1}^{n-1} h_i \otimes h_{n-i} + h_n \otimes 1, \quad \varepsilon_S(h_n) = 0, \quad n \geq 1 \tag{10.3}$$

(It is often convenient to define $h_0 = 1$ so that the multiplication can be written $\sum_{i+j=n} h_i \otimes h_j$.)

Writing

$$m: \mathbf{Symm} \otimes \mathbf{Symm} \longrightarrow \mathbf{Symm}, \quad e: \mathbf{Z} \longrightarrow \mathbf{Symm} \tag{10.4}$$

[36] The index ‘S here refers to the matter that these are the ring morphisms that give the sum part of the Witt vector functor.

for the multiplication and unit element of **Symm**, in total there is a bialgebra structure $(\mathbf{Symm}, m, e, \mu_S, \varepsilon_S)$. Finally there is an antipode

$$\iota_S\colon \mathbf{Symm} \longrightarrow \mathbf{Symm} \tag{10.5}$$

The existence of an antipode comes for free in the present case, see below. The antipode is an anti ring morphism (and hence in the present case a ring morphism because of commutativity). The antipode is determined by

$$\iota_S(h_n) = -h_n \tag{10.6}$$

All this makes the total structure $(\mathbf{Symm}, m, e, \mu_S, \varepsilon_S, \iota_S)$ a Hopf algebra over the integers. For the basic theory of Hopf algebras, see [95].

10.7. *Skew Schur functions.* (Very partial symmetric function formularium (2)). Let $\kappa, \lambda$ be two partitions and let $s_\kappa, s_\lambda$ be the corresponding Schur functions. As the Schur functions form an orthogonal basis, there are coefficients $c^{\nu}_{\kappa,\lambda}$ such that

$$s_\kappa s_\lambda = \sum_\nu c^{\nu}_{\kappa,\lambda} s_\nu \quad \text{or, equivalently,} \quad \langle s_\kappa s_\lambda, s_\nu \rangle = c^{\nu}_{\kappa,\lambda} \tag{10.8}$$

The multiplicity coeffients[37] $c^{\nu}_{\kappa,\lambda}$ have a combinatorial interpretation, see [281], Ch.1, §9, p. 143, and are, hence, nonnegative integers. Define the skew Schur function $s_{\kappa/\lambda}$ by

$$\langle s_{\kappa/\lambda}, s_\nu \rangle = \langle s_\kappa, s_\lambda s_\nu \rangle, \quad \text{or, equivalently,} \quad s_{\kappa/\lambda} = \sum_\nu c^{\kappa}_{\lambda,\nu} s_\nu \tag{10.9}$$

It turns out that $s_{\kappa/\lambda} = 0$ unless $\lambda \subset \kappa$ which by definition means that $\lambda_i \le \kappa_i$ for all $i$.[38] There is a determinantal formula for the skew Schur functions as follows

$$s_{\kappa/\lambda} = \det(h_{\kappa_i - \lambda_j - i + j})_{i,j} \tag{10.10}$$

In terms of skew Schur functions the comultiplication of **Symm** can be written

$$\mu_S(s_\kappa) = \sum_\lambda s_{\kappa/\lambda} \otimes s_\lambda \tag{10.11}$$

as follows from the duality formula $\langle xy, z \rangle = \langle x \otimes y, \mu_S(z) \rangle$ in **Symm**.

10.12. *Grading.* **Symm** is a graded Hopf algebra. That is, as an Abelian group it is a direct sum

---

[37] In the representation theoretic incarnation of **Symm**, see section 18 below, these coefficients turn up as multiplicities of irreducible representations (Littlewood-Richardson rule).

[38] $\kappa \supset \lambda$ is the case if and only if the diagram of $\kappa$ contains the diagram of $\lambda$; whence the notation.

$$\mathbf{Symm} = \bigoplus_{n=0}^{\infty} \mathbf{Symm}_n \tag{10.13}$$

where the homogenous part of degree (or weight) $n$ is spanned by the monomials in the $h_i$ of weight $n$ where $h_i$ has weight $i$; i.e. by the monomials $h_\lambda$ with $\mathrm{wt}(\lambda) = n$, and that moreover the multiplication and comultiplication and unit and co-unit are compatble with the grading, which means

$$m(\mathbf{Symm}_i \otimes \mathbf{Symm}_j) \subset \mathbf{Symm}_{i+j}, \quad e(\mathbf{Z}) \subset \mathbf{Symm}_0$$

$$\mu_S(\mathbf{Symm}_n) \subset \bigoplus_{i+j=n} \mathbf{Symm}_i \otimes \mathbf{Symm}_j, \quad \varepsilon(\mathbf{Symm}_n) = 0 \text{ for } n \geq 1$$

$$\iota_S(\mathbf{Symm}_n) \subset \mathbf{Symm}_n$$

Note that $\mathbf{Symm}_0 = \mathbf{Z}$ and from the point of view of this component $e$ and $\varepsilon$ identify with the identity on this component.

A graded Hopf algebra is connected if its degree 0 component is of rank 1. It is easy to see that for a connected graded bialgebra there is a unique antipode making it a Hopf algebra.

10.14. *Antipodes on connected graded bialgebra.* Indeed, as $H_0 = \mathbf{Z}$ (or whatever ring once is working over), an antipode must be the identity on $H_0$. Further, if $x$ is primitive, i.e. $\mu_H(x) = 1 \otimes x + x \otimes 1$, the defining requirements for an antipode say that $\iota(x) = -x$. Further, if $i$ is the smallest integer $\geq 1$ such that $H_i \neq 0$ all elements of $H_i$ must be primitive by degree considerations (and the counit properties). Finally if the antipode is known on all the $H_j$ with $j < n$ and $x \in H_n$, $\mathrm{id} \otimes \iota$ and $\iota \otimes \mathrm{id}$ are known on all the terms of $\mu_H(x)$ except $1 \otimes x$ and $x \otimes 1$. The defining requirements for an antipode then immediately give a formula for $\iota(x)$ (two formulas really). They must give the same result because anitpodes (if they exist) are unique (by the same argument that is used that inverses in groups are unique).

10.15. *Primitives of* **Symm**. For every Hopf algebra it is important to know its primitives. That is the elements $x$ that satisfy

$$\mu(x) = 1 \otimes x + x \otimes 1 \tag{10.16}$$

These form a module over the ring over which the Hopf algebra is defined (in this case **Z**) and have a Lie algebra structure under the commutator product $[x, y] = xy - yx$. Because of commutativity this Lie product is zero in the present case.

By definition $p_1 = h_1$ is a primitive. With induction, from the Newton relations (9.57) it follows immediately that all the $p_n$ are primitives. It is not difficult to see that these form a basis over $\mathbf{Z}$ for the module of primitives [39].

[39] First observe that if an element is primitive its different homogeneous summands must all be primitives. next observe that a term in the coproduct of a monomial of length $n$ can only cancel against one from another term of the same length; xontinue by using lexicographic order on exponent sequences. Finally, as $p_n \equiv (-1)^{n+1} h_1^n \mod(h_2, h_3, \cdots)$ (also from the Newton relations with induction), the $p_n$ are not divisible by any natural number other than 1 and as the different $p_n$ are of different degree no sum of them is divisible by any natural number other than 1.

So the comultiplication and counit of **Symm** satisfy

$$\mu_S(p_n) = 1 \otimes p_n + p_n \otimes 1, \ \ \varepsilon_S(p_n) = 0 \qquad (10.17)$$

and as the power sum symmetric functions form a free polynomial basis (over $\mathbf{Q}$) for the symmetric functions this is a perfectly good description of the Hopf algebra **Symm** (once it is known that it is in fact defined over the integers). Formulas (10.17) of course amount to completely the same thing as saying that at the ghost component level addition is coordinate-wise.

10.18. *Liulevicius theorem*, [269, 268]. There are precisely four graded Hopf algebra automorphisms of **Symm**. They are functorially given by

$$\text{identity}, \ a(t) \mapsto a(-t), \ a(t) \mapsto a(t)^{-1}, \ a(t) \mapsto a(-t)^{-1} \qquad (10.19)$$

and they form the Klein 4-group.

The proof is quite simple and it is surprising that this theorem was discovered so late[40]. Over the rationals of course the situation is quite different and the group of automorphisms is quite large viz a countable sum of copies of the multiplicative group of the rationals.

10.20. Here is the proof. The power sum symmetric functions are primitives for the Hopf algebra structure meaning that

$$\mu_S(p_i) = 1 \otimes p_i + p_i \otimes 1 \qquad (10.21)$$

Moreover every homogeneous primitive of weight $i$ is a scalar multiple of $p_i$. Let $\varphi$ be a graded automorphism. Automorphisms of Hopf algebras must take primitives into primitives and as a graded automorphism preserves the grading it must be the case that

$$\varphi(p_i) = ap_i$$

with $a \in \{1,-1\}$ because $\varphi$ is an invertible morphism. The last three automorphisms named in (10.19) take on $p_1, p_2$ respectively the values

$$-p_1, \ p_2; \ \ -p_1, \ -p_2; \ \ p_1, \ -p_2$$

and so composing with a suitable one from the four automorphisms (10.19) it remains to analyze the case that

$$\varphi(p_1) = p_1, \ \varphi(p_2) = p_2$$

Suppose with induction that it has been proved that $\varphi(p_i) = p_i, \ i < n \geq 3$. Look at the Newton relations (see (9.57)

[40] This has most probably to do with the unwarranted and regrettable tendency of mathematicians to think primarily of Hopf algebras, etc. as things over a field of characteristic zero.

$$nh_n = p_n + h_1 p_{n-1} + \cdots + h_{n-1} p_1 \tag{10.22}$$

Because the $h_i$ are polynomials in the $p_1, p_2, \cdots, p_i$ (albeit with rational coefficients), $\varphi(h_i) = h_i,\ i < n$. So applying $\varphi$ to (10.16) and subtracting the result from (10.16) one finds

$$n(h_n - \varphi(h_n)) = p_n \pm p_n$$

which because $n \geq 3$ and $\varphi(h_n)$ must be an integral polynomial of weight $n$, and $p_n$ is not divisible by any integer $> 1$, is only possible if on the right hand side of (10.13) the minus sign applies; i.e. only if $\varphi(p_n) = p_n$.

10.23. *Product comultiplication.* Besides the Abelian group structure on $\Lambda(A)$ and $W(A)$ there is also the product of Witt vectors and power series. These are also functorial and (hence) given by universal polynomials, i.e. by algebra morphisms

$$\mu_P:\ \mathbf{Symm} \longrightarrow \mathbf{Symm}\quad \mathbf{Symm},\quad \varepsilon_P: \mathbf{Symm} \longrightarrow \mathbf{Z}$$

and these give a second bialgebra structure $(\mathbf{Symm},\ m,\ e,\ \mu_P,\ \varepsilon_P)$ on $\mathbf{Symm}$. Because the functorial unit on the rings $\Lambda(A)$ is given by the power series $1 + t + t^2 + t^3 + \cdots$ the co-unit morphism of rings $\varepsilon_P$ is given by

$$\varepsilon_P(h_n) = 1,\quad n = 1, 2, 3, \cdots \tag{10.24}$$

The formula for the morphism governing the multiplication is less easy to describe; certainly explicitely. Various descriptions will be given later in the present section 10. There is no antipode for the second comultiplication of course (otherwise **Symm** would define a field valued functor).

Putting all these structure morphisms together there results the object

$$(\mathbf{Symm},\ m,\ e,\ \mu_S,\ \varepsilon_S,\ \iota_S,\ \mu_P,\ \varepsilon_P)$$

which is a coring object in the category **CRing** of commutative unital rings. One axiom that such an object much satisfy is distributivity on both the left and the right of the second comultiplication over the first in the category **CRing** (where the appropriate codiagonal map is the multiplication on a ring). For distributivity on the left this means that the following diagram must be commutative

$$\begin{array}{ccc}
\mathbf{Symm} & \xrightarrow{\ \mu_S\ } & \mathbf{Symm}\quad \mathbf{Symm} \\
\downarrow \mu_P & & \downarrow \mu_P \quad \mu_P \\
\mathbf{Symm}\quad \mathbf{Symm} & & \mathbf{Symm}\quad \mathbf{Symm}\quad \mathbf{Symm}\quad \mathbf{Symm} \\
\downarrow \mathrm{id} \quad \mu_S & & \downarrow \mathrm{id} \quad \tau \quad \mathrm{id} \\
\mathbf{Symm}\quad \mathbf{Symm}\quad \mathbf{Symm} & \xleftarrow{\ m\ \ \mathrm{id}\ \ \mathrm{id}\ } & \mathbf{Symm}\quad \mathbf{Symm}\quad \mathbf{Symm}\quad \mathbf{Symm}
\end{array}$$

where, as usual $\tau:\ a \quad b \mapsto b \quad a$ is the twist morphism. There is a similar diagram for distributivity on the right.

10. 25. *Various descriptions of the second comultiplication morphism on* **Symm**. At first sight the easiest description is of course in terms of the power sum symmetric functions as

$$\mu_P(p_n) = p_n \quad p_n \tag{10.20}$$

but for a variety of reasons this one is of very limited use. It is of course the definition of the multiplication in terms of ghost components.

Much better are the three descriptions that follow from the three expansions of $\prod (1-\xi_i \eta_j)^{-1}$ given in (9.62), (9.44), (9.45). These give the second comultiplication formulas

$$\mu_P(h_n) = \sum_\lambda h_\lambda \quad m_\lambda \tag{10.21}$$

$$\mu_P(h_n) = \sum_\lambda z_\lambda^{-1} p_\lambda \quad p_\lambda \tag{10.22}$$

$$\mu_P(h_n) = \sum_\lambda s_\lambda \quad s_\lambda \tag{10.23}$$

where in all three formulas the sum is over all partitions of $n$. Explicitly the fist few multiplication polynomials are

$$\begin{aligned} &\mu_P(h_1) = h_1 \quad h_1, \;\; \mu_P(h_2) = 2h_2 \quad h_2 - (h_2 \quad h_1^2 + h_1^2 \quad h_2) + h_1^2 \quad h_1^2 \\ &\mu_P(h_3) = 3h_3 \quad h_3 - 3(h_3 \quad h_2h_1 + h_2h_1 \quad h_3) + (h_3 \quad h_1^3 + h_1^3 \quad h_3) \\ &\qquad -(h_2h_1 \quad h_1^3 + h_1^3 \quad h_2h_1) + h_1^3 \quad h_1^3 \end{aligned} \tag{10.24}$$

10.25. Still other descriptions of the second comultiplication morphism and the (functorial)multiplication on $W(A)$ will be given below, in sections 15 and 11 respectively.

10.26. *Product comultiplication vs inner product.* Theorem 9.42, in particular the expansion formula (9.44), compared with the definition of the mutiplication on the Witt vectors determined by

$$(1-\xi t)^{-1} \quad (1-\eta t)^{-1} = (1-\xi\eta t)^{-1}$$

suggest that there are nontrivial interrelations between the inner product on **Symm** and the second comultiplication (= product comultiplication) on **Symm**. And so there are. For instance cocommutativity of the product comultiplication and symmetry of the inner product are equivalent. At this stage and in this setting this is just an amusing remark. There may be more to it in other contexts.

The inner aproduct is defined by formula (9.41)

$$\langle h_\lambda, m_\mu \rangle = \delta_{\lambda,\mu}$$

Write the monomial symmetric functions as (integer) linear combinations of the $h_\lambda$ :

$$m_\mu = \sum_\lambda a_\mu^\lambda h_\lambda$$

so that

$$\langle m_\mu, m_{\mu'} \rangle = a_\mu^{\mu'} \tag{10.27}$$

Now (by (9.44))

$$\mu_P(h_n) = \sum_{\mathrm{wt}(\lambda)=n} h_\lambda \otimes m_\lambda = \sum a_\lambda^{\lambda'} h_\lambda \otimes h_{\lambda'} \tag{10.28}$$

By cocommutativity of the product comultiplication, also

$$\mu_P(h_n) = \sum_{\mathrm{wt}(\lambda)=n} m_\lambda \otimes h_\lambda = \sum a_\lambda^{\lambda'} h_{\lambda'} \otimes h_\lambda \tag{10.29}$$

Comparison ot the two expressions (10.28) and (10.29) gives, see (10.27),

$$a_\lambda^{\lambda'} = a_{\lambda'}^{\lambda}$$

i.e. symmetry of the inner product. The argument also runs the other way.

## 11. QSymm, the Hopf algebras of quasisymmetric functions and NSymm, the Hopf algebra of noncommutative symmetric functions.

When looking at various universality properties of the Witt vectors and **Symm** (which is the topic of the next section) one rapidly stumbles over a (maximally) non commutative version, **NSymm**, and a (maximally) non cocommutative version, **QSymm**. This section is devoted to a brief discussion of these two objects. Somehow a good many things become easier to see and to formulate in these contexts (including certain explicit calculations). As I have said before, e.g. in [200], p. 56; [199], Ch. H1, p. 1, once one has found the *right* non commutative version, things frequently become more transparent, easier to understand, and much more elegant.

11.1. *The Hopf algebra of non commutative symmetric functions.* Let $Z = \{Z_1, Z_2, Z_3, \cdots\}$ be a countably infinite set of (noncommuting) indeterminates. As an algebra the Hopf algebra of noncommutative symmetric functions is simply the free associative algebra in these indeterminates over the integers

$$\mathbf{NSymm} = \mathbf{Z}\langle Z_1, Z_2, Z_3, \cdots \rangle = \mathbf{Z}\langle Z \rangle \tag{11.2}$$

The coalgebra structure is given on the generators $Z_n$ by

$$\mu(Z_n) = \sum_{i+j=n} Z_i \otimes Z_j = 1 \otimes Z_1 + Z_1 \otimes Z_{n-1} + Z_2 \otimes Z_{n-2} + \cdots + Z_{n-1} \otimes Z_1 + Z_n \otimes 1 \tag{11.3}$$

where the notation $Z_0 = 1$ is used in the expression in the middle. For every word over the natural numbers $\alpha = [a_1, a_2, \cdots, a_m]$, $a_i \in \mathbf{N} = \{1,2,3,\cdots\}$

$$Z_\alpha = Z_{a_1} Z_{a_2} \cdots Z_{a_m} \tag{11.4}$$

denotes the corresponding (non commutative) monomial. This incudes $Z_{[\,]} = 1$. These form a basis of the free Abelian group underlying **NSymm**. The co-unit morphism is given by

$$\varepsilon(Z_n) = 0,\ n = 1,2,3,\cdots,\quad \varepsilon(Z_{[\,]}) = 1 \tag{11.5}$$

Give $Z_n$ degree (or weight) $n$, so that the weight (or degree) of a monomial $Z_\alpha$ is $a_1 + a_2 + \cdots + a_m = \mathrm{wt}(\alpha)$. Then, obviously, **NSymm** is a graded Abelian group

$$\mathbf{NSymm} = \bigoplus_{n\geq 0} \mathbf{NSymm}_n \tag{11.6}$$

with the homogeneous part of degree $n$ spanned by the monomials (11.4) of weight $n$. Then, obviously, $(\mathbf{Symm}, m, e, \mu, \varepsilon)$ is a graded, connected (meaning that $\mathbf{NSymm}_0 = \mathbf{Z}$) bialgebra so that there is (for free, see 10.14) a suitable antipode making **NSymm** a connected, graded, non commutative, cocommutative Hopf algebra.

There is a natural surjective morphism of graded Hopf algebras

$$\mathbf{NSymm} \longrightarrow \mathbf{Symm},\quad Z_n \mapsto h_n \tag{11.7}$$

that exhibits **Symm** as the maximal commutative quotient Hopf algebra of **NSymm**. The kernel of (11.7) is the commutator ideal generated by all elements of the form $Z_i Z_j = Z_j Z_i,\ i, j \in \mathbf{N}$.

The first paper to study **NSymm** in depth was probably [162]. It was immediately followed by a slew of other publications [54, 129, 130, 128, 216, 214, 219, 249, 251, 250, 386, 384, 200, 195, 198, 193, 201][41].

11.8. *Divided power sequences. Curves.* Given a Hopf algebra $H$ a divided power sequence, or curve, in $H$ is by definition a sequence of elements

$$\gamma = (d_0 = 1, d_1, d_2, \cdots) \tag{11.9}$$

such that

$$\mu_H(d_n) = \sum_{i+j=n} d_i \otimes d_j \tag{11.10}$$

If $H$ is the covariant bialgebra of a formal group $F$ this translates into a mapping of the formal affine line into the formal scheme (variety) underlying $F$, whence the terminology curve [42]. A

[41] Many of these papers are about noncommutative symmetric functions over a field of characteristic zero; here, of course, it is **NSymm** over the integers that is really important (here and elsewhere).

[42] There is also another notion in algebra that goes by the name 'divided powers' viz the presence of a sequence of operators (usually on an ideal in a commutative ring) that behave like sending an element $x$ to $(n!)^{-1}x^n$. Hence the term curve, rather than divided power sequence (DPS) is preferred for the notion defined by (11.9), (11.10).

curve is often denoted by its generating power series

$$\gamma(t) = 1 + d_1 t + d_2 t^2 + d_3 t^3 + \cdots \tag{11.11}$$

Two such power series can be multiplied and the result is again a curve i.e a power series in $t$ of which the coefficients satsify (11.10). This defines a functor

$$Curve\colon \mathbf{Hopf} \longrightarrow \mathbf{Group} \tag{11.12}$$

from the category of Hopf algebras to the category of groups. An example of a curve is the universal curve

$$Z(t) = 1 + Z_1 t + Z_2 t^2 + \cdots \quad Curve(\mathbf{NSymm}) \tag{11.13}$$

and it is immediate from the defining property (11.10) and the freeness of $\mathbf{NSymm}$ as an associative algebra that

11.14. *Theorem* (*universality property of* **NSymm**). For every curve $\gamma(t)$ in a Hopf algebra $H$ there is a unique morphism of Hopf algebras $\mathbf{NSymm} \longrightarrow H$ that takes the universal curve $Z(t)$ to the given curve $\gamma(t)$.

That is, the functor $Curve$ is represented by **NSymm** [43]. The group operation on curves corresponds to the convolution product on $\mathbf{Hopf}(\mathbf{NSymm}, H)$. When the Hopf algebra $H$ is commutative this theorem translates into

11.15. *Theorem.* (*Second universality property of the Witt vectors*). Let $\gamma(t)$ be a curve in a commutative Hopf algebra $H$. Then there is a unique morphism of Hopf algebras $\mathbf{Symm} \longrightarrow H$ that takes the curve

$$h(t) = 1 + h_1 t + h_2 t^2 + \cdots \qquad Curve(\mathbf{Symm}) \tag{11.16}$$

into the given curve $\gamma(t)$.

In the context of commutative formal groups with $H$ the covariant bialgebra of a commutative formal group this is known as Cartier's first theorem.

11.17. *The Hopf algebra* **QSymm** *of quasi-symmetric functions.* The next sections deal with another generalization of **Symm**, dual to **NSymm** and containing **Symm**.

11.18. *Compositions and partitions.* A composition is a word (of finite length) over the natural numbers. Here a composition will be written

$$\alpha = [a_1, a_2, \cdots, a_m], \;\; a_i \quad \mathbf{N} = \{1,2,3,\cdots\} \tag{11.19}$$

---

[43] See also [149].

A composition (11.19) is said to have length $m$ and weight $\mathrm{wt}(\alpha) = a_1 + a_2 + \cdots + a_m$. Associated to a composition $\alpha$ is a partition $\lambda$ which consists of the $a_j$ arranged in nonincreasing order. Or, equivalently, $\lambda = (1^{m_1(\alpha)} 2^{m_2(\alpha)} 3^{m_3(\alpha)} \cdots)$ where $m_j(\alpha)$ is the number of entries of $\alpha$ equal to $j$.

11.19. *Monomial quasi-symmetric functions.* Consider again the ring of polynomials, i.e. power series of bounded degree[44] in a countable infinity of commuting indeterminates, i.e. an element of $\mathbf{Z}[\xi_1, \xi_2, \xi_3, \cdots]$. Such a polynomial

$$\sum c_{\alpha, i_1, \cdots, i_m} \xi_{i_1}^{a_1} \xi_{i_2}^{a_2} \cdots \xi_{i_m}^{a_m} \tag{11.20}$$

where the sum runs over runs over all compositions $\alpha = [a_1, \cdots, a_m]$ and all index sequences $i_1 < \cdots < i_m$ of the same length as $\alpha$, is said to be quasi-symmetric if

$$c_{\alpha, i_1, \cdots, i_m} = c_{\alpha, j_1, \cdots, j_m} \quad \text{for all index sequences } i_1 < \cdots < i_m \text{ and } j_1 < \cdots < j_m \tag{11.21}$$

For a given composition $\alpha = [a_1, \cdots, a_m]$ the associated monomial quasi-symmetric function is

$$[a_1, \cdots, a_m] = \sum_{i_1 < i_2 < \cdots < i_m} \xi_{i_1}^{a_1} \xi_{i_2}^{a_2} \cdots \xi_{i_m}^{a_m} \tag{11.22}$$

It is denoted by the same symbol. For example in four variables

$$\begin{aligned}
[1,2] &= \xi_1\xi_2^2 + \xi_1\xi_3^2 + \xi_1\xi_4^2 + \xi_2\xi_3^2 + \xi_2\xi_4^2 + \xi_3\xi_4^2 \\
[1,2,1] &= \xi_1\xi_2^2\xi_3 + \xi_1\xi_2^2\xi_4 + \xi_1\xi_3^2\xi_4 + \xi_2\xi_3^2\xi_4 \\
[3,2,1] &= \xi_1^3\xi_2^2\xi_3 + \xi_1^3\xi_2^2\xi_4 + \xi_1^3\xi_3^2\xi_4 + \xi_2^3\xi_3^2\xi_4
\end{aligned} \tag{11.23}$$

Note that the monomial symmetric function $m_\lambda$ defined by a partition $\lambda$ is the sum of all *distinct* monomial quasi-symmetric functions $\alpha$ whose associated partition is $\lambda$. For example

$$\begin{aligned}
&m_{(2,1)} = [2,1] + [1,2], \quad m_{(2,1,1)} = [2,1,1] + [1,2,1] + [1,1,2] \\
&m_{(3,2,1)} = [3,2,1] + [3,1,2] + [2,3,1] + [2,1,3] + [1,3,2] + [1,2,3]
\end{aligned} \tag{11.24}$$

The monomial quasi-symmetric functions form a free basis of the free Abelian group **QSymm** of all quasi-symmetric functions. This group is graded by assigning to the monomial quasi-symmetric function $\alpha = [a_1, \cdots, a_m]$ the weight $\mathrm{wt}(\alpha)$. Sum and product of quasi-symmetric functions are again quasi-symmetric, so that **QSymm** is a ring.

The ring **Symm** is a subring with the inbedding

$$h_n = \sum_{\mathrm{wt}(\alpha) = n} \alpha \tag{11.25}$$

11.26. *Overlapping shuffle product.* It is useful to have an explicit direct recipe in terms of compositions for the multiplication of monomial quasi-symmetric functions. This is given by the overlapping shuffle product.

[44] See the appendix.

Let $\alpha = [a_1, a_2, \cdots, a_m]$ and $\beta = [b_1, b_2, \cdots, b_n]$ be two compositions or words. Take a 'sofar empty' word with $n+m-r$ slots where $r$ is an integer between $0$ and $\min\{m,n\}$, $0 \le r \le \min\{m,n\}$.

Choose $m$ of the available $n+m-r$ slots and place in it the natural numbers from $\alpha$ in their original order; choose $r$ of the now filled places; together with the remaining $n+m-r-m = n-r$ empty places these form $n$ slots; in these place the entries from $\beta$ in their orginal order; finally, for those slots which have two entries, add them. The product of the two words $\alpha$ and $\beta$ is the sum (with multiplicities) of all words that can be so obtained. So, for instance

$$
\begin{aligned}
[a,b][c,d] = {} & [a,b,c,d] + [a,c,b,d] + [a,c,d,b] + [c,a,b,d] + [c,a,d,b] + [c,d,a,b] + \\
& + [a+c,b,d] + [a+c,d,b] + [c,a+d,b] + [a,b+c,d] + [a,c,b+d] + \\
& + [c,a,b+d] + [a+c,b+d]
\end{aligned}
\tag{11.27}
$$

$$
[1][1][1]=6[1,1,1]+3[1,2]+3[2,1]+[3] \tag{11.28}
$$

It is easy to see that the recipe given above gives precisely the multiplication of (the corresponding basis) quasi-symmetric functions. The shuffles with no overlap of $a_1, \cdots, a_m; b_1, \cdots, b_n$ correspond to the products of the monomials that have no $\xi_j$ in common; the other terms arise when one or more of the $\xi$ 's in the monomials making up $\alpha$ and $\beta$ do coincide. In example (11.27) the first six terms are the *shuffles*; the other terms are '*overlapping shuffles*'.

The term shuffle comes from the familiar rifle shuffle of cardplaying; an overlapping shuffle occurs when one or more cards from each deck don't slide along each other but stick edgewise together; then their values are added.

11.29. *Hopf algebra structure on* **QSymm**. Now introduce on **QSymm** the comultiplication 'cut':

$$
\begin{aligned}
\mu_{\mathbf{QSymm}}([a_1, a_2, \cdots, a_m]) = {} & 1 \quad [a_1, a_2, \cdots, a_m] \\
& + \sum_{i=1}^{n-1} [a_1, \cdots, a_i] \quad [a_{i+1}, \cdots a_m] + [a_1, a_2, \cdots, a_m] \quad 1
\end{aligned}
\tag{11.30}
$$

and the counit

$$
\varepsilon_{QSymm}(\alpha) = \begin{cases} 1 & if \quad \lg(\alpha) = 0 \\ 0 & if \quad \lg(\alpha) \ge 1 \end{cases} \tag{11.31}
$$

Here recall that the only composition of length zero is the empty composition $1 = [\ ]$ which is the unit for the overlapping shuffle multiplication.

Write $m_{\mathbf{QSymm}}$ for the overlapping shuffle multiplication and $e_{\mathbf{QSymm}}$ for the unit morphism determined by $1 \mapsto [\ ] = 1$. Then

11.32. *Theorem.* The structure $(\mathbf{QSymm}, m_{\mathbf{QSymm}}, e_{\mathbf{QSymm}}, \mu_{\mathbf{QSymm}}, \varepsilon_{\mathbf{QSymm}})$ is a commutative non cocommutative graded connected bialgebra and hence defines a Hopf algebra.

The antipode again comes for free once the the bialgebra statement has been proved because **QSymm** is connected graded. The only thing in proving the bialgebra statement which is not dead easy is the verification that the comultiplication is an algebra morphism for the overlapping shuffle product. Even that can be avoided by using graded duality as described below.

11.33. *Duality between* **NSymm** *and* **QSymm**. Introduce the nondegenerate graded pairing

$$\langle\ ,\ \rangle:\ \mathbf{NSymm}\times\mathbf{QSymm}\longrightarrow\mathbf{Z},\quad \langle Z_\alpha,\beta\rangle=\delta_{\alpha,\beta} \tag{11.34}$$

With respect to this pairing the multiplication (resp. comultiplication) on **NSymm** is dual to the comultiplication (resp. multiplication) on **QSymm**. Similarly the counit (resp. unit) on **NSymm** is dual to the unit (resp. counit) on **QSymm**. These statements amount to the following four formulas.

$$\begin{aligned}\langle m_{\mathbf{NSymm}}(Z_\alpha\otimes Z_\beta),\gamma\rangle&=\langle Z_\alpha\otimes Z_\beta,\mu_{\mathbf{QSymm}}(\gamma)\rangle\\ \langle\mu_{\mathbf{NSymm}}(Z_\alpha),\beta\otimes\gamma\rangle&=\langle Z_\alpha,m_{\mathbf{QSymm}}(\beta\otimes\gamma)\rangle\\ \langle e_{\mathbf{NSymm}}(1),\beta\rangle&=\langle 1,\varepsilon_{\mathbf{QSymm}}(\beta)\rangle\\ \langle\varepsilon_{\mathbf{NSymm}}(\alpha),1\rangle&=\langle\alpha,e_{\mathbf{QSymm}}(1)\rangle\end{aligned} \tag{11.35}$$

where $\langle\ ,\ \rangle$ on **Z** is the innerproduct $\langle r,s\rangle=rs$. Formulas (11.35) easily imply that if **NSymm** is a bialgebra, so is **QSymm** (and vice-versa). The antipodes on the two Hopf algebras are (as a consequence) also dual to each other:

$$\langle\iota_{\mathbf{NSymm}}(Z_\alpha),\beta\rangle=\langle Z_\alpha,\iota_{\mathbf{QSymm}}(\beta)\rangle$$

11.36. *Primitives of* **QSymm**. It is an easy exercise to show that the primitives for the first comultiplication, the 'sum comultiplication' , i.e. the elements for which

$$\mu_S(x)=1\otimes x+x\otimes 1$$

are precisely the linear combinations of the words of length 1. [45]

11.37. *Autoduality of* **Symm**. Under (graded) duality a sub Hopf algebra of of a Hopf algebra $H$ corresponds to a quotient Hopf algebra of the graded dual. To avoid notational confusion write temporarily **Symm'** for the algebra of symmetric functions as a subalgebra of **QSymm** and retain the notation $\mathbf{Symm}=\mathbf{Z}[h]$ for the algebra of symmetric functions as a quotient algebra of **NSymm**. The Hopf ideal $J$ in **NSymm** such that $\mathbf{NSymm}/J=\mathbf{Symm}$ is the commutator ideal spanned by all elements of the form $Z_\alpha(Z_iZ_j-Z_jZ_i)Z_\beta$. It is easy to check that a quasisymmetric function has pairing 0 (under (11.34) with all these elements if and only if it is symmetric. Thus the sub Hopf algebra of **QSymm** corresponding to the quotient **Symm** of **NSymm** is **Symm'**. Also the induced pairing is

$$\langle h_\lambda,m_\kappa\rangle=\delta_{\lambda,\kappa} \tag{11.38}$$

[45] The situation for the primitives of **NSymm** is very different. Over the rationals the primitives form the free Lie algebra in countably many generators. Over the integers the description is still more involved and was only recently found; see [Hazewinkel, 2001 #377].

which is the inner product on **Symm** as defined in 9.40. [46]

Moreover as algebras **Symm** and **Symm’** are isomorphic under the isomorphism $h_n \mapsto e_n$. To prove that **Symm** and **Symm’** are isomorphic as Hopf algebras, it remains to check that this isomorphism takes the given comultiplication on **Symm**, $h_n \mapsto \sum_{i+j=n} h_i \quad h_j$ into the comultiplication induced on **Symm’** by ‘cut’ the comultiplication on **QSymm**. But $e_n = [1,1,1,\cdots,1]$ as an element of **Symm’** **QSymm** and so ‘cut = $\mu_{\mathbf{QSymm}}$ on $e_n$ takes the value

$$\mu_{QSymm}(e_n) = \sum_{i+j=n} e_i \quad e_j$$

and things fit perfectly.

11.39. *The second comultiplication on* **QSymm**. There is a second comultiplication (morphism) on **QSymm** that is distributive over the first one on the right but not on the left. Here is the definition. Let $\alpha = [a_1, a_2, \cdots a_m]$ be a composition. A $(0,\alpha)$-matrix is a matrix whose entries are either zero or one of the $a_i$, which has no zero columns or zero rows, and in which the entries $a_1, a_2, \cdots, a_m$ occur in their original order if one orders the entries of a matrix by first going left to right through the first row, then left to right through the second row, etc.

For a matrix $M$ let $column(M)$ be the vector of column sums and $row(M)$ the vector of row sums. For instance

$$M = \begin{pmatrix} 0 & 1 & 0 & 3 \\ 1 & 2 & 0 & 1 \\ 0 & 0 & 1 & 0 \end{pmatrix} \tag{11.40}$$

is a $(0,[1,3,1,2,1,1])$- matrix with $column(M) = [1,3,1,4]$ and $row(M) = [4,4,1]$.

The second comultiplication on $QSymm$ is now given by

$$\mu_P(\alpha) = \sum_{\substack{M \text{ is a} \\ (0,\alpha)\text{-matrix}}} row(M) \quad column(M) \tag{11.41}$$

Restricted to **Symm** **QSymm**, this describes the second comultiplication on **Symm** in combinatorial terms. To specify it, it suffices to do this on the elementary symmetric pooynomials (because it is a ring morphism). I.e to use the recipe for compositions of the form $\alpha = [1,1,1,\cdots,1]$. For these compositions the matrices involved are what are usually called (0,1)-matrices.

Probably the most famous theorem on (0,1)-matrices is the Gale-Ryser theorem. Thus it would appear likely that there is some connection between Witt vectors (via symmetric functions) and the Gale Ryser theorem, [150, 345]. And indeed there is (via the additional link of the symmetric functions with the representations of the symmetric groups). It is the Snapper Liebler-Vitale Lam theorem, [196, 254]. This will be described in a little more detail in subsection 18.33 below.

---

[46] This is one main reason for working with the complete symmetric functions rather than with the elementary ones.

11.42. *Distributivity properties of* $\mu_P$ *over* $\mu_S$. The second comultiplication morphism $\mu_P$ of (11.39) is distributive on the right over the first comultiplication morphism $\mu_S$ but not on the left.

Quite generally let $(H, m, e, \mu, \varepsilon)$ be a bialgebra equipped with a second comultiplication morphism of rings $\mu_P : H \longrightarrow H \otimes H$. For this second comulltiplication to be distributive on the right over the first one (in the category of rings) one needs the following diagram (11.41) to commute. (This is precisely what is needed to ensure that the functor represented by $H$, $A \mapsto$ **CRing**$(H, A)$ have as values groups (denoted additively even when they are noncommutative) (as in the case at hand)) equipped with an additional multiplication such that $(a+b)c = ab + bc$.)

$$\begin{array}{ccc} H & \xrightarrow{\ \mu\ } & H \otimes H \\ \downarrow \mu_P & & \downarrow \mu_P \otimes \mu_P \\ H \otimes H & & H \otimes H \otimes H \otimes H \\ \downarrow \mu \otimes \mathrm{id} & & \downarrow \mathrm{id} \otimes \mathrm{tw} \otimes \mathrm{id} \\ H \otimes H \otimes H & \xleftarrow{\ \mathrm{id} \otimes \mathrm{id} \otimes m\ } & H \otimes H \otimes H \otimes H \end{array} \tag{11.43}$$

Here is a proof that this diagram is commutative in the case at hand. Take a composition $\alpha$ and let $M$ be a $(0,\alpha)$-matrix. Applying $\mu \otimes \mathrm{id}$ to $row(M) \otimes column(M)$ is the same thing as cutting the matrix M horizontally to obtain two blocks $M_1, M_2$ with the same number of columns stacked on top of each other. All terms in the lower left hand corner of diagram (11.43) coming from $\alpha$ thus are of the form $row(M_1) \otimes row(M_2) \otimes column(M)$. Now in $\mu(\alpha)$ (top horizontal arrow in the diagram) consider the term $\alpha_1 \otimes \alpha_2$ for which the length of $\alpha_1$ is the same as the number of rows of $M_1$. Let $M_1', M_2'$ be the matrices obtaine from $M_1, M_2$ by removing zero colums. These are respectively a $(0,\alpha_1)$-matrix and a $(0,\alpha_2)$-matrix. The term obtained in the lower right hand corner is now

$$\begin{aligned} & row(M_1') \otimes row(M_2') \otimes column(M_1') \otimes column(M_2') \\ & = row(M_1) \otimes row(M_2) \otimes column(M_1') \otimes column(M_2') \end{aligned} \tag{11.44}$$

Now look at the overlapping shuffle product of two compositions $\beta_1, \beta_2$ in a slightly different way. Intersperse both compositions with zeros in any way to obtain vectors of the same length $n$ and such that the two zero-interspersed compositions when put on top of each other give a $2 \times n$ matrix without zero columns. Take the column sum vector of this $2 \times n$ matrix. Doing this in all posible ways gives the overlapping shuffle product (obviously).

Now, all matrices $M$ which yield given matrices $M_1, M_2$ are obtained from $M_1', M_2'$ by interspersing them with zero columns to have the same number of columns and such that when stacked on top of each other they give a matrix without zero columns.

Combining these two remarks gives that for a given $M_1', M_2'$ the possible originating $M$ precisely give the compositions $column(M)$ that arise as terms of the shuffle product of $column(M_1')$ and $column(M_2')$, proving the commutativity of (11.43).

For left distributivity a similar diagram must be commutative. it is obtained from the one at hand by replacing the lower left hand morphism by $\mathrm{id} \otimes \mu$ and the bottom morphism by $m \otimes \mathrm{id} \otimes \mathrm{id}$. It is a trivial matter to check that left distributivity does not hold; it already fails for any $\alpha = [a_1, a_2]$ with $a_1 \neq a_2$.

If one switches 'row' and 'column' in (11.37) there results a second comultiplication morphism that is left distributive but not right distributive.

One easily sees that on a symmetric sum of quasisymmetric monomials $\mu_P$ is commutative. Indeed if $M$ is a $(0,\alpha)-$matrix then $row(M^{\mathrm{tr}}) = column(M)$ and $column(M^{\mathrm{tr}}) = row(M)$ and the transpose $M^{\mathrm{tr}}$ is a $(0,\alpha')$-matrix for some permutation $\alpha'$ of $\alpha$. Given right distributivity of this second comultiplication on **QSymm** to check that it does indeed give the functorial multiplication of Witt vectors or power series it suffies to check things on power series of the form $1+at$ which is a triviality.

11.45. *Co-unit for the second comultiplication on* **QSymm**. There is also a counit morphism $\varepsilon_P$ for $\mu_P$. It is given by

$$\varepsilon_P([a_1,a_2,\cdots,a_m]) = \begin{cases} 1 & \text{for} \quad m \le 1 \\ 0 & \text{for} \quad m \ge 2 \end{cases} \tag{11.46}$$

To be a counit for $\mu_P$ it must satisfy the condition that the following composition of morphisms is the identity

$$\mathbf{QSymm} \xrightarrow{\mu_P} \mathbf{QSymm} \quad \mathbf{QSymm} \xrightarrow{\varepsilon_P \quad \mathrm{id}} \mathbf{QSymm} \quad \mathbf{QSymm} \xrightarrow{m} \mathbf{QSymm}$$

and also the composition of morphism obtaine from this by switching $\varepsilon_P$ and id for the middle arrow.

Let $\alpha = [a_1,a_2,\cdots,a_m]$. As $\varepsilon_P$ is zero on anything of length 2 or more the only $(0,\alpha)$-matrices that need to be considered have only one row. There is just one such matrix, viz $(a_1,a_2,\cdots,a_m)$. This gives that the composed morphism in question is indeed the identity.

11.47. *Second multiplication on* **NSymm**. Dually the second comultiplication on **QSymm** give a second multiplication morphism (of coalgebras) on **NSymm**. Denoting this one by , there is the right distributivity formula (which follows from the right distributivity of $\mu_P$ over $\mu_S$ in **QSymm**).

$$(Z_\alpha Z_\beta) \quad Z_\gamma = \sum_i (Z_\alpha \quad Z_{\gamma_i'})(Z_\beta \quad Z_{\gamma_i''}) \text{ where } \mu(Z_\gamma) = \sum_i Z_{\gamma_i'} \quad Z_{\gamma_i''}$$

This formula is also in [162]. The dual of the recipe for the second comultiplication given in (11.39) is due to Solomon, [373].

11.48. *Unit for the second multiplication on* **NSymm**. The second multiplication, $m_p$ on **NSymm** defines a multiplication on each homogeneous component $\mathbf{NSymm}_n$. The unit for this second multiplication is the element $Z_n$. So the 'unit' for the second comultiplication on all of **NSymm** is the infinite sum

$$1 + Z_1 + Z_2 + Z_3 + \cdots$$

which does not live in **NSymm** but only in a suitable completion.

The same situation holds for **Symm**. Here, for each $n$ the second multiplication, sometimes called inner multiplication, turns each homogenous summand $\mathbf{Symm}_n$ into a commutative ring with unit element $h_n$, and the sum of all these elements, which does not live in **Symm** itself is the unit element for the second multiplication on all of **Symm**.

So the autoduality of **Symm** is not perfect at this level. The second multiplication and second comultiplication are nicely dual to each other but the second unit only exists in a suitable completion while the second counit is present for **Symm** itself.

This is understandable because the second multiplication and second comultiplication do not respect the grading.

## 12. Free, cofree and duality properties of Symm.

This short section concerns ways of obtaining **Symm** etc. by means of free and cofree constructions and consequences thereof. And open questions concerning these matters. All algebras and coalgebras in this ection will come with a unit resp. counit element seen as.morphisms from and to **Z**.

12.1. *Free algebras*. To set the stage consier first the well known case of free algebras. Let $M$ be an Abelian group. A free algebra on $M$ is an algebra $\mathrm{Free}(M)$ together with a morphism of Abelian groups $i_M\colon M \longrightarrow \mathrm{Free}(M)$ such that the following universality property holds

$$\begin{array}{l}\text{For every morphism of Abelian groups } \varphi\colon M \longrightarrow A \text{ to a ring } A \\ \text{there is a unique morphism of rings } \tilde{\varphi}\colon \mathrm{Free}(M) \longrightarrow A \text{ such that } \tilde{\varphi} i_M = \varphi\end{array} \qquad (12.2)$$

This defines a functor Free: **AbGroup** $\longrightarrow$ **Ring** that is left adjoint to the (forgetful) functor $U$: **Ring** $\longrightarrow$ **AbGroup** that forgets about the multipllication. Pictorially, the situation looks as follows.

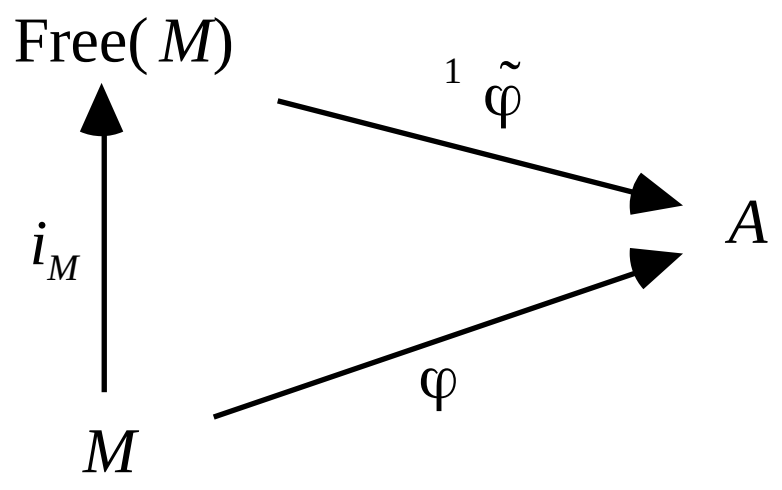

(12.3)

The symbol group 1 in the diagram above means ‘there exists a unique’.
Such a functor ‘‘Free’ exists. It is given by the tensor algebra over $M$.

12.4. *Tensor algebra and tensor coalgebra.* For an Abelian group $M$ let $T^iM = M^{\otimes i} = M \otimes M \otimes \cdots \otimes M$ (and $T^0M = \mathbf{Z}$) be the $i$-fold tensor product. There are natural isomorphisms $\psi_{n,m}\colon T^nM \otimes T^nM \longrightarrow T^{n+m}M$. Now consider

$$\begin{aligned} &TM = \bigoplus_{n=0}^{\infty} T^nM, \\ &i_M\colon M \longrightarrow TM,\quad p_M\colon TM \longrightarrow M, \\ &e_M\colon \mathbf{Z} \longrightarrow TM,\quad \varepsilon_M\colon TM \longrightarrow \mathbf{Z} \end{aligned} \qquad (12.5)$$

where the morphisms $i_M, p_M$ are obtained by identifying $M$ and $T^1M$ and $e_M, \varepsilon_M$ by identifying $\mathbf{Z}$ and $T^0M$. There is a ring structure on $TM$ given by the

$$T^mM \otimes T^nM \xrightarrow{\psi} T^{n+m}M \quad \text{and } i_M \text{ as unit morphism.} \tag{12.6}$$

This is the tensor algebra and it is the free algebra on $M$ with respect to the (canonical) morphism $i_M$.

There is also a coalgebra structure on $TM$ given by

$$(\psi^{-1}_{0,n},\psi^{-1}_{1,n-1},\psi^{-1}_{2,n-2},\cdots,\psi^{-1}_{n-1,1},\psi^{-1}_{n,0})\colon\ T^nM \longrightarrow \bigoplus_{i=0}^{n} T^iM \otimes T^{n-i}M \quad (12.7$$
$$\text{and } \varepsilon_M \text{ as counit morphism}$$

The two structures do not combine to give a bialgebra structure. Far from it.

12.8. *The graded case.* In the present context the important case is when everything in sight is graded and $M$ is a free Abelian group.

So let $M$ be a positively graded Abelian group

$$M = \bigoplus_{i=1}^{\infty} M_i \tag{12.9}$$

Then $TM$ is graded by giving a pure tensor $x_1 \otimes \cdots \otimes x_m$, $x_i \in M_{i_1}$ degree $i_1 + \cdots i_m$. If no grading on $M$ is specified it is treated as graded with all nonzero elements of degree 1. If now $M$ is moreover free with homogeneous basis

$$\{T_u\colon\ u \in U\},\ \ \text{degree}(T_u) = r_u \tag{12.10}$$

$\text{Free}(M)$ is $\mathbf{Z}\langle T_u\colon\ u \in I\rangle$. For instance if the index set is $\mathbf{N}$ and the degree of $T_i$ is $i$ $\text{Free}(M) = \mathbf{NSymm}$ (as far as the underlying graded algebra is concerned).

There is also a commutative version with $\text{Free}_{\text{comm}}(M)$ being the maximal commutative qoutient of $\text{Free}(M)$. So if $M = \bigoplus_{i=1}^{\infty} \mathbf{Z}h_i$, $\text{degree}(h_i) = i$, $\text{Free}_{\text{comm}}(M) = \mathbf{Symm}$.

12.11. *Cofree coalgebras.* All coalgebras in the following are supposed to come with a co-unit and morphisms of coalgebras are supposed to be compatible with the co-units involved. Given an Abelian group $M$ the cofree coalgebra over $M$ (if it exists) is a coalgebra $\text{CoFree}(M)$ (over the integers) together with a morphism of Abelian groups $\pi_M\colon \text{CoFree}(M) \longrightarrow M$ with the following universal property:

For every coalgebra $C$ together with a morphism of Abelian groups
$\varphi\colon\ C \longrightarrow M$ there is a unique morphism of coalgebras (12.12)
$\tilde{\varphi}\colon\ C \longrightarrow \text{CoFree}(M)$ such that $\pi_M\tilde{\varphi} = \varphi$

Pictorially this looks as follows:

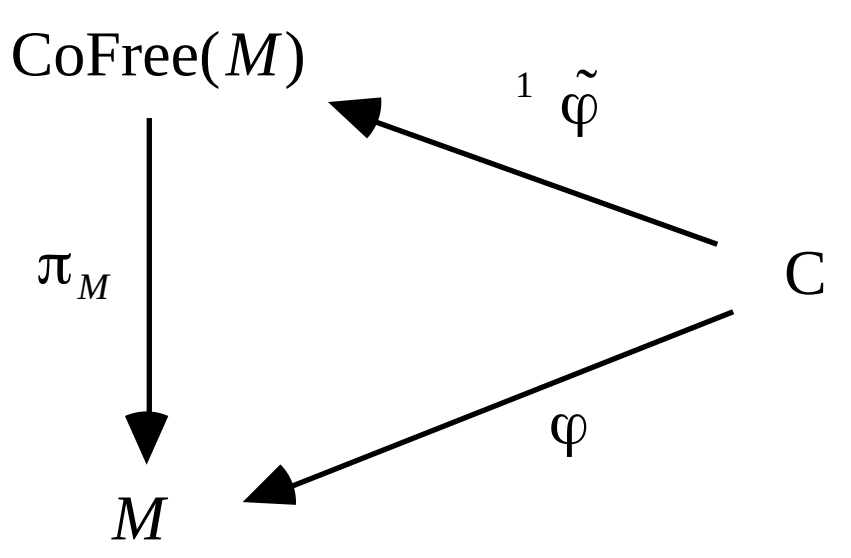

(12.13)

a picture that is completely dual to the one, see (12.3), for a free algebra over a module. The universality property (12.12), (12.13) says, given existence, that ‘CoFree’ is right adjoint to the forgetful functor which assigns to a coalgebra the underlying Abelian group. Whether the cofree coalgebra over[47] an Abelian group always exists is unkown. It does do so in the case of free Abelian groups, their duals[48] and certain other cases, see [194].

Let $M$ be (graded) free with (homogeneous) basis $\{T_u: u \in I\}$; and let $C$ be a coalgebra together with a morphism of Aberlian groups $C \xrightarrow{\varphi} M$. As an Abelian group the tensor coalgebra consists of all noncommutative polynomials in the $T_u$. This is not always large enough to serve as receiving object for the unique morphism $\tilde{\varphi}$ covering $\varphi$ as required in (12.12), (12.13). Consider for example the coalgebra $C = \mathbf{Z} \oplus \mathbf{Z}x$, $\mu(x) = x \otimes x$, $\varepsilon(x) = 1$ together with the morphism of Abelian groups that sends $\mathbf{Z}$ to zero and $x$ to one of the basis elements $T$ of $M$. It is obvious that the $n$-th componet of the covering morphism, $\tilde{\varphi}_n\colon C \longrightarrow T^i M$, must (always) be given by $\tilde{\varphi}_n = \varphi^{\otimes n}\mu_n$ where $\mu_n$ is the iterated comultiplication

$$\mu_0 = \varepsilon,\ \mu_1 = \mathrm{id},\ \mu_2 = \mu, \cdots,\ \mu_{n+1} = (\mu \otimes \mathrm{id} \otimes \cdots \otimes \mathrm{id})\mu_n,\ \cdots \qquad (12.14)$$

Thus in the example at hand $\tilde{\varphi}_n(x) = T \otimes T \otimes \cdots \otimes T$ ($n$ factors), so that the $\tilde{\varphi}(x)$ ‘is’ the power series $1 + T + T^2 + \cdots + T^n + \cdots$ which is not in $TM$. To obtain the free coalgebra on $M$ consider the completion $\hat{T}M$ of all formal noncommutative power series in the indeterminates $T_u$, $u \in I$. The free coalgebra over $M$ is now something in between $TM$ and $\hat{T}M$ of which the elements can be called for good (but varying) reasons ‘rational noncommutative power series’, ‘representative noncommutative power series’, ‘Schützenberger-recognizable noncommutative power series’, or ‘recursive noncommutative power series’, loc. cit.

Things become a good deal easier in the graded case. In this case, i.e. for graded morphisms of graded coalgebras into graded Abelian groups, the tensor coalgebra is the cofree coalgebra. To sees this first observe that a morphism of graded Abelian groups $C \longrightarrow M$ takes the degree 0 component of $C$ to zero. Next for any homogeneous $x \in C$ of degree $n$, say, all terms in its iterated coproduct $\mu_m(x)$ for $m > n$ must contain a factor from $C_0$. Hence $\tilde{\varphi}_m(x) = \varphi^{\otimes m}\mu_m(x) = 0$ for $m > n$ and things are OK.

There is also a cocommutative version where $TM$ is replaced by its maximal cocommutative sub coalgebra of all symmetric tensors.

[47] It pays to be linguistically careful here: free algebra **on** an Abelian group; cofree coalgebra **over** an Abelian group.

[48] These are not free in the infinite rank case.

In the case of graded free Abelian groups the graded free and graded cofree constructions are graded dual in that there is a natural isomorphism

$$\text{CoFree}(M)^{dual} = \text{Free}(M^{dual})$$

For the nongraded case there is also a duality but things are more complicated [49].

12.15. *Inheritance of structure.* Now let $C$ be a coalgebra and apply the free algebra construction to the augmentation module $I = \text{Ker}(\varepsilon)$. This gives (in any case) an algebra which will still be denoted Free($C$). It inherits a coalgebra structure as follows.

By functoriality the coalgebra morphism gives an algebra morphism

$$\text{Free}(C) \xrightarrow{\text{Free}(\mu)} \text{Free}(C \otimes C) \qquad (12,16)$$

Further there is the natural morphism $i_C \otimes i_C : C \otimes C \longrightarrow \text{Free}(C) \otimes \text{Free}(C)$ and by the freeness property this induces a morphism of algebras $\text{Free}(C \otimes C) \longrightarrow \text{Free}(C) \otimes \text{Free}(C)$ which when composed with (12.16) gives the desired coalgebra structure morphism

$$\text{Free}(C) \xrightarrow{\mu_{\text{Free}(C)}} \text{Free}(C) \otimes \text{Free}(C) \qquad (12.17)$$

By construction it is an algebra morphism and so Free($C$) becomes a bialgebra.

Dually, consider an algebra $A$ and apply the cofree construction to the quotient module $A/\mathbf{Z}$. Denote the result by CoFree($A$). This one inherits an algebra structure in a similar way. Applying CoFree to the multiplication morphism gives a coalgebra morphism

$$\text{CoFree}(A \otimes A) \longrightarrow \text{CoFree}(A) \qquad (12.18)$$

On the other hand there is a fairly obvious (canonical) morphism of Abelian groups $\text{CoFree}(A) \otimes \text{CoFree}(A) \longrightarrow (A \otimes A)/\mathbf{Z}$ that by the cofreeness property induces a morphism of coalgebras

$$\text{CoFree}(A) \otimes \text{CoFree}(A) \rightarrow \text{CoFree}(A \otimes A) \qquad (12.19)$$

which when composed with with (12.18) gives the desired multiplication morphism which is by construction a coalgebra morphism so that the result is again a bialgebra.

12.20. **Symm**, **NSymm** *and* **QSymm** *again*. The cofree coalgebra over the integers is

$$\text{CoFree}(\mathbf{Z}) = \mathbf{Z} \oplus \mathbf{Z}Z_1 \oplus \mathbf{Z}Z_2 \oplus \cdots$$
$$\mu(Z_n) = \sum_{i=0}^{n} Z_i \otimes Z_{n-i}, \quad \text{where } Z_0 = 1, \; \varepsilon(Z_n) = \delta_{0,n} \qquad (12.21)$$

[49] Partly this comes from the fact that while the linear algebraic dual of a coalgebra is immediately an algebra it is not true that the full dual of an algebra is a coalgebra; instead, one must take a submodle of the dual called the zero-dual.

It is now easy to see that

$$\mathbf{NSymm} = \mathrm{Free}(\mathrm{CoFree}(\mathbf{Z})), \quad \mathbf{Symm} = \mathrm{Free}_{\mathrm{comm}}(\mathrm{CoFree}(\mathbf{Z})) \tag{12.22}$$

It is far less easy (except via duality) to see this way that **QSymm** = CoFree(Free(**Z**)), but it is interesting to do the exercise as this way the overlapping shuffle product comes about naturally.

It is also not easy to see directly that **Symm** = $\mathrm{CoFree}_{\mathrm{cocomm}}(\mathrm{Free}(\mathbf{Z}))$ (except, again, via duality. The duality argument is simple.

$$\begin{aligned}\mathbf{QSymm} = \mathbf{NSymm}^{dual} &= \mathrm{Free}(\mathrm{CoFree}(Z))^{dual} \\ &= \mathrm{CoFree}(\mathrm{CoFree}(\mathbf{Z})^{dual}) = \mathrm{CoFree}(\mathrm{Free}(\mathbf{Z}^{dual})) \\ &= \mathrm{CoFree}(\mathrm{Free}(\mathbf{Z}))\end{aligned}$$

where ‘dual’ stands for graded dual.

## 13. Frobenius and Verschiebung and other endomorphisms of $\Lambda$ and the Witt vectors.

After the rather abstract stuff of the two previous sections, it is a pleasure to return in this one to such concrete things as ghost components of power series and Witt vectors.

There are a vast number of functorial oprations on the functorial rings $\Lambda(A) \quad W(A)$ which is no surpise as every ring endomorphism of **Symm** induces such a functorial operation; and, as a free ring on coutably many generators thre are very many ring endomorphisms. Things become different if one also requires such things as preservation of the Abelian group structure on the rings $W(A)$ or preservation of the full functorial ring structure or one requires other interesting properties.

Probaly the easiest to describe are the socalled Verschiebung operations.

13.1. *Verschiebung*. The Verschiebung [50] operators are defined on the $\Lambda(A) = 1 + tA[[t]]$ by

$$\mathbf{V}_n\colon \Lambda(A) \longrightarrow \Lambda(A), \quad a(t) \mapsto a(t^n) \tag{13.2}$$

To see what this means at the ghost component level, apply $t\dfrac{d}{dt}\log$. Now if

$$\frac{d}{dt}\log a(t) = p_1 + p_2 t + p_3 t^2 + \cdots, \quad t(\frac{d}{dt}\log a(t^n)) = t(nt^{n-1}(p_1 + p_2 t^n + p_3 t + \cdots))$$

and so , on the ghost component level Verschiebung is given by

$$\mathbf{V}_n p_r = \begin{cases} np_{r/n} & \text{if } r \text{ is divisible by } n \\ 0 & \text{if } r \text{ is not divisible by } n \end{cases} \tag{13.3}$$

[50] The word means ‘shift’ and comes for the German . So it should indeed be written with an upper case first letter. In the case of the $p$-adic Witt vectors it was introduced by Witt.

This corresponds to the ring endomorphism of **Symm** (also denoted by $\mathbf{V}_n$)

$$\mathbf{V}_n\colon \mathbf{Symm} \longrightarrow \mathbf{Symm},\quad h_r \mapsto \begin{cases} h_{r/n} & \text{if } r \text{ is divisible by } n \\ 0 & \text{if } r \text{ is not divisible by } n \end{cases} \tag{13.4}$$

13.5. *Proposition.* The Verschiebung morphisms define additive functorial endomorphisms of the $\Lambda(A)$. That is they are Hopf algebra endomorphisms of **Symm**.

This is obvious from the definition. It can also be checked by the ussual ghost component argument using(13.3).

13.6. *Frobenius operations.* For the Frobenius operations on the $\Lambda(A)$ use again the symmetric function point of view. So again write the universal power series $h(t)$ as

$$h(t) = 1 + h_1 t + h_2 t^2 + \cdots = \prod_i (1-\xi_i t)^{-1} \tag{13.7}$$

and define

$$\mathbf{f}_n h(t) = \prod_i (1-\xi_i^n t)^{-1} \tag{13.8}$$

Formally, of course, one observes that the right hand side of (13.7) is symmetric in $\xi$. So there are universal polynomials $Q_{n,1}(h), Q_{n,2}(h), \cdots$ $\mathbf{Symm} = \mathbf{Z}[h_1, h_2, \cdots]$ for the coefficients in the power series (13.8) and now define

$$\mathbf{f}_n a(t) = 1 + Q_{n,1}(a)t + Q_{n,2}(a)t^2 + \cdots \tag{13.9}$$

The sequences of polynomials $Q_n(h) = (Q_{n,1}(h), Q_{n,2}(h), \cdots)$ that define the Frobenius operations are determined by

$$w_r(Q_n(h)) = w_{nr}(h) \tag{13.10}$$

In particular

$$Q_{n,1}(h) = w_n(h) \tag{13.11}$$

and

$$Q_{p,p^r}(h) \text{ only involves the variables } h_1, h_p, \cdots, h_{p^{r+1}} \tag{13.12}$$

The higher $Q_{n,r}(h)$ are difficult to write down explicitely.

The Frobenius Hopf algebra endomorphisms of **Symm** coresponding to the Frobenius operations are (of course) given by $h_r \mapsto Q_{n,r}(h)$.

For later purposes it is useful to have a little explicit knowledge about these polynomials.

13.13. *Lemma.* The polynomials $Q_{n,r}(h)$ determining the Frobenius operations satisfy

$$Q_{p,p^r}(h) \equiv (h_{p^r})^p \ \bmod p, \ \text{ for } p \text{ a prime number} \tag{13.14}$$

$$Q_{n,r}(h) \equiv nh_{rn} \ \bmod (h_1, h_2, \cdots, h_{rn-1}) \tag{13.15}$$

Both statements are proved by induction. From (13.11) it follows that (13.15) holds for $r = 1$. Suppose with induction that it holds for all $1 \le i < r$. By the defining property and this induction assumption $w_r(Q_n) \equiv rQ_{n,r} \ \bmod (h_1, h_2, \cdots, h_{nr-1})$.

On the other hand $w_{nr}(h) \equiv nrh_{nr} \ \bmod (h_1, h_2, \cdots, h_{nr-1})$. Comparing these two congruences proves (13.15).

Further (13.14) holds for $r = 0$, again by (13.11). Suppose (13.14) holds for $0 \le i < r$. It follows that

$$p^i (Q_{p,p^i})^{p^{r-i}} \equiv p^i (h_{p^{i-1}}) p^{r-i+1} \tag{13.16}$$

and consequently the first $r$ terms of

$$w_{p^r}(Q_p) = (Q_{p,1})^{p^r} + p(Q_{p,p})^{p^{r-1}} + \cdots + p^{r-1}(Q_{p,p^{r-1}})^p + p^r Q_{p,p^r}$$

Cancel $\bmod p^r$ with the first $r$ terms of

$$w_{p^{r+1}}(h) = (h_1)^{p^{r+1}} + p(h_p)^{p^r} + \cdots \ p^{r-1}(h_{p^{r-1}})^{p^r} + p^r (h_{p^r})^p + p^{r+1} h_{p^{r+1}}$$

leaving

$$p^r Q_{p,p^r} \equiv p^r (h_{p^r})^p + p^{r+1} h_{p^{r+1}}$$

proving (13.14).

In the same formal vein let $\zeta_n$ be a primitive $n$-th root of unity. Then always

$$\prod_{j=1}^{n} (1 - \zeta_n^j t^{1/n}) = (1 - t)$$

and so the Frobenius operation can also be written as

$$\mathbf{f}_n a(t) = \prod_{j=1}^{n} a(\zeta_n^j t^{1/n}) \tag{13.17}$$

The ghost components of $h(t)$ as in (13.7) are the power sums in the $\xi$ and so on the level of ghost components the Frobenius operations are characterized by

$$\mathbf{f}_n p_r = p_{nr} \tag{13.18}$$

The Frobenius operations $\mathbf{f}_p$ for $p$ a prime number have a Frobenius like property:

$$\mathbf{f}_p a(t) \equiv a(t)^{\ p} \ \text{mod}\ p \tag{13.19}$$ [51]

Here $^{p}$ means taking the $p$-th power in the ring $\Lambda(A)$. The congruence (13.19) takes place in $\Lambda(A)$, so it means that there is a power series $b(t)$ $\Lambda(A)$ such that $b(t)^p(\mathbf{f}_p a(t)) = a(t)^{\ p}$.

To prove this first observe that for power series of the form $(1-xt)^{-1}$ by the definition of Frobenius and product

$$\mathbf{f}_p(1-xt)^{-1} = (1-x^p t)^{-1} = ((1-xt)^{-1})^{\ p} \tag{13.20}$$

Now in any commutative ring, including $\Lambda(A)$, $(y+z)^p \equiv y^p + z^p \ \text{mod}\, p$ and so from (13.20) it follows that (13.19) holds for any finite product power series

$$(1-\xi_1 t)^{-1} \cdots (1-\xi t)^{-1}$$

(and in fact under the right circumstances (such as here) also for infinite products).

Alternatively, for an arbitrary power series take its Witt vector coordinate expression

$$a(t) = \prod_{d=1}^{\infty} (1-x_d t^d)^{-1} \tag{13.21}$$

This stabilizes in that the first $n$ coefficients of $a(t)$ only depend on the first $n$ $x$-coordinates. As (13.19) already holds for the all factors on the right of (13.21) one gets (13.19) in full generality by a simple limit argument.

13.22. *Theorem.* The $\mathbf{f}_n$ : **Symm** ¾→ **Symm** induce functorial ring endomorphisms on the rings $\Lambda(A)$. They are also the only ring endomorphisms of **Symm** that do so.

13.23. *Corollary.* The only ring automorphism of **Symm** that preserves its coring object structure (as an object in the category of rings) is the identity.

I.e. the object **Symm** is rigid as a coring object in the category **CRing.** Note that in contrast to the Liulevicius theorem 10.8, homogeneity is not needed.

13.24. *Proof of the corollary.* The Frobenius morphisms $\mathbf{f}_n$ satisfy $\mathbf{f}_n p_r = p_{nr}$. So the only one that induces an operation that is always injective is $\mathbf{f}_1 = \text{id}$ (as is seen by looking at the

[51] Let $K/\mathbf{Q}_p$ be an unramified extension of the $p$-adic numbers with ring of integers $A$. Then the Galois group $Gal(K/\mathbf{Q}_p)$ is cyclic with a generator $\sigma$. Let $k$ be the residue field of $K$. Then, as extensions of fields $A/\mathbf{Z}_p = W_{p^\infty}(k)/\ W_{p^\infty}(\mathbf{F}_p)$ and $\sigma$ is the Frobenius endomorphism of $W_{p^\infty}(k)$. It is characterized by $\sigma(a) \equiv a^p \,\text{mod}\ p$. This has to do with the terminology empl;oyed here.

operation on any $\mathbf{Q}$-algebra. [52]

13.25. *Proof of theorem* 13.22. For the first statement of the theorem is suffices to verify this on the ghost components, that is the the power sums $p_r$. On these the sum comultiplication, sum counit, product comultiplication ,product counit are respectively given by

$$\mu_S(p_r) = 1 \quad p_r + p_r \quad 1,\ \ \varepsilon_S(p_r) = 0,\ \ \mu_P(p_r) = p_r \quad p_r,\ \ \varepsilon_P(p_r) = 1 \tag{13.26}$$

So the characterizing property (13.18) of the Frobenius morphism proves that they preserve the functorial ring structure on the $\Lambda(A)$.

Now let $\varphi$ be a ring endomorphism of **Symm** that respects the structures (13.26). In particular that means that $\varphi$ must take primitives of the Hopf algebra **Symm** into primitives. The space of primitives of **Symm** consists of the linear (integer) combinations of the $p_r$. So the image of $p_r$ under $\varphi$ is of the form

$$\varphi(p_r) = \sum c_{rj} p_j \tag{13.27}$$

Now use $\mu_P(\varphi(p_r)) = (\varphi \quad \varphi)\mu_P(p_r) = \varphi(p_r) \quad \varphi(p_r)$ and $1 = \varphi(1) = \varphi(\varepsilon_P(p_r)) = \varepsilon_P(\varphi(p_r))$ to find that the coefficients $c_{rj}$ must satisfy

$$c_{rj}c_{rj} = c_{rj},\ \ c_{rj}c_{rj'} = 0 \text{ if } j \neq j',\ \ \sum_j c_{rj} = 1$$

and this is only possible if all but one of these coefficients are zero and that last one is equal to 1. Thus there is a mapping $\sigma\colon \mathbf{N} \longrightarrow \mathbf{N}$ such that

$$\varphi(p_r) = p_{\sigma(r)} \tag{13.28}$$

Now consider the Newton relations (9.57):

$$p_r = rh_r - (p_1h_{r-1} + p_2h_{r-2} + \cdots + p_{r=1}h_1) \tag{13.29}$$

It easily follows with induction that

$$p_r \equiv (-1)^{r+1}h_1^r \ \mathrm{mod}(h_2, h_3, \cdots) \tag{13.30}$$

Let $n = \sigma(1)$ and suppose with induction that it has been shown that

$$\varphi(p_u) = p_{nu},\ \ \varphi(h_u) \equiv 0 \ \mathrm{mod}(h_2, h_3, \cdots) \text{ for } r-1 \geq u \geq 2 \tag{13.31}$$

[52] Incidentally, all Frobenius operations except the identiy are also not injective on the integers. Indeed let $p$ be a prime number and consider the vector $a = (a_1, a_2, a_3, \cdots)$ $\mathbf{Z}^{\mathbf{N}}$ with $a_n = p$ if $n$ is not a multiple of $p$, and $a_n = 0$ if $n$ is a multiple of $p$. Then by the ghost-vector criterium in 9.93 (with $A = \mathbf{Z}$ and all the endomorphisms the identity) this vector is ghost-Witt vector i.e. there is a Witt vector $x$ $W(\mathbf{Z})$ such that $w(x) = a$. For this vector $\mathbf{f}_p(x) = 0$.

Now apply $\varphi$ to the Newton relation (13.29) and use (13.28), (13.30) and the induction hypothesis (13.31) to find that

$$(-1)^{\sigma(r)+1}h_1^{\sigma(r)} \equiv \varphi(p_r) \equiv r\varphi(h_r) - ((-1)^{(r-1)n+1}h_1^{(r-1)n})(-1)^{n+1}h_1^n) \bmod (h_2, h_3, \cdots)$$

and so $\sigma(r) = nr$ and $\varphi(h_r) \equiv 0 \bmod (h_2, h_3, \cdots)$. It follows that $\varphi = \mathbf{f}_n$.

13.32. Note that the proof does not use any integrality of coefficients statement and so the theorem still holds over the rationals. Indeed it works over any ring and so the theorem is true for the coring object algebras $\mathbf{Symm}_R = \mathbf{Symm} \otimes R$ over any ring $R$.

13.33. *Frobenius endomorphisms of* **QSymm**. The Frobenius morphisms $\mathbf{f}_n$ on **Symm** extend to Frobenius endomorphisms on $\mathbf{QSymm} \supset \mathbf{Symm}$. The formula is simple

$$\mathbf{f}_n([a_1, a_2, \cdots, a_m]) = [na_1, na_2, \cdots, na_m] \tag{13.34}$$

It is ovious from the original definitions of the first and second comultiplication, and the description of the overlapping shuffle product that these are Hopf algebra endomorphisms that preserve the second comultiplication and second co-unit. They are also the only ring endomorphisms of **QSymm** that do so.

To prove this consider a Hopf algebra endomorphism $\varphi$ of **QSymm** that preserves the second comultiplication. First it must take primitives into primitives. The primitives are spanned by the compositions of length 1, see 11.36. Further

$$\mu_P([a]) = [a] \otimes [a], \quad \varepsilon_P([a]) = 1$$

and so, exactly as above, there is a mapping $\sigma: \mathbf{N} \longrightarrow \mathbf{N}$ such that

$$\varphi([a]) = [\sigma(a)] \tag{13.35}$$

Let $\sigma(1) = n$ and let $J_2$ be the ideal in **QSymm** spanned by all compositions of length 2 or more. This is indeed an ideal. Now also $\varepsilon_P\varphi = \varepsilon_P$ and it follows that $\varphi(J_2) \subset J_2$. Further $[r] \equiv [1]^r \bmod J_2$ and so

$$\varphi([r]) \equiv \varphi([1]^r) = \varphi([1])^r = [n]^r \equiv [nr]$$

and so $\varphi$ is equal to $\mathbf{f}_n$ on the submodule spanned by the [a]. Now consider a composition of length 2. As $\varphi$ respects the first comultiplication

$$\begin{aligned}\mu_S(\varphi([a,b])) &= (\varphi \otimes \varphi)(1 \otimes [a,b] + [a] \otimes [b] + [a,b] \otimes 1) \\ &= 1 \otimes \varphi([a,b]) + [na] \otimes [nb] + \varphi([a,b]) \otimes 1\end{aligned} \tag{13.36}$$

The term $[na] \otimes [nb]$ on the right of (13.36) can only come under $\mu_S$ from the composition $[na, nb]$ and if there were any other compositions involved in $\varphi([a,b])$ that would show up after applying $\mu_S$. So $\varphi([a,b]) = [na, nb]$. This argument easily continues.

13.37. *Corollary*. The Frobenius operators on **Symm** are given on the monomial symmetric functions by

$$\mathbf{f}_n m_{(\lambda_1,\lambda_2,\cdots,\lambda_m)} = m_{(n\lambda_1, n\lambda_2,\cdots,n\lambda_m)} \tag{13.38}$$

There are (of course) also other ways of seeing this, but this is a particularly elegant way. Thus on the monomial symmetric functions the Frobenius operators are easy to describe. It is far more difficult to describe them on the $h_\lambda$, though of course there are the messy formulas giving the $h$'s in terms of the $p$'s. Dually the Verschiebung operators are easy to describe on the $h$'s and difficult to describe on the monomial symmetric functions.

13.39. *Corollary*. There are no second multiplication preserving Hopf algebra endomorphisms of **NSymm** that descend to a Frobenius endomorphism of **Symm** as a quotient of **NSymm** except the identity.

Indeed such a morphism would involve some degree increasing part. Its dual would therefore involve some degree increasing part and be a second comultiplication preserving Hopf algebra endomorphism of **QSymm**, which is impossible by 13.25 above.

This gives another negative answer, as conjectured, to a question posed in [200, 201].

13.40. The Frobenius morphisms $\mathbf{f}_p$, for $p$ a prime number also satisfy the Frobenius like property

$$[a_1,a_2,\cdots,a_m]^p \equiv [pa_1,pa_2,\cdots,pa_m] \bmod p \tag{13.41}$$

This is more or less immediate from the definition of the overlapping shuffle product. Indeed, the terms of this $p$-th power are the column sums of all matrices without zero columns of which the rows are obtained from $a_1,a_2,\cdots,a_m$ by interspersing zeros. If all the rows are equal the column sum $[pa_1,pa_2,\cdots,pa_m]$ results. If not all rows are equal any (non identity) cyclic permutation yields another different matrix of the same type (this uses that $p$ is prime) and so all the other terms in the shuffle power occur with a coefficient divisible by $p$.

13.42. *'Multiplication by n' operator*. Consider the operator of adding (in $\Lambda(A)$ or $W(A)$) an element to itself $n$ times, i.e in $\Lambda(A)$ taking the $n$-th power of a power series. This operator is denoted

$$[n]:\ \Lambda(A) \longrightarrow \Lambda(A).\quad a(t) \mapsto a(t)^n \tag{13.43}$$

This operator as an endomorphism of **Symm** is given by the composition of maps

$$\mathbf{Symm} \xrightarrow{\mu_n} \mathbf{Symm}^{\otimes n} \xrightarrow{m_n} \mathbf{Symm} \tag{13.44}$$

where $\mu_n$ is the $n$-fold first comultiplication as defined in (12.14) and $m_n$ is the $n$-fold multiplication. This composed map can be written down for any Hopf algebra $H$. In general it is not a Hopf algebra endomorphism, nor even an algebra or coalgebra morphism unless $H$ is commutative and cocommutative. Then it is a Hopf algebra endomorphism. In spite of not being a Hopf algebra endomorphism in general these maps have proved to be most useful in various investigations in Hopf algebra theory, see e.g. [323].

13.45. *Homothety operations*. Consider the Hopf algebra $\mathbf{Symm}_R = \mathbf{Symm} \otimes R$ over a ring $R$. For every $R$-algebra $A$ and every $u \in R$ consider the operation

$$\langle u \rangle a(t) = a(ut), \quad a(t) \in \Lambda(A) \tag{13.46}$$

These are the homothety operators and they clearly define additive functorial Abelian group endomorphisms of the $\Lambda(A)$. The associated Hopf algebra endomorphism of $\mathbf{Symm}_R$ is

$$\langle u \rangle (h_n) = u^n h_n \tag{13.47}$$

As will be seen later the homothety, Verschiebung, and Frobenius endomorphisms together, in a very precise sense, generate all the $R$-Hopf algebra endomorphisms of $\mathbf{Symm}_R$.

There are quite a good many relations amoung all these operators. They are summed up in the following theorem.

13.48. *Theorem.* There are the following identies between operations on the functors $\Lambda_R(-)$, $W_R(-)$, respectively the Hopf algebra $\mathbf{Symm}_R$.

$$\langle u \rangle \langle u' \rangle = \langle uu' \rangle \tag{13.49}$$

$$\langle 1 \rangle = \mathbf{f}_1 = \mathbf{V}_1 = \mathrm{id} \tag{13.50}$$

$$\mathbf{V}_m \mathbf{V}_n = \mathbf{V}_{nm} \tag{13.51}$$

$$\mathbf{f}_m \mathbf{f}_n = \mathbf{f}_{mn} \tag{13.52}$$

$$\text{if } \gcd(m,n) = 1 \text{, then } \mathbf{f}_m \mathbf{V}_n = \mathbf{V}_n \mathbf{f}_m \tag{13.53}$$

$$\mathbf{f}_n \mathbf{V}_n = [n] \tag{13.54}$$

$$\mathbf{V}_n (a(t) \ast \mathbf{f}_n b(t)) = (\mathbf{V}_n a(t)) \ast b(t), \quad a(t), b(t) \in \Lambda(A) \tag{13.55}$$

$$\langle \mathbf{f}_n x, y \rangle = \langle x, \mathbf{V}_n y \rangle, \quad x, y \in \mathbf{Symm} \tag{13.56}$$

$$\langle u \rangle \mathbf{V}_n = \mathbf{V}_n \langle u^n \rangle \tag{13.57}$$

$$\mathbf{f}_n \langle u \rangle = \langle u^n \rangle \mathbf{f}_n \tag{13.58}$$

$$\langle u \rangle + \langle v \rangle = \sum_{n=1}^{\infty} \mathbf{V}_n \langle r_n(u,v) \rangle \mathbf{f}_n \tag{13.59}$$

where the $r_d(X,Y)$ are the integer coefficient polynomials in two variables determined by

$$(1-Xt)(1-Yt) = \prod_{d=1}^{\infty}(1-r_d(X,Y)t^d) \tag{13.60}$$

This, (13.59), (13.60) is of course the formula for the addition of Teichmüller representatives in Witt coordinates. Applying $t\frac{d}{dt}\log$ to the two sides of (13.60) it follows that the polynomials $r_j(X,Y)$ are determined by the relation

$$X^n + Y^n = \sum_{d|n} dr_d(X,Y)^{n/d} \tag{13.61}$$

For example

$$r_1(X,Y) = X+Y$$
$$r_2(X,Y) = -XY$$
$$r_3(X,Y) = -(X^2Y + XY^2)$$
$$r_4(X,Y) = -(X^3Y + 2X^2Y^2 + XY^3)$$
$$r_5(X,Y) = -(X^4Y + 2X^3Y^2 + 2X^2Y^3 + XY^4)$$
$$r_6(X,Y) = -(X^5Y + 3X^4Y^2 + 4X^3Y^3 + 3X^2Y^4 + XY^5)$$

It immediately strikes one's attention that the coefficients in $r_2, r_3, \cdots, r_6$ are all negative. This is true for all $n \geq 2$ as an immediate consequence of the Reutenauer-Scharf-Thibon result 9.71. Another striking fact is that all monomials that possibly can occur do in fact occur with nonzero coefficient. I know of no proof for that but definitely believe it to be true.

13.62. ***Caveat.*** All these formulas are written from the functorial operations point of view (except (13.56) where the only possible interpretation is in terms of endomorphisms). So, for instance, (13.54) means that for an element $a(t)\ \ \Lambda(A)$ there is the equality

$$\mathbf{f}_n(\mathbf{V}_n(a(t))) = [n](a(t)) \tag{13.63}$$

If $U_1, U_2, U$ are functorial operations and $u_1, u_2, u$ are the endomorphism that induce them, then seeing an element $a(t)\ \ \Lambda(A)$ as a morphism $\mathbf{Symm} \longrightarrow A$ (as must be done to have this correpondence between endomorphisms and operators), the relation beteen the operation $U$ and the endomorphism $u$ is

$$U(a(t)) = a(t)\circ u \tag{13.64}$$

where the small circle denotes composition. It follows that under the correspondence 'functorial operation' $\Leftrightarrow$ 'endomorphism' the order of composition reverses. Indeed

$$U_1(U_2(a(t))) = U_1(a(t)\circ u_2) = (a(t)\circ u_2)\circ u_1 = a(t)\circ(u_2\circ u_1)$$

13.65. *Proof of theorem* 13.48. Most of these are pretty trivial using the ghost component formalism (see also 5.16). As before for any natural number $r$ and element $a(t)\ \ \Lambda(A)$ let

$$s_r(a(t)) = \text{ coefficient of } t^r \text{ in } t\frac{d}{dt}\log(a(t)) \tag{13.66}$$

Then the characterizing ghost component description of the various functorial operations are

$$s_r(\mathbf{V}_n a(t)) = \begin{cases} ns_{r/n}(a(t)) & \text{if } n \text{ divides } r \\ 0 & \text{if } n \text{ does not divide } r \end{cases} \tag{13.67}$$

$$s_r(\mathbf{f}_n a(t)) = s_{nr}(a(t)) \tag{13.68}$$

$$s_r(\langle u \rangle a(t)) = u^r s_r(a(t)) \tag{13.69}$$

$$s_r([n]a(t)) = ns_r(a(t)) \tag{13.70}$$

So, using this, here is a proof of (13.54):

$$s_r(\mathbf{f}_n\mathbf{V}_n a(t)) = s_{rn}(\mathbf{V}_n a(t)) = ns_r(a(t)) = s_r([n]a(t))$$

using (13.68), (13.67), (13.70) in the order named. And here is how (13.59) is tackled with the ghost component formalism

$$s_r(\mathbf{V}_n\langle r_n(u,v)\rangle \mathbf{f}_n a(t) = \begin{cases} ns_{r/n}(\langle r_n(u,v)\rangle \mathbf{f}_n a(t)) & \text{if } n \text{ divides } r \\ 0 & \text{if } n \text{ does not divide } r \end{cases}$$

$$= \begin{cases} nr_n(u,v)^{r/n} s_{r/n}(\mathbf{f}_n a(t)) & \text{if } n \text{ divides } r \\ 0 & \text{if } n \text{ does not divide } r \end{cases}$$

$$= \begin{cases} nr_n(u,v)^{r/n} s_r(a(t)) & \text{if } n \text{ divides } r \\ 0 & \text{if } n \text{ does not divide } r \end{cases}$$

and so

$$s_r\left(\sum_{n=1}^{\infty}{}^{\wedge}\, \mathbf{V}_n\langle r_n(u,v)\rangle \mathbf{f}_n(a(t))\right) = \sum_{n|r} nr_n(u,v)^{r/n} s_r(a(t)) = (u^n + v^n)s_r(a(t)) = s_r(\langle u\rangle a(t) +_{\Lambda} \langle v\rangle a(t))$$

where of course, as indicated, the sum of elements in $\Lambda(A)$ is taken according to the addition law in the Abelian group $\Lambda(A)$.

Here is the proof of (13.53) as a further illustration of ghost component techniques

$$s_r(\mathbf{f}_m\mathbf{V}_n a(t)) = s_{rm}(\mathbf{V}_n a(t)) = \begin{cases} s_{rm/n}(a(t)) & \text{if } n \text{ divides } rm \\ 0 & \text{otherwise} \end{cases}$$

and as $(m,n) = 1$ the divisibility condition in the formula above is the same as '$n$ divides $r$', which is whar turns up when calculating $s_r(\mathbf{V}_n\mathbf{f}_m a(t))$.

All the other statements of theorem 13.48 are proved the same way, except of course the duality statement (13.56). For this statement use the orthonormal dual pair of bases $\{h_\lambda\}$, $\{m_\lambda\}$, see (9.40), (9.42). Let $\kappa = (\kappa_1, \kappa_2, \cdots, \kappa_m)$, $\lambda = (\lambda_1, \lambda_2, \cdots, \lambda_n)$ be two partitions. Recall that

$$\mathbf{f}_r m_\lambda = m_{(r\lambda_1, r\lambda_2, \cdots, r\lambda_n)}, \quad \mathbf{V}_r h_\kappa = \begin{cases} h_{(\kappa_1/r, \kappa_2/r, \cdots, \kappa_m/r)} & \text{if all } \kappa_j \text{ are divisible by } r \\ 0 & \text{otherwise} \end{cases}$$

Thus

$$\langle \mathbf{f}_r m_\lambda, h_\kappa \rangle = \begin{cases} 1 & \text{if } \lg(\lambda) = \lg(\kappa) \text{ and } \kappa_i = r\lambda_i \text{ for all } i \\ 0 & \text{otherwise} \end{cases} \tag{13.71}$$

$$\langle m_\lambda, \mathbf{V}_r h_\kappa \rangle = \begin{cases} 1 & \text{if } \lg(\lambda) = \lg(\kappa),\ r \text{ divides all } \kappa_i \text{ and } \kappa_i / r = \lambda_i \text{ for all } i \\ 0 & \text{otherwise} \end{cases} \tag{13.72}$$

comparing (13.71) and (13.72) now gives the duality formula (13.56)

13.73. *Convention.* Thus from (13.56) the Frobenius and Verschiebung endomorphisms of **Symm** = **Z**[$h$] are (graded) dual to each other. But **Symm** is also isomorphic to its graded dual. That makes it a trifle difficult to agree what should be called Frobenius and what Verschiebung. The convention is that the degree nondecreasing endomorphisms (of these) are always called Frobenius and the degree nonincreasing ones Verschiebung. Both are Hopf algebra endomorphisms.

Then the Frobenius endomorphisms always respect the second comultiplication (and the corresponding co-unit), but not the second multiplication, and thus induce functorial ring endomorphisms of the $W(A) \cong \Lambda(A) = \mathbf{CRing}(\mathbf{Symm}, A)$, while the Verschiebung endomorphisms respect the second multiplication (and corresponding unit), but not the second comultiplication, and thus induce ring endomorphisms on the rings of coalgebra morphims $\mathbf{CoAlg}(C, \mathbf{Symm})$.

13.74. *Ring of Hopf endomorphisms of* **Symm**. Naturally the next step is to try to determine the full ring of Hopf algebra endomorphisms of **Symm**, which is the same as the determination of all additive functorial operations on the functor $\Lambda$. Consider again, see 12.11 ff, the cofree coalgebra over **Z**, i.e. the coalgebra

$$\mathrm{CoFree}(\mathbf{Z}) = \bigoplus_{n=0}^{\infty} \mathbf{Z}h_n, \quad h_0 = 1, \quad \mu(h_n) = \sum_{i+j=n} h_i \otimes h_j, \quad \varepsilon(h_n) = 0 \text{ for } n \geq 1 \tag{13.75}$$

Now **Symm** = **Z**[$h$] is the free algebra on CoFree(**Z**), more precisely it is the free algebra over the augmentation submodule $\bigoplus_{n=1}^{\infty} \mathbf{Z}h_n \subset \mathrm{CoFree}(\mathbf{Z})$ and hence a Hopf algebra endomorphism of **Symm** is the same thing as a morphism of coalgebras from CoFree(**Z**) to **Symm** that takes $h_0$ to 1.

$$\mathbf{Hopf}(\mathbf{Symm}, \mathbf{Symm}) \cong \mathbf{CoAlg}'(\mathrm{CoFree}(\mathbf{Z}), \mathbf{Symm}) \tag{13.76}$$

where the prime indicates that only morphisms that take $h_0$ to 1 are permitted.
In more pedestrian terms, as **Symm** = **Z**[$h$] is free on the generators $h_n$ an algebra endomorphism $\varphi$ of **Symm** is the same thing as specifying a sequence of polynomials $\varphi(h_n)$ and this sequence yields an endomorphism of Hopf algebras if and only if the sequence $1, \varphi(h_1), \varphi(h_2), \varphi(h_3), \cdots$ is a curve (= divided power sequence).

The graded dual of $\mathrm{CoFree}(\mathbf{Z})$ is $\mathbf{Z}[T]$ the ring of polyomials in one variable $T$.

Thus, taking graded duals, and using that **Symm** and its graded dual are isomorphic as Hopf algebras (and even as coring objects in the category of rings and as ring objects in the category of coalgebras) one sees that

$$\mathbf{CoAlg}(\mathrm{CoFree}(\mathbf{Z}), \mathbf{Symm}) \quad \mathbf{CRing}(\mathbf{Symm}, \mathbf{Z}[T]) \tag{13.77}$$

and thus, tracing out what the prime in (13.68) means for the right hand side of (13.69,

$$\mathbf{Hopf}(\mathbf{Symm}, \mathbf{Symm}) \quad \mathbf{CRing}'(\mathbf{Symm}, \mathbf{Z}[t]) \tag{13.78}$$

where this time the prime means that only morphisms are allowed that take the degee $\geq 1$ part of **Symm** into the degree $\geq 1$ part of $\mathbf{Z}[T]$.

But his last object is easy to describe: an element of $\mathbf{CRing}'(\mathbf{Symm}, \mathbf{Z}[t])$ is simply an infinite sequence of polynomials in $T$ with constant terms zero. Things go exactly the same way for $\mathbf{Symm}_A = \mathbf{Symm} \quad A$.

As the isomorphism between **Symm** and its graded dual takes the $h_n$ into the $e_n$ it is best to see this sequence of polynomials as the images of the $e_n$.

Thus in this way a Hopf algebra endomorphism $\varphi$ of $\mathbf{Symm}_A$ is exactly the same thing as an $\text{infinite} \times \text{infinite}$ matrix with coefficients in $A$

$$M_\varphi = \begin{pmatrix} a_{11} & a_{12} & a_{13} & \cdots \\ a_{21} & a_{22} & a_{23} & \cdots \\ a_{31} & a_{32} & a_{33} & \cdots \\ \vdots & \vdots & \vdots & \ddots \end{pmatrix}$$

such that in each row there are but finitely many coeffients that are nonzero.

It remains to figure out how these matrices are to be added and multiplied (composed) when interpreted as encodings of Hopf endomorphisms of the Hopf algebra **Symm**. That will be the business of section 15 below.

To do something similar for the characteristic $p$ case as well (and the $p$-adic case) one needs to figure out how the $p$-adic Witt vectors fit with the power series point of view. That is the subject matter of the next section.

## 14. Supernatural and other quotients of the big Witt vectors

It is when trying to figure out how the $p$-adic Witt vectors fit with the power series point of view of sections 9-13 that the Witt vector coordinates have a decided advantage. Recall that in terms of Witt vector coordinates, 9.63, as a set

$$W(A) = \{(x_1, x_2, x_3, \cdots) \colon \ x_n \quad A\} \tag{14.1}$$

and that two Witt vectors are added and multiplied by means of universal polynomials with integer coordinates

$$x +_W y = (\mu_{S,1}(x,y), \mu_{S,2}(x,y), \mu_{S,3}(x,y), \cdots)$$
$$x \times_W y = (\mu_{P,1}(x,y), \mu_{P,2}(x,y), \mu_{P,3}(x,y), \cdots) \tag{14.2}$$

where $x = (x_1, x_2, x_3, \cdots)$, $y = (y_1, y_2, y_3, \cdots$ , and where the polynomials $\mu_{S,i}$, $\mu_{P,i}$ are recursively given by

$$w_n(\mu_{S,1}(X;Y), \mu_{S,2}(X;Y), \mu_{S,3}(X;Y), \cdots) = w_n(X) + w_n(Y)$$
$$w_n(\mu_{P,1}(X;Y), \mu_{P,2}(X;Y), \mu_{P,3}(X;Y), \cdots) = w_n(X) w_n(Y) \tag{14.3}$$

with

$$w_n(X) = \sum_{d|n} d X_d^{n/d} \tag{14.4}$$

The important fact to notice is now that $w_n(X)$ depends only on the $X_d$ for $d$ a divisor of $n$, and hence that the $n$-th addition and multiplication polynomials $\mu_{S,n}$, $\mu_{P,n}$ are polynomials that only involve the $X_d$ and $Y_d$ with $d$ a divisor of $n$. Thus, for instance, $\mathbf{Z}[X_1, X_p, X_{p^2}, X_{p^3}, \cdots]$, where $p$ is a prime number, is a sub Hopf algebra and sub coring object of $\mathbf{Z}[X] = \mathbf{Z}[X_1, X_2, X_3, \cdots]$ , which means that it defines a quotient functor, which is manifestly $W_{p^\infty}$, of $W$ so that the $p$-adic Witt vectors are a functorial quotient of the big Witt vectors. There are lots more such quotient functors.

14.5. *Nests*. A nest is a nonempty subset $N$ of the natural numbers $\mathbf{N}$ such that together with any $n \in N$ all the divisor of $n$ are also in $N$. Some simple examples of nests are

$$\text{for any given fixed } n \in \mathbf{N}, \text{ the set } \{1,2,3,\cdots,n\}$$
$$p^\infty = \{1, p, p^2, p^3, \cdots\} \tag{14.6}$$
$$\text{for any fixed } n \in \mathbf{N}, \text{ the set } \{d: d \text{ divides } n\}$$

Unions and intersections of nests are nests. Every nest contains the natural number 1.

14.7. *Supernatural numbers.* A supernatural number is a formal expression of the form

$$\mathfrak{n} = \prod_p p^{\alpha_p} \quad \alpha_p \in \mathbf{N} \cup \{0, \infty\}$$

where the product is over all prime numbers $p$. This is exactly the same as specifying for every prime number $p$ an element of the extended natural numbers $\{0\} \cup \mathbf{N} \cup \{\infty\}$. Given a supernatural number $\mathfrak{n}$ there is a nest associated to it (often denoted with the same symbol), viz

$$\mathfrak{n} = N_{\mathfrak{n}} = \{m \in \mathbf{N}: v_p(m) \le \alpha_p \text{ for all prime numbers } p\}$$

where $v_p$ is the $p$-adic valuation on the integers. The last two examples in (14.6) are nests coming from a supernatural number. The first example from (14.6) is not of this form (if $n \geq 3$).

14.8. *Supernatural quotients of the bigWitt vectors.* Let $N$ be a nest. Let $\mathbf{Symm}(N) = \mathbf{Z}[X_n : n \in N] \subset \mathbf{Symm}$. By the remarks made in the beginning of this section these form sub Hopf algebras of **Symm** and hence define quotient Wiit vector functors $W_N(A)$ of the big Witt vectors, called the *N*-adic Witt vectors. The most important ones are

$W_{p^\infty}(A)$, the $p$-adic Witt vectors

$W_{p^n}(A)$, the $p$-adic Witt vectors of length $n+1$

$W_n(A)$, the Witt vectors of length $n$.

The last named quotient functor is the one defined by the nest $\{1,2,\cdots,n\}$. This notation is not entirely consistent with the other two. However, the quotient Witt vectors defined by the nest $\{d\colon\ d \mid n\}$ are so seldom used (if ever) that this seems justified.

Sometimes the curious appellation "nested Witt vectors" is used for the elements of a $W_N(A)$.

14. 9. *Operations on nested Witt vectors.* A first question now is now which of the functorial additive operations $\langle a \rangle$, $\mathbf{V}_r$, $\mathbf{f}_r$ survive to define operations on the various quotients $W_N$. Here of course these functorial operations are the ones defined by transferring the operations denoted by the same symbols on $\Lambda(A)$. They are therefore characterized by

$$w_r \mathbf{f}_n = w_{rn}, \quad w_r \mathbf{V}_n = \begin{cases} w_{r/n} & \text{if } n \text{ divides } r \\ 0 & \text{otherwise} \end{cases}, \quad w_r \langle a \rangle = a^r w_r \tag{14.10}$$

The answer to the above question of which operations 'descend' to operations on the various $W_N$ is fairly simple.

The homothety operations $\langle a \rangle$ always define a quotient operation.

For the Verschiebung operations the situation is simple. The polynomials defining $\mathbf{V}_n$ are either zero or of the form $X_{r/n}$ for $r$ a multiple of $n$. So by the nest property they always exist on the quotients $W_N$. Of course many of them may be zero.

The situation with the Frobenius operators is different and a slightly more complicated. Let $p$ be a prime number. The polynomials defining $\mathbf{f}_p$ are degree increasing by a factor of $p$ and because $\mathbf{f}_p$ is characterized by $\mathbf{f}_p(p_n) = p_{pn}$ at the ghost component level, the *n*-th polynomial defining $\mathbf{f}_p$ involves the indeterminate $X_{np}$. See also lemma 13.13 above. So for $\mathbf{f}_p$ to descend to a nest quotient in full generality, i.e. to exist for the functor $W_N$ on all of **CRing,** it is necessary and sufficient that the nest contain together with any $n \in N$ also all its $p$-power multiples $p^i n$. In particular $\mathbf{f}_p$ exists on the $p$-adic Witt vectors.

However, if one restricts attention to the functor on rings of characteristic $p$ only $\mathbf{f}_p$ always exists on the quotients $W_N$ because in that case the Frobenius is given by raising each coordinate of a Witt vector to its $p$-th power. See again lemma 13.13.

14.11. *Möbius function.* The next question is whether there exists lifts, that is an additive functor morphism $W_N \longrightarrow W$ which composed with the quotient map $W \longrightarrow W_N$ gives the identity on $W_N$. This will involve the Möbius function from number theory and combinatorics. Explicitely this function is defined by

$$\mu(1) = 1$$
$$\mu(n) = (-1)^r \quad \text{if } n \text{ is the product of } r \text{ different prime numbers} \tag{14.12}$$
$$\mu(n) = 0 \quad \text{if } n \text{ is divisible by the square of a nontrivial prime number}$$

The function has the characterizing property

$$\sum_{d|n} \mu(d) = 1 \tag{14.13}$$

(which also defines it recursively). Let $\mathbf{N}(p) = \{n \in \mathbf{N}: (p,n) = 1\}$ be the set of all natural numbers relatively prime to a given prime number $p$. It immediately follows from the above that

$$\sum_{d|n,\ d\in \mathbf{N}(p)} \mu(d) = \begin{cases} 1 & \text{if } n \text{ is a power of } p \\ 0 & \text{otherwise} \end{cases} \tag{14.14}$$

14.15. *Sectioning the projection* $W(A) \longrightarrow W_{p^\infty}(A)$. To set the stage here is the abstract situation. Suppose there is an Abelian group $M$ together with a surjective projection $\pi\colon M \longrightarrow M_p$ to another group $M_p$. Suppose there is a section, i.e. a morphism of Abelian groups $s\colon M_p \longrightarrow M$ such that $\pi s = id_{M_p}$. Then $\Phi = s\pi$ is an idempotent endomorphism of $M$ with image $s(M_p)$ and kernel $\mathrm{Ker}(\pi)$. Further the projection $\pi$ induces an isomorphism from $\Phi(M) = s(M_p)$ to $M_p$. It is such an idempotent endomorphism $\Phi$ that will now be constructed in the case of the canonical projection

$$\pi_p\colon W(A) \longrightarrow W_{p^\infty}(A) \tag{14.16}$$

in the case of $\mathbf{Z}_{(p)}$ algebras, i.e. rings $A$ in which all prime numbers except the give prime number $p$ are invertible.

To this end consider the operation

$$\Phi_{p\text{-typ}} = \sum [\frac{\mu(n)}{n}]\mathbf{V}_n\mathbf{f}_n \tag{14.17}$$

This functorial operation satisfies

$$w_r\Phi_{p\text{-typ}} = \begin{cases} w_r & \text{if } r \text{ is a power of } p \\ 0 & \text{otherwise} \end{cases} \tag{14.18}$$

(and it is characterized by this property, of course). This is an immediate consequence of (14.14) and (14.10). Note that it follows that $\Phi_{p\text{-typ}}$ is idempotent. (First for rings $A$ of characteristic zero and then for all by functoriality).

14.19. *Remark.* $\Phi_{p\text{-typ}}$ commutes with the operations $\langle a\rangle$, $\mathbf{f}_p$, $\mathbf{V}_p$.

14.20. *Comments.* $\Phi_{p\text{-typ}}$ is not easy to write down explicitely, say, in terms of its defining polynomials. It is definitely not something like murdering each coordinate that is not at a $p$-th power position.

Easy explicit calculations show that it does not preserve the unit element and that it is not multiplicative. This already shows up at the second coordinate for all prime numbers larger than 2 and in the third coordinate for the prime number 2.

14.21. *Proposition.* For each $\mathbf{Z}_{(p)}$-algebra $A$, a Witt vector $x = (x_1, x_2, x_3, \cdots) \in W(A)$ is in the image of $\Phi_{p\text{-typ}}(A)\colon W(A) \longrightarrow W(A)$ if and only if $\mathbf{f}_{p'}x = 0$ for all prime numbers $p'$ different from $p$.

First, if $\mathbf{f}_{p'}x = 0$ for all prime numbers different from $p$, then $\mathbf{f}_n x = 0$ for all $n \in \mathbf{N}(p)$, $n \geq 2$ and so $\Phi_{p\text{-typ}}\, x = [1]x = x$. Second, if $A$ is of characteristic zero and $x \in W(A)$ is in the image of $\Phi_{p\text{-typ}}(A)$ then the characterizing property (14.18) shows that $\mathbf{f}_{p'}x = 0$. Third suppose that $x \in W(A)$ is in the image of $\Phi_{p\text{-typ}}$. Choose any lift $\tilde{x} \in W(\tilde{A})$ of $x \in W(A)$ for some characteristic zero $\mathbf{Z}_{(p)}$-algebra $\tilde{A}$ covering Â. Then, by idempotency, $\Phi_{p\text{-typ}}(\tilde{A})\tilde{x}$ is also a lift of $x$, and so $\mathbf{f}_{p'}\Phi_{p\text{-typ}}\tilde{x} = 0$ implies $\mathbf{f}_{p'}x = 0$.

14.21. *Theorem.* For $\mathbf{Z}_{(p)}$-algebras the canonical projection $\pi_p\colon W(A) \longrightarrow W_{p^\infty}(A)$ induces an isomorphism of Abelian groups $\Phi_{p\text{-typ}}(W(A)) \xrightarrow{\cong} W_{p^\infty}(A)$.

Pictorially things are described by the following commutative diagram..

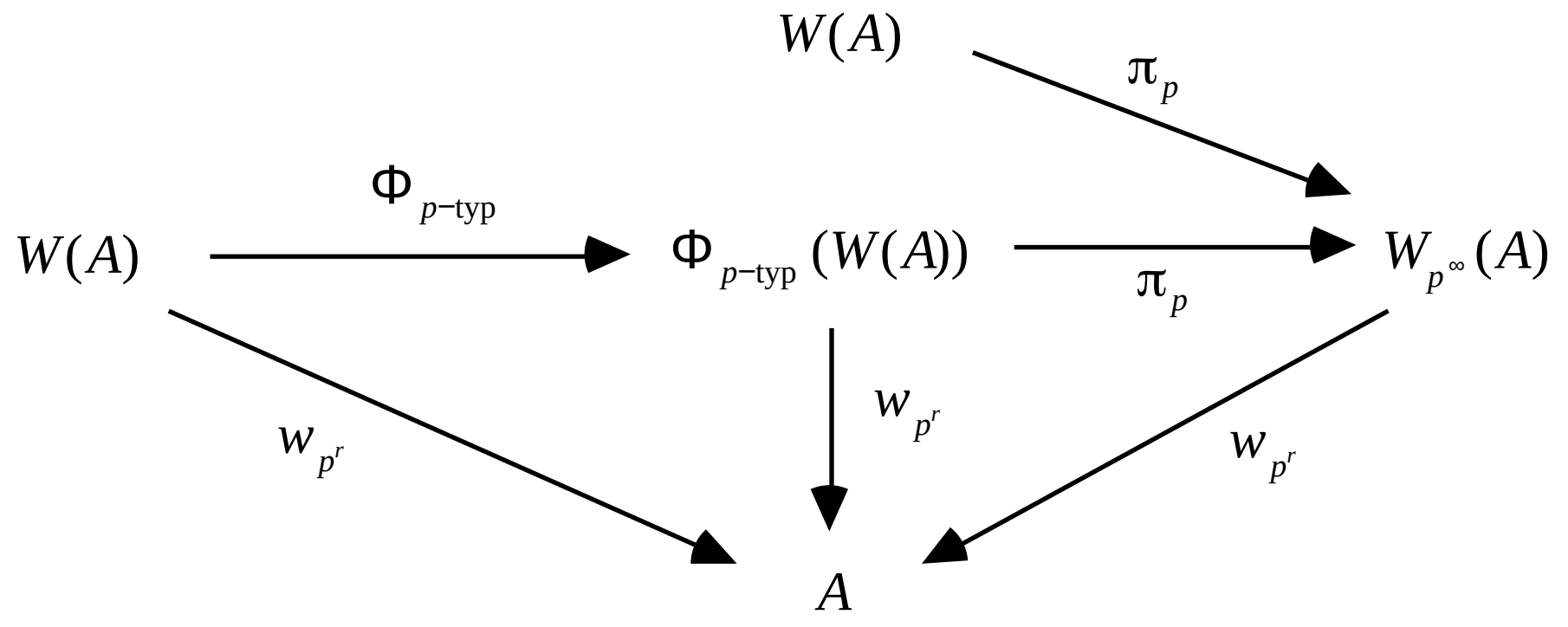

(14.22)

It follows just about immediately that the theorem holds for $A$ of charateristic zero. Now let $A$

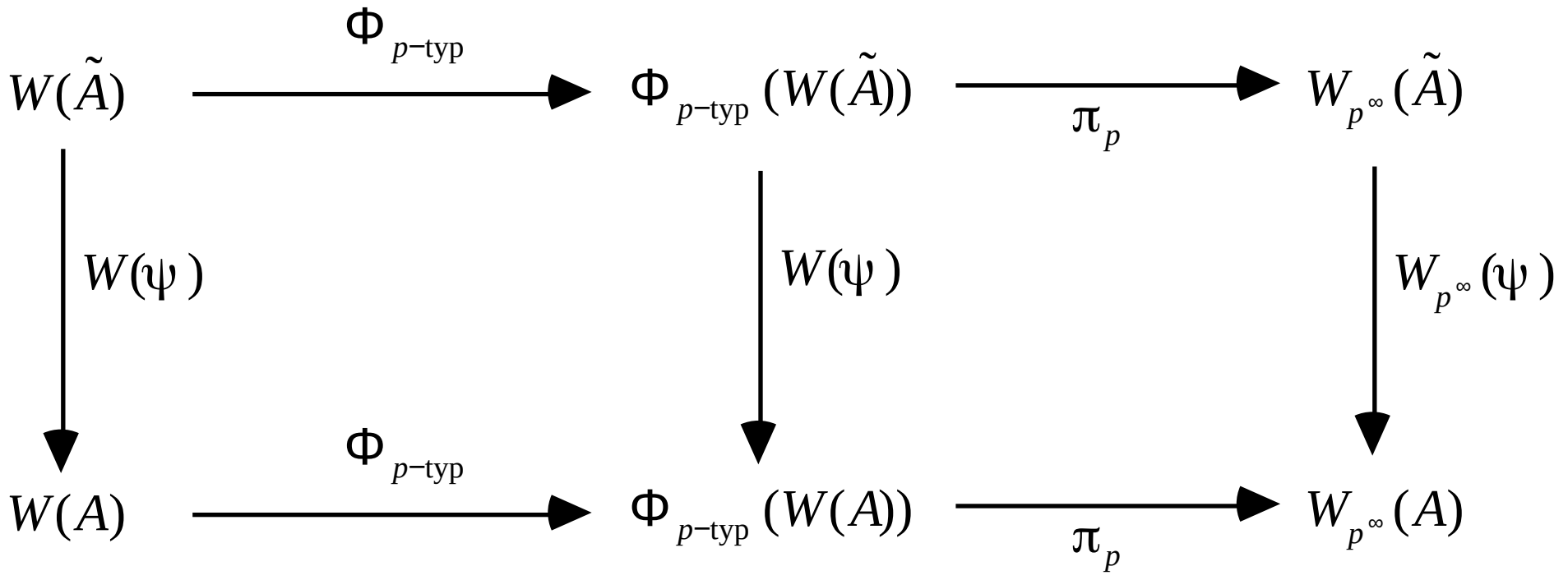

(14.23)

be any $\mathbf{Z}_{(p)}$-algebra and choose a charateristic zero $\mathbf{Z}_{(p)}$-algebra $\tilde{A}$ that covers it with corresponding projection $\psi: \tilde{A} \longrightarrow A$ and set $J = \mathrm{Ker}(\psi)$. Consider the commutative diagram (14.23) above.

It follows immediately from the fact that the upper right $\pi_p$ is an isomorphism that the lower right $\pi_p$ is surjective. Now let $x \in \Phi_{p-\mathrm{typ}}W(A)$, and take a $y \in W(A)$ that maps to $x$ and take a $\tilde{y} \in W(\tilde{A})$ that lifts it. Then $\tilde{x} = \Phi_{p-\mathrm{typ}}(\tilde{y})$ is a lift of $x$. Now suppose that $\pi_p(x) = 0$ (lower right hand corner in the diagram). Then $\pi_p(\tilde{x})$ has all its coordinates in $J$. By the commutativity of diagram (14.22) another, possibly different, lift of $x$ can be obtained by applying $\Phi_{p-\mathrm{typ}}$ to the vector in $W(\tilde{A})$ that has the same coordinates as $\pi_p(\tilde{x})$ at the $p$ power spots and zeros elsewhere. Let this vector be $\tilde{x}'$. But the morphism $\Phi_{p-\mathrm{typ}}$ is given by some universal polynomials with coefficients in $\mathbf{Z}_{(p)}$ and so all coordinates of $\tilde{x}'$ are in $J$, and so $x = W(\psi)\tilde{x}' = 0$. This proves that the lower $\pi_p$ in diagram (14.23) is also injective.

14.24. The morphism $\Phi_{p-\mathrm{typ}}$ (transferred to the power series case) is an example of '$p$-typification'. This is a general notion, due to Pierre Cartier, defined for curves in any commutative formal group. The one involved here is the second simplest example, namely the one for the one dimensional multiplicative formal group $\hat{\mathbf{G}}_m$. See [85] and [192], section 15.2.

14.25. *Decomposition.* It follows pretty directly from the preceding that over a $\mathbf{Z}_{(p)}$-algebra each Witt vector can be written uniquely as a sum

$$x = \sum_{n \in \mathbf{N}(p)} \mathbf{V}_n y_n, \quad y_n \in \Phi_{p-\mathrm{typ}}W(A) \tag{14.26}$$

of shifted $p$-typical Witt vectors. And thus for $\mathbf{Z}_{(p)}$-algebras the functor of the big Witt vectors decomposes (as an Abelian group) into a direct product of $\mathbf{N}(p)$ copies of the $p$-adic Witt vectors.

14.27. *Iterated nested Witt vectors.* Let $N$ and $M$ be two nests. These can be multiplied to form a new nest

$$MN = \{mn: m \in M, n \in N\} \subset \mathbf{N} \tag{14.28}$$

There is now, in [28], the theorem (observation) that if $M \cap N = \{1\}$ there is a canonical functorial isomorphism

$$W_{MN}(A) \longrightarrow W_M(W_N(A))$$

which essentially come from the Artin-Hasse exponential $W(-) \longrightarrow W(W(-))$ which is the essential part of the cotriple structure on the functor of the big Witt vectors, see subsections 16.43 and 16.59 below.

The case $M = \{1, p, p^2, p^3, \cdots\}$, $N = \mathbf{N}(p) = \{n \in \mathbf{N}: (n,p) = 1\}$, so that $MN = \mathbf{N}$ can already be found, more or less, in [336].

This is a good, actually better, substitute for the decomposition (14.26), and it works also for rings that are not $\mathbf{Z}_{(p)}$-algebras

## 15. Cartier algebra and Dieudonné algebra

Consider again the matter of Hopf algebra endomorphisms of **Symm**, or, equivalently, the matter of determining all functorial additive operations on the Witt vectors $W(A) \quad \Lambda(A)$. More generally consider the Hopf algebra $\mathbf{Symm}_R$ over a ring $R$ and consider the isomorphic functors $\Lambda$ and $W$ on the category of $R$-algebras and the additive operations on them, or, equivalently, the Hopf $R$-algebra endomorphisms of $\mathbf{Symm}_R$. As was already indicated these can be desribed by some infinity$\times$ infinity matrix with only finitely many nonzero entries in each row. See 13.74 above. There are other ways to see such a thing and to describe all Hopf algebra endomorphisms over $R$.

15.1. $(h-h)$-matrix of a Hopf $R$-algebra endomorphism of $\mathbf{Symm}_R$. One other way is as follows. Because $\mathbf{Symm}_R$ is free as an algebra over the ring $R$ an algebra endomorphism $\varphi$ of it is uniquely specified by a sequence of polynomials, e.g. the images $\varphi(h_n)$ of the free polynomial generators $h_n$. To be a Hopf algebra endomorphism the sequence of polynomials

$$d_0 = 1, d_1 = \varphi(h_1), d_2 = \varphi(h_2), d_3 = \varphi(h_3), \cdots \tag{15.2}$$

must form a curve (divided power sequence) in $\mathbf{Symm}_R$, which means that

$$\varepsilon_S(d_n) = 0 \quad \text{for } n \geq 1 \quad \text{and} \quad \mu_S(d_n) = \sum_{i+j=n} d_i \quad d_j \tag{15.3}$$

where $\mu_S$ and $\varepsilon_S$ are the sum comultiplication and cooresponding counit of $\mathbf{Symm}_R$. This means that the $d_n$ have constant terms zero and, I claim, are uniquely specified by the matrix

$$M_{(h-h)}(\varphi) = \left(\langle \varphi(h_i), h_j \rangle\right)_{i,j} \tag{15.4}$$

I like to call this matrix the $(h-h)$-matrix of the endomorphism $\varphi$. There are only finitely many nonzero entries in each row because the $\varphi(h_i)$ are polynomials (and hence bounded in weight).

The proof of the claim is simplicity itself. By the autoduality of $\mathbf{Symm}_R$ one has for example

$$\langle \varphi(h_n), h_k h_l \rangle = \left\langle \sum_{i+j=n} \varphi(h_i) \quad \varphi(h_j), h_k \quad h_l \right\rangle = \sum_{i+j=n} \langle \varphi(h_i), h_k \rangle \langle \varphi(h_j), h_l \rangle$$

and thus, with induction, the $(h-h)$-matrix determines the innerproducts of the polynomials $\varphi(h_n)$ with every monomial in the $h$'s and hence specifies the polynomials themselves completely. One can also of course take any other free polynomial basis for $\mathbf{Symm}_R$, for instance the elementary symmetric functions, which gives the $(h-e)$-matrix of an endomorphism.

It is perhaps interesting to calculate these matrices for some of the more important endomorphisms. For instance the Frobenius and Verschiebung endomorphisms. By definiton $\langle h_1, e_1 \rangle = 1$ and $\langle h_n, e_n \rangle = 0$ for $n \geq 2$. It follows with induction using the Newton relations, the fact that $\mathbf{f}_n$ is a ring endomorphism, and duality, that $\langle p_n, e_n \rangle = (-1)^{n+1}$ and hence that $\langle \mathbf{f}_n h_1, e_n \rangle = (-1)^{n+1}$. Using again the Newton relations, ring morphism, and duality and induction one further finds $\langle \mathbf{f}_n h_r, e_{nr} \rangle = 0$. Or use $\langle \mathbf{f}_r h_n, e_{nr} \rangle = \langle h_n, \mathbf{V}_r e_{nr} \rangle = \langle h_n, e_n \rangle$ to see this. As all other entries of the matrix must be zero by degree considerations it follows that the $(h-e)$-matrix of the Frobenius Hopf algebra endomorphism $\mathbf{f}_n$ of **Symm** has an entry $(-1)^{n+1}$ at spot $(1,n)$ and zeros everywhere else.

Even easier one finds that the $(h-e)$-matrix of the $n$-th Verschiebung operator has a 1 at spot $(n,1)$ and zeros everywhere else. This is encouraging and suggests that every matrix of the tye specified can arise. And this is indeed the case and can be proved this way. But this is not the easiest or most elegant way to see this.

15.5 *The DE-matrix of a Hopf algebra endomorphism of* $\mathbf{Symm}_R$. Let $R[T]$ be the algebra of polynomials in a single indeterminate over the ring $R$. The determining example for an additive operation on the functor $\Lambda_R$ (the functor $\Lambda$ restricted to the category of $R$-algebras) , and hence for the corresponding Hopf algebra endomorphism of $\mathbf{Symm}_R$, is the element

$$(1-Tt)^{-1} \quad \Lambda(R[T]) \tag{15.6}$$

in the sense that what a functorial operation does to this one example determines it completely. This claim will be examined a bit more closely further on. For the moment the argument is as follows. Let $\varphi$ be an additive operation and suppose that

$$\varphi((1-Tt)^{-1}) = 1 + a_1(T)t + a_2(T)t^2 + a_3(T)t^3 + \cdots \tag{15.7}$$

where the $a_n(T)$ are polynomials in $T$ with coefficients in $R$. Because the operation is supposed additive it must take the zero of the Abelian group $\Lambda(R[T])$, which is 1, into itself and so $a_i(0) = 0$ for all $i$. Then the naive argument goes at follows.

Take any element $b(t) \quad \Lambda(A)$ for an $R$-algebra $A$. Write it formally as a product

$$b(t) = 1 + b_1 t + b_{21} t^2 + b_3 t^3 + \cdots = \prod_{i=1}^{\infty} (1-\eta_i t)^{-1} \tag{15.8}$$

Then

$$\varphi(b(t)) = \prod_{i=1}^{\infty} (1 + a_1(\eta_i)t + a_2(\eta_i)t^2 + a_3(\eta_i)t^3 + \cdots) \tag{15.9}$$

This is more than a bit shaky. For one thing it is far from clear whether the product on the right hand side exists in a suitable sense. To make sense of things and also for other purposes rewrite the right hand side of (15.7 as a a product

$$\varphi((1-Tt)^{-1}) = \prod_{m,n\geq 1} (1 - c_{mn} T^n t^m)^{-1} \tag{15.10}$$

where for any $m$ there are only finitely many $n$ such that $c_{mn} \neq 0$. This can be done in precisely one way. The finiteness condition comes about precisely because the $a_i(T)$ in (15.7) are polynomials (not power series for instance). The matrix

$$M_{DE}(\varphi) = (c_{mn})_{m,n} \tag{15.11}$$

is the $DE$-matrix of $\varphi$ . [53].

15.12. *Proof of the DE principle* (*splitting principle*). Manifestly, by the definition of the Frobenius, Verschiebung and homothety operations the right hand side of (15.10) is the operation

$$\sum_{m,n} \mathbf{V}_m \langle c_{mn} \rangle \mathbf{f}_n \tag{15.13}$$

applied to the determining example power series $(1-Tt)^{-1}$. It remains to show that (15.13) makes unique sense when applied to an arbitrary element $b(t)$. To this end write $b(t)$ in Witt vector coordinates as

$$b(t) = \prod_d (1 - x_d t^d)^{-1} \tag{15.14}$$

Fix for the moment a power $m_0$ of $t$. Because of the finiteness condition there is an $n_0$ such that $c_{nm} = 0$ for all $m \le m_0$, $n \ge n_0$. Because of the definition of $\mathbf{V}_m$ the first $m_0$ coefficients of (15.13) applied to a finite product like (15.14) are determined by the part of (15.13) with $m \le m_0$. Now an $\mathbf{f}_n$ is given by a series of polynomials of weights $n, 2n, 3n, \cdots$. It follows that $\mathbf{f}_n$ applied to an $(1-xt^r)^{-1}$ with $r > m_0 n$ gives a result that has the first $m_0$ coefficients (except the constant term 1) equal to zero. Thus to calculate the first $m_0$ coefficients of (15.13) applied to (15.14) it suffices to apply it to the finite product

$$\prod_{d \le m_0 n_0} (1 - x_d t^d)^{-1} = \prod_i^{<\infty} (1 - \eta_i t)^{-1} \tag{15.15}$$

where the right hand side of (15.15) is also a finite product. And this is perfectly well defined by additivity and functoriality. This proves the $DE$ principle. This principle, apart form the infinity questions just taken care of, is an algebraic analogue of the splitting principle in algebraic topology as used in, say, topological K-theory [54].

15.16. *Cartier algebra.* By now it is clear that the ring of all additive operations on the functor is the set of all expressions (15.13) (with the stated finiteness condition on the coefficients involved). This ring, which is also the ring of Hopf $R$-algebra endomorphisms of the Hopf algebra $\mathbf{Symm}_R$, is called the Cartier algebra on $R$ and denoted $\mathrm{Cart}(R)$ [55].

---

[53] $DE$ stands for 'determining example'.

[54] Suppose one has an additive (in some suitable sense) operation in, say, topological $K$-theory. To find out what it does to a given bundle over a space take a suitable covering space such that the induced buncle splits into a sum of line bundles. This can always be done. Thus, by functoriality and additivity it suffices to know what the operation does to line bundles. And that in turn is governed (again by functoriality) by what it does to the universal line bundle over infinite dimensional projective space. From this point of view I should have called (15.6) the universal example (overworking the word 'universal').

The 'principle' that under suitable circomstances it suffices to verify properties just for the case of line bundles is (or at least used to be) called the 'verification principle'; see e.g. [27].

[55] When seen as endomorphisms of the Hopf algebra $\mathbf{Symm}_R$ the elements of $\mathrm{Cart}(R)$ need to be

There are also already a number of calculating rules, as given by the formulas from theorem 13.48 viz formulas (13.49)-(13.54), (13.57)-(13.59). The remaining question is whether these suffice. The most troublesome one (potentially) seems to be adding two expressions like (15.13). So consider a sum

$$\sum_{m,n} \mathbf{V}_m \langle b_{mn} \rangle \mathbf{f}_n + \sum_{m,n} \mathbf{V}_m \langle c_{mn} \rangle \mathbf{f}_n \tag{15.17}$$

Manipulating these directly according to the rules given, which needs especially (13.51), just seems to lead to more and more sums at first sight. More care is needed. To see how things work calculate what the sum (15.17) perscribes for the first few coeficients, say 2. That means that all terms with a $\mathbf{V}_m,\ m \geq 3$ can be discarded leaving

$$\sum_n \langle b_{1n} \rangle \mathbf{f}_n + \sum_n \mathbf{V}_2 \langle b_{2n} \rangle \mathbf{f}_n + \sum_n \langle c_{1n} \rangle \mathbf{f}_n + \sum_n \mathbf{V}_2 \langle c_{1n} \rangle \mathbf{f}_n \tag{15.18}$$

Now $\langle b_{1n} \rangle + \langle c_{in} \rangle = \sum_{i=1}^{\infty} \mathbf{V}_i \langle r_i(b_{1n}, c_{1n} \rangle \mathbf{f}_i$ and so the $\mathbf{V}_1$ part of the sum (15.17) is given by

$$\sum_n \langle b_{1n} + c_{1n} \rangle \mathbf{f}_n$$

(which is a finite sum). What is left is a sum

$$\mathbf{V}_2 \langle r_2(b_{1n}, c_{1n}) \rangle \mathbf{f}_2 + \sum_n \mathbf{V}_2 \langle b_{2n} \rangle \mathbf{f}_n + \sum_n \mathbf{V}_2 \langle c_{2n} \rangle \mathbf{f}_n$$

and so the $\mathbf{V}_2$ part of (15.17) is given by

$$\mathbf{V}_2 \langle b_{21} + c_{21} \rangle \mathbf{f}_1 + \mathbf{V}_2 \langle (r_2(b_{12}, c_{12}) + b_{22} + c_{22}) \rangle \mathbf{f}_2 + \sum_{n \geq 3} \mathbf{V}_2 \langle (b_{2n} + c_{2n}) \rangle \mathbf{f}_n$$

Continuing this way (which is tedious) one sees that the calculating rules given suffice to deal with sums.

Composing operations means dealing with products like

$$\mathbf{V}_m \langle b \rangle \mathbf{f}_n \mathbf{V}_r \langle c \rangle \mathbf{f}_s \tag{15.19}$$

Let $d$ be the greatest common divisor of $n$ and $r$. Then the product (15.19) is equal to

$$\begin{aligned} \mathbf{V}_m \langle b \rangle \mathbf{f}_n \mathbf{V}_r \langle c \rangle \mathbf{f}_s &= \mathbf{V}_m \langle b \rangle \mathbf{f}_{n/d} \mathbf{f}_d \mathbf{V}_d \mathbf{V}_{r/d} \langle c \rangle \mathbf{f}_s = \mathbf{V}_m \langle b \rangle \mathbf{f}_{n/d} [d] \mathbf{V}_{r/d} \langle c \rangle \mathbf{f}_s \\ &= [d] \mathbf{V}_m \langle b \rangle \mathbf{f}_{n/d} \mathbf{V}_{r/d} \langle c \rangle \mathbf{f}_s = [d] \mathbf{V}_m \langle b \rangle \mathbf{V}_{r/d} \mathbf{f}_{n/d} \langle c \rangle \mathbf{f}_s = [d] \mathbf{V}_{mr/d} \langle b^{r/d} c^{n/d} \rangle \mathbf{f}_{sn/d} \end{aligned}$$

So, back to (finite) sums again. The only remaining thing to check is that no infinite sums turn up when multiplying (composing) two things like (15.13) which is immediate.

So the given calculating rules suffice.

15.20. *Witt vectors as endomorphisms.* Quite generally if $\Omega$ is a unital-commutatve-ring-valued functor on **CRing** to itself its restriction to $R$-algebras has $\Omega(R)$ as part of its ring of additive endo operations. Indeed an $R$-algebra structure on a ring $A$ is the same thing as a morphism of rings $R \longrightarrow A$, which on applying $\Omega$ yields a morphism of rings $\Omega(R) \longrightarrow \Omega(A)$ which, in turn, for each $x \quad \Omega(R)$ an associated additive operation $(x, a) \mapsto xa$. [56]

So in the present case of the big Witt vectors $W(R)$ must be part of the ring af additive endo operations of $W_R \quad \Lambda_R$. It remains to identify these Witt vectors among the $\sum_{m,n} \mathbf{V}_m \langle c_{mn} \rangle \mathbf{f}_n$.

This is not difficult. Inside the algebras of operations on $\Lambda_R$ there are the operations of the form

$$\sum_{n=1}^{\infty} \mathbf{V}_n \langle x_n \rangle \mathbf{f}_n, \quad x_n \quad R \tag{15.21}$$

which take the determining example $(1 - Tt)^{-1}$ into the power series

$$\prod_{n=1}^{\infty} (1 - x_n T^n t^n)^{-1} \tag{15.22}$$

which is a power series with constant term 1 in the variable $Tt$ in Witt coordinate form. As addition of operations goes pointwise, the power series attached to the sum of two elements of the form (15.21) is the product of the corresponding power power series (15.22). It follows that the special operations (15.21) constitute a subgroup of all operations that is (isomorphic to) the additve group of the Witt vectors $W(R)$.

Now consider the composition of an operation (15.21) with an operation of the form $\langle \eta \rangle$ with associated power series $(1 - \eta Tt)^{-1} = \langle \eta \rangle (1 - Tt)^{-1}$. Applying $\langle \eta \rangle$ to (15.22) by definition gives

$$\prod_{n=1}^{\infty} (1 - x_n \eta^n T^n t^n)^{-1}$$

By the definition of the multiplication of Witt vectors, see (9.16) taking a product of a Witt vector $a(tT)$ with one of the form $(1 - \eta Tt)^{-1}$ is the same as substituting $\eta tT$ for $tT$. Thus composition of operations and Witt vector multiplication agree in this case. As should be the case for the operation of $\langle \eta \rangle$ to fit with the algebra structure $\Lambda(R) \longrightarrow \Lambda(A)$. By additivity and distributivity (both of Witt vector multiplication over Witt vector adition and operation composition over operation addition) it folllows that a product of Witt vectors of the form

$$(1 - x t^n T^n)^{-1} \quad \left(\prod_{n=1}^{\infty} (1 - x_n T^n t^n)^{-1}\right)$$

corresponds exactly to the compostion of the corresponding operations. A degree argument now finishes the proof that multiplication of Witt vectors corresponds to the compostion of operations of the form (15.21).

written in the form $\sum \mathbf{f}_n \langle c_{mn} \rangle \mathbf{V}_m$. See the 'caveat' 16.62 below. The need for the finiteness condition is perhaps even clearer in this interpretation.

[56] If the functor takes noncommutative rings as values, or more generally groups with an extra multiplication structure as values, there are two such operations, one on each side and these are additive if and only if resp. left distributivity or right distributivity holds for the extra multiplication structre over the group structure.

15.23. *Theorem.* The Witt vectors $W(R)$ under the correpondence (15.21)-(15.22) form a functorial subring of the rings of operations $\mathrm{End}(\Lambda_R) \quad \mathrm{End}_R(\mathbf{Symm}_R)$.

A somewhat remarkable fact in this case is that in the description 15.16 the subring of operations corresponding to $\Lambda(R)$ comes in Witt vector coordinates.

Still another way to see that the Witt vectors are part of $\mathrm{Cart}(R)$ will be briefly discussed further below.

15.24. *Further description of the Cartier algebra.* It is now natural to describe the ring of operations $\mathrm{Cart}(R)$ as an overring of $W(R)$ ([84]). Let $H$ be the skew polynomial ring over $W(R)$ in commuting indeterminates $\mathbf{f}_p$ , $p$ a prime number, with the commutation relation with Witt vectors given by $\mathbf{f}_p x = x^{(p)}\mathbf{f}_p$ where $x^{\mathbf{f}_p}$ de notes the Witt vector arising from applying the $p$-th Frobenius operation to the Witt vector $x$. Next take power series in the commuting indeterminates $\mathbf{V}_p$ with commutation relations

$$x\mathbf{V}_p = \mathbf{V}_p x^{(p)} \quad \text{and} \quad \mathbf{f}_p\mathbf{V}_{p'} = \mathbf{V}_{p'}\mathbf{f}_p \quad \text{for} \quad p \neq p'$$

Let this ring be S. In it take the ideal $J$ generated by the elements

$$\mathbf{f}_p\mathbf{V}_p - p, \quad \mathbf{V}_p x\mathbf{f}_p - x^{\mathbf{V}_p}, \quad x \quad W(R), \; p \text{ a prime number}$$

Then $\mathrm{Cart}(R) = S / J$ .

15.25. *The ring of additive endo operations of* $W_{p^\infty}$. The situation in the case of the ring valued functor of the $p$-adic Witt vectors is a bit different. There is no ‘determining example’ (as far as I can see) like in the case of the big Witt vectors and it is not true that the operations $\mathbf{V} = \mathbf{V}_p$, $\mathbf{f} = \mathbf{f}_p$ (the only surviving Frobenius and Verschiebung operations) together with the homothety operations suffice to generate all of them. In fact for $p \neq 2$ these are not even enough to produce the operation $[2]$ of adding a $p$-adic Witt vector to itself.

However, as above in 15.20 there are the $p$-adic Witt vectors $W_{p^\infty}(\mathbf{Z})$ as operations on $W_{p^\infty}$ and together with $\mathbf{V}$, $\mathbf{f}$ these do suffice to generate all additive operations. This is the topic of the next few pages. In this case I find it more convenient to work in terms of Hopf algebra endomorphisms of the representing Hopf algebra of the $p$-adic Witt vectors. This is the Hopf algebra

$$WH^{(p)} = \mathbf{Z}[X_0, X_1, X_2, \cdots] \qquad (15.26)$$ [57]

with addition and multipication comultiplications determined by the $p$-adic Witt polynomials as in section 5 above. For this part of this section 15, temporarily, the notations of section 5 will be used again.

15.27. *Endomorphisms of the infinite dimensional additive group.* To start with consider the infinite dimensional additive group. This is the functor

---

This applies for instance to the case of the functor defined by the ring of quasi-symmetric functions **QSymm**, see [200].

$$A \mapsto A^{\mathbf{N}}, \text{ represented by the ring } \mathbf{Z}[U] = \mathbf{Z}[U_0, U_1, U_2, \cdots] \tag{15.28}$$

with sum comultiplication and product comultiplication and corresponding counits

$$\mu_S(U_n) = 1 \quad U_n + U_n \quad 1, \ \varepsilon_S(U_n) = 0, \ \ \mu_P(U_n) = U_n \quad U_n, \ \varepsilon_P(U_n) = 1 \tag{15.29}$$

A Hopf algebra endomorphism of this Hopf algebra is given by a series of polynomials with coefficients in the integers

$$G(U) = (G_0(U), G_1(U), G_2(U), \cdots)$$

which are additive. This means that

$$G_n(U_1 + U_1', U_2 + U_2', U_3 + U_3', \cdots) = G_n(U_1, U_2, U_3, \cdots) + G_n(U_1', U_2', U_3', \cdots) \tag{15.30}$$

There are over a ring of characteristic zero, such as **Z**, not many such polynomials. The only ones are the linear ones, i.e. polynomials of the form

$$G_n(U) = \sum_i^{<\infty} a_{n,i} U_i \ , \ \ a_{n,i} \quad \mathbf{Z} \tag{15.31}$$

where, as indicated, the sums are finite.

15.32. *Determination of the Hopf algebra endomorphisms of* $WH^{(p)}$. Now let $\varphi$ be a Hopg algebra endomorphism of $WH^{(p)}$. It is given by a series of polynomials $(\varphi_0, \varphi_1, \varphi_0, \cdots)$. These, when substituted in the polynomials $w_0, w_1, w_2, \cdots$ define an Hopf algebra endomorphism of the infinite additive group and so are linear in the $w_i(X)$. Thus a series of polynomials $(\varphi_0, \varphi_1, \varphi_0, \cdots)$ defines a Hopf algebra endomorphism of $WH^{(p)}$ if and only if

$$w_n(\varphi(X)) = \sum_{i=0}^{\infty} a_{n,i} w_i(X) \tag{15.33}$$

for suitable integers $a_{n,i}$ of which only finitely many are nonzero for each $n$. Note that the integrality of the coefficients $a_{n,i}$ follows from (15.33). Just look at the coefficients of the powers of $X_0$ on the left and right hand side of (15.33). It also follows that the right hand side 'comes from' an integral vector over $WH^{(p)}$. This ring is of the type described in the 'ghost component integrality lemma' 9.93. The ring endomorphism $\varphi_p$ is the one that takes each $X_i$ to $X_i^p$. Here we are dealing with $p$-adic Witt vectors so the only denominators involved are powers of the prime number $p$. So the integrality criterion simplifies to the single sequence for the prime number $p$. So in the case at hand

$$\sum_{i=0}^{\infty} a_{n,i} w_i(X^p) \equiv \sum_{i=0}^{\infty} a_{n+1,i} w_i(X) \bmod p^{n+1} \tag{15.34}$$

is necessary and sufficient for the sequence $\varphi(X)$ as determined by (15.33) to define a Hopf algebra endomorphism of $WH^{(p)}$ and all Hopf algebra endomorphisms of $WH^{(p)}$ arise (uniquely) in this way.

Looking at the coefficients of the powers $X_0^{p^j}$ it follows from (15.34) that the integers $a_{n,i}$ satisfy

$$\begin{aligned} &a_{m,0} \equiv 0 \text{ mod } (p^m), \quad m = 1,2,3,\cdots \\ &a_{m,i} \equiv a_{m+1,i+1} \text{ mod } (p^{m+1}), \quad m,i \quad \mathbf{N} \quad \{0\} \end{aligned} \tag{15.35}$$

(and these congruences are also sufficient for (15.34) to hold). Now consider the following vectors over the integers

$$\begin{aligned} &b_m = (p^{-m}a_{m,0}, p^{-m}a_{m+1,1}, p^{-m}a_{m+2,2}, \cdots), \ m = 1,2,3,\cdots \\ &c_n = (a_{0,n}, a_{1,n+1}, a_{2,n+2}, a_{3,n+3}, \cdots), \quad n = 0,1,2,3,\cdots \end{aligned} \tag{15.36}$$

By (15.35) these satisfy the criterium of lemma (9.93) for being a $p$-adic ghost Witt vector. So there are Witt vectors

$$x_m, y_n \quad W_{p^\infty}(\mathbf{Z}) \quad \text{such that} \quad w(x_m) = b_m, \ w(y_n) = c_n \tag{15.37}$$

By definition the endomorphism $\mathbf{f}$ takes $w_r(X)$ into $w_{r+1}(X)$ and the endomorphism $\mathbf{V}$ takes $w_r(X)$ into $pw_{r-1}(X)$ and by the definition of the multiplication of Witt vectors

$$w(x \ {}_{W_{p^\infty}} \ (X_0, X_1, X_2, \cdots) = (w_0(x)w_0(X), w_1(x)w_1(X), w_2(x)w_2(X), \cdots)$$

for any $x \quad W_{p^\infty}(\mathbf{Z})$. So the endomorphism

$$\sum_{n=1}^{\infty} \mathbf{f}^n c_n \ + \ c_0 \ + \sum_{m=1}^{\infty} b_m \mathbf{V}^m \tag{15.38}$$

takes $w_r(X)$ into

$$\sum_{n=1}^{\infty} \mathbf{f}^n a_{r,n+r} w_r(X) + \ a_{r,r} w_r(X) \ + \sum_{m=1}^{r} (p^{-m} a_{r,r-m}) p^m w_{r-m}(X) \ = \ \sum_{i=1}^{\infty} a_{r,i} w_i(X)$$

exactly as in (15.33). Because a Witt vector $x \quad W_{p^\infty}(\mathbf{Z})$ starts with $j$ zeroes, i.e lies in the ideal $\mathbf{V}^j W_{p^\infty}(\mathbf{Z})$, if and only if its associated ghost vector starts with $j$ zeroes (as also follows from the ghost vector criterium of lemma 9.93), the finiteness condition of 15.33) translates into the statement that for each $r$ all but finitely many of the Witt vectors $c_n$ in (15.38) lie in the ideal $\mathbf{V}^r W_{p^\infty}(\mathbf{Z})$. Thus the ring of Hopf algebra endomorphisms of $WH^{(p)}$ consists of all expressions (15.38) with the finiteness condition just stated and with the calculating rules

$$\mathbf{Vf} = p, \quad \mathbf{V}x = x^{(p)}\mathbf{V}, \quad x\mathbf{f} = \mathbf{f}x^{(p)} \tag{15.39}$$

where $x^{(p)}$ denotes the result of applying the operation $\mathbf{f}$ to the $p$-adic Witt vector $x$ and where $p = [p]$ is the $p$-fold sum of the unit element in $W_{p^\infty}(\mathbf{Z})$.

15.40. The same technique can be used to determine the ring of endomorphisms of the Hopf algebra **Symm** as the representing Hopf algebra of the big Witt vectors, giving another way to obtain the ring $\mathrm{Cart}_{\mathbf{Z}}$.

15.41. The techniques of (15.32) as a method to determine the endomorphisms of the Hopf algebra $WH_R^{(p)}$ break down when $R$ is not a ring of characteristic zero for two reasons. One is that ghost component calulations as used are then not determining and second because in characteristic $p > 0$ there are other endomorphism of the additive groups than the linear ones; e.g. raising to the $p$-th power. Still I believe that the same picture holds. One reason is that the picture holds for the $WH_k^{(p)}$ when $k$ is a perfect field of characteristic $p > 0$. this will be described below. However, the proofs are quite different. Still every endomorphism in this case does come from one defined in characteristic zero (after the fact, i.e. after they have been determined in characteristic $p$). So there still appears to be a bit of work to do.

15.42. *Dieudonné algebras.* Let $k$ be a perfect field of characteristic $p > 0$. The Dieudonné algebra $\mathbf{D}_k$ consists of all expressions

$$\sum_{n=1}^{<\infty} \mathbf{f}^n c_n \; + \; c_0 \; + \; \sum_{m=1}^{<\infty} b_m \mathbf{V}^m, \quad c_n, b_m \quad W_{p^\infty}(k) \tag{15.43}$$

subject to the calculating rules

$$\mathbf{Vf} = \mathbf{fV} = p, \quad \mathbf{V}x = x^{(p)}\mathbf{V}, \quad x\mathbf{f} = \mathbf{f}x^{(p)} \tag{15.44}$$

Recall that in this case the Frobenius operation on $W_{p^\infty}(k)$ is given by

$$x^{(p)} = (x_0^p, x_1^p, x_2^p, x_3^p, \cdots)$$

and hence is an isomorphism, so that in (15.43) the Frobenius and Verschiebung symbols can be written on the left or the right as convenient. Also recall that (see 13.13) the sequence of polynomials that define the Frobenius endomorphism on $WK_k^{(p)}$ is

$$X_0^p, X_1^p, X_2^p, \cdots$$

so that these endomorphisms in characteristic $p$ take care of the 'raising to the power $p$' additive operations. One result in the present case is now (see [105], Ch. V, §1, No 3, p. 550):

15.45. *Theorem.* The endomorphim ring (over $k$) of the sub Hopf algebra $k[X_0, X_1, \cdots X_{n-1}] \quad WH_k^{(p)}$ is isomorphic to the ring $\mathbf{D}_k / \mathbf{D}_k \mathbf{V}^n$.

15.46. *Hirzebruch polynomials.* (Very partial symmetric function formularium (3)). Consider again an additive (functorial) operation on the groups of Witt vectors $\Lambda(A)$. Such an operation is given uniquely by a series of polynomials

$$1 = K_0, K_1, K_2, \cdots \quad \mathbf{Z}[e] = \mathbf{Z}[e_1, e_2, e_3, \cdots] \tag{15.47}$$

such that

$$\mu_S(K_n) = \sum_{i=0}^{n} K_i \quad K_{n-i}\,, \quad \text{where} \quad \mu_S(e_n) = \sum_{i=0}^{n} e_i \quad e_{n-i} \tag{15.48}$$

is the sum comultiplication on **Symm** $= \mathbf{Z}[h] = \mathbf{Z}[e]$[58]. If the $K_n$ are also homogeneous of degree $n$ (where $e_i$ has degree $i$ this is precisely what is defined in [213], Ch. I, §1, p. 9 as a multiplicative sequence.

The operation $\psi_K$ defined by the mulriplicative sequence (15.47) is of course

$$\psi_K\colon\ a(t) = 1 + a_1 t + a_2 t^2 + \cdots + \ \mapsto\ 1 + K_1(a)t + K_2(a)t^2 + \cdots \tag{15.49}$$

A multiplicative sequence is uniquely determined by what in loc. cit. is called its characteristic power series

$$Q(z) = 1 + K_1(1,0,0,\cdots)z + K_2(1,0,0,\cdots)z^2 + \tag{15.50}$$

which can be seen as the value of the operation on the determining example $1+zt$ because by the homogeneity condition on the $K_n$

$$\psi_K(1+zt) = Q(zt) \tag{15.51}$$

The homogeneity condition is also clearly equivalent to the property that in the DE formalism of (15.5) the DE-matrix of the operation is diagonal and so the operation on $\Lambda(A)$ in question is multiplication by a single Witt vector from $\Lambda(\mathbf{Z})$ (or $\Lambda(\mathbf{Q})$ if rational coefficients are allowed in the K (and in the characteristic power series). This single particular Witt vector is $Q(-t)^{-1}$.

The multiplicative sequence itself is recovered by considering the product

$$\prod_{i=1}^{\infty} Q(\xi_i t)$$

For each $n$ the coefficient of $t^n$ is a homogeneous symmetric polynomial in the $\xi$ and hence a homogeneous polynomial of weight $n$ in the $e_i$. These are the original multiplicative sequence polynomials by multiplicativity of the sequence or, equivalently, additivity of the operation defined by them.

Inversely, start with any power series $Q(t)$ over the integers or rationals and form the sum and product

$$s_Q(t) = \sum_{i=1}^{\infty} Q(\xi_i t), \quad p_Q(t) = \prod_{i=1}^{\infty} Q(\xi_i t) \tag{15.52}$$

where for the sum case it is assumed that $Q(0) = 0$ and for the product case $Q(0) = 1$. The coefficients of $t^n$ in (15.52) are symmetric in the $\xi$ and hence are polynomials in the elementary symmetric functions. These are the $n$-th additive and multiplicative Hirzebruch polynomials defined by the power series $Q(t)$.

[58] Elementary symmetric functions are used here in order to conform with usage in algebraic topology such as in [213] where these things were first introduced.

For example the Todd genus operation is defined by the power series

$$Q(z) = \frac{z}{1-e^{-z}} = 1 + \frac{1}{2}z + \sum_{n=1}^{\infty} (-1)^n \frac{B_n}{(2n)!} z^{2n} \tag{15.53}$$

and the Todd genus of a vectorbundle is obtained by applying this operation of the total Chern class $1 + c_1(V)t + c_2(V)t^2 + \cdots + c_n(V)t^n$ of a complex vectorbundle of dimension $n$.

The additive Hirzebruch polynomials serve to define functorial additive operations from $\Lambda(A)$ to $A^{\mathbf{N}}$. An example is the Chern character of a vector bundle defined by applying the operation associated to the power series $e^z$ to the total Chern class of a vector bundle [59].

15.54. *Open problem.* In subsection 15.5 (on the DE principle), formula 15.7 to be precise, an endomorphism of the Hopf algebra is given by an element of $\Lambda(\mathbf{Z}[T])$, the Witt vectors over the ring of polynomials in one indeterminate. One wonders how the multiplication of these as Witt vectors fits with the other bits of structure (composition and addition of endomorphisms).

## 16. More operations on the $\Lambda$ and $W$ functors: $\lambda$-rings.

There can be other operations than additive ones on an (Abelian) group valued functor that are important. This happens e.g. in algebraic topology and algebraic geometry with the operations on various functors induced by exterior powers (or, equivalently, symmetric powers) of vectorbundles. There are such functorial operations on the big Witt vectors. Indeed the Witt vectors $\Lambda(A)$ are functorial $\lambda$-rings [60]. In this setting the Witt vectors have three more universality properties (that are interrelated).

16.1. *Third universality property of the Witt vectors.* The functor $\Lambda$: $\mathbf{CRing} \longrightarrow \lambda-\mathbf{Ring}$ is right adjoint to the forgetful functor the other way. This means the following. First there is a canonical projection $\Lambda(A) \xrightarrow{\pi} A$ (which is in fact the functorial morphism $s_1$ (see (9.8), the first ghost component morphism) and second there is the universal property (that says that $\Lambda(A)$ is the cofree object over $A$) to the effect that for every morphism of rings $\varphi$: $S \longrightarrow A$ from a $\lambda$-ring to $A$ there is a unique lift, that is a morphism of $\lambda$-rings $\tilde{\varphi}$: $S \longrightarrow \Lambda(A)$ such that $\pi\tilde{\varphi} = \varphi$ as illustrated in the diagram below [61].

$$\begin{array}{ccc} & & \Lambda(A) \\ & \nearrow^{\tilde{\varphi}} & \downarrow \pi \\ S & \xrightarrow{\ \varphi\ } & A \end{array}$$

[59] Here $Q(0) = 1 \neq 0$ and so a finite factorization $1 + c_1(V)t + \cdots + c_n(V)t^n = \prod_{i=1}^{n} (1 + \beta_i t)$ is used.

[60] I do not know whether the notation '$\Lambda$' for the functor $A \mapsto 1 + tA[[t]]$ has anything to do with the fact that the $1 + tA[[t]]$ are functorial $\lambda$-rings. I believe not. However, $\Lambda$ is also a standard notation for **Symm** the representing object of this functor.

[61] This property rather nicely takes care of the complaint on page 1 of [263].

16.2. *Fourth universality property of the Witt vectors.* **Symm**, the representing object of the functor $\Lambda$, is the free $\lambda$-ring on one generator.

16.3. *Fifth universality property of the Witt vectors.* There is a comonad (= cotriple) structure on $\Lambda$ and the coalgebras for this comonad are precisely the $\lambda$-rings.

16.4. *Definition of $\lambda$-rings*. A $\lambda$-ring is a unital commutative ring that comes equipped with a an extra collection of operations

$$\lambda^i\colon A \longrightarrow A, \quad i = 1,2,3,\cdots \tag{16.5}$$

which 'behave like exterior powers'. In particular they are not additive for $i \geq 2$. The phrase 'behave like exterior powers' means like the exterior powers of vectorbundles over a space or like the exterior powers of representations of a group. Thus the first requirement is that

$$\lambda^n(a+b) = \lambda^n(a) + \sum_{i=1}^{n-1} \lambda^i(a)\lambda^{n-i}(b) + \lambda^n(b) = \sum_{i=0}^{n} \lambda^i(a)\lambda^{n-i}(b) \tag{16.6}$$

where by definition $\lambda^0(a) = 1$ for all $a \in A$. In addition the exterior product operations are supposed to satisfy some properties on products and when iterated in the sense that there are some universal polynomials which give $\lambda^n(xy)$ and $\lambda^m(\lambda^n(x))$ in terms of the $\lambda^i(x), \lambda^j(y)$. These polynomials are specified by defining exterior products on the rings $\Lambda(A)$ and declaring these particular rings with these particular exterior products to be $\lambda$-rings [62]. That goes as follows.

A morphism of $\lambda$-rings is a morphism of rings that respects the exterior product structure, i.e. a unital morphism of rings $\varphi : A \longrightarrow B$ such that

$$\varphi(\lambda_A^i(x)) = \lambda_B^i(\varphi(x)) \tag{16.7}$$

The exterior product structure on the $\Lambda(A)$ is defined by

---

[62] Lambda rings were introduced by Alexandre Grothendieck in the course of his investigations into Riemann-Roch type theorems, [48, 178, 179]. In spite of its attractive sounding title the reader is advised to steer clear of [148]. This book contains a great deal of most interesting material; it also contains nasty mistakes (as in the second paragraph of p. 15) and is sloppily written; see also the review by K R Coombes in Math. Rev. (88h:14011).

In the older literature, see e.g. loc. cit. and [27, 239] a lambda ring is a ring with exterior product operations satisfying just the additivity condition (16.6) and the term 'special lambda ring' was used for what here is called a lambda ring.

$$\lambda^m a(t) = \prod_{i_1 < i_2 < \cdots < i_m} (1 - \xi_{i_1}\xi_{i_2}\cdots\xi_{i_m} t)^{-1} \quad \text{when} \quad a(t) = \prod_{i=1}^{\infty} (1 - \xi_i t)^{-1} \tag{16.8}$$

By definition a $\lambda$-ring is a unital commutative ring $A$ equipped with exterior product operations $\lambda^i$ such that

$$\begin{array}{l} \lambda^1 = \mathrm{id} \text{ and the mapping } A \longrightarrow \Lambda(A) \text{ given by} \\ x \mapsto (1 - \lambda^1(x)t + \lambda^2(x)t^2 - \lambda^3(x)t^3 + \lambda^4(x)t^4 \cdots)^{-1} \end{array} \tag{16.9}$$

is a morphism of $\lambda$-rings (where $\Lambda(A)$ is given the exterior product structure (16.8)). Writing

$$\lambda_t(x) = 1 + \lambda^1(x)t + \lambda^2(x)t^2 + \lambda^3(x)t^3 + \cdots \tag{16.10}$$

the second part of(16.9) says that

$$x \mapsto \lambda_{-t}(x)^{-1}, \quad A \longrightarrow \Lambda(A) \tag{16.11}$$ [63]

is a morphism of $\lambda$-rings,while the first part says that $s_1(\lambda_{-t}^{-1}) = \mathrm{id}$.

Of course to make the whole thing consistent it has to be proved that $\Lambda(A)$ with the exterior product as specified by (16.8) is a $\lambda$-ring as defined via (16.9). That is it has to be proved that

$$\Lambda(A) \longrightarrow \Lambda(\Lambda(A)) \tag{16.12}$$

is a morphism of $\lambda$-rings. This will be done below, mixing in a number of lemmas which will also be useful later.

16.14. *$\lambda$-rings vs $\sigma$-rings*. Given a ring with exterior product operations $\lambda^n$ define the corresponding symmetric power operations $\sigma^n$ by the formula

$$\sum_{i=0}^{n} (-1)^i \lambda^i(x)\sigma^{n-1}(x) = 0, \quad \text{where} \quad \sigma^0(x) = 1 = \lambda^0(x) \tag{16.15}$$

If the exterior product operations satisfy (16.6) then so do the symmetric product operations. Also, obviously a morphism of rings respects exterior products if and only if it respects symmetric products. Write

$$\sigma_t(x) = 1 + \sigma_1(x)t + \sigma_2(x)t^2 + \sigma_3(x)t^3 + \cdots \tag{16.16}$$

---

[63] There are a couple minus signs here compared to the definition as given in e.g. [239, 324]. The reason for that is that in the present text (for good reasons) the ring structure on $\Lambda(A)$ is specified by the 'Hadamard like product' $(1-xt)^{-1} \quad (1-yt)^{-1} = (1-xyt)^{-1}$ with unit element $(1-t)^{-1}$ instead of $(1+xt)\ '(1+yt) = 1 + xyt$ with unit element $1+t$ (as in [239, 324). There are more good reasons for the particular choice made here (Which is also the choice made in the first treatments of this subject, [84, 86]) For instance, as already mentioned, it works out nicer for the autoduality of **Symm**). Still another good reason is that the multiplication $(1-xt)^{-1} \quad (1-yt)^{-1} = (1-xyt)^{-1}$ is how zeta functions of varieties over a finite field multiply in the sense that $\varsigma(X \times Y) = \varsigma(X) \quad \varsigma(Y)$, see [263], p. 2.

then by definition a ring is a $\sigma$-ring if and only if

$$\sigma_t:\ x \mapsto \sigma_t(x),\quad A \longrightarrow \Lambda(A) \tag{16.17}$$

is a morphism of $\sigma$-rings. The relation between the $\lambda^n$ and $\sigma^m$ can be succintly written

$$\sigma_t(x) = (\lambda_{-t}(x))^{-1} \tag{16.18}$$

and so $\sigma_t$ is precisely the mapping that figures in condition (16.11). Thus a $\sigma$-ring is exactly the same as a $\lambda$-ring. But the formulation of the property is just a tiny bit more elegant in terms of $\sigma$-rings (no minus signs).

It is a nice little exercise to show that the $\sigma$-ring structure on $\Lambda(A)$ is given more explicitely by

$$\sigma^m\Big(\prod_i (1-\xi_i t)^{-1}\Big) = \prod_{i_1 \le i_2 \le \cdots \le i_m} (1-\xi_{i_1}\xi_{i_1}\cdots\xi_{i_1} t)^{-1} \tag{16.19}$$

This is in fact the Wronski relation again between elementary and complete symmetric functions (when the complete symmetric functions are written out as the sum of all monomial symmetric functions of the same weight).

16.20. *Adams operations*. Given a ring with exterior product structure operations $\lambda^n$ the associated Adams operations [64] [65] are defined by

$$\sum_{n=1}^{\infty} \Psi^n(x)t^n = -t\frac{d}{dt}\log(\lambda_{-t}(x)) = t\frac{d}{dt}\log(\sigma_t(x)) \tag{16.21}$$

where the last equality of course follows from (16.18). Here also the formulation is just a bit more elegant in terms of the $\sigma$-structure.

16.22. *Theorem.* The Adams operations on the rings $\Lambda(A)$ with exterior product operations given by (16.8) or, equivalently, with symmetric product operations given by (16.19), are the Frobenius operations.

I like to call this observation the first 'Adams = Frobenius' theorem. The proof is easy. First take a power series of the form $x = (1-at)^{-1}$. Then $\lambda^1(x) = x,\ \lambda^i(x) = 0$ for $i \ge 2$ and hence

[64] The operations are named after J Frank Adams who first defined them in algebraic topology in the context of complex vector bundles and topological *K*-theory, [14].

[65] Consider a Hopf algebra $H$ and let $\mu_n$ and $m_n$ be the $n$-th iterates of its comultiplication and multiplication. Then the composite is the additive map (in general nothing more) [*n*], the n-fold sum of the identity (under convolution). In some of the literature, e.g. [141, 167, 270, 323] these maps are called Adams operations. That is a bit unfortunate as it does not fit in e.g. the case of **Symm**. For one thing if the Hopf algebra is graded these maps [*n*] are homogeneous (of degree zero) while the Addams operations are degree increasing. Also it does not fit with the case where the Adams operations came from, the cohomology of the classifying space **BU.**. This does not mean that these maps are not important. They are, see loc. cit.

$$\lambda_{-u}(x) = 1 - xu, \quad -u\frac{d}{du}\log(\lambda_{-u}(x)) = \frac{xu}{1-xu} = xu + x^2u^2 + x^3u^3 + \cdots$$

And so

$$\Psi^n((1-at)^{-1}) = ((1-at)^{-1})^n = (1-a^nt)^{-1} = \mathbf{f}_n((1-at)^{-1})$$

by the definition of product and Frobenius operations on $\Lambda(A)$. Now write an arbitrary element from $\Lambda(A)$ (formally) in the form

$$x = \prod_i (1-\xi_i t)^{-1}, \text{ i.e. } x = x_1 +_\Lambda x_2 +_\Lambda x_3 +_\Lambda \cdots, \quad x_i = (1-\xi_i t)^{-1}$$

in terms of the addition in a $\Lambda(\tilde{A})$. The same calulation as above gives

$$\Psi^n(x) = x_1^n +_\Lambda x_2^n +_\Lambda x_3^n +_\Lambda \cdots = \prod_i ((1-\xi_i t)^{-1})^n = \prod_i (1-\xi_i^n t)^{-1} = \mathbf{f}_n x$$ [66]

16.23. Also observe that by the definitions the Frobenius operations on the $\Lambda(A)$, and hence the Adams operations, commute with the exterior and symmetric power operations.

16.34. *Existence lemma* (*for exterior power* (*and other*) *structures*), [192], lemma 17.6.8, p.138. Let $A$ be a characteristic zero ring with ring endomorphisms $\varphi_n\colon A \to A$ for all $n \in \mathbf{N}$ such that $\varphi_1 = \mathrm{id}$, $\varphi_m\varphi_n = \varphi_{mn}$ and such that $\varphi_p(a) \equiv a^p \bmod (pA)$ for all prime numbers $p$ and $a \in A$. then there exists a unique mapping $D_A\colon A \to \Lambda(A)$ (resp. $D_A\colon A \to W(A)$) such that $s_nD_A = \varphi_n$ (resp. $w_nD_A = \varphi_n$). Moreover, this mapping is a ring morphism.

This is an immediate consequence of the ghost Witt vector integrality lemma 9.93.

16.35. *Lemma.* Let $\alpha\colon A \to B$ be a ring morphism and let both rings have exterior products $\lambda_A^n$, $\lambda_B^n$ with associated Adams operators $\Psi_A^n$, $\Psi_B^n$. Then if $\alpha\Psi_A^n = \Psi_B^n\alpha$ for all $n$, also $\alpha\lambda_A^n = \lambda_B^n\alpha$.

This comes about because the Adams operations and the exterior products are related by (universal) polynomials with rational coefficients.

Indeed these polynomials had better be given explicitely as they will also be needed further on. They are as follows

---

[66] Here, again, the splitting principle is used together with what has been called the ‘verification principle’, which says that under suitable circumstances it suffices to verify things for line bundles.

$$\Psi^n = \det\begin{pmatrix} \lambda^1 & 1 & 0 & \cdots & 0 \\ 2\lambda^2 & \lambda^1 & 1 & \ddots & \vdots \\ 3\lambda^3 & \lambda^2 & \ddots & \ddots & 0 \\ \vdots & \vdots & \ddots & \lambda^1 & 1 \\ n\lambda^n & \lambda^{n-1} & \cdots & \lambda^2 & \lambda^1 \end{pmatrix} \qquad (16.40)$$

$$n!\lambda^n = \det\begin{pmatrix} \Psi^1 & 1 & 0 & \cdots & 0 \\ \Psi^2 & \Psi^1 & 2 & \ddots & \vdots \\ \Psi^3 & \Psi^2 & \ddots & \ddots & 0 \\ \vdots & \vdots & \ddots & \Psi^1 & n-1 \\ \Psi^n & \Psi^{n-1} & \cdots & \Psi^2 & \Psi^1 \end{pmatrix} \qquad (16.41)$$

Some care must be taken in reading these formulas. E.g. in $\lambda^2 = \frac{1}{2}(\Psi^1)^2 - \frac{1}{2}\Psi^2$, $(\Psi^1)^2(x)$ must not be read as $\Psi^1(\Psi^1(x))$ but as $(\Psi^1(x))^2$.

These determinental formulas are exactly the same as those linking power sums and elementary symmetric functions in symmetric function theory ([281], p. 28), which is as must be because the defining formulas are the same in the two cases [67].

16.42. *Clarence Wilkerson theorem*. Let $A$ be as above in lemma 16.34. Then there is a unique $\lambda$-ring structure on $A$ such that the associated Adams operators are the given ring morphisms $\varphi_n$. Moreover the thus defined exterior powers commute with the given morphisms $\varphi_n$.

Proof. By the existence lemma 16.34 there is a unique ring morphism $\sigma_t\colon A \longrightarrow \Lambda(A)$ such that $s_n\sigma_t = \varphi_n$ for all $n$. This defines the exterior powers. Now observe that for all $a \in A$

$$s_m\sigma_t(\varphi_n(a)) = \varphi_m(\varphi_n(a)) = \varphi_{mn}(a)$$

and on the other hand, using the defining property of the Frobenius operation on the Witt vectors

$$s_m\mathbf{f}_n(\sigma_t(a)) = s_{mn}(\sigma_t(a)) = \varphi_{mn}(a)$$

Because $A$ is of characteristic zero (so that $s\colon \Lambda(A) \longrightarrow A^{\mathbf{N}}$ is injective) this implies that $\mathbf{f}_n\sigma_t = \sigma_t\varphi_n$. Thus by lemma 16.35 $\sigma_t$ respects the exterior powers and thus is a morphism of $\lambda$-rings as required.

16.43. *Corollary*. The rings of Witt vectors $\Lambda(A)$ are functorial $\lambda$-rings.

Proof. First suppose that A is of characteristic zero. Then so is $\Lambda(A)$. Now let $AH$ be the mapping $\Lambda(A) \longrightarrow \Lambda(\Lambda(A))$ defined by the given exterior product structure on $\Lambda(A)$. At this stage almost nothing is known about what properties it has. However, by definition of the Adams operations these satisfy $s_{n,\Lambda(A)} \circ AH = \Psi^n$. The operations $\Psi^n$ on $\Lambda(A)$ have been shown to be the Frobenius operations (the first Adams = Frobenius theorem) and these Frobenius operations satisfy the conditions of the existence lemma 16.34. So the uniqueness part says that $AH$ is equal to the morphism whose existence is guaranteed by that lemma and that one is a ring morphism by that lemma and respects exterior powers by theorem 16.42. This proves the corollary for characteristic zero rings. For an arbitrary ring take any characteristic zero cover. Everything in sight is functorial and so a little diagram chasing gives the result also in this case.

[67] It is because of this that the Adams operations are often called ‘power operations’.

16.44. *Remark.* The following diagram commutes for any $\lambda$-ring.

$$\begin{array}{ccc} \Lambda(A) & \xrightarrow{\ \mathbf{f}_n\ } & \Lambda(A) \\ \uparrow \sigma_t & & \uparrow \sigma_t \\ A & \xrightarrow{\ \Psi^n\ } & A \end{array}$$

For instance because $\mathbf{f}_n$ and $\Psi^n$ are given by the same polynomials in the exterior products (see (16.40) and the $\sigma_t$ are morphisms of $\lambda$-rings. In fact more generally if $A$ and $B$ are $\lambda$-rings and $A \xrightarrow{\ \alpha\ } B$ is a morphism of $\lambda$-rings then the morphism commutes with the respective Adams operations.

It follows from the commutativity of the above diagram that $\Psi^m\Psi^n = \Psi^{mn}$ and that the $\Psi^n$ respect the exterior powers.

16.45. *Remark.* The notation *'AH'* for the functorial morphism of $\lambda$-rings $\Lambda(A) \longrightarrow \Lambda(\Lambda(A))$, which could also be properly denoted $\sigma_{t,\Lambda(A)}$, stands for 'Artin-Hasse'. It is in fact a variant of the Artin Hasse exponential [68] in algebraic number theory. More precisely if one takes a finite field $k$, identifies $\Lambda(k) = W(k)$, takes the quotient $W(k) \longrightarrow W_{p^\infty}$, and sees the latter ring as the unramified mixed characteristic complete discrete valuation ring with residue field $k$ one finds a morphism of rings $W_{p^\infty}(k) \longrightarrow \Lambda(W_{p^\infty}(k))$ which is in fact the classical Artin-Hasse exponential, see [192], sections 17.5, 17.6 and E2.

By its definition the Artin-Hasse exponential $AH$ is such that the left hand diagram below commutes. It follows that the right hand diagram also commutes.

$$\begin{array}{ccc} & & \Lambda(\Lambda(A)) \\ & \overset{AH}{\nearrow} & \downarrow s_{n,\Lambda(A)} \\ \Lambda(A) & \underset{\mathbf{f}_n}{\searrow} & \Lambda(A) \end{array} \qquad\qquad \begin{array}{ccc} & & \Lambda(\Lambda(A)) \\ & \overset{AH}{\nearrow} & \downarrow \Lambda(s_{n,A}) \\ \Lambda(A) & \underset{\mathbf{f}_n}{\searrow} & \Lambda(A) \end{array}$$

16.46. *Remark.* The simplest example of a $\lambda$-ring is probably the ring of integers with the exterior power structure $\lambda_t(x) = (1+t)^x$. So in this case $\lambda_t$ is a kind of exponential, and that is a good way to think about it. In fact there situations where the presence of an exterior power structure does lead to exponential type isomorphims between the underlying ‘additive’ group of a ring and a multiplicative group of units of it, see e.g. [27], [388, 389].

[68] For the definition of the Artin-Hasse exponential in algebraic number theory and some of its uses see e.g. [412] to start with.

16.47. *The essence of a* $\lambda$*-ring structure*. Let $A$ be a $\lambda$-ring. Then it comes with with a morphism of $\lambda$-rings $\sigma_t\colon A \longrightarrow \Lambda(A)$, $\sigma_t(x) = (1-\lambda^1(x)t+\lambda^2(x)t^2-\lambda^3(x)t^3+\cdots)^{-1}$. That means that there are formulas for the exterior product $\lambda^n(x+y)$ of a sum of elements of A because $\sigma_t$ is additive; there are formulas for the exterior product $\lambda^n(xy)$ of a product of elements of $A$ because $\sigma_t$ is multiplicative; and there are formulas for the iterations of exterior products $\lambda^m(\lambda^n(x))$ because $\sigma_t$ is respects exterior products (as a morphism of $\lambda$-rings). These formulas are given in terms of the addition, multiplication, and exterior products on the Witt vector ring $\Lambda(A)$. But these formulas are universal [69]; they do not in any way depend on the ring $A$; they are certain polynomials with integer coefficients determined by certain manipulations with symmetric functions over the integers.

Thus there are formulas for $\lambda^n(x+y)$, $\lambda^n(xy)$, $\lambda^m(\lambda^n(x))$ as polynomials in the $\lambda^i(x)$ and $\lambda^j(y)$ and, which is the real essence of the story, these polynomials have their coefficients in the integers and are the same for all $\lambda$-rings $A$.

These polynomials can be calculated. Either by working directly with the one universal example of the Witt vectors $1+h_1(\xi)t+h_2(\xi)t^2+h_3(\xi)t^3+\cdots$ and $1+h_1(\eta)t+h_2(\eta)t^2+h_3(\eta)t^3+\cdots$ over $\mathbf{Z}[h(\xi);h(\eta)] \quad \mathbf{Z}[\xi;\eta]$, or, which I find easier, by using the determinantal formulas (16.40) and (16.41).

For example

$$\begin{aligned}
\lambda^3(xy) &= \frac{1}{6}(\Psi^1)^3(xy) - \frac{1}{2}(\Psi^1\Psi^2)(xy) + \frac{1}{3}(\Psi^3)(xy) \\
&= \frac{1}{6}\Psi^1(x)^3\Psi^1(y)^3 - \frac{1}{2}\Psi^1(x)\Psi^1(y)\Psi^2(x)\Psi^2(y) + \frac{1}{3}\Psi^3(x)\Psi^3(y) \\
&= \lambda^1(x)\lambda^2(x)\lambda^1(y)\lambda^2(y) + \lambda^1(x)^3\lambda^3(y) + \lambda^3(x)\lambda^1(y)^3 \\
&\quad -3\lambda^1(x)\lambda^2(x)\lambda^1(y)^3 - 3\lambda^1(x)^3\lambda^1(y)\lambda^2(y) + 3\lambda^3(x)\lambda^3(y)
\end{aligned}$$

$$\begin{aligned}
\lambda^2(\lambda^2(x)) &= \frac{1}{2}(\Psi^1)^2(\lambda^2(x)) - \frac{1}{2}\Psi^2(\lambda^2(x)) \\
&= \frac{1}{2}\lambda^2(x)^2 - \frac{1}{2}\Psi^2(\frac{1}{2}\Psi^1(x)^2 - \frac{1}{2}\Psi^2(x)) \\
&= \frac{1}{2}\lambda^2(x)^2 - \frac{1}{4}\Psi^2(x)^2 + \frac{1}{4}\Psi^4(x) \\
&= \lambda^1(x)\lambda^3(x) - \lambda^4(x)
\end{aligned}$$

16.48. $\Psi$*-rings*. A $\Psi$-ring is a ring together with family of ring endomorphisms $(\Psi^n)_{n\ \mathbf{N}}$ such that

$$\Psi^1 = \text{id},\quad \Psi^m\Psi^n = \Psi^{mn} \tag{16.49}$$

I shall say that $A$ is a ‘$\Psi$-ring with Frobenius morphism like property’ if moreover

[69] Overworking that unhappy word again; but I know of none other that meets the case.

$$\Psi^p(a) = a^p \bmod p \quad \text{for all prime numbers } p \tag{16.50}$$

Then if $A$ is a $\lambda$-ring the associated Adams operations turn it into a $\Psi$-ring by what was remarked in 16.44 above. It is moreover a $\Psi$-ring with Frobenius morphism like property. To see this it suffices to remark that that the Adams operations are given in terms of the symmetric powers (resp. the exterior powers) by the same polynomials that relate the power sum symmetric functions to the complete symmetric functions (resp. the elementary symmetric functions). Now $h_p \equiv h_1^p \bmod p$ and the statement follows. Another way to get it is to look at the defining equation for $\Psi^p$ in Witt coordinates:

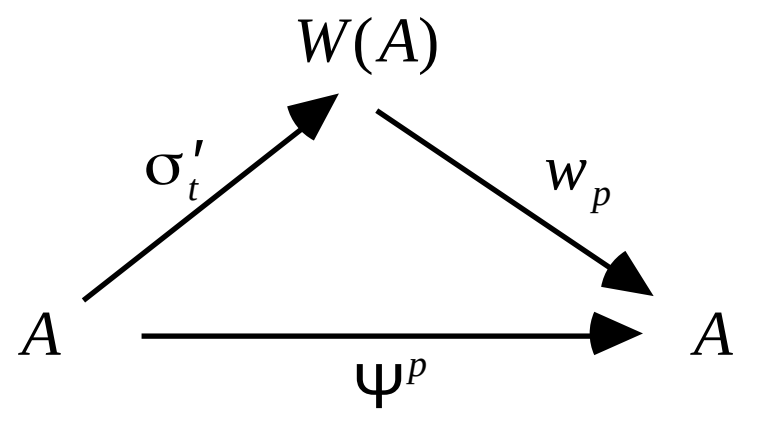

where $\sigma'_t$ is $\sigma_t$ followed by the isomorphism $e_A\colon \Lambda(A) \longrightarrow W(A)$. Because $w_1\sigma'_t = \mathrm{id}$, $\sigma'_t(a)$ is of the form $(a, x_2, x_3, \cdots)$ and so

$$\Psi^p(a) = w_p(a, x_2, x_3, \cdots) = a^p + px_p \equiv a^p \bmod p$$

16.51. Thus for a characteristic zero ring the statements "$A$ is a $\lambda$-ring' and '$A$ is a $\Psi$-ring with Frobenius morphism like property' are equivalent. If $A$ is not of characteristic zero this need not be the case.

16.52. Another result in this setting is that if $A$ is a characteristic zero ring with exterior power operations defined on it such that the associated Adams operations turn it into a $\Psi$-ring then it is a $\lambda$-ring.

Note the differences between the two theorems 16.51 and 16.52.

Thus the difference between $\Psi$-rings and $\lambda$-rings is an integrality matter (and such matters can be important). The idea of a $\Psi$-ring can easily be extended to a noncommutative setting (see below) and dualized. In view of the autoduality of **Symm** there should be some interesting notion of what a 'co-$\lambda$-coring' could be. Some sort of co-$\Psi$-coring with an extra integrality property.

16.53. There are a number of noncommutative rings in various parts of algebra that behave almost like $\Psi$-rings and $\lambda$-rings. Thus it makes sense to develop a theory of noncommutative $\Psi$-rings just like it was (and is) important to have a theory of noncommutative symmetric functions; the theory of **NSymm** briefly alluded to above in section 11.

The beginnings of a theory of noncommutative $\Psi$-rings have been established, see [325] (and [324] for a brief account). And it turns out that this theory relates nicely to the theory of noncommutative symmetric functions, loc. cit.

The next topic is the comonad structure on the functor of the Witt vectors.

16.54. *Monads and comonads.* A monad[70] in a category $\mathcal{C}$ is a triple $(T,m,e)$ consisting of an endofunctor $T\colon \mathcal{C} \longrightarrow \mathcal{C}$, a natural transformation $m\colon TT \longrightarrow T$ (where $TT$ stand for the iterate of $T$, i.e. $TT(A) = T(T(A))$) and a natural transformation $e\colon \mathrm{id} \longrightarrow T$ such that the following diagrams commute.

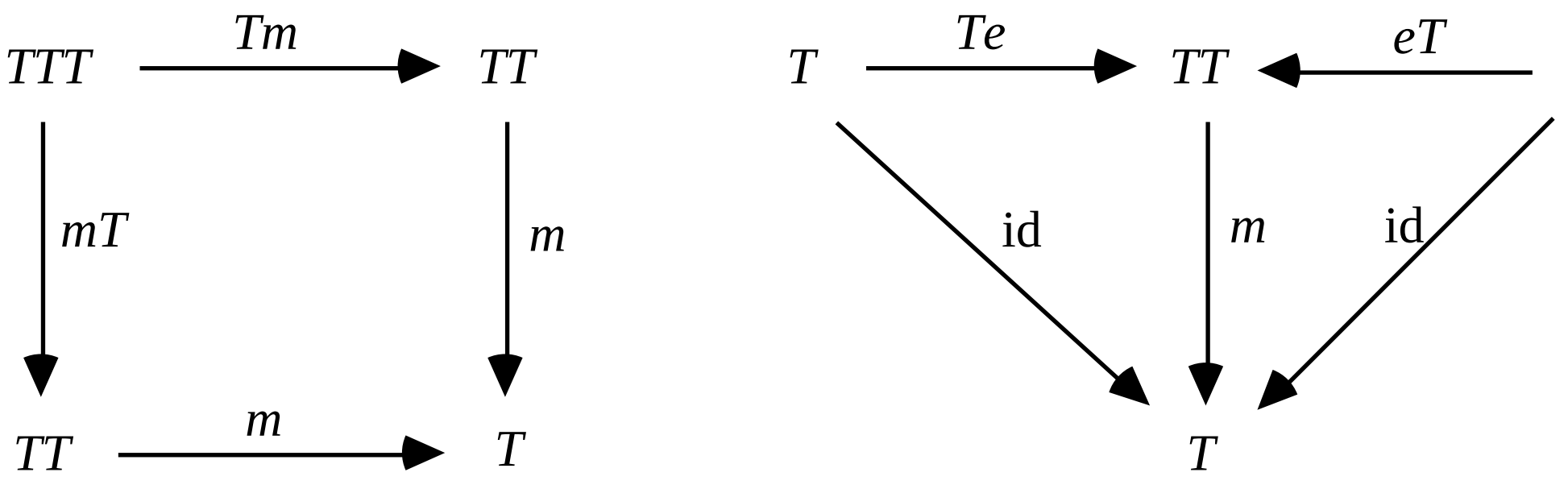

(16.55)

Here '*Tm*' means '*m*' first and then '*T*', i.e. at an object $C$ take the morphism $m_C\colon T(C) \longrightarrow TT(C) = T(T(C))$ and then apply the functor $T$ to it to obtain a morphism that is also often denoted $T(m_C)$; on the other hand '*mT*' at an object $C$, i.e. '*T*' first and then '*m*', is '*m*' at the object $T(C)$, often written $m_{T(C)}$.

An algebra (in the category $\mathcal{C}$) for a monad $(T,m,e)$, also called a $T$-algebra, is an object $A \in \mathcal{C}$ together with a morphism $\alpha\colon TA \longrightarrow A$ such that the following diagram commutes

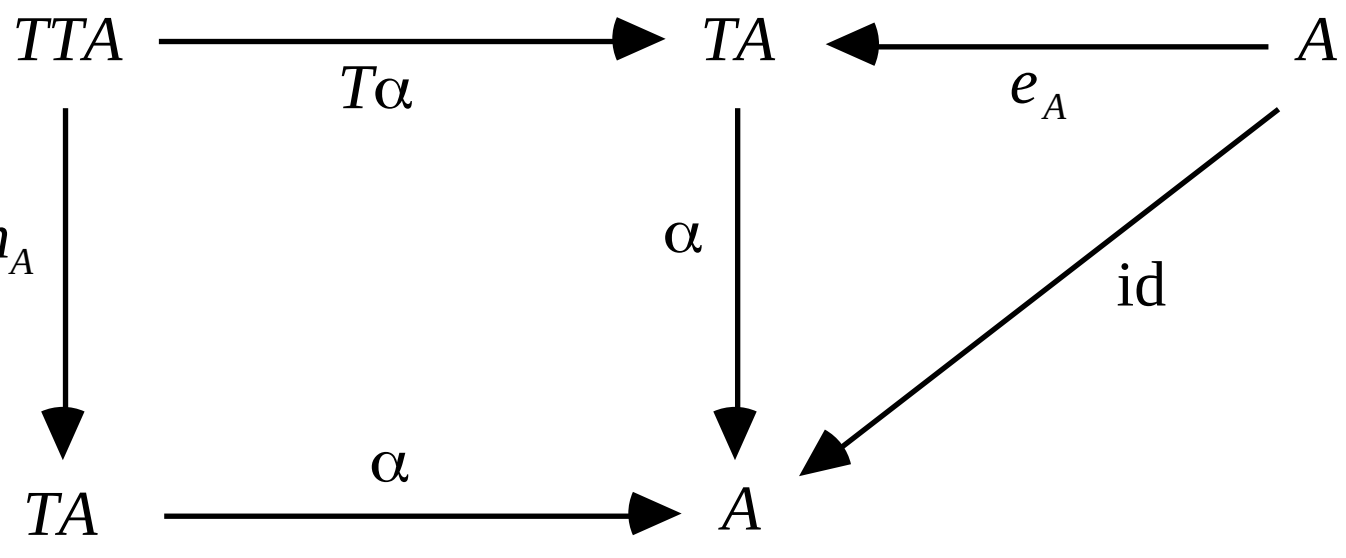

(16.56)

The opposite notions, i.e. the same notions in the opposite category, are those of a comonad and a coalgebra. Explicitely: a comonad $(T,\mu,\varepsilon)$ in a category $\mathcal{C}$ is an endo functor $T$ of $\mathcal{C}$ together with a morphism of functors $\mu\colon T \longrightarrow TT$ and a morphism of functors $\varepsilon\colon T \longrightarrow \mathrm{id}$ such that

$$(T\mu)\mu = (\mu T)\mu, \quad (\varepsilon T)\mu = \mathrm{id} = (T\varepsilon)\mu \qquad (16.57)$$

And a coalgebra for the comonad $(T,\mu,\varepsilon)$ is an object in the category $\mathcal{C}$ together with a morphism $\sigma\colon C \longrightarrow TC$ such that

[70] Other names for monad that are used (or have been used) in the literature are 'triad', standard construction', 'fundamental construction', 'triple'.

Monads first arose in work of Roger Godement, [171], in connection with the construction of simplicial objects and (standard) resolutions in connection with sheaf cohomology.

$$\varepsilon_C \sigma = \mathrm{id}, \quad (T(\sigma)\sigma = (\mu_{TC})\sigma \qquad (16.58)$$

The equations (16.57) and (16.58) can of course be written out in diagram form, the diagrams involved being the ones one gets by reversing all arrows in (16.55) and (16.56). This seems hardly worth the ink and paper in view of the fact that they will shortly appear explicitely in the case of the particular comonad of the Witt vectors and its coalgebras, which are precisely the $\lambda$-rings (or, better, $\sigma$-rings).

A standard text on monads and comonads is [33]; for a thorough up-to-date account see [284] [71].

16.59. *Comonad structure on the Witt vectors.* The claim is that the Witt vector functor $\Lambda$ (or $W$) is a comonad in the category **CRing**, with the comonad morphism given by the Artin-Hasse exponential and the first ghost component morphism and that the coalgebras for this comonad are precisely the $\lambda$-rings.

To prove this it must be shown that for each ring $A$ the following diagrams (16.60) are commutative and that an exterior product structure $\sigma_t\colon A \longrightarrow \Lambda(A)$ gives a $\lambda$-ring structure if and only if the diagrams (16.61) below are commutative.

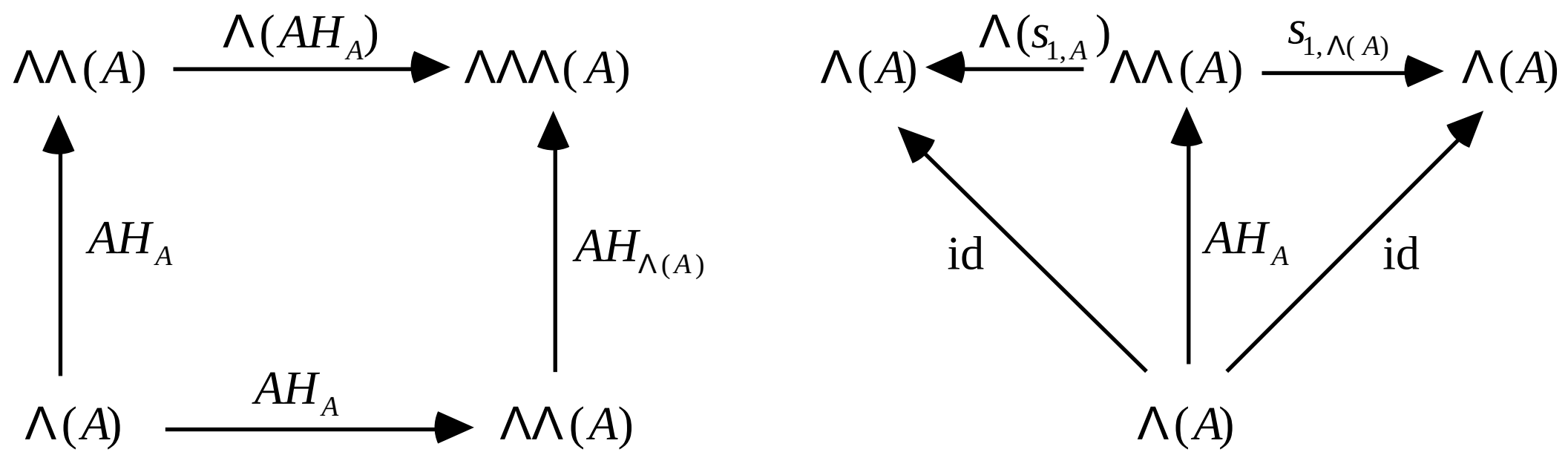

(16.60)

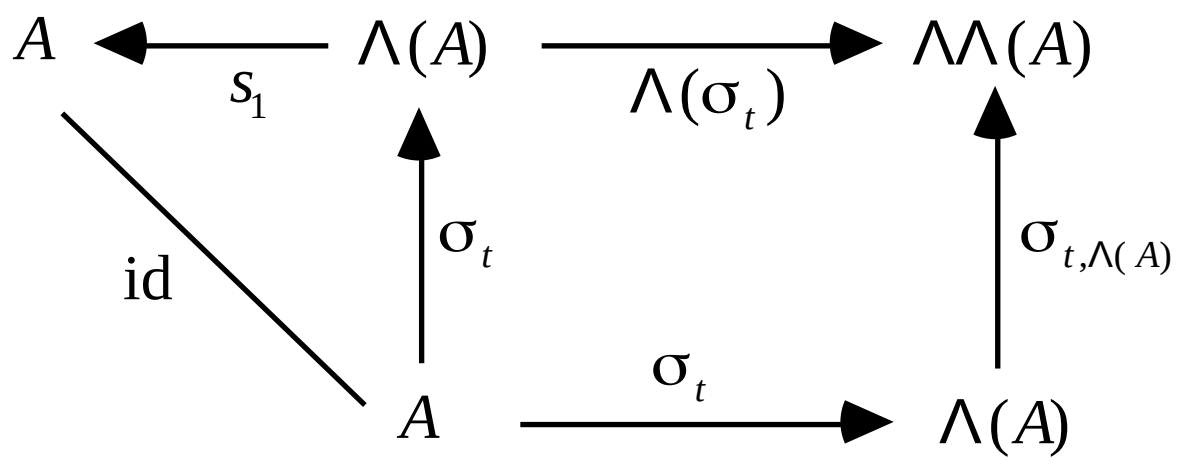

(16.61)

This looks like there is a fair amount of work to do. Especially the triple iterate at the upper right hand corner of the left diagram in (16.60) could seem a bit intimidating. Actually everything has already been proved. Start with the commutativity of the right hand part of diagram (16.61). This is just saying (in a slightly fancy way) that $\sigma_t$ is a morphism of $\lambda$-rings (with the fact that it is a ring morphism coming from the requirement that things take place in the category of rings). The left hand part of (16.61) just says that $\lambda^1 = \mathrm{id}$.

As $AH$ is the morphism determined by the exterior power structure on $\Lambda(A)$ (see the proof of 16.43) the commutativity of the left diagram of (16.60) is just proving that $\Lambda(A)$ is a $\lambda$-ring. The remainder of (16.60) is part of the definition of $AH$ (right half) and a consequence of that (see 16.45). The right half of the right diagram of (16.60) is also the first part of the requirement that $\Lambda(A)$ be a $\lambda$-ring.

16.62. *Theorem. Cofreeness of the Witt vectors.* The $\lambda$-ring $\Lambda(A)$ of Witt vectors over $A$,

[71] And a, in my opinion, inspiring account.

together with the first ghost component morphism $s_1 : \Lambda(A) \longrightarrow A$ is the cofree $\lambda$-ring over $A$.

This means that if $S$ is a $\lambda$-ring and $\alpha: S \longrightarrow A$ is a morphisms of rings (i.e forget the exterior product structure on $S$ for the moment), then there is a unique morphism of $\lambda$-rings $\tilde{\alpha}: S \longrightarrow \Lambda(A)$ that lifts $\alpha$. This is illustrated in the following diagram.

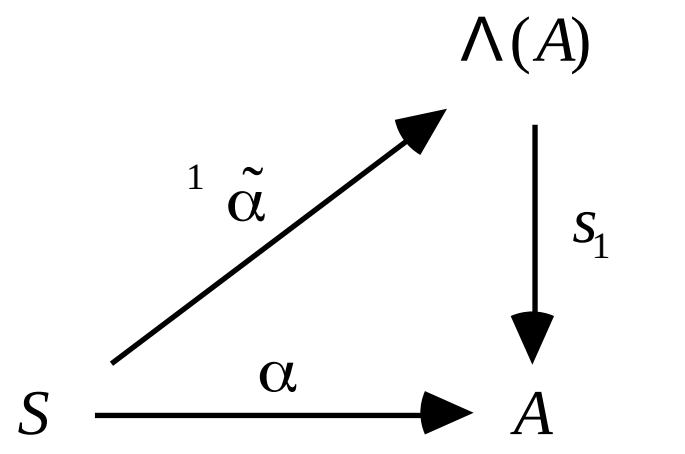

(16.63)

By now this is simple to prove. First remark that if

$$a(t) = 1 + a_1 t + a_2 t^2 + a_3 t^3 + \cdots = \prod_i (1 - \xi_i t)^{-1} \qquad \Lambda(A)$$

then by the definition of the $\lambda$-ring structure on $\Lambda(A)$ (see (16.8)

$$\lambda^n a(t) \equiv 1 + e_n(\xi)t \equiv 1 + ((-1)^{n+1}(h_n(\xi) + P(h_1(\xi), \cdots, h_{n-1}(\xi)))t \mod t^2$$

for some polynomial $P$ with integer coefficients. And so

$$s_1(\lambda^n a(t)) = (-1)^{n+1} a_n + P(a_1, \cdots, a_{n-1}) \tag{16.64}$$

Now let $\alpha: S \longrightarrow A$ be a morphism of rings. Define $\tilde{\alpha}: S \longrightarrow \Lambda(A)$ as the composition $S \longrightarrow \Lambda(S) \xrightarrow{\Lambda(\alpha)} \Lambda(A)$. Then this $\tilde{\alpha}$ does the job, and, with induction, it is also the only possibility because it is required to be a morphism of $\lambda$-rings and (16.64) [72].

16.65. $\lambda$-ring structure on **Symm** itself. There is a very simple and obvious $\lambda$-ring structure on the ring $\mathbf{Z}[\xi] = \mathbf{Z}[\xi_1, \xi_2, \xi_3, \cdots]$ of polynomials over the integers, viz the one determined by

$$\lambda^1(\xi_i) = \xi_i, \ \lambda^j(\xi_i) = 0 \text{ if } j \geq 2 \tag{16.66}$$

The easiest way to prove that there is indeed such a $\lambda$-ring structure is to consider the ring endomorphisms given by

$$\Psi^n(\xi_i) = \xi_i^n \tag{16.67}$$

These ring endomorphisms satisfy all the requirements of the Wilkerson theorem 16.42 and so there is a $\lambda$-ring structure on $\mathbf{Z}[\xi]$ for which the ring endomorphisms (16.67) are the associated

[72] Again the proof would be more elegant if things were formulated for $\sigma$-rings.

Adams operations. The relations between the Adams operations and the exterior product operations are the usual ones linking power sums and elementary symmetric functions. In particular

$$n!\lambda^n = \det\begin{pmatrix} \Psi^1 & 1 & 0 & \cdots & 0 \\ \Psi^2 & \Psi^1 & 2 & \ddots & \vdots \\ \Psi^3 & \Psi^2 & \ddots & \ddots & 0 \\ \vdots & \vdots & \ddots & \Psi^1 & n-1 \\ \Psi^n & \Psi^{n-1} & \cdots & \Psi^2 & \Psi^1 \end{pmatrix} \tag{16.68}$$

and it follows immediately that the $\lambda$-operations (exterior powers) in casu are indeed those of (16.66).

Now consider the subring of the symmetric functions $\mathbf{Z}[h]$   $\mathbf{Z}[\xi]$. Obviously this subring is stable under the ring endomorphisms (16.67), so applying the Wilkerson theorem16.42 again there is a $\lambda$-ring structure on $\mathbf{Z}[h]$ for which the morphisms (16.66) are the associated Adams opererations which is obviously the restriction of the $\lambda$-ring structure (16.65) to the symmetric functions.

16.69. $\lambda$*-ring structure on* **QSymm**. In between the symmetric functions and $\mathbf{Z}[\xi]$ there sit the quasi-symmetric functions $\mathbf{Z}[h] = \mathbf{Symm}$   $\mathbf{QSymm}$   $\mathbf{Z}[\xi]$. Obviously, **QSymm** is stable under the ring morphisms (16.67) and so there is a corresponding $\lambda$-ring structure on **QSymm**, which is the restriction of the one on $\mathbf{Z}[\xi]$ and which extends the one just discussed on **Symm**. The associated Adams operations are the ring morphisms (16.67) which are the same as the Frobenius endomorphisms of 13.33.

16.70. *Lyndon words.* The proper tails of a word $\alpha = [a_1, a_2, \cdots, a_m]$, $a_i$   $\mathbf{N}$, over the natural numbers are the words $[a_i, a_{i+1}, \cdots, a_m]$, $1 < i \le m$ . A word is Lyndon if it is lexicographically smaller than each of its proper tails. For instance $[1,3,2]$ is Lyndon and $[2,1,1]$ and $[1,1]$ are not Lyndon. A word $\alpha$ is primitive if the greatest common divisor of its entries is 1.

16.71. *Theorem,* [195, 198] (*Free polynomial generators for* **QSymm** *over the integers*). The elements $\lambda^n\alpha$, $n$   $\mathbf{N}$, $\alpha$ primitive word, form a free polynomial basis over the integers of the ring **QSymm**.

This has been a vexing problem since 1972 which was finally solved in 1999, see [193]. **QSymm** is a strong candidate for a theory of noncommutative Witt vectors (as a tool for classifying noncommutative formal groups amoung other things).

16.72. *Discussion.* The many different (?) operations on **Symm**. There are by now some five potentially different unary operations on **Symm** and it is perhaps wise to list them.

(a) The exterior product operations that define the $\lambda$-ring structure on **Symm**. These are not additive of course.

(b) The corresponding Adams operations as described by (16.66)

(c) The rings $\Lambda(A)$ are functorial $\lambda$-rings. The functor $A \mapsto \Lambda(A)$ is represented by **Symm**: $\Lambda(A)$ = **CRing**(**Symm**, $A$). The functorial exterior product operations on $\Lambda(A)$ must therefore come from ring endomorphisms of **Symm**. These endomorphism cannot be Hopf algebra endomorphisms (because otherwise the functorial exterior product operations on $W(A)$ would be additive.

(d) The rings $\Lambda(A)$ have functorial Adams operations. These also must come from ring endomorphisms of **Symm**. Moreover these must be Hopf algebra endomorphisms and even coring object morphisms in the category of rings.

(e) The Frobenius endomorphisms as defined in subsection 13.6. These are the ones that induce the functorial Frobenius operarions on the $\Lambda(A)$.

The first Adams = Frobenius theorem, 16.22, says that the morphisms (d) and (e) are the same. By the remarks already made the only one which could conceivable also be the same is (b). And this is indeed the case

16.73. *Theorem.* (Second Adams = Frobenius theorem). The Adams operations $\Psi^n$ coming from the $\lambda$-ring structure on **Symm** are the same as the endomorphisms $\mathbf{f}_n$ of **Symm** that induce the Frobenius endomorphisms of the functor $\Lambda(A)$ of the big Witt vectors.

Proof. This follows immediately from the fact that $\mathbf{f}_n(p_m) = p_{nm}$ compared with (16.66).

There is of course still more structure on **Symm**. **Symm** being autodual there are also all the duals of (a)-(e). It is not evident what all these dual operations are.

16.74. *Theorem.* Universal $\lambda$-ring on one generator. The ring of symmetric polynomials with the $\lambda$-ring structure defined above is the universal $\lambda$-ring on one generator [73].

This means the following. For each $\lambda$-ring $A$ and element $a \in A$ there is a unique morphism of $\lambda$-rings $\varphi : Symm \longrightarrow A$ such that $\varphi(e_1) = a$.

Proof. As is easily verified from e.g. (16.58) the $\lambda$-ring structure on **Symm** satisfies $\lambda^n(e_1) = e_n$. It follows that the only ring morphism that could possibly work is defined by $\varphi$: $e_n \mapsto \lambda^n(a)$. And there is such a ring morphism because **Symm** is free on the $e_n$. That this is actually a morphism of $\lambda$-rings requires a bit more work, as follows. Let $x$ be an element of **Symm**, i.e. a polynomial in the elementary symmetric functions

$$x = P_x(e_1, e_2, \cdots) = P_x(\lambda^1(e_1), \lambda^2(e_1), \cdots)$$

Now consider $\lambda^n(x)$. Because composition of lammda operations, a lambda operation applied to a product, and a lambda operation applied to a sum, are given by 'universal polynomials', that means the same polynomials for any $\lambda$-ring, see 16.47, there is a universal polynomial $Q_{n,P_x}$ (with coefficients in the integers) [74] such that for any $\lambda$-ring and any element $a$ in it

[73] There is a far reaching generalization. As will be discussed below **Symm** is isomorphic to $RS$ the direct sum of the representation rings of the symmetric groups. For a fixed finite group $G$ let $G \times_{wr} S_n$ be the wreath product of $G$ and $S_n$. The direct sum of the representation rings $R(G \times_{wr} S_n)$ is the free lambda ring on the irreducible representations of $G$ as generators, [280].

[74] This polynomial is in fact the plethysm $e_n \circ P_x$, see subsection 16.76 below.

$$\lambda^n(P_x(\lambda^1(a),\lambda^2(a),\cdots)) = Q_{n,P_x}(\lambda^1(a),\lambda^2(a),\cdots)$$

Also $\varphi$ is a ring homomorphism and so commutes with polynomials

$$\varphi(Q(x_1,x_2,\cdots)) = Q(\varphi(x_1),\varphi(x_2),\cdots)$$

It follows that $\varphi$ commutes with the lambda operations, so that it is a morphism of $\lambda$-rings.

16.75. *Comonadability vs cofreeness vs representability vs freeness.* The various notions that passed under review in the previous few subsections are far from being unrelated. The following matters were discussed.

(a) The functor $\Lambda: \mathbf{CRing} \longrightarrow \lambda-\mathbf{Ring}$ takes a ring $A$ into the cofree $\lambda$-ring over $A$.

(b) The functor $\Lambda$ comes with a comonad structure and the coalgebras af this comad are precisely the $\lambda$-rings.

(c) The functor $\Lambda$ is representable. The representing ring is a $\lambda$-ring and is the free $\lambda$-ring on one generator.

The relation (a) - (b) is a well known part of monad and comonad theory (from the early days of this theory). Every monad or comonad comes from an adjoint pair of functors, see [284, 33], and of course cofreeness of $\Lambda(A) \longrightarrow A$ translates into the observation that $\Lambda$ is right adjoint to the functor $\lambda$-**Ring** $\longrightarrow$ **CRing** that forgets about the exterior product structure.

Here is how the relation (b) - (c) goes[75]. Let $\mathcal{C}$ be a category and let $(T,\mu,\varepsilon)$ be a comonad in $\mathcal{C}$. Now let $(Z, z \in T(Z))$ represent the functor $T$. That is, there is a functorial bijection $\mathcal{C}(Z,A) \longrightarrow T(A)$, $f \mapsto T(f)(z)$. The monad structure gives in particular a morphism $\sigma: Z \longrightarrow TZ$, viz the image of $\mathrm{id}_Z$ under $\mu_Z: T(Z) = \mathcal{C}(Z,Z) \longrightarrow T(T(Z)) = \mathcal{C}(Z,T(Z))$. This defines a 'coalgebra for $T$' structure on Z. Now let $(A,\sigma)$ be a coalgebra for the comonad $T$ and let $a$ be an element of $A$. Consider the element $\sigma(a) \in T(A) = \mathcal{C}(Z,A)$. This gives a unique morphism of $T$-coalgebras that takes $z$ into $a$. There are of course a number of things to verify both at this categorical level and to check that these categorical considerations fit with the expliit constructions carried out in the previous subsections. This is straightforward.

16.76. *Plethysm.* (Very partial symmetric function formularium (4)). Given an element $a \in A$ in a $\lambda$-ring A and a polynomial $f \in$ **Symm** define

$$\alpha_a(f) = \beta_f(a) = f(\lambda^1 a, \lambda^2 a, \cdots) \qquad (16.77)$$

where in the last part of (16.77) $f$ is seen as a polynomial in the elementary symmetric functions; i.e. in the expression for $f$ as a polynomial in the elementary symmetric functions these elementary symmetric functions are replaced by the exterior powers of the element $a$. If $f$ is seen as a polynomial in the complete symmetric functions the same result is obtained by repacing these with the symmetric powers of the element $a$. This can lead to confusion so it is better to have a description that does not depend on what free polynomial basis is used for **Symm**. That goes as follows. Given $a \in A$ and $f \in$ **Symm**

[75] I have not found this relation in the monad literature (which does not mean it is not there somewhere). There should be a general theorem to this effect also in a more general context than 'categories of sets with structure'.

Let $\alpha_a$: **Symm** $\longrightarrow A$ be the unique morphism of $\lambda$ - rings such that $\phi_a(h_1) = a$. Then $\alpha_a(f) = \beta_f(a)$ (16.78)

This looks at $a \in A$ as something that defines (and is defined by) a morphism of $\lambda$-rings from the free $\lambda$-ring on one generator to $A$.

Another way to look at (16.77) is as the definition of a functorial operation on $\lambda$-rings defined by the polynomial $f$. This is the most used way to look at (16.77), (16.78).

In particular **Symm** itself is a $\lambda$-ring, and so, taking a $g \in$ **Symm** for $a$ in (16.77), (16.78) there results a new kind of composition of polynomials

$$f \circ g = \alpha_g(f) = \beta_f(g) = f(\lambda^1 g, \lambda^2 g, \cdots) \tag{16.79}$$

This composition law on **Symm** is called plethysm. More precisely it is called ‘outer plethysm’. There is also something called ‘inner plethysm’ which has to do with the $\lambda$-ring structures on the homogeneous summands **Symm**$^{(n)}$ of **Symm** [76]. In the representation theoretic incarnation of **Symm** in terms of representations of the general linear groups, outer plethysm corresponds to composition of representations.

In terms of the $\lambda$-ring interpretation mentioned above it also coresponds to composition: the functorial operation on $\lambda$-rings defined by $f \circ g$ on $\lambda$-rings is the composition of the functorial operations defined by $f$ and $g$:

$$\beta_{f \circ g}(a) = \beta_f(\beta_g(a)) \tag{16.80}$$

Plethysm is associative

$$(f \circ g) \circ h = f \circ (g \circ h) \tag{16.81}$$

This is a special case of (16.80).

16.80. *Calculation of plethysms.* Write $g \in$ **Symm** **Z** as a sum of monomials in the $\xi$ $g = \sum_\tau c_\tau \xi^\tau$ and write (formally)

$$\prod (1 + y_i t) = \prod_\tau (1 + \xi^\tau t)^{c_\tau} \tag{16.81}$$

then

$$f \circ g = f(y_1, y_2, \cdots) \tag{16.82}$$

where $f$ is seen as a symmetric polynomial in the $\xi$.

Thus for instance $f \circ p_n = f(\xi_1^n, \xi_2^n, \cdots)$, a plethysm that was used in 9.93 above when discussing the Witt symmetric functions which give the Witt coordinates of a power series.

[76] These homogenous summnds are rings under the second multiplication on **Symm**. In the representation theoretic incarnation of **Symm** they correspond with the ring of (virtual) representations $R(S_n)$ of the symmetric group $S_n$. The exterior powers there are the exterior powers of representations.

To see that (16.82)-(16.83) fits with the definition, first note that for the $\lambda$-ring structure on $\mathbf{Z}[\xi]$, $\lambda^1(\xi_i)=\xi_i$, $\lambda^n(\xi_i)=0$ for $n\geq 2$. Thus $\sigma_t(\xi_i)=(1-\xi_i t)^{-1}$ and because $\sigma_t$ is a morphism of rings

$$\sigma_t(g)=\prod_{\tau}(1-\xi^{\tau}t)^{-c_{\tau}} \tag{16.83}$$

and so

$$\lambda_t(g)=\sigma_{-t}(g)^{-1}=\prod_{\tau}(1+\xi^{\tau}t)^{c_{\tau}}=\prod(1+y_i t)$$

so that the elementary symmetric functions in the $y$ are the $\lambda^n g$.

16.84. *Distributivity of the functorial plethysm operations.* The plethysm operations $\beta_f$ on $\lambda$-rings $A$ are not additive. But they are distributive over addition and multiplication in a Hopf algebra way as follows. Let

$$\mu_S(f)=\sum_i f'_{S,i}\quad f''_{S,i}\ ,\quad \mu_P(f)=\sum_i f'_{P,i}\quad f''_{P,i} \tag{16.85}$$

Then

$$\beta_f(a+b)=\sum_i f'_{S,i}(a)f''_{S,i}(b)\ ,\quad \beta_f(ab)=\sum_i f'_{P,i}(a)f''_{P,i}(b) \tag{16.86}$$

This is seen as follows. given $a,b\in A$, they define Witt vectors $\sigma_t(a),\sigma_t(b)\in\Lambda(A)$ which, in turn, are morphisms of rings **Symm** $\longrightarrow A$. Because $\sigma_t$ is a morphism of rings the Witt vectors of $a+b$ and $ab$ are the Witt vector sum and Witt vector product of the Witt vectors of $a$ and $b$. As morphisms of rings **Symm** $\longrightarrow A$. Witt vectors are summed and multiplied as follows. The sum Witt vector corresponds to the composite

$$\mathbf{Symm}\longrightarrow\mathbf{Symm}\otimes\mathbf{Symm}\longrightarrow A\otimes A\longrightarrow A$$

and the product Witt vector corresponds to the composite

$$\mathbf{Symm}\longrightarrow\mathbf{Symm}\otimes\mathbf{Symm}\longrightarrow A\otimes A\longrightarrow A$$

Formula (16.86) follows.

For those polynomials for which there are nice formulas for the sum comultiplication morphism $\mu_S$ or for the product comultiplication $\mu_p$ (16.85) and (16.86) can give useful formulas for calculating plethysm operations. For instance, using 10.7

$$s_{\kappa}\circ(f+g)=\sum_{\lambda\subset\kappa}(s_{\lambda}\circ f)(s_{\kappa/\lambda}\circ g) \tag{16.87}$$

It also follows (again) from (16.85), (16.86) that the plethysm operations $\beta_{p_n}$ on $\lambda$-rings defined by the power sum symmetric functions $p_n$ are both additive and multiplicative; they are of course the Adams operations (explaining once more why these are often called "power operations'.

Rather few of the plethysm operations defined by a polynomial are additive: only those that come from a sum of power sum symmetric functions. On the other hand there are all the additive operations generated by the Frobenius, Verschiebung and homothety operations (the Cartier ring). It is unknown whether all these together with the plethysm operations generate all operations on the functor $\Lambda$.

There is quite a bit of literature on plethysm, mostly in a representation theoretic context and much is on how to calulate it in special cases. A sampling is [8, 66, 82, 282, 313, 394, 395], [312, 306, 43, 42, 429].

## 17. Necklace rings.

There is (for suitable rings) a third coordinatization of the unital power series over a ring $A$, i.e. the elements of $\Lambda(A) = 1 + tA[[t]]$, besides the power series coordinates and Witt vector coordinates considered so far. These coordinates go by the name necklace coordinates. The three systems of coordinates are related by

$$1 + a_1 t + a_2 t^2 + a_2 t^2 + \cdots = \prod_{i=1}^{\infty} (1 - x_i t^i)^{-1} = \prod_{i=1}^{\infty} (1 - t^i)^{-c_i} \qquad (17.1)$$

Unlike in the case of Witt vector coordinates it is not always possible to find necklace coordinates. It can certainly be done for the case that $A$ is a $\mathbf{Q}$-algebra, but there is but little interest in that as in that case

$$W(A) \quad \Lambda(A) \quad Gh(A) = A^{\mathbf{N}}$$

But there are quite a few rings $A$ that are not $\mathbf{Q}$-algebras for which such a representation in necklace coordinates can always be found with the $c_i \quad A$; in particular the rather important case $A = \mathbf{Z}$. This will be discussed a bit further below. On the other hand this can not be done for the case that $A$ is a ring of polynomials, in particular not for $A = \mathbf{Symm}$ which is the context in which the universal example of a one power series lives; that is the power series

$$h(t) = 1 + h_1 t + h_2 t^2 + h_3 t^3 + \cdots \qquad (17.2)$$

from which all others are obtained

Why I like to call the coordinates on the right hand side of (17.1) by the name necklace coordinates will become apparent in a minute.

17.3. *Necklace polynomials.* Consider a totally ordered alphabet of $\alpha$ letters, say $\{1,2,\cdots,\alpha\}$. A word over this alphabet is primitive (also called aperiodic) if it is not a concatenation power of a strictly smaller word. A word is Lyndon if it is strictly smaller in te lexicographic order than any of its (non-identity) cyclic permutations. So in particular it is primitive. For the equivalence of this definition of Lyndon word with the one used in 16.71 see e.g. [273], section 4.4.

A necklace (also called circular word) is an equivalence class under cyclic permutations of words. A necklace is primitive if the words in the equivalence class are primitive. Thus the Lyndon words can be regarded as a systematic choice of representatives of primitive necklaces.

The formula for the number of primitive necklaces of length $n$ on $\alpha$ letters, i.e. with $\alpha$ different colors of beads, has been known since Colonel Moreau, [298], in 1872. It is

$$M(\alpha;n) = n^{-1}\sum_{d|n}\mu(d)\alpha^{n/d} \tag{17.4}$$

These expressions, seen as polynomials in an indeterminate $\alpha$ are known as the necklace polynomials. Note that though these polynomials are integer valued for every integer argument they do not have integer coefficients. For instance for a prime number $p$

$$M(\alpha;p) = p^{-1}(\alpha^p - \alpha) \tag{17.5}$$

The same expression (17.4) turns up in other contexts. For instance in the theory of free Lie algebras. Consider the free Lie algebra generated by $\alpha$ symbols. Give each symbol weight one. The free Lie algebra (over the integers or any field) is then graded and the graded part of weight $n$ has rank $M(\alpha,n)$. This can be seen for instance by the socalled Hall set construction of a basis (as an Abelian group) for the free Lie algebra, using the set of Lyndon words as a Hall set. See Ch. 5 in [334] for details. In this context of free Lie algebras (17.4) is known as the Witt formula, [418].

In loc. cit., p. 153, Witt writes: “Es ist merkwürdig, dass diese Rangformel übereinstimmt mit der bekannten Gausschen Formel für die Anzahl der Primpolynome $x^n + a_1x^{n-1} + \cdots + a_n$ im Galoisfeld von $q$ Elementen.” Later, Solomon Golomb, [174], found indeed a bijection between primitive necklaces and irreducible polynomials. This correspondence is not yet entirely satisfactory in that it depends on the choice of a primitive element for the Galois field. See also [335] for a discussion of this correspondence and other occurences of expression (17.4).

Consider the power series in Witt coordinates

$$\sum_{n\geq 0}\beta_n t^n = \prod_{n\geq 1}(1-M(\alpha,n)t^n)^{-1} \tag{17.6}$$

Then $\beta_n$ is the rank of the homogeneous component of degree $n$ of the free associative algebra over $\mathbf{Z}$ in $\alpha$ symbols. This results from the Witt formulas by using the Poincaré-Birkhoff-Witt theorem combined with the Milnor-Moore type result that says that the free associative algebra is the universal enveloping algebra of the free Lie algebra on the same set of symbols. For some further manifestations of the necklace polynomials and related expressions in such varied fields of inquiry as the Feyman identity, the elliptic modular function, multiple zeta values, (symbolic) dynamics and fixed points, formal groups, see also [70, 101, 113, 218, 234, 261, 260, 292, 299, 316, 317].

17.7. *Necklace polynomial formulas.* Here are two interesting formulas from [291].

$$M(\alpha\beta;n) = \sum_{[i,j]=n}(i,j)M(\alpha;i)M(\beta;j) \tag{17.8}$$

where $(i,j)$ denotes the greatest common divisor of $i$ and $j$, and $[i,j]$ denotes their least common multiple.

$$M(\beta^r;n) = \sum\frac{j}{n}M(\beta;j) \tag{17.9}$$

where the sum ranges over all $j$ such that $[j,r]=nr$ (which implies that $j$ is a multiple of $n$).

17.10. *Cyclotomic identity.* The socalled cyclotomic identity is

$$\frac{1}{(1-\alpha t)} = \prod_{n=1}^{\infty} \left( \frac{1}{1-t^n} \right)^{M(\alpha;n)} \tag{17.11}$$

(It is a tradition to write it precisely this way.) All three identities (17.8), (17.8), (17.11) are identities about polynomials. As such it suffices to prove them for integer values of the variables $\alpha, \beta$ which is (usually) done by combinatorial means. Cf e.g. [290, 291].

There is a very elegant 'symmetric' generalization of the cyclotomic identity due to Volker Strehl, [375]

$$\prod_{n=1}^{\infty} \left( \frac{1}{1-\alpha t^n} \right)^{M(\beta;n)} = \prod_{n=1}^{\infty} \left( \frac{1}{1-\beta t^n} \right)^{M(\alpha;n)} \tag{17.12}$$ [77]

There are some other generalizations which will be described later in a section devoted to generalized Witt vectors.

17.13. *Motivational remarks regarding the necklace algebra functor.* Now consider power series expressed by means of necklace coordinates as in (the right hand side of) 17.1. Multiplying two such expression means adding their necklace coordinates. The multiplication defined on the Witt vector ring is determined by

$$\frac{1}{(1-\alpha t)} \quad \frac{1}{(1-\beta\, t)} = \frac{1}{(1-\alpha\beta\, t)} \tag{17.14}$$

Taking into account distributivity, the cyclotomic identity (17.11) and formula (17.8) this dictates the following definition of the necklace algebra (and its ghost components morphism).

17.15. *Definition of the necklace algebra functor.* As an Abelian group the necklace ring over a ring $A$ is the infinite product $A^{\mathbf{N}}$ of all sequences $(a_1, a_2, a_3, \cdots)$, $a_i \in A$ with component wise addition. Two such sequences are multiplied according to the rule

$$\begin{aligned} &(a_1, a_2, a_3, \cdots) \quad (b_1, b_2, b_3, \cdots) = (c_1, c_2, c_3, \cdots) \\ &c_n = \sum_{[i,j]=n} (i,j) a_i b_j \end{aligned} \tag{17.16}$$

This is clearly functorial and this functor will be denoted $Nr$: **CRing** $\longrightarrow$ **CRing**. There is a functorial ghost components ring morphism $u\colon Nr(A) \longrightarrow Gh(A) = A^{\mathbf{N}}$ to the direct product of $\mathbf{N}$ copies of the ring $A$ defined by

$$u(a) = (u_1(a), u_2(a), u_3(a), \cdots), \quad u_n(a) = \sum_{d|n} d a_d \tag{17.17}$$

Given a power series for which necklace coordinates exist, so that (17.1) holds, one then has for these coordinates

---

[77] This has the flavour of a 'reciprocity formua' and one wonders if it relates to some other 'reciprocity results' in mathematics.

$$w(x) = s(a) = u(c) \tag{17.18}$$

in the ghost ring $Gh(A) = A^{\mathbf{N}}$ (with component-wise addition and multiplication).

Necklace rings (in the sense defined above [78] are a special kind of convolution rings as defined in [403, 402, 405]. For some first investigations in the algebraic theory of necklace rings see [404, 322].

17.19. *Binomial rings*. Obviously, from (17.1), necklace coordinates are hardly compatible with the presence of torsion.

Further as

$$(1-t^n)^x = 1 - xt^n + \frac{x(x-1)}{2!}t^{2n} - \frac{x(x-1)(x-2)}{3!}t^{3n} + \cdots \tag{17.20}$$

it seems necessary for necklace coordinates to exist over a ring $A$ to require that together with an $x \in A$ it also contains the binomial coefficients $(n!)^{-1}x(x-1)\cdots(x-n+1)$. This leads to the idea of a binomial ring.

A (commutative unital) ring $A$ is said to be binomial if it is torsion free (as a **Z**-module) and if together with any $x$ in it it also contains the binomial coefficients

$$\binom{x}{n} = \frac{x(x-1)(x-2)\cdots(x-n+1)}{n!} \tag{17.21}$$

These rings have made their appearance before in 1958, in the fundamental work of Philip Hall on the theory of nilpotent groups, [186], often referred to as the 'Edmonton notes".

Obviously any **Q**-algebra is binomial. But there are many rings that are binomial that are not algebras over the rationals. First and foremost certain rings of functions with values in the integers (under pointwise addition and multiplication), notably polynomials.

17.22. *Free binomial rings*. Let $X = \{X_i : i \in I\}$ be a set (of indeterminates). By definition $IVal[X]$ is the ring of all polynomials with rational coefficients that take integer values on integers; i.e integer valued polynomials. These are also sometimes called numerical polynomials.

Examples are the binomial coefficients polynomials

$$\binom{X_i}{n} = \frac{X_i(X_i-1)(X_i-2)\cdots(X_i-n+1)}{n!} \tag{17.22}$$

and it is a theorem that the monomials in these form a basis of $IVal[X]$ over the integers (as an Abelian group), see [45], §45, p. 240ff.

These are indeed the free binomial rings in the technical sense that the functor

[78] There is another kind of object in this more or less same corner of algebra that has 'necklace' in its name. Viz necklace Lie algebras, necklace Hopf algebra, necklace Lie coalgebra, see [Bocklandt, 2002 #402; Gan, 2007 #401; Ginzburg, 2006 #400; Schedler, 2004 #399]. These have little to do with the necklace approach to Witt vectors.

*IVal*: **Set** $\longrightarrow$ **BinRing** (where **BinRing** stands for the sub category of binomial rings in **CRing**) is left adjoint to the forgetful functor the other way, see [134], p. 168, propositiom 2:

$$\mathbf{BinRing}(IVal[X],A) \quad \mathbf{Set}(X,A) \tag{17.23}$$

A monograph on integer valued polynomials is [80] [79].

Other examples of binomial rings are the $p$-adic integers, the profinite completion $\hat{\mathbf{Z}}$ of the integers and localizations of the integers, and a special kind of $\lambda$-rings.

17.24. *Binomial rings vs $\lambda$-rings*. Let $A$ be a torsion free ring such that $x^p \equiv x \bmod p$ for all prime numbers $p$. Then , taking $\varphi_n = \mathrm{id}$ for all $n$ the conditions of the Wilkerson theorem 16.42 are satisfied and so there is a $\lambda$-ring structure on $A$ for which all the Adams operations are the identity. For this $\lambda$-ring structure

$$\lambda^n(x) = \binom{x}{n}, \quad \lambda_t(x) = (1+t)^x, \quad \sigma_t(x) = (1-t)^{-x} \tag{17.25}$$

and thus $A$ is a binomial ring. To see (17.25) simply calculate the determinant of the matrix $M_x = (m_{i,j})$, $m_{i,j} = x$ for $i \geq j$, $m_{i,i+1} = i$, $m_{i,j} = 0$ for $j \geq i+2$, see (16.41). Actually one can do a little better and prove that a $\lambda$-ring for which all the Adams operations are equal to the identity is automatically torsion free, see [134], section 5.

Inversely let $A$ be a binomial ring then taking the (tentative) exterior product operations to be the binomial coeffients one gets a ring with exterior products with associated Adams operations that are all the identity and hence a $\lambda$-ring structure.

17.26. *Frobenius and Verschiebung on necklace rings*. From the above it is clear that for binomial rings necklace coordinates always exist and that the necklace ring over a binomial ring is isomorphic (functorially) to the ring of Witt vectors over that ring.

Using this isomorphism one can of course transfer the Frobenius and Verschiebung operations to necklace rings. Actually these are always defined. Here are the explicit formulas.

$$\begin{aligned} &\text{For } c = (c_1, c_2, c_3, \cdots) \quad Nr(A), \quad \mathbf{f}_r c = (b_1, b_2, b_3, \cdots), \quad b_n = \sum_{[j,r]=nr} n^{-1} j c_j \\ &\mathbf{V}_r c = (\underbrace{0,0,\cdots,0,c_1}_{r}, \underbrace{0,0,\cdots,0,c_2}_{r}, \cdots) \end{aligned} \tag{17.27}$$

and one easily proves the usual formulas; either directly or via the ghost components or via ..., such as

$$\mathbf{V}_r\mathbf{V}_s = \mathbf{V}_{rs}, \ \mathbf{f}_r\mathbf{f}_s = \mathbf{f}_{rs}, \ \mathbf{f}_r\mathbf{V}_r = [r], \ \mathbf{f}_r\mathbf{V}_s = \mathbf{V}_s\mathbf{f}_r \ \text{ if } \ (r,s) = 1 \tag{17.28}$$

---

[79] As far as I know this is also the only monograph on the topic.

17.29. *Witt vectors vs necklaces in general.* Let $A$ be a torsion free ring [80]. Consider the mapping

$$\psi_A \colon Nr(A) \longrightarrow \Lambda(A),\quad (c_1, c_2, c_3, \cdots) \mapsto \prod_{n=1}^{\infty} (1 - t^n)^{-c_n} \tag{17.30}$$

This is by the very definition of addition and multiplication on $Nr(A)$ a morphism of rings and it is compatible with the Frobenius and Verschiebung morphims on the two sides.

When $A$ is a binomial ring it is an isomorphism. When $A$ is not binomial this is not the case. However, it is an isomorphism for $A_{\mathbf{Q}} = A \otimes \mathbf{Q}$ and it is (of course) dead easy to describe the subring of $Nr(A_{\mathbf{Q}})$ which corresponds to the subring $\Lambda(A) \subset \Lambda(A_{\mathbf{Q}})$.

For each $\alpha \in A$ let $M(\alpha)$ be the necklace vector

$$M(\alpha) = (M(\alpha;1), M(\alpha;2), M(\alpha;3), \cdots) \in Nr(A_{\mathbf{Q}}) \tag{17.31}$$

Then by the cyclotomic identity and the definition of the Verschiebung operations

$$\psi_{A_{\mathbf{Q}}}(M(\alpha)) = (1 - \alpha t)^{-1},\quad \psi_{A_{\mathbf{Q}}}(\mathbf{V}_r M(\alpha)) = (1 - \alpha t^r)^{-1} \tag{17.32}$$

and it follows that $\psi_{A_{\mathbf{Q}}}$ induces an isomorphism

$$\{\sum_{r=1}^{\infty} \mathbf{V}_r M(\alpha_r) \colon \alpha_i \in A\} \xrightarrow{\ \psi_A\ } \Lambda(A) \subset \Lambda(A_{\mathbf{Q}}) \tag{17.33}$$

By the left half of (17.32) it is clear that the necklace vectors $M(\alpha)$ are the analogs of Teichmüller representatives.

For a binomial ring, like the integers, it follows from (17.31)-(17.33) that each necklace vector can be written uniquely as an infinite sum

$$\sum_{r=1}^{\infty} \mathbf{V}_r M(\alpha_r),\quad M(\alpha_r) = (M(\alpha_r;1).M(\alpha_r;2), M(\alpha_r;3), \cdots) \in Nr(A) \tag{17.34}$$

17.35. *Modified necklace rings*. The description (17.33) of the ring of Witt vectors in terms of necklace vectors is not very satisfactory or elegant. One can in fact do a good deal better using modified necklace polynomials, [316]. In loc. cit. Young-Tak Oh defines modified necklace polynomials over a $\lambda$-ring $A$ which explicitely involve the Adams operations as follows

$$M(r;n) = n^{-1} \sum_{d|n} \mu(d) \Psi^d (r^{n/d}) \tag{17.36}$$

He also introduces much related [81] multivariable versions of the necklace polynomials

---

[80] As far as I know no investigations have been carried out on the necklace ring over $A$ when $A$ is not torsion free, for instance when $A$ is a finite field. Yet it is in that setting that something interesting may happen that is different from what happens in the case of the Witt vectors. Also the analogue of the $p$-typification endomorphism defined by the Frobenius and Verschiebung operations has not been investigated for necklace rings.

[81] This again involves the splitting principle technique/philosophy.

$$M(X;n) = n^{-1}\sum_{d|n}\mu(d)p_d(X) \tag{17.37}$$

where the $p_d$ 's are the power sums in the $X$'s. These polynomials satisfy very similar properties to (17.7) and (17.8) and lead to a definition of necklace ring for $\lambda$-rings which are isomorphic to the rings of Witt vectors.

It would be very nice if there were combinatorial interpretations of these modified necklace polynomials.

17.38. *Adjoints of the inclusion* $\mathbf{BinRing} \quad \mathbf{CRing}$. Now consider the inclusion of the binomial rings into the category of rings. This inclusion functor has both a left adjoint (free objects), $Bin^U$, and a right adjoint (cofree objects), $Bin_U$, characterized and defined by the functorial properties

$$\mathbf{BinRing}(Bin^U(A),B) \quad \mathbf{CRing}(A,B), \quad \mathbf{BinRing}(B,Bin_U(A)) \quad \mathbf{CRing}(B,A)$$

for binomial rings $B$ and rings $A$.

The construction of $Bin^U(A)$ is much like $IVal[X]$ compared to $\mathbf{Z}[X]$. For $Bin_U(A)$ take $\Lambda(A)$ and take the subring of all elements on which all the Frobenius operations are the identity. For more details and proofs see [134].

17.39. The functors $BinW$. The sum comultiplication $\mu_S$ and the product comultiplication $\mu_P$ on $\mathbf{Z}[h]$ that define the ring valued functor of the big Witt vectors, extend uniquely to $IVal[h_1,h_2,h_3,\cdots]$ and $IVal[h_1,h_2,h_3,\cdots,h_n]$ to define coring objects and hence ring-valued sub functors of the Witt vectors and the truncated Witt vectors. These would seem to merit some investigation. Even the simplest one $IVal[h_1]$ already has thought provoking properties, see [45], p.241 ff.

17.40. *Carryless Witt vectors.* Like the real numbers in decimal notation the Witt vectors, when added or multiplied, involve 'carry-overs'. The necklace approach to them can be used to describe a carry-less version, [Metropolis, 1983 #18]. The price, however, is high, too high in my opinion, in that Witt vectors in this approach are seen as equivalence classes rather than single objects.

17.41. *To ponder and muse about.* A nice theory of necklace rings as in [316] works precisely for $\lambda$-rings, i.e. the very objects which are coalgebras for the main functor involved, the one of the Witt vectors. There is something circular about this; almost incestuous, of much the same flavour as one often meets in various parts of category theory. To me this is something that needs to be pondered and mulled over a bit.

17.42. *Apology.* Within the published literature on necklace algebras one finds the notion of what are called aperiodic rings. The formulas involved are, in my opinion, essentially empty given the necklace formulas and, better, Witt vector formulas. So in this chapter I will not discuss these aperiodic rings.

17.43. A selection of references on necklaces and necklace algebras is [21, 22, 55, 70, 115, 120, 144, 145, 174, 260 , 261, 291, 294, 303, 314, 315, 316, 322, 329, 399, 401, 404].

## 18. Symm vs $\bigoplus_n R(S_n)$ and vs $R_{\text{rat}}(\mathbf{GL}_\infty)$

That symmetric functions and representations of the symmetric groups and general linear groups have much to do with one another has been known since the early days of the previous century (Alfred Young, Issai Schur, Georg F Frobenius, ...). The realization that these things become more elegant and better understandable from the Hopf algebra point of view is of a more recent date, [269, 268, 239, 425]. Here, mostly, the latter approach will be outlined. For more details see [199], Ch. 4. First, however, here is a streamlined version of what might be called classical Schur-Frobenius theory as presented in [281], pages 112-114

18.1. *The ring* $R(S)$. Let $S_n$ be the group of permutations on $n$ letters, usually taken to be $\{1,2,\cdots,n\}$. For each $n$ let $R(S_n)$ be the free Abelian (Grothendieck) group spanned by the irreducible representations of the symmetric group $S_n$; or, what is the same, the Grothendieck group of isomorphism classes of (complex) representations of $S_n$, or, what is again the same, the free Abelian group with as basis the irreducible characters of $S_n$. An element from $R(S_n)$ that comes from an actual representation will be called real (as opposed to virtual (not as opposed to complex)). These are the ones that are nonnegative integral sums of irrrducible representation. Other elements of $R(S_n)$ are sometimes referred to as virtual representations.

As an Abelian group $R(S)$ is the direct sum of all these groups of representations:

$$R(S) = \bigoplus_{n=0}^{\infty} R(S_n) \tag{18.2}$$

where, by decree, $R(S_0) = \mathbf{Z}$. A product is defined on $R(S)$ as follows. Take a representation $\rho$ of $S_p$ and a representation $\sigma$ of $S_q$. Taking the (outer) tensor product gives a representation $\rho \otimes \sigma$ of $S_p \times S_q$ on the tensor product of the representation spaces of $\rho$ and $\sigma$. Now consider $S_p \times S_q$ as the subgroup of $S_n = S_{p+q}$ of those permutations that take the set of the first $p$ elements into itself and that take the set of the last $q$ elements into itself. Now induce the representation $\rho \otimes \sigma$ up to $S_n$. Further take $1 \in \mathbf{Z} = R(S_0)$ to be the unit element. This defines an associative and commutative multiplication on $R(S)$ that makes it into a graded ring. As a rule this outer product of two irreducible representations is not irreducible and determining the multiplicities of the irreducible representations that occur in it is and always has been a major part of the representation theory of the symmetric groups (Littlewood-Richardson rule [82]).

18.3. *Scalar product* [83]. If $f$, $g$ are functions on a finite group $G$ their scalar product is defined by

$$\langle f,g\rangle_G = \frac{1}{\#G} \sum_{x \in G} f(x)g(x^{-1}) \tag{18.4}$$

This scalar product is used to define a scalar product on all of $R(S)$ by

$$\langle f,g\rangle = \sum_{n\geq 0} \langle f_n, g_n\rangle_{S_n} \tag{18.5}$$

where $f = \sum f_n,\ g = \sum g_n \in R(S)$ (and $\langle 1,1\rangle = 1$ for $1 \in R_0(S) = \mathbf{Z}$).

[82] [James, 1981 #459], theorem 2.8.13, p. 93.

[83] This inner product is sometimes called 'Hall inner product'.

18.6. *Characteristic map.* Let $w \in S_n$ be a permutation on $n$ letters. It decomposes as a product of disjoint cycles and the lengths of these cycle define a partition $\lambda(w)$ of $n$ called the cycle type of $w$. Define a mapping

$$\psi : S_n \longrightarrow \mathbf{Symm}_n, \quad w \mapsto p_{\lambda(w)} \tag{18.7}$$

Where $p_{\lambda(w)}$ is the power sum monomial, see (9.61), defined by the cycle type of $w$. The next step is the definition of a morphism of Abelian groups called the characteristic map

$$\mathrm{ch}: R(S) \longrightarrow \mathbf{Symm}_{\mathbf{C}} = \mathbf{Symm} \otimes_{\mathbf{Z}} \mathbf{C} \tag{18.8}$$

as follows. If $f$ is a (virtual) character of $S_n$

$$\mathrm{ch}(f) = \langle f, \psi \rangle_{S_n} = \frac{1}{n!} \sum_{w \in S_n} f(w)\psi(w^{-1}) = \frac{1}{n!} \sum_{w \in S_n} f(w)\psi(w) = \sum_{\mathrm{wt}(\lambda)=n} z_\lambda^{-1} f_\lambda p_\lambda \tag{18.9}$$ [84]

where $f_\lambda$ is the value of $f$ on the cycle class $\lambda$ and $z_\lambda$ is the number from (9.62) that gives the scalar product of the power sum monomials with themselves and where it is used that the cycle class of $w \in S_n$ is the same as that of $w^{-1}$. (Recall that $f$ as a character has the same value on each element of a conjugacy class). It follows that, using also $\langle p_\lambda, p_\mu \rangle = z_\lambda \delta_{\lambda,\mu}$,

$$\langle \mathrm{ch}(f), \mathrm{ch}(g) \rangle = \sum_{\mathrm{wt}(\lambda)=n} z_\lambda^{-1} f_\lambda g_\lambda = \langle f, g \rangle_{S_n} \tag{18.10}$$

and hence that ch is an isometry. The basic theorem is now that

18.11. *Theorem.* The characteristic map is (induces) an isometric isomorphism from $R(S)$ onto **Symm**.

There are a certain number of things to prove. But this is not a text on representation theory, let alone a text on the representation theory of the symmetric groups, that vast and fascinating subject; so I will restrict myself to a brief sketch. First, ch is multiplicative. This is handled by Frobenius reciprocity (which will also turn up later in the Hopf algebra approach).

Next one shows that the identity character $\eta_n$, i.e. the character of the trivial representation of $S_n$ corresponds under ch to the complete symmetric function $h_n$. (The sign character corresponds to the elementary symmetric function $e_n$.)

Further for each partition $\lambda$ of $n$ *define*

$$\chi_\lambda = \det(\eta_{\lambda_i - i + j})_{1 \le i,j \le n} \in R(S_n) \tag{18.12}$$

These are (possibly virtual) characters of $S_n$. Their images under ch are the Schur functions and so $\langle \chi_\kappa, \chi_\lambda \rangle = \delta_{\kappa,\lambda}$ and as ch is an isometry they are up to sign irreducible characters. As the number of conjugacy classes of $S_n$ is equal to the number of partitions of $n$ they must form a basis for $R(S_n)$ and so ch indeed induces an isomorphism of $R(S_n)$ onto **Symm**.

---

[84] Where for the second expression in (18.9) definition (18.4) is extended a bit in that $g$ has its values in **Symm**.

It remains to verify that the $\chi_\lambda$ are in fact real characters (as opposed to virtual) which is done by checking their value on the trivial permutations.

18.13. *Enter Hopf algebras.* The topologist knows **Symm** as $H(\mathbf{BU};\mathbf{Z})$, the cohomology of the classifying space **BU** of the unitary group. And he remembers that this is a Hopf algebra. So wouldn't it be nice if $R(S)$ were a Hopf algebra too and if ch were an isomorphism of Hopf algebras. This is indeed the case. The first to notice that $R(S)$ is a Hopf algebra would appear to have been Burroughs, [79]. It was, however, Arunas Liulevicius, [269, 268], who first systematically exploited this point of view. To quote a bit more from the second of his two papers: "The aim of this paper is to present a ridiculously simple proof of a theorem on representations rings of the symmetric groups which concisely presents theorems of Frobenius, Atiyah, and Knutson" (and Schur in my view of things).

Whether this approach is really simpler than the very streamlined presentation of Macdonald outlined above is debatable. A Hopf algebra is a heavy structure. But it certainly brings in new and more functorial and universal points of view, which are, I believe, important.

18.14. *The Hopf algebra structure on* $R(S)$. For a composition $\alpha = [a_1, a_2, \cdots, a_m]$ of weight $n$ the corresponding Young subgroup is $S_\alpha = S_{a_1} \times S_{a_2} \times \cdots \times S_{a_m}$. It consists of all permutations that map all the sets of letters $\{1, \cdots a_1\}$, $\{a_1 + 1, \cdots, a_1 + a_2\}$, $\cdots$ $\{a_1 + \cdots a_{m-1} + 1, \cdots, a_1 + \cdots + a_m = n\}$ to themselves. Define a coproduct structure by

$$\mu_S(\rho) = 1 \otimes \rho + \sum_{i=1}^{n-1} \mathrm{Res}^{S_n}_{S_i \times S_{n-i}}(\rho) + \rho \otimes 1 \in \bigoplus_{i=0}^{n} R(S_i) \otimes R(S_{n-i}) \tag{18.15}$$

and a counit (augmentation)

$$\varepsilon_S \colon RS \longrightarrow \mathbf{Z}, \ \ \varepsilon = \begin{cases} \text{id on } R(S_0) = \mathbf{Z} \\ 0 \ \text{ on } R(S_n), \ n \geq 1 \end{cases} \tag{18.16}$$

The basic fact is now that together with the ring structure already defined and used this makes $R(S)$ a Hopf algebra. The harder part of this is to prove the compatibility of the multiplication with the comultiplication. This is taken care of by the Mackey double coset theorem which describes what happens when a representation is first induced up from a subgroup and then restircted to another subgroup, see [199], Chapt. 4 for details. Actually it is the projection formula which does the job (together with a description of double cosets of Young subgroups of the type $S_i \times S_{n-i} \subset S_n$. This 'induction restriction projection formula' is the following [85]:

18.17. *Theorem.* Projection formula, Frobenius axiom. Let $G$ be a finite group with subgroup $H$, and let $V$ be an $H$-module and $W$ a $G$-module, then

$$\mathrm{Ind}^G_H(V \otimes \mathrm{Res}^G_H(W)) = \mathrm{Ind}^G_H(V) \otimes W \tag{18.18}$$

18.19. *Hopf algebras continue their insidious work.* For an account of how to use the Hopf algebraic structure so far described in the representation theory of the symmetric groups, see the two papers of Liulevicius already quoted. Here I will now continue to describe some more Hopf algebraic structure, culminating in the Zelevinsky structure theorem. This involves further fairly heavy machinery, and certainly does not give the easiest way to get at the representation theory of the symmetric groups. But this way is, in my view, quite important.

[85] Caveat: this is **not** the same instance of a projection formula as occurred in theorem 13.48.

As to the various bits of structure on $R(S)$ the situation is now as follows:

(i) $R(S)$ is a connected graded $\mathbf{Z}$ module.

(ii) There is a preferred basis consisting of $1$ and the irreducible representations with a corresponding inner product for which this basis is orthonomal.

(iii) Multiplication is positive (because taking the outer tensor product of true (real) representations and than inducing up gives a true representation).

(iv) Comultiplication is positive (because restricting a true representation yields true representations.

(v) Comultiplication is multiplication preserving (or, equivalently, the multiplication is comultiplication preserving). This is an immediate consequence of the Mackey double coset theorem combined with the description of double cosets of Young subgroups.

(vi) Comultiplication and multiplication are dual to each other with respect to the inner product. This follows from Frobenius reciprocity in the form that induction is adjoint to restriction (both on the left and on the right) combined with the fact that the scalar product of two characters on a group counts the number of irreducible representations that they have in common.

(vii) The counit morphism is a morphism of algebras.

(viii) There is just one element of the preferred basis that is primitive. This is the unique irreducible representation of $S_1$. Indeed for all $n \geq 2$ a real representation of $S_n$ must restrict to some real (as opposed to virtual) representation of $S_1 \times S_{n-1}$.

18.20. *PSH algebras*. The acronym ‘PSH’ stands for ‘positive selfadjoint Hopf’. This implies sort of that there is also a positive definite inner product and that we are working over something like the integers or the reals, where positive makes sense. What is not mentioned is that these Hopf algebras are also supposed to be connected and graded. As will be seen the assumptions ‘positive, selfadjoint, graded’ are all three very strong; so strong that these algebras can actually be classified. Indeed they are all products of one example [86] and that example is the Hopf algebra of the symmetric functions. The notion is due to Zelevinsky, [425], and the classification theorem is his.

The precise definition is as follows. A PSH algebra is a connected graded Hopf algebra over the integers, so that

$$H = \bigoplus_n H_n, \quad H_0 = \mathbf{Z}, \quad rk(H_n) < \infty \tag{18.21}$$

which is free as an Abelian group and which comes with a given, ‘preferred’ homogeneous basis $\{\omega_i : i \in I\} = \mathcal{B}$. Define an inner product $\langle\ ,\ \rangle$ on H by declaring this basis to be orthonormal. Then the further requirements are:

$$\text{Selfadjointness:} \quad \langle xy,z\rangle = \langle x \otimes y, \mu(z)\rangle,\ \langle \varepsilon(x),a\rangle = \langle x,e(a)\rangle,\ x,y,z \in H, a \in \mathbf{Z} \tag{18.22}$$

$$\text{Positivity: let } \omega_i\omega_j = \sum_r a^r_{i,j}\omega_r,\ \ \mu(\omega_r) = \sum_{i,j} b^{i,j}_r \omega_i \otimes \omega_j \text{, then } a^r_{i,j},\ b^{i,j}_r \geq 0 \tag{18.23}$$

The Zelevinsky classification theorem now says that a PSH algebra with just one primitive amoung the preferred basis elements is isomorphic to the Hopf algebra of the symmetric functions (as a Hopf algebra). The proof proceeds by the inductive construction in any PSH algebra with just one primitive preferred basis element $p$ of a series of elements that behave just like the $h_n, e_n, p_n$ of symmetric function theory. The key observation here is that the powers of

[86] More precisely they are all products of that one example, but the factors are possibly degree shifted.

$p$ must involve all preferred basis elements. In the representations theoretic case $R(S)$ this corresponds to the fact that by definition of the outer multiplication $p^n$ is the regular representation of $S_n$. Here $p$ is the identity representation of $S_1$. For the details of the proof see [425], [199] Chapter 4 [87].

18.24. *Bernstein morphism* [88]. The final step of the proof involves a construction that I think of particular interest, that is not yet well understood, and that deserves more study.

Let $H$ be any graded commutative and associative Hopf algebra. Let

$$\mu_n(x) = \sum_i x_{i,1} \otimes x_{i,2} \otimes \cdots \otimes x_{i,n} \tag{18.25}$$

be the $n$-fold comultiplication written as a sum of tensor products of homogeous components. Now define

$$\beta_n : H \longrightarrow H[\xi_1, \xi_2, \cdots] = H \otimes \mathbf{Z}[\xi_1, \xi_2, \cdots]$$

by

$$\beta_n(x) = \sum_i x_{i,1} x_{i,2} \cdots x_{i,n} \xi_1^{\deg(x_{i,1})} \xi_2^{\deg(x_{i,2})} \cdots \xi_n^{\deg(x_{i,n})} \tag{18.26}$$

Because $H$ is coassociative and cocommutative this is symmetric in the variables $\xi_1, \cdots, \xi_n$ so that there is an induced algebra morphism [89]

$$\beta_n : H \longrightarrow H \otimes \mathbf{Z}[h_1, \cdots, h_n]$$

and because $H$ is graded this stabilizes in $n$ giving the Bernstein morphism [90]

$$\beta : H \longrightarrow H \otimes Symm \tag{18.27}$$

In the case that $H$ is **Symm** this turns out to be the product comultiplication morphism!. The final step of the proof of the Zelevinsky classification theorem is now to compose the Bernstein morphism with a morphism $H \longrightarrow \mathbf{Z}$ that takes the inductively constructed anlogues of the $h_n$ all to one. In the case $H = \mathbf{Symm}$ this is the counit morphism of the product comultiplication morphism $\mu_P$. [91]

---

[87] Th proof as written down in loc. cit. is just my own attempt to write down Zelevinsky's proof in a way that I could understand it.

[88] The construction is due to Joseph N Bernstein.

[89] The fact that this is a morphism of algebras uses commutativity; otherwise multiplication is not an algebra morphism

[90] For a general graded commutative Hopf algebra the Bernstein morphism defines a coaction of **Symm** on $H$ and then by duality also an action of **Symm** on the graded dual of $H$.

[91] Note that without the characteristic map isomorphism or the Zelevinsky theorem intself, it is not yet clear that **Symm** is PSH (but $R(S)$ is). The trouble is positivity; specifically the fact that the product of two Schur symmetric functions is a nonnegative linear combination of Schur functions. This seems a fact that is not so easy to

As promised in 11.39 above, there now follow a few remarks on (0,1)-matrices and their links with Witt vectors and representation theory. This requires some preparation.

18.28. *Majorization ordering.* Let $\alpha = (a_1, a_2, \cdots, a_n)$ and $\beta = (b_1, b_2, \cdots b_n)$ be two vectors of non-negative real numbers of the same $l_1$-norm, i.e. $a_1 + a_2 + \cdots + a_n = b_1 + b_2 + \cdots + b_n$. Denote by $\bar{\alpha} = (\bar{a}_1, \bar{a}_2, \cdots, \bar{a}_n)$ a reordering (rearrangement, permutation) of $\alpha$ such that $\bar{a}_1 \geq \bar{a}_2 \geq \cdots \geq a_n$.

The majorization ordering is now defined by

$$\alpha \geq_{maj} \beta \quad \Leftrightarrow \quad \sum_{i=1}^{r} \bar{a}_i \geq \sum_{i=1}^{r} \bar{b}_i, \ r = 1, 2, \cdots, n \tag{18.29}$$

This ordering occurs is many parts of mathematics under many different names: majority ordering, dominance ordering, natural ordering, specialization ordering, Snapper ordering, Ehresmann ordering, mixing ordering. Parts of mathematics where it plays an important role include: families of algebraic geometric vectorbundles, families of representations, families of nilpotent matrices, Grassmann manifolds, control theory, representation theory, thermodynamics, convex function theory (Schur convex functions [92]), doubly stochastic matrices, (0,1) matrices, inequality theory (Muirhead inequalities, a far reaching generalization of the geometric mean - arithmetic mean inequality), representation theory, ... . See [3]. Many of these uses of the majorization ordering are related, see [202].

18.30. *Conjugate partition*. Let $\alpha = [a_1, a_2, \cdots, a_m]$ be a partition of $n$. Then the conjugate partition $\alpha^{conj} = [a'_1, a'_2, \cdots, a'_{m'}]$ is defined by

$$a'_i = \#\{j \colon a_j \geq i\} \tag{18.31}$$

So, for instance, $[4,4,3,2,1,1,1]^{conj} = [7,4,3,2]$. (If the partition is displayed as a diagram (either in the French or Anglo-Saxon manner), the conjugate partition looks at columns instead of rows.)

There are (at least) three [93] applications of the majorization ordering in the representation theory of the symmetric groups. In addition there is the Gale-Ryser theorem which is of immediate relvance here as it deals with the existence of (0,1)-matrices and thus has things to say about the second comultiplication on **Symm** (and hence the multiplication of Witt vectors.

18.32. *Gale-Ryser theorem*. Let $\alpha$ and $\beta$ be two partitions. then there is a $(0,1)-$matrix with row sum $\alpha$ and column sum $\beta$ if and only if $\alpha^{conj} \geq_{maj} \beta$.

Of course a similar theorem holds for compositions instead of partitions; this amounts to taking permutations of columns and permutations of rows.

For example take the $(0,1)-$matrix

establish directly (without going to representation theory) (and, hence, is a bit of a blemish on symmetric function theory).

[92] Same Schur; totally different topic.

[93] Anothe (related one), due to Kraft and de Concini, is too far from the present topic to discuss.

$$M = \begin{pmatrix} 0 & 1 & 0 & 0 & 1 \\ 1 & 0 & 1 & 0 & 1 \\ 0 & 1 & 0 & 1 & 0 \end{pmatrix}$$

This one has as row sum the composition [2,3,2] with associated partition [3,2,2]; and the column sum is the composition [1,2,1,1,2] with associated partition [2,2,1,1,1]. Also $[3,2,2]^{conj} = [3,3,1]$. And, indeed $[3,3,1] \geq_{maj} [2,2,1,1,1]$.

Of the three theorems in the representation theory of the symmetric groups that involve the majorization ordering the one closest to the present concerns (Witt vectors) is the Snapper Liebler-Vitale Lam theorem.

18.33. *Snapper Liebler-Vitale Lam theorem.* Let $\alpha$ and $\beta$ be partitions and let $S_\alpha$ be the Young subgroup defined by $\alpha$. Then $\rho_\beta$, the irreducible representation corresponding to $\beta$, occurs in $Ind_{S_\alpha}^{S_n}(I)$ if and only if $\alpha \leq_{maj} \beta$.

Here $I$ stands for the trivial representation of $S_\alpha$. For a proof see [203] and/or the references therein.

For completeness sake and possible future applications, below there is the Ruch-Schönhofer theorem. Under the isomorphism between **Symm** and $R(S)$ the trivial representation of $S_n$ corresponds to the complete symmetric function $h_n$ and the sign representation $A_{S_n}$ corresponds to the elementary symmetric function $e_n$. So this theorem could certainly be relevant.

18.34. *Ruch-Schönhofer theorem.* The representations $\mathrm{Ind}_{S_\alpha}^{S_n}(I_{S_\alpha})$ and $Ind_{S_\beta}^{S_n}(A_{S_\beta})$ have an irreducible representation in common (which is the same as saying that their inner product is nonzero), if and only if $\alpha \leq_{maj} \beta$.

18.35. *Outer plethysm and inner plethysm.* In section 16.67 the composition operator 'outer plethysm' on **Symm** made its appearance. But **Symm** is isomorphic to $R(S)$ as Hopf algebras (with extra product multiplication and extra coproduct multiplication with counit). Thus there is an outer plethysm operation on $R(S)$. The problem is to describe this operation in representation theoretic terms. This has been called the'outer plethysm problem'.

The other way: each of the summands of $R(S) = \bigoplus_{n=0}^{\infty} R(S_n)$ is a $\lambda-$ ring in its own right, giving rise to plethym operations called inner plethysm, and the 'inner plethysm problem' is to describe these in symmetric function terms when transferred to **Symm**.

18.36. *Outer plethysm problem.* Outer plethysm gives in particular a composition of $I_{S_n}$, the trivial representation of $S_n$ with any $\sigma \in R(S_k)$. The result is denoted $h_n(\alpha)$ in [239], p. 135. The trivial (or identity) representation of $S_n$ corresponds to $h_n \in$ **Symm**, whence the notation.

In loc. cit. , p.135 Donald Ivar Knutson guesses [94] the following representation theoretic description of $h_n(\alpha)$. Given a representation $\alpha$ of $S_k$ on $V$ first construct the induced representation of the wreath product $S_n[S_k]$ on $V^{n}$, then induce up the representation obtained to $S_{nk}$ using the natural inclusion $S_n[S_k] \subset S_{nk}$.

[94] His word.

This 'guess' has since turned out to be correct, see [221, 313, 394].

18.37. *Inner plethysm problem.* The $\lambda-$ ring structure on $R(S_n)$ gives of course rise to Frobenius operators which via the isomorphism give rise to (degree preserving) Adams-Frobenius type operators on the homogenous component of weight (degree) $n$ of **Symm**. Using the inner product there are corresponding adjoint operators. These are described nicely in section 3 of [349].

Besides the five references just quoted here is a selection of further references on plethysm: [8, 43, 42, 59, 66, 82, 267, 265, 266, 282, 306, 312, 313, 341, 351, 385, 394, 395, 429]

## 19. Burnside rings.

Let $G$ be a finite group. The Burnside ring $B(G)$ is one of the fundamental (representation like) rings attached to $G$. It is the Grothendieck ring of finite $G$-sets with sum and product induced by disjoint union and Cartesian product respectively. Or, equivalently, the ring of permutation representations.

So as a group $B(G)$ is a finitely generated Abelian group with as (canonical) basis the transitive $G$-sets (orbits). According to some the Burnside ring was introduced by Andreas Dress in [117], where he proved a somewhat unexpected kind of result, viz. that a finite group is solvable if and only if the spectrum of its Burnside ring is Zariski connected, i.e. if and only if $0$ and $1$ are the only idempotents in $B(G)$.

Others attribute the introduction of the concept to Louis Solomon [372].

By now there is a substantial literature on Burnside rings: the ZMath database lists at the moment [95] 305 publications with "Burnside ring" in title or abstract. For a fair selection see the bibliography of [61].

"[The Burnside ring] is in many ways the universal object to consider when looking at the category of G-sets. It can be viewed an an anlogue of the ring **Z** of integers for this category."

"The ring $B(G)$ is also functorial with respect to $G$ and subgroups of $G$ and this leads to the Mackey functor or Green functor point of view. ... The Burnside Mackey functor is a typical example of projective Mackey functor. It is also a universal object in the category of Green functors."

These are two quotes from the introduction of [61]. So if the Burnside ring of an arbitrary finite group is already as nice as all that, what about the Burnside ring of a really nice group like the integers. That one must then be super-nice. And it is. It turns out to be the ring of Witt vectors of the integers, $W(\mathbf{Z})$, a discovery due to Andreas Dress and Christian Siebeneicher, [115, 120]. [96]

19.1. *G-sets*. Let $G$ be a group. A $G$-set is a set $X$ together with an action (on the left) of $G$ on it; that is a mapping $G \times X \longrightarrow X, (g,x) \mapsto gx$ such that $g(h(x)) = (gh)(x)$, $g,h \in G, x \in X$ and $e(x) = x$ where $e$ is the identity element of $G$. A morphism of $G$-sets is a mapping $f: X \longrightarrow Y$ such that $f(gx) = g(f(x))$. This defines the category $G$-**Set** of $G$-sets.

[95] 19 April 2008.

[96] This also of course, given the 'universality remark' just quoted, provides a seventh universality property of the Witt vectros. Working out what this one really means is still an open matter.

For $x \in X$ the subgroup

$$G_x = \{g \in G: gx = x\} \tag{19.2}$$

is the stabilizer of $x \in X$ and

$$Gx = \{gx: g \in G\} \tag{19.3}$$

is the orbit of $x$ (or through $x$). If $H$ is a subgroup multiplication on the left induces a $G$-set structure on the set of left cosets $G/H$ which has a single orbit (a transitive $G$-set). The orbit $Gx$ of an element is isomorphic as a $G$-set to $G/G_x$.

A $G$-set $X$ gives rise to a permutation representation $\rho(X)$ over any ring of which the underlying module is the free module with basis $X$ and with the action of $G$ given by the the permutations of basis elements determined by the $G$-set structure of $X$.

19.4. *Induction and restriction for $G$-sets.* Let $H$ be a subgroup of a group $G$. Then if $X$ is a $G$-set retricting the action to $H$ gives an $H$-ser, defining a functor $\mathrm{Res}_H^G: G-\mathbf{Set} \longrightarrow H-\mathbf{Set}$.

As is to be expected there is also a functor the other way called induction. Let $Y$ be an $H$-set. Consider the set $G \times Y$ with the equivalence relation determined by

$$(g,y) \sim (gh^{-1}, hy),\ g \in G, h \in H, y \in Y \tag{19.5}$$

The set of equivalence classes is suggestively denoted $G \times_H Y$. (Note that $(g,y) \mapsto (gh^{-1}, hy)$ defines an action of $H$ on $G \times Y$ and that the equivalence classes are the orbits for this action.) Multiplication on the left $(g',(g,y)) \mapsto (g'g, y)$ induces an action of $G$ on $G \times_H Y$ giving a $G$-set denoted $\mathrm{Ind}_H^G(Y)$. It is obvious how to make this into a functor $\mathrm{Ind}_H^G: H-\mathbf{Set} \longrightarrow G-\mathbf{Set}$. Induction is left adjoint to restriction (but not right adjoint in general). There is also a product formula (projection formula, Frobenius identity)

$$X \times \mathrm{Ind}_H^G(Y) \cong \mathrm{Ind}_H^G(\mathrm{Res}_H^G(X) \times Y) \tag{19.6}$$

For more details on this and some more related material see e.g. [409] p. 811, [61] pp 744-745.

Restriction and induction are compatible with the notions of the same name for representations; i.e. they are compatible with the mapping $\rho$ that assigns to a $G$-set its associated permutation representation.

19.7. *Almost finite G-sets.* Given a group $G$ its profinite completion $\hat{G}$ is its completion for the topology of normal subgroups of finite index. Or, equivalently it is the projective limit

$$\hat{G} = \varprojlim G/N\ , \quad N \text{ runs over the normal subgroups of finite index} \tag{19.8}$$

A $G$-set $X$ for a group $G$ is almost finite if for each $x \in X$ the stabilizer subgroup $G_x$ is a subgroup of finite index (so that all orbits are finite) and such that subsets of invariants

$$X^U = \{x \quad X : G_x \quad U\} \tag{19.9}$$

are finite for every subgroup of finite index. (So that there are only finitely many orbits of isomorphism type $G/U$ for each such $U$.)

It does not matter whether one works in this with the group $G$ or its profinite completion. In the case of a profinite group one works with $G$-spaces. These are Hausdorff spaces with a continuous action [97] of $G$. See [115], p. 5. Most papers that deal with the present subject (Burnside rings and Witt-Burnside rings) take the profinite point of view.

For a subgroup of finite index $H \quad G$ restriction takes almost finite $G$-sets into almost finite $H$-sets and induction takes almost finite $H$-sets into almost finite $G$-sets.

Finite disjoint unions and and finite Cartesian products of almost finite $G$-sets are almost finite.

19.10. *Burnside theorem.* [78] Chapter XII theorem 1. Let $G$ be a finite group and $X$ and $Y$ finite G-sets. Then the following are equivalent

(i) The $G$-sets $X$ and $Y$ are isomorphic

(ii) For any subgroup $H$ of $G$ the sets of invariants $X^H$ and $Y^H$ have the same cardinality.

This still holds for almost finite $G$-sets where one only need consider subgroups of finite index.

19.11. *Burnside ring.* The (completed) Burnside ring $\hat{B}(G)$ of a group is now defined as the Grothendieck ring of the category of almost finite $G$-sets.

It can also be defined as the projective limit of the usual Burnside rings (of finite groups) $B(G/N)$ where $N$ ranges over allnormal subgroups (resp. closed normal subgroups) of finite index [98].

19.12. *Almost finite cyclic sets*. In the case of the group of integers all this specializes as follows. A cyclic set is simply a set with a left action of the group of integers on it. That is it is a set with a specified bijection. A cyclic set is almost finite if every orbit is of finite length and if for every $n \quad \mathbf{N}$ there are only finitely many orbits of length $n$. Obviously (finite) disjoint unions and (finite) Cartesian products of almost finite cyclic sets are again almost finite and so there is the Grothendieck ring $\hat{B}(\mathbf{Z})$ of almost finite cyclic sets.

If $B(\mathbf{Z})$ denotes the Grothendieck ring of finite cyclic sets, $\hat{B}(\mathbf{Z})$ is the completion of $B(\mathbf{Z})$ under a suitable natural topology on $B(\mathbf{Z})$, see [120], p. 3 and below.

19.13. *A remarkable commutative diagram. The* $W(\mathbf{Z})$, $\hat{B}(\mathbf{Z})$, $\Lambda(\mathbf{Z})$, $Nr(\mathbf{Z})$ *isomorphisms.* . The ring $\hat{B}(\mathbf{Z})$ fits in the following commutative diagram.

---

[97] Here $G$ is given the topology defined by the normal subgroups of finite index (which define the profinite structure) and the set $X$ is given the discrete topology.

[98] This explains the notation $\hat{B}$ and the terminology ‘completed’.

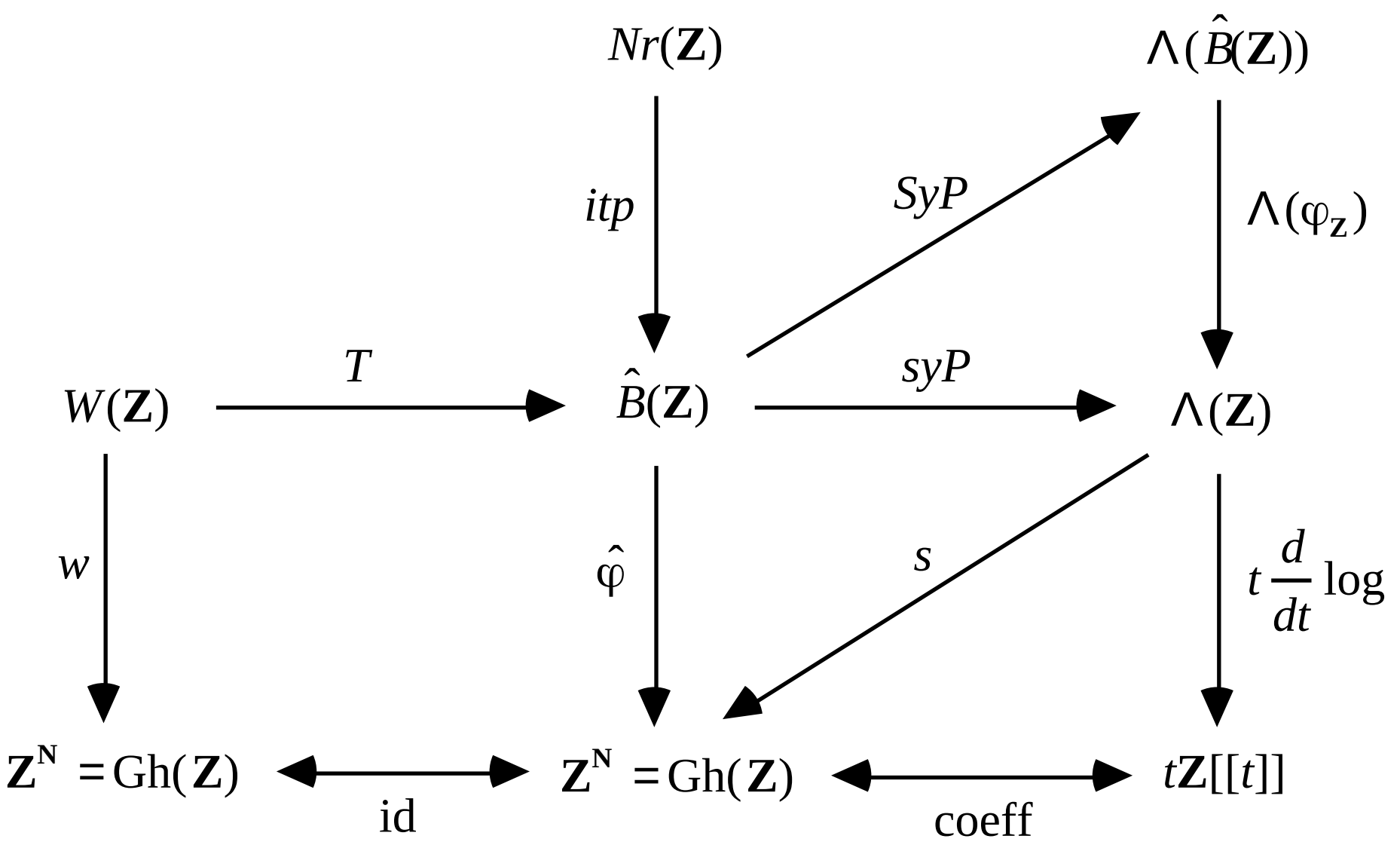

(19.14)

All the horizontal arrows in this diagram are isomorphisms (of rings with operators) and so is the morphism denoted '*itp*' (which stands for 'interpretation').. The morphisms w, $\hat{\varphi}$, s, $t\frac{d}{dt}\log$ are ghost component morphisms and injective. Finally, $\varphi_{\mathbf{Z}}$ and $\Lambda(\varphi_{\mathbf{Z}})$ are surjections and *SyP* is an injection. Those morphisms which have not already occurred in earlier sections, such as w, s and the logarithmic derivative, will be elucidated below of course.

One of the more remarkable aspects of this diagram is the fact that the composite $syP \circ T$ precisely embodies the very nasty coordinate change from Witt vector coordinates to power series coordinates encoded by the power series identity

$$\prod_n \frac{1}{(1-x_n t^n)} = 1 + a_1 t + a_2 t^2 + \cdots \qquad (19.15)$$

I consider this a main contribution from [120]. (The necklace coordinates also fit in as suggested by the diagram.)

19.16. *The ghost component morphism* $\hat{\varphi}: \hat{B}(\mathbf{Z}) \longrightarrow \mathrm{Gh}(\mathbf{Z}) = \mathbf{Z}^{\mathbf{N}}$. Given an almost finite cyclic set X and an element $n \quad \mathbf{N}$, consider the subgroup of finite index $n\mathbf{Z} \quad \mathbf{Z}$ and define

$$\varphi_{n\mathbf{Z}}(X) = \#\{x \quad X : x \text{ is invariant under } n\mathbf{Z}\} \qquad (19.17)$$

This induces ring morphisms $\hat{B}(\mathbf{Z}) \longrightarrow \mathbf{Z}$ and these combine to define the ghost component ring morphism $\hat{\varphi}: \hat{B} \longrightarrow Gh(\mathbf{Z})$. The morphism $\hat{\varphi}$ is injective because of the (extension to almost finite sets of the) Burnside theorem 19.10.

Take the restrictions of the $\varphi_{n\mathbf{Z}}$ to $B(\mathbf{Z})$ and give $B(\mathbf{Z})$ the coarsest topology for which all these restrictions are continuous (with $\mathbf{Z}$ discrete). Then $\hat{B}(\mathbf{Z})$ is the completion of $B(\mathbf{Z})$ for this topology.

The morphism $\hat{\varphi}$ is injective; it is not surjective; an element $y \quad \mathrm{Gh}(\mathbf{Z})$ is in its image if and only if

$$\sum_{d|n} \mu(d^{-1}n)y(d) \equiv 0 \mod n \tag{19.18}$$

for all $n \in \mathbf{N}$. The same condition turned up in connection with the necklace ring in section 17 above ([113, 118, 126], [391] pp 11-12).

19.19. *The isomorphism itp*: $Nr(\mathbf{Z}) \longrightarrow \hat{B}(\mathbf{Z})$. The only transitive almost finite cyclic sets are the coset spaces $C_n = \mathbf{Z}/n\mathbf{Z}$. So every element of $\hat{B}(\mathbf{Z})$ is a (possibly infinite) sum

$$\sum_n b_n C_n\,,\quad b_n \in \mathbf{Z} \tag{19.20}$$

As a group $Nr(\mathbf{Z}) = \mathbf{Z}^{\mathbf{N}}$ with coordinate wise addition. So assigning to an element $\beta = (b_1, b_2, b_3, \cdots) \in \mathbf{Z}^{\mathbf{N}}$ the formal difference of almost finite cyclic sets

$$\sum_{b_n>0} b_n C_n - \sum_{b_n<0} (-b_n) C_n$$

defines an isomorphism '*itp*' of Abelian groups $Nr(\mathbf{Z}) \longrightarrow \hat{B}(\mathbf{Z})$. It remains to see how *itp* behaves with respect to multiplication.

Now observe [99] that the product of two transitive cyclic sets $C_r$ and $C_s$ decomposes as the sum of $(r,s)$ copies of $C_{[r,s]}$. Here $(r,s)$ is the greatest common divisor of $r$ and $s$ and $[r,s]$ is their least common multiple. Given the definition of the multiplication on the necklace ring, see 17.15, it follows that *itp* is an isomorphism of rings.

19.21. *The isomorphism* $T$: $W(\mathbf{Z}) \longrightarrow \hat{B}(\mathbf{Z})$. This isomorphism of rings is denoted $\tau$ in [120] and called 'Teichmüller" there. This may be (indeed, is) appropriate in some sense but is also potentially confusing in view of the Teichmüller representatives mapping $\tau: A \longrightarrow W(A)$ in Witt vector theory.

The definition of the isomorphism $T$ involves induction. Let $X$ be a $\mathbf{Z}$-set. Via the isomorphism $\mathbf{Z} \longrightarrow n\mathbf{Z} \subset \mathbf{Z}$, $r \mapsto nr$ it can be seen as an $n\mathbf{Z}$-set. Now induce [100] this one up to $\mathbf{Z}$. The result is denoted $\mathrm{ind}_n(X)$. In concreto this works out as follows. Consider the product $\mathbf{Z} \times X$ and the equivalence relation on it defined by $(z,x) \sim (z - nu, ux)$, $u \in \mathbf{Z}$. Then $\mathrm{ind}_n(X)$ is the set of equivalence classes of this equivalence relation with the action induced by the left action of $\mathbf{Z}$ on itself, $(z',(z,x)) \mapsto (z'+z, x)$.

It readily follows that

$$\mathrm{ind}_n(C_r) = C_{rn}$$

and hence that

$$\mathrm{ind}_n\Big(\sum_r b_r C_r\Big) = \sum_r b'_r C_r \tag{19.22}$$

[99] It appears, see [120] p. 5, that it was this observation that lead to the investigations of Andreas Dress and Christian Siebeneicher that culminated in [115, 120, 126].

[100] See 19.4 above.

where

$$b'_r = \begin{cases} b_{r/n} & \text{if } n \text{ divides } r \\ 0 & \text{otherwise} \end{cases} \tag{19.23}$$

and hence that

$$\varphi_{r\mathbf{Z}}(\mathrm{ind}_n(X)) = \begin{cases} n\varphi_{r/n}(X) & \text{if } n \text{ divides } r \\ 0 & \text{otherwise} \end{cases} \tag{19.24}$$

So, at the ghost component level $\mathrm{ind}_n$ behaves just like Verschiebung in the case of the big Witt vectors.

Now let $q \quad \mathbf{N} \quad \{0\}$. With $q^{(\mathbf{Z})}$ denote the set of maps $\mathbf{Z} \longrightarrow \{1,2,\cdots,q\}$ that factor through some set of cosets $\mathbf{Z}/n\mathbf{Z}$. That is the continous maps $\mathbf{Z} \longrightarrow \{1,2,\cdots,q\}$ where $\{1,2,\cdots,q\}$ has the discrete topology and $\mathbf{Z}$ the topology of subgroups of finite index. The action of $\mathbf{Z}$ is given by

$$(zf)(z') = f(z' - z) \tag{19.25}$$

Note that $0^{(\mathbf{Z})} = \quad$. It is fairly immediate that

$$(qq')^{(\mathbf{Z})} = q^{(\mathbf{Z})} \times q'^{(\mathbf{Z})}, \quad \varphi_{n\mathbf{Z}}(q^{(\mathbf{Z})}) = q^n \tag{19.26}$$

$$\hat{\varphi}(q^{(\mathbf{Z})}) = (1, q, q^2, q^3, \cdots) = \mathrm{coeff}(1-q)^{-1} \tag{19.27}$$

Finally define

$$T\colon W(\mathbf{Z}) \longrightarrow \hat{B}(\mathbf{Z}),\ (x_1, x_2, x_3, \cdots) \mapsto \sum_{n=1}^{\infty} \mathrm{ind}_n(x_n^{(\mathbf{Z})}) \tag{19.28}$$

Then it follows from what has been noted just above that

$$\varphi_{n\mathbf{Z}}(x_1, x_2, x_3, \cdots) = \sum_{d|n} dx_d^{n/d} \tag{19.29}$$

and that

$$\prod_{n=1}^{\infty}(1 - x_n t^n)^{-1} = \prod_{n=1}^{\infty}(1 - t^n)^{-b_n} \quad \Leftrightarrow \quad T(x_1, x_2, x_3, \cdots) = \sum_{i=1}^{\infty} b_i C_i \tag{19.30}$$

Formula (19.29) shows that the lower left square in the diagram (19.14) is commutative. Using the characterizations of the image of the ghost components morphisms w and $\hat{\varphi}$, see (9.95) and (19.20) above, and their injectivity, it follows that $T$ is an isomorphism.

Further (19.30) says that the composite morphism $(itp)^{-1} \circ T$ exactly embodies the coordinate transformation between Witt vector coordinates and necklace coordinates encoded by the power series identity

$$\prod_{n=1}^{\infty}(1-x_n t^n)^{-1} = \prod_{n=1}^{\infty}(1-t^n)^{-b_n}$$

19.31. *Symmetric powers of* $G$-sets. Given a $G$-set $X$ its $n$-th symmetric power is obtained by first taking the $n$-fold Cartesian product $X^n$ with diagonal action

$$g(x_1,x_2,\cdots,x_n)=(gx_1,gx_2,\cdots,gx_n) \tag{19.32}$$

The symmetric group $S_n$ acts on this by permuting coordinates

$$\sigma(x_1,x_2,\cdots,x_n)=(x_{\sigma^{-1}(1)},x_{\sigma^{-1}(2)},\cdots,x_{\sigma^{-1}(n)})$$

and this action commutes with the diagonal action of $G$. So there is an induced action of $G$ on the set of orbits $S^nX$. This is the $n$-th symmetric power of the $G$-set $X$. Alternatively, and better, an element from $S^nX$ can be described as a multiset of size $n$, i.e. a function $f: X \longrightarrow \mathbf{N} \cup \{0\}$ such that $\sum_{x\in X} f(x) = n$. In this picture the action is $(gf)(x) = f(g^{-1}x)$. The symmetric powers of almost finite $G$-sets are almost finite.

It follows rather immediately from the construction that

$$S^n(X \coprod Y) = \coprod_{i+j=n} S^iX \times S^jY \tag{19.33}$$

and it follows that the map

$$X \mapsto \sum_{n=1}^{\infty}(S^nX)t^n \tag{19.34}$$

induces a morphism of Abelian groups

$$SyP: \hat{B}(\mathbf{Z}) \longrightarrow \Lambda(\hat{B}(\mathbf{Z})) \tag{19.35}$$

Now compose this with $\Lambda(\varphi_{\mathbf{Z}}): \Lambda(\hat{B}(\mathbf{Z})) \longrightarrow \Lambda(\mathbf{Z})$ to obtain a morphism of Abelian groups

$$syP: \hat{B}(\mathbf{Z}) \longrightarrow \Lambda(\mathbf{Z}) \tag{19.36}$$

(and the commutative upper triangle of diagram (19.14)). Thus the coefficient of $t^n$ in $syP(X)$ is the number of invariant elements in the $n$-th symmetric power $S^nX$.

It now turns out that in fact $syP$ is an isomorphism of rings and even an isomorphism of $\lambda-$ rings. The key to that is the easy observation that

$$\varphi_{\mathbf{Z}}(S^nC_r) = \begin{cases} 1 & \text{if } r \text{ divides } n \\ 0 & \text{otherwise} \end{cases} \tag{19.37}$$

(Indeed, a function $f:\{0,1,\cdots,r-1\}\longrightarrow \mathbf{N}\quad \{0\}$ is invariant under ‘left shift modulo $r$ of its argument’ if and only if all its values are equal.)

Formula (19.37) serves to prove that the lower right square of diagram (19.14) is commutative (as well as the two triangles there, that $syP$ is an isomorphism of rings, and that the composed morphism $syP\circ itp$ embodies the coordinate change encoded by the power series identity

$$\prod_{n=1}^{\infty}(1-t^n)^{-b_n} = 1+a_1t+a_2t^2+a_3t^3+ \qquad (19.38)$$

When combined with (19.30) this shows that the composite morphism $syP\circ T: W(\mathbf{Z})\longrightarrow \Lambda(\mathbf{Z})$ precisely gives the nasty coordinate change formulas between Wiit vector coordinates and power series coordinates encoded by the power series identity

$$\prod_{n=1}^{\infty}(1-x_nt^n)^{-1} = 1+a_1t+a_2t^2+a_3t^3+\cdots$$

which I feel is a major insight from [120].

19.39. *Frobenius and Verschiebung on* $W(\mathbf{Z})$, $\hat{B}(\mathbf{Z})$, $\Lambda(\mathbf{Z})$. On all the rings in diagram (19.14) there are Frobenius and Verschiebung operators. All have been defined before except the ones on $\hat{B}(\mathbf{Z})$. On this Burnside ring ‘Verschiebung’ is ‘ind’ and Frobenius is restriction (followed by the identification coming from the obvious (canonical) isomorphism of Abelian groups $\mathbf{Z}\longrightarrow n\mathbf{Z}$.

It now turns out (and is verified with no great difficulty) that all morphism indicated in the diagram are compatible with the Frobenius operators and all except $SyP$ are compatible with Verschiebung. ($SyP$ is most definitely not compatible with Verschiebung.)

As a matter of fact, identifying $\hat{B}(\mathbf{Z})$ with $\Lambda(\mathbf{Z})$ via syP turns SyP into the morphism of $\lambda$-rings $\sigma_t:\Lambda(\mathbf{Z})\longrightarrow \Lambda(\Lambda(\mathbf{Z}))$ of $\lambda$-ring theory [101], see section 16 above.

19.40. This concludes the treatment here of the Burnside ring of the integers (or its profinite completion $\hat{\mathbf{Z}}$). I have pretty much followed [120] apart from some interspersed remarks. I hope and believe that the outline above is sufficient that the reader can fill in all details. But if needed they can be found in loc. cit.

However, there is more to the Burnside ring story in connection with Witt vectors. For every profinite group there is a Witt vector like functor $\mathbf{CRing}\longrightarrow \mathbf{CRing}$ called the Witt Burnside functor. Below there is the main theorem from [126] about them. First some definitions and notation.

19.41. Let $G$ be a profinite group. That is a projective limit $G=\varprojlim G/N$ of finite quotients group. Give $G$ the topology defined by the collection of normal subgroups $N$ of finite index. A $G$-space is a Hausdorff spcae on which $G$ acts continuously from the left. Such a $G$-space is almost finite if it is discrete and if for any open subgroup $U\quad G$ the number of invariants under $U$ is finite. Set

$$\varphi^U(X)=\#\{x\quad X: ux=x \text{ for all } u\quad U\} \qquad (19.42)$$

[101] This is one more indication that it is better to work with symmetric powers than with exterior powers.

For any open subgroup $U$ there is the transitive almost finite (in fact finite) $G$-space of left cosets $G/U$. Denote with $osg(G)$ the set of all open subgroups of $G$ and with $cosg(G)$ the quotient set of conjugacy classes of open subgroups.

Finally let $\hat{B}(G)$ be the "completed Burnside ring of $G$, that is the Grothendieck ring of (isomorphism classes) of almost finite $G$-spaces with addition induced by disjoint union and multiplication induced by the Cartesian product of $G$-spaces.

19.43. *Existence theorem of the Witt-Burnside functors* [126]. Let $G$ be a profinite grouup. There exists a unique covariant functor $W_G$ from the category of unital commutative rings to itself such that as a set $W_G(A)$ is the set $A^{cosg(G)}$ of all functions from $cosg(G)$ to $A$. that is all functions from $osg(G)$ to $A$ that are constant on conjugacy classes of subgroups, and with $W_G(h)\colon W_G(A) \longrightarrow W_G(B)$ given by composition $\alpha \mapsto \alpha \circ h, \alpha \in W_G(A)$, such that for all open subgroups the map

$$\psi_U^A\colon W_G(A) \longrightarrow A \tag{19.44}$$

defined by

$$\sum_{U \le_{scjg} V \le G} \varphi^U(G/V)\alpha(V)^{(V:U)} \tag{19.45}$$

is a natural transformation of functors from $W_G$ to the identity. Here the sum (9.45) is over all $V$ such that $U$ is subconjugate to $U$ (denoted $U \le_{scjg} V$), which means that there is a $g \in G$ such that $U$ is a subgroup of $gVg^{-1}$, the 'index' $(V{:}U) = (G{:}U)/(G{:}V)$,. and in the sum (9.45) exactly one summand is taken for each conjugacy class of subgroups $V$ with $U \le_{scjg} V$.

Moreover, $W_G(\mathbf{Z}) = \hat{B}(G)$ and $W_{\mathbf{Z}} = W$ the functor of the big Witt vectors.

There are also Frobenius and Verschiebung like functorial endomorphisms coming respectively from restriction and induction. The Frobenius morphism are ring morphisms, the Verschiebung morphism are morphisms of Abelian groups. These morphisms have the usual kinds of properties.

Since the appearance of the foundational papers [126, 120] a number of other papers have appeared on the topic of Witt-Burnside functors, giving refinements, further developments, applications and interrelations, simplifications and complications; see [73, 116, 122, 135, 156, 175, 318, 319, 315].

19.46. $\beta$*-rings*. It seems clear from [370] that there is no good way to define a $\lambda$-ring structure on Burnside rings, see also [158]. There are (at least) two different choices giving pre-$\lambda$-rings but neither is guaranteed to yield a $\lambda$-ring. Of the two the symmetric power construction seems to work best.

Instead one needs what are called $\beta$-operations [102] (power operations), first introduced in [56]., and $\beta$-rings. There is a $\beta$-operation for each conjugacy class of subgroups of the symmetric groups $S_n$ and they are constructed by means of symmetric powers of $G$-sets. Every

[102] I have tried to track down where the appellation 'beta' comes from. Unsuccessfully. But it seems likely that it is some sort of philosophical mix between the 'B' from Burnside and the '$\lambda$ ' from $\lambda$-ring.

$\lambda$-ring is a $\beta$-ring (but not vice versa).

There is some vagueness about what is precisely the right definition of a $\beta$-ring, see e.g. the second paragraph on page 2 of [183] and the second comment at the end of §3 of that preprint.

However, it seems clear that the free $\beta$-ring on one generator must be the ring

$$B(S) \;=\; \bigoplus_{n=0}^{\infty} B(S_n) \tag{19.47}$$

This is the direct sum of the Burnside rings of the symmetric groups (with $B(S_0) = \mathbf{Z}$ by decree) equipped with an outer product defined completely analogously as in the case of $R(S)$, see section 18 above. I.e.

$$XY = \mathrm{Ind}_{S_i \times S_j}^{S_{i+j}} (X \times Y) \tag{19.48}$$

There is also a coproduct making it a Hopf algebra and the underlying ring is again a ring of polynomials in countably many indeterminates over the integers.

The natural morphism

$$\rho \colon B(S) \longrightarrow R(S) \tag{19.49}$$

that assings to an $S_n$-set the corresponding permutation representation is a surjective (but not injective) morphism of Hopf algebras and also a morphism of $\lambda$-rings, see [183] p. 11.

A selection of papers on $\beta$-rings is [56, 158, 183, 300, 301, 344, 397, 396, 433].

**Coda**.
Some people write about the Witt polynomials as mysterious polynomials that come out of nowhere [103]. To me they are so elegant and natural that they simply cry out for deep study. This is my main reason for presenting things here as I have done above, even though, as has been shown, they are not really needed.

Should a reader inadvertently get really interested in the Witt vector ring functors he/she is recommended to work through the 57 exercises on the subject in [63], 26 pp worth, just for the formulations only, mostly contributed, I have been told, by Pierre Cartier. He/she can then continue with the exercises on Cohen rings in same volume which also involve a fair amount of Witt vector stuff.

**Appendix. The algebra of symmetric functions in infinitely many indeterminates**.
In several places in the sections above there occur expressions like

$$\prod_{i=1}^{\infty} (1 + \xi_i t + \xi_i^2 t^2 + \xi_i^3 t^3 + \cdots ) \tag{A.1}$$

and statements that the coefficients of each power in $t$ in (A.1) are symmetric functions in the infinity of commuting variables $\xi_1, \xi_2, \xi_3, \cdots$ and that they are, hence, polynomials in the elementary (or complete) symmetric functions in the infinity of commuting indeterminates $\xi_1, \xi_2, \xi_3, \cdots$ . The few pages in this appendix are meant for those who feel (rightly) a bit nervous about statements like this even though the meaning seems intuitively clear [104].

[103] Also, see section 3 above, they definitely do not come out of nowhere.

A.2. *Power series in infinitely many variables*. Let $I$ be any set and let $\{\xi_i : i \in I\}$ be a corresponding set of commuting indeterminates. An exponent sequence for $I$ is a map $e\colon I \longrightarrow \mathbf{N} \cup \{0\}$ of finite support. I.e. there are only finitely many $i \in I$ for which $e(i) \neq 0$. Let $E(I)$ be the set of all exponent sequences. For each $e \in E(I)$ introduce a symbol $\xi^e$. Giving $I$ a total order, one can think of $\xi^e$ as

$$\xi^e = \prod_{e(i)\neq 0} \xi_i^{e(i)} \tag{A.3}$$

where the product (monomial) on the right hand side is written down in the order specified by the ordering of $I$ [105]. Two monomials, i.e. symbols, $\xi^e$ are multiplied by the rule $\xi^e \xi^{e'} = \xi^{e+e'}$ where $e+e'$ is the point-wise sum $(e+e')(i) = e(i) + e'(i)$. The ring of formal power series over a base ring $k$ is now defined as

$$k\langle\langle \xi_i : i \in I \rangle\rangle = \{ \sum_{e \in E(I)} a_e \xi^e : a_e \in k \} \tag{A.4}$$

Two such expressions as on the right of (A.3) are multiplied by the rule

$$(\sum_{e \in E(I)} a_e \xi^e)(\sum_{e \in E(I)} b_e \xi^e) = (\sum_{e \in E(I)} c_e \xi^e), \quad c_e = \sum_{e'+e''=e} a_{e'} b_{e''} \tag{A.5}$$ [106]

This product rule is well defined because of the finite support condition on exponent sequences.

A product like (A.1) is now by definition interpreted as the formal power series one obtains by multiplying any finite number of factors from it. That is

$$\begin{gathered}\prod_{i=1}^{\infty} (1 + \sum_{\mathrm{wt}(e_i)\geq 1} a_{i,e_i} \xi^{e_i} t^{\mathrm{wt}(e_i)}) = 1 + \sum_{\mathrm{wt}(e)\geq 1} c_e \xi^e t^{\mathrm{wt}(e)}, \\ \mathrm{wt}(e) = \sum_i e(i), \quad c_e = \sum_{\substack{e_1+\cdots+e_r=e \\ \mathrm{wt}(e_j)\geq 1}} a_{i_1,e_1} \cdots a_{i_r,e_r}\end{gathered} \tag{A.6}$$

and similar products with the ‘counting variable’ $t$ left out.

A.7. *Symmetric and quasisymmetric power series*. Probably every one (who is likely to get his fingers on this chapter) knows, or can easily guess at, what is a symmetric polynomial (over the integers, or any other base ring $k$).

Say, the polynomial is in $n$ (commuting) variables $\xi_1, \xi_2, \cdots, \xi_n$. Then the polynomial $f = f(\xi_1, \xi_2, \cdots, \xi_n)$ is symmetric iff for every permutation $\sigma$ of $\{1, 2, \cdots, n\}$

---

[104] To my mind the matters touched upon here illustrate well issues in foundational mathematics (intuitionism) having to do with realized infinities vs potential infinities.

[105] But this is not needed.

[106] All this just formalizes what everyone knows intuitively. One can also do these things for noncommuting indeterminates, and things are (curiously enough?) actually easier in that case. For totally ordered index sets one can go much further and make sense of infinite ordered products and sums by using injective limit ideas to make sense of infinite ordered sums and products of elements of the base field and the integers by treating them as sequences with two sequences equal if they eventually agree.

$$f(\xi_{\sigma(1)}, \xi_{\sigma(2)}, \cdots, \xi_{\sigma(n)}) = f(\xi_1, \xi_2, \cdots, \xi_n) \tag{A.8}$$

And, in that case, the main theorem of symmetric functions says that: the polynomial $f$ is a polynomial $p_f$ in the complete symmetric functions $h_1, h_2, \cdots, h_n$ or, equivalently, in the elementary symmetric functions, and "the polynomial $p_f$ is independent of the number of variables involved provided there are enough of them" (meaning more than or equal to the degree of $f$). The last phrase needs explaining. Also this strongly suggests that the best way to work with symmetric polynomials is to take an infinity of variables.

A little bit more notation is useful. Let $S_I$ be the group (under composition) of all bijections $\sigma\colon I \longrightarrow I$ such that there are only finitely many $i \in I$ with $\sigma(i) \neq i$. For $S_{\mathbf{N}}$ the notation $S_\infty$ is usually used. For an exponent sequence $e$ let $e^\sigma$ be the exponent sequence $e^\sigma(i) = e(\sigma^{-1}(i))$[107].

A power series $\sum_{e \in E(I)} a_e \xi^e$ is now symmetric if and only if $a_e = a_{e^\sigma}$ for all $\sigma \in S_I$. A polynomial in an infinite set of commutating indeterminates is a finite sum $\sum_{e \in E(I)} a_e \xi^e$ i.e. a power series with all but finitely many of the coefficients $a_e$ unequal to zero. A power series $\sum_{e \in E(I)} a_e \xi^e$ is of bounded degree if and only if there is a natural number $n$ such that $a_e = 0$ for all $e$ with $\mathrm{wt}(e) > n$.. The complete symmetric functions in the $\xi$ are by definition

$$h_n = \sum_{\mathrm{wt}(e)=n} \xi^e \tag{A.9}$$

The main theorem for symmetric power series in an infinity of indeterminates now takes the following form.

A.10. *Theorem*. Every symmetric power series can be uniquely written as a power series in the complete symmetric functions. A bounded degree symmetric power series is a polynomial in the complete symmetric functions.

There is of course a similar theorem in terms of the elementary symmetric functions.

The algebra **Symm** is the algebra of polynomials in the countable set of complete symmetric functions $h_n$ for the case of the index set $I = \mathbf{N}$.

A.11. *Projective limit description.* For the finicky (or pernickety) the following projective limit construction is perhaps more congenial. For each $n$ consider the algebra morphism

$$\pi_{n+1,n}\colon k[\xi_1, \cdots, \xi_n, \xi_{n+1}] \longrightarrow k[\xi_1, \cdots, \xi_n],\ \xi_i \mapsto \xi_i \ \text{ for } i \in \{1, \cdots, n\},\ \xi_{n+1} \mapsto 0$$

These are graded algebra morphisms

Let $\mathbf{Symm}_k^{(n)}$ be the subalgebra of symmetric polynomials $k[x_1, \cdots, x_n]$. This gives a graded projective systems of graded algebras

$$\pi_{n+1,n}\colon \mathbf{Symm}_k^{(n+1)} \longrightarrow \mathbf{Symm}_k^{(n)}$$

and **Symm** is the graded projective limit of this system. So, for instance, the symmetric power series $\sum_{e \in E(\mathbf{N})} \xi^e$ is not in **Symm**. In this picture theorem A.10 is obtained by the usual theorem for symmetric functions in a finite number of variables and using that for each given degree the coefficients involved in that degree stabilize as $n \to \infty$.

[107] The exponent '-1' is there to ensure that $e^{\sigma\tau} = (e^\sigma)^\tau$ ; not that that is important in the present context.